\PassOptionsToPackage{unicode}{hyperref}
\PassOptionsToPackage{hyphens}{url}
\PassOptionsToPackage{dvipsnames,svgnames,x11names}{xcolor}
\documentclass[
  11pt]{article}
\usepackage{xcolor}
\usepackage[margin=1in]{geometry}
\usepackage{amsmath,amssymb}
\usepackage{iftex}
\ifPDFTeX
  \usepackage[T1]{fontenc}
  \usepackage[utf8]{inputenc}
  \usepackage{textcomp} 
\else 
  \usepackage{unicode-math} 
  \defaultfontfeatures{Scale=MatchLowercase}
  \defaultfontfeatures[\rmfamily]{Ligatures=TeX,Scale=1}
\fi
\usepackage[]{lmodern}
\ifPDFTeX\else
\fi
\IfFileExists{upquote.sty}{\usepackage{upquote}}{}
\IfFileExists{microtype.sty}{
  \usepackage[]{microtype}
  \UseMicrotypeSet[protrusion]{basicmath} 
}{}
\makeatletter
\@ifundefined{KOMAClassName}{
  \IfFileExists{parskip.sty}{%
    \usepackage{parskip}
  }{
    \setlength{\parindent}{0pt}
    \setlength{\parskip}{6pt plus 2pt minus 1pt}}
}{
  \KOMAoptions{parskip=half}}
\makeatother
\usepackage{longtable,booktabs,array}
\usepackage{calc} 
\usepackage{etoolbox}
\makeatletter
\patchcmd\longtable{\par}{\if@noskipsec\mbox{}\fi\par}{}{}
\makeatother
\IfFileExists{footnotehyper.sty}{\usepackage{footnotehyper}}{\usepackage{footnote}}
\makesavenoteenv{longtable}
\usepackage{graphicx}
\makeatletter
\newsavebox\pandoc@box
\newcommand*\pandocbounded[1]{
  \sbox\pandoc@box{#1}%
  \Gscale@div\@tempa{\textheight}{\dimexpr\ht\pandoc@box+\dp\pandoc@box\relax}%
  \Gscale@div\@tempb{\linewidth}{\wd\pandoc@box}%
  \ifdim\@tempb\p@<\@tempa\p@\let\@tempa\@tempb\fi
  \ifdim\@tempa\p@<\p@\scalebox{\@tempa}{\usebox\pandoc@box}%
  \else\usebox{\pandoc@box}%
  \fi%
}
\def\fps@figure{htbp}
\makeatother
\NewDocumentCommand\citeproctext{}{}

\makeatletter
 \let\@cite@ofmt\@firstofone
 \def\@biblabel#1{}
 \def\@cite#1#2{{#1\if@tempswa , #2\fi}}
\makeatother
\newlength{\cslhangindent}
\newlength{\csllabelwidth}
\newenvironment{CSLReferences}[2] 
 {\begin{list}{}{%
  \setlength{\itemindent}{0pt}
  \setlength{\leftmargin}{0pt}
  \setlength{\parsep}{0pt}
  \ifodd #1
   \setlength{\leftmargin}{\cslhangindent}
   \setlength{\itemindent}{-1\cslhangindent}
  \fi
  \setlength{\itemsep}{#2\baselineskip}}}
 {\end{list}}
\usepackage{calc}

\providecommand{\tightlist}{%
  \setlength{\itemsep}{0pt}\setlength{\parskip}{0pt}}
\usepackage{bookmark}
\IfFileExists{xurl.sty}{\usepackage{xurl}}{} 
\hypersetup{
  pdftitle={Guaranteed indirect delivery in directed networks},
  pdfauthor={Angel Ivanov Raychev angel.ivanov.raychev@gmail.com},
  colorlinks=true,
  linkcolor={blue},
  filecolor={Maroon},
  citecolor={Blue},
  urlcolor={blue},
  pdfcreator={LaTeX via pandoc}}

\ifXeTeX

\fi
\ifPDFTeX
\pdftrailerid{}
\pdftrailer{/ID [ <90d4b0e3ee7715a3525aa795d9cb9f92> <90d4b0e3ee7715a3525aa795d9cb9f92> ]}
\fi
\ifLuaTeX
\pdfvariable trailerid {[ <90d4b0e3ee7715a3525aa795d9cb9f92> <90d4b0e3ee7715a3525aa795d9cb9f92> ]}
\fi

\title{Guaranteed indirect delivery in directed networks}
\usepackage{etoolbox}
\makeatletter
\providecommand{\subtitle}[1]{
  \apptocmd{\@title}{\par {\large #1 \par}}{}{}
}
\makeatother
\subtitle{The joint extremal function and explicit attaining
constructions}
\author{Angel Ivanov Raychev\\
\href{mailto:angel.ivanov.raychev@gmail.com}{\nolinkurl{angel.ivanov.raychev@gmail.com}}}
\date{24 October 2026}

\usepackage{needspace}
\begin{document}
\maketitle
\begin{abstract}
A messenger and a pursuer move alternately on a finite loopless directed
graph, with complete information and waiting permitted. The messenger
moves first and wins upon reaching a designated recipient, before
interception there. We determine the maximum number \(M(n,m)\) of
ordered nonadjacent source--recipient pairs for which delivery is
guaranteed against every initial pursuer position, among graphs with
\(n\) vertices and \(m\) arcs. We give the exact value and an explicit
attaining graph for every admissible pair \((n,m)\). For \(n\ge12\), a
piecewise formula describes sparse growth, the interval
\(2n\le m\le n(n-3)\) in which every missing pair can be guaranteed, a
dense staircase, and the zero range above \(n(n-2)\). Complete finite
tables cover smaller orders. The proofs combine degree and support
bounds, local strategy composition, protected routing, and constructions
in the missing-arc graph. Finite exceptional cases use explicit strategy
certificates and complete computational exclusions. Reproducible
programs and certificates accompany the paper; a Lean development
verifies the complete formula for \(n\ge12\) and \(m\ge24\), the claimed
scores of prescribed constructions at every admissible pair, and further
finite upper bounds. Its remaining obligations are stated separately.
\end{abstract}

{
\setcounter{tocdepth}{1}
\small
\makeatletter
\renewcommand{\l@section}{\@dottedtocline{1}{0em}{2.3em}}
\makeatother
\tableofcontents
}
\clearpage
\section{Results and scope}\label{results-and-scope}

\subsection{Introduction}\label{introduction}

A messenger has just obtained a message at a station \(s\) of a directed
network. A pursuer occupies another station, and both positions are
visible. They move alternately, the messenger first, using one directed
connection or waiting. The messenger wants to deliver to a designated
recipient \(t\) before interception. Receipt completes the task
immediately, even if the pursuer is at \(t\): the pursuer can intercept
the messenger but cannot disable the recipient.

A direct connection \(s\to t\) makes delivery immediate. We study the
additional communication supported by the network: source--recipient
pairs without a direct connection for which delivery is nevertheless
guaranteed against every allowed initial pursuer position and policy.
Write \(I(G)\) for their number.

We determine \[M(n,m)=\max_{|V(G)|=n,\ |E(G)|=m} I(G)\] for
\textbf{every admissible pair of integers} \(n\ge1\) and
\(0\le m\le n(n-1)\), and give explicit adjacency rules or a specified
composition procedure attaining each value, with parameters determined
from \((n,m)\) and a proof of guaranteed delivery for its counted pairs.
The vertex-only and arc-only maxima follow by taking the corresponding
maxima: \[
M_V(n)=\max_{0\le m\le n(n-1)}M(n,m),\qquad
M_E(m)=\max_{n\ge1:\,m\le n(n-1)}M(n,m).
\]

Isolated vertices are allowed. Removing them preserves \(I\), so the
arc-only maximum is finite: a graph with \(m>0\) arcs has at most \(2m\)
nonisolated vertices. Our constructions prescribe adjacency or
composition directly from \((n,m)\); finite exceptions are specified by
certified lists of arcs. The results also include separate probabilistic
and deterministic derandomization statements. Time bounds establish
termination and, in several families, stronger guarantees. Optimization
of delivery time is a separate problem, discussed in the
\href{https://angelraychev.com/cop-prisoner/delivery-time/}{companion}.

The original problem and a layered subset construction arose in
Alexandra Ignatova's 2023 project, mentored by Angel Raychev. A public
2023 abstract records that provenance (Ignatova 2023); the retained
manuscript is dated February 2024 (Ignatova 2024). Ignatova and Raychev
had already established the upper bound \(I(G)\le n(n-3)\) and
investigated two-in/two-out-regular constructions. This account
reconstructs that bound, distinguishes the original construction from
its corrected analysis, and supplies new attaining families among the
September 2026 extensions. The present account completes the joint
extremal classification with explicit attaining constructions.

\subsubsection{Relation to pursuit games and
reachability}\label{relation-to-pursuit-games-and-reachability}

Destination-based pursuit has a classical antecedent in the
cat-and-mouse game. Chandra and Stockmeyer's work on alternation
introduced a directed version (Chandra and Stockmeyer 1976). Greenlaw,
Hoover, and Ruzzo record a variant in which the mouse moves first, both
players may wait, and the cat is forbidden to enter the goal (Greenlaw
et al. 1995, Appendix A.11.2). Huq supplies P-completeness proofs for
directed and undirected variants and compares the rules (Huq 2015). In
his formulation the cat moves first, repeated positions give a draw, and
moves traverse edges. Goal priority also has an earlier precedent:
Downey and Thilikos describe a directed cat-and-mice game in which a
mouse wins upon reaching the goal even if the cat is already there
(Downey and Thilikos 2011, sec. 3.3). Thus this receipt rule itself is
not claimed as new.

Here there is one messenger, it moves first, both players may wait, and
receipt takes priority over interception at the recipient. For each
source--recipient pair we quantify over every permitted initial pursuer
position, then maximize the number of guaranteed nonadjacent ordered
pairs among graphs with prescribed \((n,m)\). This joint extremal
problem, together with explicit attaining graphs, is the subject of the
paper. We do not transfer a complexity classification between the cited
rule sets.

Guarding games give another destination-based comparison. Fomin et al.
(Fomin, Golovach, Hall, et al. 2011) and Fomin, Golovach, and Lokshtanov
(Fomin, Golovach, and Lokshtanov 2011) study an intruder entering a
protected region. Their guards remain inside that region and occupied
entry vertices are blocked. Here the pursuer can move throughout the
graph and receipt succeeds even at an occupied recipient. Their
guard-number and decision problems therefore have different rules and
objectives from the present fixed-resource pair count.

The traditional Cops and Robbers game studies capture or indefinite
evasion; see Nowakowski and Winkler (Nowakowski and Winkler 1983) and
the monograph of Bonato and Nowakowski (Bonato and Nowakowski 2011).
Delivery requires reaching a specified vertex, so an evasion strategy or
a usual cop-win characterization alone does not resolve our question.
Once the recipient is fixed, the present game is a finite reachability
game on positions of both players and the turn. Its attractor and
rank-certificate arguments are standard finite-game methods (Berwanger
2009). Our extremal question concerns their simultaneous consequences
for all source--recipient pairs and all graphs with prescribed vertex
and arc counts. State-digraph labels and relational characterizations
for generalized pursuit games also appear in Bonato and MacGillivray
(Bonato and MacGillivray 2017). We use this established finite-game
viewpoint; the classification is not a new general method for solving
reachability games.

The difficulty is not a shortage of paths alone. An additional arc can
provide a route to the messenger while also enabling an interception by
the pursuer. Accordingly, neither \(I(G)\) nor the validity of an
attaining construction is monotone under arbitrary addition of arcs. We
address this by proving strategies against a larger permitted pursuer
graph while restricting the messenger to a compulsory subgraph. Every
graph between those two graphs then inherits the same guarantee. The
separating set viewpoint also connects the two-move constructions to
earlier work on separating systems and oriented diameter-two graphs
(Bollobás and Scott 2007); the directed safety condition remains an
additional requirement in the present game.

Some attaining examples are already consequences of classical tournament
adjacency properties. A tournament is \textbf{3-existentially closed} if
every orientation pattern between three distinct vertices and a new
vertex is realized. Such tournaments and their constructions were
studied by Bonato, Gordinowicz, and Prałat (Bonato et al. 2010); Javier,
Llano, and Zuazua (Javier et al. 2026, Corollary 4.3) prove this
property for Paley tournaments of prime order \(p\equiv3\pmod4\),
\(p\ge19\). For distinct \(s,t,c\), choose a witness \(u\) with
\(s\to u\), \(u\to t\), and \(u\to c\). The last arrow excludes
\(c\to u\), so \(s\to u\to t\) is safe against the pursuer's response.
If \(c=t\), the implied two-vertex extension property gives
\(s\to u\to t\) and again excludes \(t\to u\). Thus every missing pair
in a 3-existentially closed tournament is guaranteed within two
messenger moves. Since a tournament has \(n(n-1)/2\) arcs, these known
graphs already attain \(M(n,n(n-1)/2)=n(n-1)/2\) at their orders. This
conclusion is an elementary consequence of their adjacency property
under our rules; it is not a claim that the cited papers formulate this
extremal function. The protected-routing theorem below requires only
specified extension patterns and additionally controls the exact budget
and arbitrary arrows outside the routing set.

Other ingredients are likewise classical: incomparable hub subsets are
Sperner families (Sperner 1928), the finite-field parabola is a Sidon
set (Erdős and Turán 1941; Eberhard and Manners 2023), and the routing
derandomization uses conditional expectations (Alon and Spencer 2016).
Their delivery-specific hypotheses and applications are proved below.
The contribution claimed here is the complete joint extremal
classification, with matching upper bounds and prescribed attaining
graphs at every admissible \((n,m)\), including the sparse and dense
endpoint regimes and the finite exceptions. The game, these general
tools, the classical Sidon envelope, and the tournament-density examples
just described are not claimed as new.

Our computational proof method follows the principle of certifying
algorithms: a producer supplies a witness and a separate checker
verifies its meaning (McConnell et al. 2011; Alkassar et al. 2014). An
attaining graph certifies a lower bound; an upper bound also needs a
complete exclusion of better graphs. The latter requires proved coverage
of degree profiles and every branch discarded by pruning. Lean 4 (Moura
and Ullrich 2021) checks the formal checker theorems and the concrete
certificates. The verification section states the current coverage; an
untrusted search report is never promoted to a Lean theorem by its
numerical agreement.

\subsubsection{Proof roadmap}\label{proof-roadmap}

The complete formula and finite tables appear first. Source/sink
deletion and live-row/live-column counting provide the sparse upper
bounds, while counting failed pairs in the missing-arc graph gives the
dense staircase. Periodic local strategies attain the sparse endpoints;
robust circles and hub/guide constructions cover the intervening
budgets. Missing-arc cores and functional attachments attain the dense
band. Finally, degree-profile reductions and complete finite exclusions
establish the exceptional small values. The verification section
specifies which steps have ordinary proofs, which use finite
computation, and which have Lean verification.

\subsubsection{Results and proof
boundaries}\label{results-and-proof-boundaries}

\textbf{Complete classification.} For every \(n\ge1\) and
\(0\le m\le n(n-1)\), we determine \(M(n,m)\) and supply an explicit
attaining graph. The combined formula below covers every order
\(n\ge12\); the complete finite tables cover orders one through eleven.
Ordinary parameterized constructions cover the infinite ranges.
Certified listed graphs cover the genuine finite exceptions.

The classification has two ingredients:

\begin{enumerate}
\def\labelenumi{\arabic{enumi}.}
\tightlist
\item
  \textbf{A formula at every order from twelve.} Put \(A_m=M_E(m)\) for
  \(0\le m\le23\). The small-budget table specifies these 24 constants.
  For larger budgets, the sparse parity formulas, missing-pair main
  interval, final dense staircase, and zero range determine every value.
  The sparse graph is padded with isolates; the other cases use the
  stated circle, stage, hub, code, and dense-block recipes. The main
  formula starts sharply at twelve because \(M(11,22)=79<88\); the dense
  staircase starts sharply at eleven. The full main interval has no
  uniform two-move deadline.
\item
  \textbf{Complete small-host values.} The finite sparse-budget theorem
  settles the global values at twelve through twenty-one arrows.
  Source/sink deletion, live-row and live-column capacities, and
  complete joint-degree exclusions resolve the exceptional smaller
  hosts. The finite tables give every remaining value, with literal
  adjacency and full winning/losing certificates. In particular,
  \(M_V(7)=21\), \(M_V(8)=32\), \(M_V(9)=48\), \(M_V(10)=65\), and
  \(M_V(11)=85\). From twelve onward, \(M_V(n)=n(n-3)\).
\end{enumerate}

Both parts include constructions with parameters fixed by \((n,m)\).
Separate reproducibility checks verify the numerical coverage and the
availability of direct adjacency, prescribed composition, or a saved
finite graph for each case. These checks supplement the coverage and
strategy proofs; numerical equality alone is not an explicit
construction.

The complete mathematical classification uses ordinary proofs and finite
computer-assisted proofs. Independent reviews check the finite domain
reductions, pruning, and the actual alternating game; their precise
replay scopes are stated with the evidence. The Lean development
formalizes selected structural implications and examples, not this
entire classification.

The proof uses a shared source/sink deletion and positive-support
induction for the sparse upper bounds, separated pair codes for a
shorter main-interval construction, and a functional attachment identity
for dense graphs. The following construction and support results retain
their full parameter ranges and stronger guarantees:

\begin{itemize}
\tightlist
\item
  \textbf{Protected routing (Theorem 54).} For a specified routing set
  of size \(b\), put \(Q=b(2n-b-1)\) and \(p=a/Q\). If
  \((b-3)p^2(1-p)>4\log n+\log(8/3)\), exactly \(a\) incident arrows can
  be chosen deterministically so that every assignment of exterior
  arrows preserves universal two-move delivery. This gives every exact
  budget \(a\le m\le a+(n-b)(n-b-1)\) and subsumes the former
  full-vertex conditioning theorem, including its random-sampling
  guarantee.
\item
  \textbf{Sparse support (Theorems 56--57).} With \(2k\) arrows,
  \(k\ge6\), a graph outside the \(k\)-vertex two-regular shape has
  \(I\le(k-1)(k-3)\). With \(2k+1\) arrows, \(k\ge3\), a graph with
  \(I>k^2-4k+6\) reduces to a \(k\)-vertex induced core with at least
  the same score, or to three specified \(k+1\)-vertex degree families.
  A stronger original-graph reduction shows that, for \(k\ge7\), a graph
  with \(2k\) arrows and \(I>k^2-5k+8\) already has exactly \(k\)
  nonisolated vertices, allowing arbitrary degrees in that core. These
  quantitative statements retain the regular equality information and do
  not assume that deleting arrows preserves strategies.
\item
  \textbf{Direct dense construction (Theorem 58).} Its two-move endpoint
  works at every \(n\ge10\). Its arbitrary-subset bridge covers
  \(2n\le h\le3n\) at every \(n\ge16\). Two explicit residue
  constructions enlarge the dense envelope to \(n^2/4-O(n)\) missing
  arrows, preserving two-move delivery and every optional subset. Its
  first-step two-move domain is \(n\ge16\), and its original full-band
  domain is \(n\ge20\). Combining its blocks with root-colored cycles
  extends the full-band two-move domain to \(n\ge16\).
  Arbitrary-duration block addition and the fixed attachments give the
  uniform staircase from order eleven; safe-root exclusions and literal
  graphs determine every smaller dense-band case. The smaller-order
  delivery bounds remain separate.
\item
  \textbf{Antichain composition (Theorem 59).} Incomparable hub subsets
  attach outside vertices to a fixed two-part core, preserving two-move
  delivery under every choice of outside arrows. The generalized theorem
  includes separated two-element codes: for \(k\ge6\) and
  \(2k\le n\le k(k-3)\) they give every budget from \((4k-2)n-2k(2k+1)\)
  through \(n(n-3)\). Their overlap with circles covers
  \(3n\ldots n(n-3)\) for every \(n\ge79\). The supplied larger and
  mixed codes retain their independent ranges; the two-move sparsity
  consequence is retained in the companion.
\item
  \textbf{Support concentration (Theorem 60).} With
  \(\eta=m-2\sqrt{I(G)}\), the nonisolated support has at most
  \(m/2+4\eta\) vertices; all but \(12\eta\) arrows lie in a core of
  size between \(m/2-3\eta\) and \(m/2\). This controls near-extremizers
  below the earlier exact-support cutoffs without asserting that
  deleting the exceptional arrows preserves delivery.
\end{itemize}

The general arguments have ordinary mathematical proofs, with the
main-interval finite bridge additionally using complete independently
replayed state certificates. The accompanying Lean development verifies
actual-game semantics, selected bounds and exact cells, finite
certificate interfaces, and strategy implications for routing and
envelope constructions. The dense staircase bound and explicit attainers
are formalized throughout the dense region from order eleven in
\texttt{DenseSharpClassification.\allowbreak{}lean}; dense cells at smaller orders
retain their separate proof boundaries. Its precise formal coverage is
stated below; the complete classification is not a Lean theorem.

\subsection{The combined joint
maximum}\label{the-combined-joint-maximum}

The sparse constructions can be placed inside a larger network by adding
isolated vertices. This observation joins the arc-budget extrema to the
two dense regimes and gives a single eventual formula.

\textbf{Corollary 47 (complete joint classification).} Let \(n\ge12\)
and \(0\le m\le n(n-1)\). Put \(A_m\) equal to the following listed
value:

{\def\LTcaptype{none} 
\Needspace{5\baselineskip}
\begin{longtable}[]{@{}lrrrrrrrrrrrr@{}}
\toprule\noalign{}
\(m\) & 0 & 1 & 2 & 3 & 4 & 5 & 6 & 7 & 8 & 9 & 10 & 11 \\
\midrule\noalign{}
\endhead
\bottomrule\noalign{}
\endlastfoot
\(A_m\) & 0 & 0 & 0 & 0 & 1 & 1 & 2 & 2 & 4 & 5 & 6 & 8 \\
\end{longtable}
}

{\def\LTcaptype{none} 
\Needspace{5\baselineskip}
\begin{longtable}[]{@{}lrrrrrrrrrrrr@{}}
\toprule\noalign{}
\(m\) & 12 & 13 & 14 & 15 & 16 & 17 & 18 & 19 & 20 & 21 & 22 & 23 \\
\midrule\noalign{}
\endhead
\bottomrule\noalign{}
\endlastfoot
\(A_m\) & 9 & 11 & 13 & 15 & 19 & 28 & 33 & 40 & 53 & 60 & 79 & 83 \\
\end{longtable}
}

Then \[
M(n,m)=
\begin{cases}
A_m,&0\le m\le23,\\
k(k-3),&m=2k,\quad24\le m\le2n,\\
k(k-3)+2,&m=2k+1,\quad25\le m<2n,\\
n(n-1)-m,&2n\le m\le n(n-3),\\
2n-r-\lceil r/2\rceil-1,
&m=n(n-3)+r,\quad1\le r\le n,\\
0,&n(n-2)<m\le n(n-1).
\end{cases}
\] Every value has an explicit attaining construction. For
\(1\le n\le11\), the complete tables in the finite-values section supply
every admissible value and its certified listed graph or specified
composition.

The two formulas at \(m=2n\) agree. For even \(m\), the first formula
also holds for every \(24\le m\le2n\) on every host order \(n\ge m/2\).

\begin{figure}
\centering
\includegraphics[width=1\linewidth,height=\textheight,keepaspectratio,alt={The positive exact regimes at n=115 for m\textbackslash ge24: sparse growth, the main interval, and the final dense band. Separate axis scales show the two narrow endpoint regimes. The dotted continuation in the last panel leads to zero above n(n-2).}]{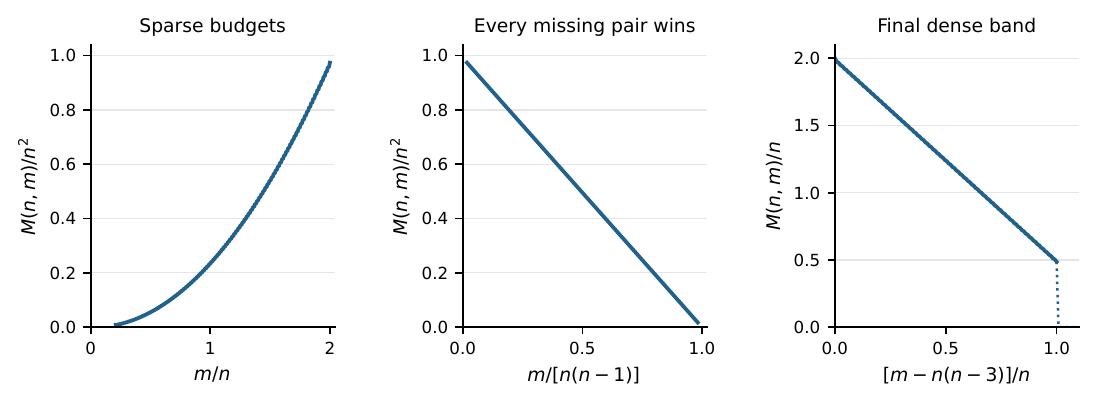}
\caption{The positive exact regimes at \(n=115\) for \(m\ge24\): sparse
growth, the main interval, and the final dense band. Separate axis
scales show the two narrow endpoint regimes. The dotted continuation in
the last panel leads to zero above \(n(n-2)\).}
\end{figure}

\textbf{Proof.} For \(0\le m\le11\), use the complete small-budget
theorem and its listed witnesses. The finite sparse-budget theorem
supplies \(12\le m\le21\); its largest least attaining host is twelve.
Thus every \(A_m\) through21 is attained at every \(n\ge12\) by adding
isolates.

The explicit eleven-vertex, 22-arrow extremizer and its matching
arc-only bound give \(M(n,22)=79\) after adding isolates. Theorem 53
gives \(M(n,23)=83\) on every host of order at least twelve.

For even \(m=2k\ge24\), Theorem 40 constructs an attaining \(k\)-vertex
graph. If \(k\le n\), add \(n-k\) isolated vertices. This does not
change \(I\): pairs between components have no delivery path; a pursuer
starting at an added isolate stays there and cannot obstruct an existing
delivery strategy. The arc-only upper bound supplies the matching upper
bound for the joint problem.

For odd \(m=2k+1\ge25\), Corollary 45.1 constructs an attaining graph on
\(k+1=(m+1)/2\) vertices. When \(m<2n\), its support has at most \(n\)
vertices. Add isolates and apply the odd arc-only upper bound. This
gives \(k(k-3)+2\).

The complete main-interval theorem supplies the third case from order
twelve, preserving Theorem 44 and all its earlier constructions. The
complete-band corollary in the root-colored cycle section supplies the
fourth from order eleven. Four moves suffice for that whole dense band;
its stronger two-move corollary applies from order sixteen. Finally
Theorem 13 gives zero above \(n(n-2)\). The ranges cover every stated
integer budget. The finite-values section and its complete source-pinned
registry cover all smaller orders. \(\square\)

At order eleven the three initial values are
\[M(11,22)=79,\qquad M(11,23)=78,\qquad M(11,24)=79.\] For
\(25\le m\le88\), the saved finite stages give \(M(11,m)=110-m\). The
displayed staircase and zero branches then cover all remaining budgets.
Thus every \(n\ge11\), \(22\le m\le n(n-1)\) has an exact value and an
explicit attaining graph. The special order-eleven values are proved in
the small-order section; they are not obtained by extending the sparse
parity formula below its hypotheses.

The two parts together determine every admissible integer pair,
including the small host exceptions. Stronger optional-subset and
delivery-time conclusions remain attached to the construction theorems
that prove them; the classification does not transfer those properties
to every other attainer of the same value.

\subsection{Complete values and explicit
constructions}\label{the-remaining-classification-domain}

Every admissible pair \(n\ge1\), \(0\le m\le n(n-1)\) is determined,
with an explicit attaining graph. Corollary 47 gives the full formula
for \(n\ge12\); the complete finite tables cover \(1\le n\le11\). No
numerical or explicit-construction gap remains within the stated
objective.

The numerical and construction obligations were checked separately. The
independent closure in \texttt{research/\allowbreak{}finish-\allowbreak{}ledger-\allowbreak{}review.\allowbreak{}json}
preserves every previous exact value and closes all 101 former
representatives. All 493,924 admissible cells in the retained audit box
through order 114 have matching bounds. This finite check is
supplemented by the ordinary parameterized construction proofs and host
reductions; it is not used to infer an infinite theorem merely from
finite success.

For budgets \(m=12,13,\ldots,21\), the least attaining host orders are
\[7,8,8,8,8,10,9,11,10,12.\] Adding isolates gives every larger host at
each corresponding budget. The remaining smaller-host values are
explicitly tabulated. Above those budgets, the established sparse, main,
dense, and zero constructions cover the infinite ranges. The exact
23-arrow value 83 first fits at order 12.

The separate construction audit finds a concrete adjacency rule,
prescribed composition, or certified listed graph for every value. It
needs no generic graph search or orientation derandomization to
instantiate an attainer. Fixed stages retain every-optional-subset
guarantees inside their own envelopes; other constructions retain their
separately proved delivery bounds. Those stronger properties are not
inferred for unrelated graphs of the same score.

Finite upper proofs combine ordinary degree and source/sink arguments
with complete residual graph enumerations. Their reviews distinguish
full independent game replay from a second full graph enumeration and
from sampled pruning tests. The exact dependencies are in the evidence
guide and approved registry. The new concrete classification is not
fully instantiated in Lean.

Delivery-time optimization, classification of all maximizing graphs, and
a structural characterization of all guaranteed pairs remain outside
this paper's scope. Further simplification and fuller formalization do
not represent unknown values or absent explicit attainers.

\subsection{Complete tables through eleven
vertices}\label{complete-tables-through-eleven-vertices}

These tables and interval rules determine \(M(n,m)\) for every
admissible pair with \(1\le n\le11\). An em dash means \(m>n(n-1)\), so
no such loopless graph exists. At every order, \(M(n,m)=0\) for
\(m>n(n-2)\); this rule covers any final budgets omitted below.

\subsubsection{Orders one through seven}\label{orders-one-through-seven}

The three blocks list budgets zero through thirty-five. Together with
the zero-tail rule, they cover every admissible budget at these orders.

{\def\LTcaptype{none} 
\Needspace{11\baselineskip}
\begin{longtable}[]{@{}lllllllllllll@{}}
\toprule\noalign{}
\(n\backslash m\) & 0 & 1 & 2 & 3 & 4 & 5 & 6 & 7 & 8 & 9 & 10 & 11 \\
\midrule\noalign{}
\endhead
\bottomrule\noalign{}
\endlastfoot
1 & 0 & --- & --- & --- & --- & --- & --- & --- & --- & --- & --- &
--- \\
2 & 0 & 0 & 0 & --- & --- & --- & --- & --- & --- & --- & --- & --- \\
3 & 0 & 0 & 0 & 0 & 0 & 0 & 0 & --- & --- & --- & --- & --- \\
4 & 0 & 0 & 0 & 0 & 1 & 1 & 2 & 1 & 1 & 0 & 0 & 0 \\
5 & 0 & 0 & 0 & 0 & 1 & 1 & 2 & 2 & 4 & 3 & 3 & 3 \\
6 & 0 & 0 & 0 & 0 & 1 & 1 & 2 & 2 & 4 & 5 & 6 & 7 \\
7 & 0 & 0 & 0 & 0 & 1 & 1 & 2 & 2 & 4 & 5 & 6 & 8 \\
\end{longtable}
}

{\def\LTcaptype{none} 
\Needspace{11\baselineskip}
\begin{longtable}[]{@{}
  >{\raggedright\arraybackslash}p{(\linewidth - 24\tabcolsep) * \real{0.0769}}
  >{\raggedright\arraybackslash}p{(\linewidth - 24\tabcolsep) * \real{0.0769}}
  >{\raggedright\arraybackslash}p{(\linewidth - 24\tabcolsep) * \real{0.0769}}
  >{\raggedright\arraybackslash}p{(\linewidth - 24\tabcolsep) * \real{0.0769}}
  >{\raggedright\arraybackslash}p{(\linewidth - 24\tabcolsep) * \real{0.0769}}
  >{\raggedright\arraybackslash}p{(\linewidth - 24\tabcolsep) * \real{0.0769}}
  >{\raggedright\arraybackslash}p{(\linewidth - 24\tabcolsep) * \real{0.0769}}
  >{\raggedright\arraybackslash}p{(\linewidth - 24\tabcolsep) * \real{0.0769}}
  >{\raggedright\arraybackslash}p{(\linewidth - 24\tabcolsep) * \real{0.0769}}
  >{\raggedright\arraybackslash}p{(\linewidth - 24\tabcolsep) * \real{0.0769}}
  >{\raggedright\arraybackslash}p{(\linewidth - 24\tabcolsep) * \real{0.0769}}
  >{\raggedright\arraybackslash}p{(\linewidth - 24\tabcolsep) * \real{0.0769}}
  >{\raggedright\arraybackslash}p{(\linewidth - 24\tabcolsep) * \real{0.0769}}@{}}
\toprule\noalign{}
\begin{minipage}[b]{\linewidth}\raggedright
\(n\backslash m\)
\end{minipage} & \begin{minipage}[b]{\linewidth}\raggedright
12
\end{minipage} & \begin{minipage}[b]{\linewidth}\raggedright
13
\end{minipage} & \begin{minipage}[b]{\linewidth}\raggedright
14
\end{minipage} & \begin{minipage}[b]{\linewidth}\raggedright
15
\end{minipage} & \begin{minipage}[b]{\linewidth}\raggedright
16
\end{minipage} & \begin{minipage}[b]{\linewidth}\raggedright
17
\end{minipage} & \begin{minipage}[b]{\linewidth}\raggedright
18
\end{minipage} & \begin{minipage}[b]{\linewidth}\raggedright
19
\end{minipage} & \begin{minipage}[b]{\linewidth}\raggedright
20
\end{minipage} & \begin{minipage}[b]{\linewidth}\raggedright
21
\end{minipage} & \begin{minipage}[b]{\linewidth}\raggedright
22
\end{minipage} & \begin{minipage}[b]{\linewidth}\raggedright
23
\end{minipage} \\
\midrule\noalign{}
\endhead
\bottomrule\noalign{}
\endlastfoot
1 & --- & --- & --- & --- & --- & --- & --- & --- & --- & --- & --- &
--- \\
2 & --- & --- & --- & --- & --- & --- & --- & --- & --- & --- & --- &
--- \\
3 & --- & --- & --- & --- & --- & --- & --- & --- & --- & --- & --- &
--- \\
4 & 0 & --- & --- & --- & --- & --- & --- & --- & --- & --- & --- &
--- \\
5 & 2 & 2 & 2 & 1 & 0 & 0 & 0 & 0 & 0 & --- & --- & --- \\
6 & 7 & 6 & 7 & 6 & 6 & 6 & 6 & 4 & 4 & 4 & 4 & 3 \\
7 & 9 & 10 & 10 & 11 & 12 & 12 & 14 & 15 & 18 & 21 & 14 & 12 \\
\end{longtable}
}

{\def\LTcaptype{none} 
\Needspace{11\baselineskip}
\begin{longtable}[]{@{}
  >{\raggedright\arraybackslash}p{(\linewidth - 24\tabcolsep) * \real{0.0769}}
  >{\raggedright\arraybackslash}p{(\linewidth - 24\tabcolsep) * \real{0.0769}}
  >{\raggedright\arraybackslash}p{(\linewidth - 24\tabcolsep) * \real{0.0769}}
  >{\raggedright\arraybackslash}p{(\linewidth - 24\tabcolsep) * \real{0.0769}}
  >{\raggedright\arraybackslash}p{(\linewidth - 24\tabcolsep) * \real{0.0769}}
  >{\raggedright\arraybackslash}p{(\linewidth - 24\tabcolsep) * \real{0.0769}}
  >{\raggedright\arraybackslash}p{(\linewidth - 24\tabcolsep) * \real{0.0769}}
  >{\raggedright\arraybackslash}p{(\linewidth - 24\tabcolsep) * \real{0.0769}}
  >{\raggedright\arraybackslash}p{(\linewidth - 24\tabcolsep) * \real{0.0769}}
  >{\raggedright\arraybackslash}p{(\linewidth - 24\tabcolsep) * \real{0.0769}}
  >{\raggedright\arraybackslash}p{(\linewidth - 24\tabcolsep) * \real{0.0769}}
  >{\raggedright\arraybackslash}p{(\linewidth - 24\tabcolsep) * \real{0.0769}}
  >{\raggedright\arraybackslash}p{(\linewidth - 24\tabcolsep) * \real{0.0769}}@{}}
\toprule\noalign{}
\begin{minipage}[b]{\linewidth}\raggedright
\(n\backslash m\)
\end{minipage} & \begin{minipage}[b]{\linewidth}\raggedright
24
\end{minipage} & \begin{minipage}[b]{\linewidth}\raggedright
25
\end{minipage} & \begin{minipage}[b]{\linewidth}\raggedright
26
\end{minipage} & \begin{minipage}[b]{\linewidth}\raggedright
27
\end{minipage} & \begin{minipage}[b]{\linewidth}\raggedright
28
\end{minipage} & \begin{minipage}[b]{\linewidth}\raggedright
29
\end{minipage} & \begin{minipage}[b]{\linewidth}\raggedright
30
\end{minipage} & \begin{minipage}[b]{\linewidth}\raggedright
31
\end{minipage} & \begin{minipage}[b]{\linewidth}\raggedright
32
\end{minipage} & \begin{minipage}[b]{\linewidth}\raggedright
33
\end{minipage} & \begin{minipage}[b]{\linewidth}\raggedright
34
\end{minipage} & \begin{minipage}[b]{\linewidth}\raggedright
35
\end{minipage} \\
\midrule\noalign{}
\endhead
\bottomrule\noalign{}
\endlastfoot
1 & --- & --- & --- & --- & --- & --- & --- & --- & --- & --- & --- &
--- \\
2 & --- & --- & --- & --- & --- & --- & --- & --- & --- & --- & --- &
--- \\
3 & --- & --- & --- & --- & --- & --- & --- & --- & --- & --- & --- &
--- \\
4 & --- & --- & --- & --- & --- & --- & --- & --- & --- & --- & --- &
--- \\
5 & --- & --- & --- & --- & --- & --- & --- & --- & --- & --- & --- &
--- \\
6 & 2 & 0 & 0 & 0 & 0 & 0 & 0 & --- & --- & --- & --- & --- \\
7 & 12 & 11 & 10 & 10 & 9 & 8 & 8 & 7 & 6 & 5 & 4 & 2 \\
\end{longtable}
}

\subsubsection{Orders eight through eleven: the low
budgets}\label{orders-eight-through-eleven-the-low-budgets}

The next two blocks list budgets zero through twenty-three. Budget
twenty-four appears with the interval data immediately afterward.

{\def\LTcaptype{none} 
\Needspace{8\baselineskip}
\begin{longtable}[]{@{}lllllllllllll@{}}
\toprule\noalign{}
\(n\backslash m\) & 0 & 1 & 2 & 3 & 4 & 5 & 6 & 7 & 8 & 9 & 10 & 11 \\
\midrule\noalign{}
\endhead
\bottomrule\noalign{}
\endlastfoot
8 & 0 & 0 & 0 & 0 & 1 & 1 & 2 & 2 & 4 & 5 & 6 & 8 \\
9 & 0 & 0 & 0 & 0 & 1 & 1 & 2 & 2 & 4 & 5 & 6 & 8 \\
10 & 0 & 0 & 0 & 0 & 1 & 1 & 2 & 2 & 4 & 5 & 6 & 8 \\
11 & 0 & 0 & 0 & 0 & 1 & 1 & 2 & 2 & 4 & 5 & 6 & 8 \\
\end{longtable}
}

{\def\LTcaptype{none} 
\Needspace{8\baselineskip}
\begin{longtable}[]{@{}
  >{\raggedright\arraybackslash}p{(\linewidth - 24\tabcolsep) * \real{0.0769}}
  >{\raggedright\arraybackslash}p{(\linewidth - 24\tabcolsep) * \real{0.0769}}
  >{\raggedright\arraybackslash}p{(\linewidth - 24\tabcolsep) * \real{0.0769}}
  >{\raggedright\arraybackslash}p{(\linewidth - 24\tabcolsep) * \real{0.0769}}
  >{\raggedright\arraybackslash}p{(\linewidth - 24\tabcolsep) * \real{0.0769}}
  >{\raggedright\arraybackslash}p{(\linewidth - 24\tabcolsep) * \real{0.0769}}
  >{\raggedright\arraybackslash}p{(\linewidth - 24\tabcolsep) * \real{0.0769}}
  >{\raggedright\arraybackslash}p{(\linewidth - 24\tabcolsep) * \real{0.0769}}
  >{\raggedright\arraybackslash}p{(\linewidth - 24\tabcolsep) * \real{0.0769}}
  >{\raggedright\arraybackslash}p{(\linewidth - 24\tabcolsep) * \real{0.0769}}
  >{\raggedright\arraybackslash}p{(\linewidth - 24\tabcolsep) * \real{0.0769}}
  >{\raggedright\arraybackslash}p{(\linewidth - 24\tabcolsep) * \real{0.0769}}
  >{\raggedright\arraybackslash}p{(\linewidth - 24\tabcolsep) * \real{0.0769}}@{}}
\toprule\noalign{}
\begin{minipage}[b]{\linewidth}\raggedright
\(n\backslash m\)
\end{minipage} & \begin{minipage}[b]{\linewidth}\raggedright
12
\end{minipage} & \begin{minipage}[b]{\linewidth}\raggedright
13
\end{minipage} & \begin{minipage}[b]{\linewidth}\raggedright
14
\end{minipage} & \begin{minipage}[b]{\linewidth}\raggedright
15
\end{minipage} & \begin{minipage}[b]{\linewidth}\raggedright
16
\end{minipage} & \begin{minipage}[b]{\linewidth}\raggedright
17
\end{minipage} & \begin{minipage}[b]{\linewidth}\raggedright
18
\end{minipage} & \begin{minipage}[b]{\linewidth}\raggedright
19
\end{minipage} & \begin{minipage}[b]{\linewidth}\raggedright
20
\end{minipage} & \begin{minipage}[b]{\linewidth}\raggedright
21
\end{minipage} & \begin{minipage}[b]{\linewidth}\raggedright
22
\end{minipage} & \begin{minipage}[b]{\linewidth}\raggedright
23
\end{minipage} \\
\midrule\noalign{}
\endhead
\bottomrule\noalign{}
\endlastfoot
8 & 9 & 11 & 13 & 15 & 19 & 19 & 21 & 22 & 24 & 25 & 28 & 29 \\
9 & 9 & 11 & 13 & 15 & 19 & 24 & 33 & 35 & 37 & 39 & 40 & 44 \\
10 & 9 & 11 & 13 & 15 & 19 & 28 & 33 & 38 & 53 & 56 & 60 & 61 \\
11 & 9 & 11 & 13 & 15 & 19 & 28 & 33 & 40 & 53 & 58 & 79 & 78 \\
\end{longtable}
}

{\def\LTcaptype{none} 
\Needspace{8\baselineskip}
\begin{longtable}[]{@{}llll@{}}
\toprule\noalign{}
\(n\) & \(M(n,24)\) & Interval of budgets \(m\) & Value on that
interval \\
\midrule\noalign{}
\endhead
\bottomrule\noalign{}
\endlastfoot
8 & 32 & \(24\le m\le34\) & \(56-m\) \\
9 & 48 & \(24\le m\le54\) & \(72-m\) \\
10 & 61 & \(25\le m\le70\) & \(90-m\) \\
11 & 79 & \(25\le m\le88\) & \(110-m\) \\
\end{longtable}
}

\subsubsection{The remaining eight-vertex middle
budgets}\label{the-remaining-eight-vertex-middle-budgets}

{\def\LTcaptype{none} 
\Needspace{5\baselineskip}
\begin{longtable}[]{@{}lllllll@{}}
\toprule\noalign{}
\(n\backslash m\) & 35 & 36 & 37 & 38 & 39 & 40 \\
\midrule\noalign{}
\endhead
\bottomrule\noalign{}
\endlastfoot
8 & 19 & 18 & 16 & 15 & 14 & 12 \\
\end{longtable}
}

\subsubsection{The small dense tails}\label{the-small-dense-tails}

Put \(m=n(n-3)+r\). The entries below give \(M(n,m)\); the preceding
intervals already include \(r=0\) at orders nine and ten, and the
eight-vertex endpoint \(m=40\) appears above.

{\def\LTcaptype{none} 
\Needspace{7\baselineskip}
\begin{longtable}[]{@{}lllllllllll@{}}
\toprule\noalign{}
\(n\backslash r\) & 1 & 2 & 3 & 4 & 5 & 6 & 7 & 8 & 9 & 10 \\
\midrule\noalign{}
\endhead
\bottomrule\noalign{}
\endlastfoot
8 & 11 & 10 & 9 & 8 & 7 & 6 & 4 & 3 & 0 & 0 \\
9 & 14 & 12 & 11 & 10 & 9 & 8 & 6 & 5 & 3 & 0 \\
10 & 16 & 15 & 14 & 12 & 11 & 10 & 8 & 7 & 5 & 4 \\
\end{longtable}
}

At order eleven the uniform staircase applies throughout the remaining
nonzero tail. For \(m=88+r\) and \(1\le r\le11\),

\[M(11,88+r)=21-r-\lceil r/2\rceil.\]

For \(100\le m\le110\), the value is zero.

\subsubsection{The uniform classification from twelve
vertices}\label{the-uniform-classification-from-twelve-vertices}

Every entry of the small-arrow table below has an explicit attainer on
at most twelve vertices. Consequently it equals \(M(n,m)\) for every
\(n\ge12\) at that budget, as well as the unrestricted-host maximum
\(M_E(m)\).

{\def\LTcaptype{none} 
\Needspace{5\baselineskip}
\begin{longtable}[]{@{}lllllllllllll@{}}
\toprule\noalign{}
\(m\) & 0 & 1 & 2 & 3 & 4 & 5 & 6 & 7 & 8 & 9 & 10 & 11 \\
\midrule\noalign{}
\endhead
\bottomrule\noalign{}
\endlastfoot
\(M_E(m)\) & 0 & 0 & 0 & 0 & 1 & 1 & 2 & 2 & 4 & 5 & 6 & 8 \\
\end{longtable}
}

{\def\LTcaptype{none} 
\Needspace{5\baselineskip}
\begin{longtable}[]{@{}lllllllllllll@{}}
\toprule\noalign{}
\(m\) & 12 & 13 & 14 & 15 & 16 & 17 & 18 & 19 & 20 & 21 & 22 & 23 \\
\midrule\noalign{}
\endhead
\bottomrule\noalign{}
\endlastfoot
\(M_E(m)\) & 9 & 11 & 13 & 15 & 19 & 28 & 33 & 40 & 53 & 60 & 79 & 83 \\
\end{longtable}
}

For \(n\ge12\) and \(m\ge24\), put \(L=n(n-1)\),
\(k=\lfloor m/2\rfloor\) and, in the dense band, \(r=m-n(n-3)\). Then

\[
M(n,m)=\begin{cases}
k(k-3)+2(m\bmod2),&24\le m\le2n,\\
L-m,&2n<m\le n(n-3),\\
2n-r-\lceil r/2\rceil-1,&n(n-3)<m\le n(n-2),\\
0,&n(n-2)<m\le L.
\end{cases}
\]

The sparse formulas have fitting even- and odd-budget cores: \(k\le n\)
in the even case and \(k+1\le n\) in the odd case. The complete
main-interval theorem and the staircase theorem from order eleven cover
the next two bands. The small-arrow attainers extend by isolated
vertices. Conversely, removing isolated vertices leaves at most
\(2m\le46\) vertices when \(1\le m\le23\), so the finite ledger and
isolate invariance also establish the unrestricted-host upper bounds in
the small-arrow table. At \(m=0\) every graph has value zero.

\subsubsection{Provenance and
attainment}\label{provenance-and-attainment}

The values are generated from the complete classification ledger,
\texttt{publication/\allowbreak{}classification-\allowbreak{}coverage.\allowbreak{}json}, by
\texttt{src/\allowbreak{}finite\_\allowbreak{}value\_\allowbreak{}tables.\allowbreak{}py}, using an independent
affine-segment decoder. The generator refuses unresolved cells, checks
that the displayed tables and intervals cover all 451 admissible cells
through order eleven, and checks the uniform formulas against every
stored order from twelve onward. It also requires a fresh
explicit-construction audit,
\texttt{research/\allowbreak{}explicit-\allowbreak{}construction-\allowbreak{}coverage.\allowbreak{}json}, with zero
known-value construction gaps. These paths identify files in the
downloadable research-source archive.

The infinite statement uses the ordinary sparse constructions in
\texttt{paper/\allowbreak{}sparse-\allowbreak{}all-\allowbreak{}orders.\allowbreak{}md}, the complete main interval in
\texttt{paper/\allowbreak{}main-\allowbreak{}interval-\allowbreak{}complete.\allowbreak{}md}, the dense staircase in
\texttt{paper/\allowbreak{}dense-\allowbreak{}third.\allowbreak{}md}, and the zero regime. The finite
exclusions retain their individual source reviews and complete literal
winning/losing rank checks in the small-case registry,
\texttt{research/\allowbreak{}small-\allowbreak{}case-\allowbreak{}results.\allowbreak{}json}. This table readback adds no
graph search and makes no new Lean claim.

\section{Game and structural bounds}\label{game-and-structural-bounds}

\subsection{The game and its elementary
structure}\label{the-game-and-its-elementary-structure}

Let \(G=(V,E)\) be a finite loopless digraph. Opposite arcs are allowed.
A legal move from \(x\) goes to \(x\) or to an outgoing neighbor; write
\[C(x)=\{x\}\cup N^+(x).\] Fix \(s,t\in V\) and \(c\ne s\). If \(s=t\),
receipt has already occurred. Otherwise the messenger moves first. On
its turn it chooses \(v\in C(s)\). If \(v=t\), delivery wins
immediately. If \(v\ne t\) and \(v=c\), it is captured. Otherwise the
pursuer chooses \(c'\in C(c)\); equality \(c'=v\) is capture. The
process repeats. Infinite play without delivery is a messenger loss.

Let \(W_G(s,t)\) mean that, for every initial \(c\ne s\), there exists a
messenger strategy guaranteeing delivery against every pursuer strategy.
The messenger observes \(c\); its strategy may depend on it. Then \[
I(G)=\#\{(s,t):s\ne t,\ (s,t)\notin E,\ W_G(s,t)\}.
\] If \(w(G)\) counts all guaranteed pairs including the diagonal, then
exactly \[w(G)=n+m+I(G).\]

The terminal rule matters. In the directed diamond \(s\to a\to t\),
\(s\to b\to t\), the messenger chooses the intermediate vertex not
occupied by the pursuer and delivers on its next turn. But reaching a
vertex that would be safe as a \emph{recipient} does not justify using
it as an intermediate stop. The pursuer gets an intervening move before
the message can be forwarded again.

\begin{figure}
\centering
\pandocbounded{\includegraphics[keepaspectratio,alt={Two routes allow a guaranteed indirect delivery. Adding the arc from a to b lets a pursuer at a block both first moves.}]{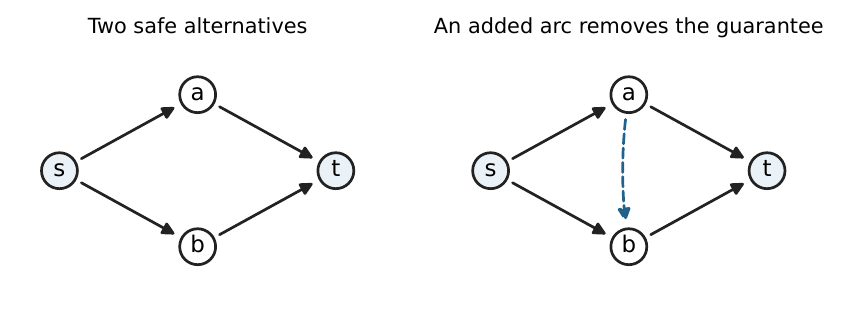}}
\caption{Two routes allow a guaranteed indirect delivery. Adding the arc
from a to b lets a pursuer at a block both first moves.}
\end{figure}

\subsubsection{Camping obstructions}\label{camping-obstructions}

\textbf{Lemma 1 (source domination).} If \(c\ne s\) and
\(N^+(s)\subseteq C(c)\), there is no guaranteed indirect pair with
source \(s\).

\textbf{Proof.} The pursuer starts at \(c\) and waits while the
messenger waits. Every nonterminal departure from \(s\) is either an
immediate collision or is captured on the following pursuer turn.
Waiting forever is not delivery. \(\square\)

\textbf{Lemma 2 (predecessor domination).} Put \(D_t=N^-(t)\). If
\(D_t\subseteq C(c)\), an indirect guaranteed source for \(t\) must
equal \(c\).

\textbf{Proof.} For any other source the pursuer may start at \(c\).
Every indirect delivery must first enter \(D_t\), where the pursuer can
capture before the next messenger turn. \(\square\)

In particular, a source counted by \(I\) has outdegree at least two, and
a recipient counted by \(I\) has indegree at least two. An indegree-one
recipient is blocked by its unique predecessor. In an undirected
network, interpreted as two opposite arcs per edge, the pursuer can
start at \(t\) and catch the messenger upon its first arrival at a
neighbor of \(t\). Thus \(I=0\) for every such network.

The two degree conditions are necessary, not sufficient. Neither the
absence of a one-vertex directed separator nor the existence of two
internally disjoint routes characterizes this game. Adding \(a\to b\) to
the diamond destroys its only indirect guarantee, although both routes
remain. Consequently \(I\) is not monotone under adding arcs.

\textbf{Lemma 3 (two-move criterion).} For an indirect pair \((s,t)\),
put \(B=N^+(s)\cap N^-(t)\). Delivery is guaranteed within two messenger
moves exactly when
\[\text{for every }c\ne s,\qquad B\not\subseteq C(c).\]

\textbf{Proof.} Choose an intermediate vertex of \(B\setminus C(c)\),
survive every pursuer response, and deliver. Conversely, any successful
two-move indirect delivery must use such an intermediate. Waiting on the
first move cannot reach a nonadjacent recipient on the second.
\(\square\)

\subsubsection{Components, zero minima, and the vertex
bound}\label{components-zero-minima-and-the-vertex-bound}

Removing isolated vertices preserves \(I\). More generally,
\[I(G\sqcup H)=I(G)+I(H).\] There are no cross-component routes; a
pursuer in the other component cannot interfere with a route in the
messenger's component.

\textbf{Proposition 4.} At every admissible \((n,m)\) the minimum of
\(I\) is zero.

\textbf{Proof.} For \(n\ge2\), list all nonloop arcs row by row and
retain the first \(m\). There are full outgoing rows, at most one
partial row, and empty rows. Full rows have only direct recipients;
empty rows cannot move. If a full row exists, its vertex blocks every
indirect attempt from the partial row by Lemma 1. Otherwise successors
of the partial row have no outgoing arcs, so no indirect route exists.
The empty and one-vertex cases are immediate. \(\square\)

\textbf{Theorem 5 (original upper bound).} For every \(n\ge3\),
\[I(G)\le n(n-3).\]

\textbf{Proof.} A recipient of indegree zero or one has no indirect
guaranteed sources. A recipient of indegree \(d\ge2\) has at most
\(n-1-d\le n-3\), after excluding itself and its direct predecessors.
Sum over recipients. \(\square\)

For \(n\ge4\), equality forces every recipient to have indegree two and
every missing pair to be guaranteed. Every source must then have
outdegree at least two; the degree sum forces outdegree exactly two.
Conversely, a two-in/two-out-regular graph with every pair guaranteed
attains equality. Thus attaining the vertex bound is exactly a sparse
universal-delivery construction problem.

\subsection{Sharp arc-budget upper
bounds}\label{sharp-arc-budget-upper-bounds}

The source/sink and positive-support lemma below bounds an actual
source/sink graph by preceding arc budgets and a positive-degree graph
by its eligible rectangle. Its ordinary induction proves the following
bounds and equality shapes before deriving the stronger support results.
The induction starts from the elementary results through six arrows; it
uses no finite extremum as a premise.

\textbf{Theorem 6.} If \(m=2k\) and \(k\ge5\), then \[I(G)\le k(k-3).\]

\textbf{Theorem 7.} If \(m=2k+1\) and \(k\ge5\), then
\[I(G)\le k(k-3)+2.\]

\textbf{Proof of Theorems 6 and 7.} Apply the ordinary sparse parity
proposition, whose proof is the increasing-budget induction
(56.2)--(56.6). \(\square\)

\textbf{Theorem 8.} For \(k\ge6\), equality in Theorem 6 holds exactly
when deleting isolated vertices leaves a two-in/two-out-regular graph on
\(k\) vertices with every pair guaranteed.

\textbf{Theorem 9.} For \(k\ge6\), equality in Theorem 7 holds exactly
when the nonisolated support has \(k+1\) vertices; one vertex \(u\) has
outdegree one and indegree two; a different vertex \(v\) has indegree
one and outdegree two; every other degree is two; \(u\to v\) is absent;
and every missing pair from an outdegree-two source to an indegree-two
recipient is guaranteed.

\textbf{Proof of Theorems 8 and 9.} The same parity proposition proves
the equality shapes. In even support the number of missing pairs is
\(k(k-3)\). In the odd absent-arrow family the eligible count is
\(k(k-3)+2\); every other shape has a strictly smaller upper bound. The
stated guarantees therefore characterize equality. \(\square\)

Adding the absent \(u\to v\) in Theorem 9 gives a two-regular digraph.
It does not follow that the completed graph is universally winning.

\subsection{A common upper-bound argument for sparse
graphs}\label{a-common-upper-bound-argument-for-sparse-graphs}

The same source/sink deletion and positive-degree count prove the sharp
arc bounds, their equality shapes, and the quantitative support
reductions. All graphs are simple and loopless; opposite arrows are
allowed. An induced subgraph retains every arrow between its vertices.

Put \(A=\{d^+\ge2\}\), \(B=\{d^-\ge2\}\), \(q=|A|\), \(r=|B|\),
\(b=|A\cap B|\), and \(e=e(A,B)\). Camping excludes every counted pair
outside \(A\times B\). Direct and diagonal pairs are disjoint, so
\[I(G)\le qr-b-e,\qquad e\ge\max(0,2q+2r-m). \tag{56.1}\]

\subsubsection{Deletion and positive
support}\label{deletion-and-positive-support}

\textbf{Lemma (source/sink and positive-support bounds).} Write
\(A_m=M_E(m)\). An actual source of outgoing degree \(d>0\), or an
actual sink of incoming degree \(d>0\), satisfies \[
I(G)\le
\begin{cases}
A_{m-1},&d=1,\\
A_{m-d}+\lfloor(m-d)/2\rfloor,&d\ge2.
\end{cases} \tag{56.2}
\] On a fixed support of order \(N\), the added term may also be
replaced by \(N-1-d\). Consequently every graph with a nonisolated
source or sink obeys
\[I(G)\le D_m:=\max\{A_{m-1},A_{m-2}+\lfloor(m-2)/2\rfloor\}. \tag{56.3}\]
If every indegree and outdegree is positive, then \[
q,r\le\min(N,\lfloor m/2\rfloor,m-N),\qquad
I(G)\le B_m(N,q,r):=qr-(q+r-N)_+-(m-2N+q+r)_+.
\tag{56.4}
\] Here \(x_+=\max(0,x)\).

\textbf{Proof.} An old-source messenger cannot enter a deleted source. A
messenger delivering to an old recipient cannot enter a deleted sink and
subsequently leave it. Restrict the cop to the remaining graph: every
counted pair with both endpoints surviving therefore remains a guarantee
there. Only the deleted source row or sink column can add pairs.

Every live recipient column needs two predecessors of positive indegree.
An indegree-zero predecessor cannot be reached after the initial
position, and starting there would make the pair direct. If there is
only one usable predecessor, a cop may camp there. Separately, every
live source row needs two nonsink successors: a cop camps at a sole
nonsink successor, while every alternative departure traps the messenger
at a sink. Thus the deleted source's arrows supply none of the first
capacity, and the deleted sink's arrows supply none of the second. At
most \(\lfloor(m-d)/2\rfloor\) columns or rows can contribute. For
\(d=1\) the deleted row or column is empty; \(N-1-d\) simply counts its
nonadjacent distinct pairs.

Adding a fresh sink with one incoming arrow preserves \(I\). An old
winning strategy works while the cop stays in the old graph; if the cop
enters or starts at the sink, it is trapped and any old route to the
recipient suffices. Old losing pairs retain a cop strategy confined to
the old graph, and the new sink contributes no indirect row or column.
Hence \(A_m\) is nondecreasing, which gives (56.3). This is a budget
comparison, not arbitrary arrow deletion.

For positive degrees, the vertices outside each active set have the
corresponding degree one. This gives \(q,r\le m-N\) and
\[e(A,B)=m-2N+q+r+w,\] where \(w\ge0\) counts arrows from \(A\)'s
complement to \(B\)'s complement. Subtract this lower bound and
\(b\ge(q+r-N)_+\) from the eligible rectangle. The two deletion
arguments used the original arrow directions throughout; no
graph-reversal assumption is needed. \(\square\)

\textbf{Arithmetic lemma.} In positive support the following bounds
hold: \[
\begin{array}{c|c|c}
m&\text{support}&\text{bound on }I\\ \hline
2k,\ k\ge5&N=k,\ (q,r)\ne(k,k)&k^2-4k+2\\
2k,\ k\ge5&N\ge k+1&k^2-5k+8\\
2k+1,\ k\ge5&N=k+1,\ (q,r)\ne(k,k)&k^2-4k+4\\
2k+1,\ k\ge4&N\ge k+2&k^2-5k+10
\end{array} \tag{56.5}
\] For the odd support reduction, the first odd bound remains valid at
\(k=3,4\), and at \(N=k\) a coordinate below \(k\) gives \(k^2-4k+1\).
For \(k=3\), positive support \(N\ge5\) at seven arrows has \(I\le3\).

\textbf{Proof.} In (56.4), increasing either coordinate when the other
is at least two cannot decrease the expression: its product increases by
at least two, while the two positive-part subtractions increase by at
most two. At the fixed supports, substitute the largest permitted
coordinates. Coordinates zero or one give \(I\le\min(q,r)\max(q,r)\) and
are bounded separately. For \(k=3,4\), infeasible small coordinates can
be discarded using \(m\le N+q(N-2)\) and its incoming counterpart. At
\(k=3\), \(N=3\) is impossible and the nonmaximal-coordinate bound at
\(N=4\) is one. At \(k=4\), the corresponding bounds at \(N=4,5\) are
one and four. These are the only small fixed-support cases outside the
uniform substitution.

For the large-support cases put \(p=m-N\). Then \(q,r\le p\) and
\[B_m(N,p,p)=f_m(p):=p^2-(3p-m)_+-(4p-m)_+.\] If both coordinates are at
least two, their maximum is bounded by \(f_m(p)\); otherwise \(I\le p\).
Moreover \(f_m(p+1)-f_m(p)\ge2p-6\), so \(f_m\) is nondecreasing for
\(p\ge3\). The values at \(p=0,1,2,3\) are immediate. Since \(p\le k-1\)
in both large-support cases, substitution at \(p=k-1\) gives (56.5); its
value also dominates \(p\) in the stated domains. At \(k=3,m=7\), the
possibilities \(p\le2\) give the separate bound three. \(\square\)

\subsubsection{Ordinary induction for the parity
bounds}\label{ordinary-induction-for-the-parity-bounds}

\textbf{Proposition (sparse parity bounds and equality).} For every
\(k\ge5\), \[A_{2k}\le k(k-3),\qquad A_{2k+1}\le k(k-3)+2. \tag{56.6}\]
Even equality forces a two-in/two-out-regular nonisolated support of
order \(k\), with every missing pair guaranteed. Odd equality forces
support \(k+1\), with distinct vertices \(u,v\) of respective degree
pairs \((d^+,d^-)=(1,2),(2,1)\), every other degree two, \(u\to v\)
absent, and every eligible missing pair guaranteed. These conditions
also suffice.

\textbf{Proof.} The elementary ordinary results through six arrows give
\(A_5=1,A_6=2\). Applying (56.3)--(56.4) successively gives \[
\begin{array}{c|rrrrr}
m&7&8&9&10&11\\\hline
D_m\text{ using preceding bounds}&3&5&6&9&10\\
\max B_m\text{ on positive support}&3&4&6&10&12\\
A_m\text{ upper bound}&3&5&6&10&12 .
\end{array}
\] These small integer maxima require only \(2\le N\le m\),
\(m\le N(N-1)\), and (56.4); they enumerate no graphs. Thus the seeds at
ten and eleven arrows are ordinary and do not use their sharper finite
values.

Induct in the arrow budget. At \(m=2k\), (56.3) gives
\(D_{2k}\le(k-1)(k-3)<k(k-3)\). Positive support \(N<k\) has fewer than
\(k(k-3)\) missing pairs; \(N\ge k+1\) is bounded by (56.5). At \(N=k\),
the missing-pair bound is \(k(k-3)\), and a coordinate below \(k\) gives
the smaller first bound in (56.5). Hence equality requires \(q=r=k\),
which forces every degree two, and all its missing pairs to win.

At \(m=2k+1\), the deletion bound is at most
\[\max\{k(k-3),k^2-4k+5\}=k(k-3).\] Positive \(N\le k\) has at most
\(k(k-3)-1\) missing pairs, and \(N\ge k+2\) is bounded by (56.5). At
\(N=k+1\), a coordinate below \(k\) again gives a strict bound. For
\(q=r=k\), the degree sums force exactly one degree-one vertex on each
side. If these coincide, \(b=k\); otherwise \(b=k-1\). Since
\(e=2k-1+w\), the eligible count is \(k(k-3)+1-w\) in the first case and
\(k(k-3)+2-w\) in the second. Here \(w\) is precisely the possible arrow
from the low-outdegree vertex to the low-indegree vertex. This proves
the bound and equality statement. The seed cases \(k=5\) have the same
strict support and deletion comparisons. \(\square\)

\subsubsection{A gap before another even support shape can
compete}\label{a-gap-before-another-even-support-shape-can-compete}

\textbf{Theorem 56 (even sparse stability).} Let \(k\ge6\) and let \(G\)
have exactly \(2k\) arcs. At least one of the following holds:

\begin{enumerate}
\def\labelenumi{\arabic{enumi}.}
\tightlist
\item
  Deleting isolated vertices leaves a graph on exactly \(k\) vertices,
  with indegree and outdegree two at every vertex.
\item
  \(I(G)\le(k-1)(k-3)\).
\end{enumerate}

\textbf{Proof.} Put \(C=(k-1)(k-3)\). The parity induction and (56.3)
bound every source/sink graph by \(C\). In positive support, \(N<k\) has
at most \((k-1)(k-2)-2k<C\) missing pairs. At \(N=k\), every nonregular
shape has \(I\le k^2-4k+2<C\); at \(N\ge k+1\), (56.5) gives
\(k^2-5k+8\le C\). Only the stated regular support remains. \(\square\)

Thus \(I>k(k-3)-(k-3)\) forces that shape; it does not itself assert
universal delivery. The same dichotomy holds at \(k=5\) using the finite
value \(A_{10}=6\le8\). It fails at \(k=4\): the six-vertex, eight-arrow
witness has \(I=4>3\). The theorem for \(k\ge6\) is entirely ordinary.

\textbf{Corollary 56.1.} For \(k\ge6\),
\[A_{2k}\le\max\{M(k,2k),(k-1)(k-3)\}.\] For example the retained finite
core bounds give \(A_{12}\le15\), \(A_{14}\le24\), \(A_{16}\le35\), and
\(A_{22}\le81\). These weaker consequences are superseded numerically by
the finite-budget classification, but retain their separate support
interpretation.

\subsubsection{Three odd exception families and an induced
core}\label{three-odd-exception-families-and-an-induced-core}

\textbf{Theorem 57 (odd sparse support reduction).} Let \(k\ge3\), let
\(G\) have \(2k+1\) arcs, and suppose \(I(G)>k^2-4k+6\). Then at least
one holds:

\begin{enumerate}
\def\labelenumi{\arabic{enumi}.}
\tightlist
\item
  \(G\) has an induced \(k\)-vertex subgraph \(H\) with
  \(I(H)\ge I(G)\). Either \(H\) is two-in/two-out regular with \(2k\)
  arcs, or it has \(2k+1\) arcs, one outgoing degree three and one
  incoming degree three, and all other corresponding degrees two. The
  exceptional vertices may coincide.
\item
  Removing isolates leaves \(k+1\) vertices, with distinct vertices
  \(u,v\) of degree pairs \((1,2),(2,1)\) and all other degrees two. The
  two possibilities \(u\to v\) absent or present are separate families.
\item
  Removing isolates leaves \(k+1\) vertices, with one vertex of indegree
  and outdegree one and all other degrees two.
\end{enumerate}

\textbf{Proof.} Put \(T=k^2-4k+6\). A source/sink of other degree at
least two gives \(I\le A_{2k-1}+k-1\le k^2-4k+5<T\) for \(k\ge5\). At
\(k=4,3\) the ordinary seeds give respectively \(A_7+3\le6=T\) and
\(A_5+2=3=T\). Thus only other degree one is possible. Delete such a
vertex. Its row or column contributes zero, so the remaining graph has
\(2k\) arcs and score at least \(I(G)>T\). For \(k\ge6\), Theorem 56
applies because \(T=(k-1)(k-3)+3\): deleting its isolates leaves a
regular \(k\)-vertex induced core with at least the original score. It
remains induced in \(G\), proving case 1. For \(k=3,4,5\), degree-one
deletion is impossible since \(A_6\le2\), \(A_8\le5\), \(A_{10}\le10\),
respectively, are below \(T\).

Otherwise all degrees are positive. The missing-pair count excludes
\(N<k\), and (56.5) excludes \(N\ge k+2\) (including its seven-arrow
case). At \(N=k\), a coordinate below \(k\) gives \(k^2-4k+1<T\), so
\(q=r=k\) and the one-excess degree pattern in case 1 follows. At
\(N=k+1\), a coordinate below \(k\) gives \(k^2-4k+4<T\). Thus
\(q=r=k\), and the degree-one vertices either differ or coincide, giving
cases 2 and 3. These exhaust the graph itself, without reversing arrows
or contracting routing vertices. \(\square\)

For \(k=0,1,2\), the elementary values \(A_1=A_3=0,A_5=1\) make the
conditional statement vacuous.

\textbf{Corollary 57.1 (odd bound and equality shape).} For \(k\ge5\),
\(A_{2k+1}\le k(k-3)+2\), with equality exactly in the absent-arrow
distinct-exception family when every eligible missing pair is
guaranteed.

\textbf{Proof.} This is the already proved parity proposition.
Equivalently, the families in Theorem 57 have respectively at most
\(k(k-3)\), \(k(k-3)-1\), \(k(k-3)+2\), \(k(k-3)+1\), and \(k(k-3)+1\)
eligible missing pairs. Only the absent-arrow family can attain the
bound. \(\square\)

The arithmetic checker and its separate audit record verify the finite
substitutions and inequality instances used here. They supplement these
ordinary proofs; they neither enumerate graphs nor supply a new Lean
proof.

\subsection{Quantitative concentration near the sparse
optimum}\label{quantitative-concentration-near-the-sparse-optimum}

The exact support alternatives in Theorems 56 and 57 apply above
specific cutoffs. The following estimate continues below those cutoffs.
It says that losing only a linear number of pairs from the leading
quadratic bound forces only a bounded number of exceptional vertices and
arrows. It makes no claim that deleting those exceptions preserves a
delivery strategy.

\textbf{Theorem 60.} Let \(G\) have \(m\) arrows and \(I=I(G)\). Put
\[A=\{v:d^+(v)\ge2\},\qquad B=\{v:d^-(v)\ge2\},\]
\[q=|A|,\quad r=|B|,\quad b=|A\cap B|,\quad a=m-2q,\quad d=m-2r.\] If
\(R\) and \(C\) count the source rows and recipient columns containing
at least one guaranteed indirect pair, then
\[R\le\min\{q,b+2a+\lfloor d/2\rfloor\},\qquad
C\le\min\{r,b+\lfloor a/2\rfloor+2d\},\qquad I\le RC.\] Define
\(\eta=m-2\sqrt I\ge0\). The induced core \(K=G[A\cap B]\) satisfies
\[m/2-3\eta\le |V(K)|\le m/2.\] The graph has at most \(m/2+4\eta\)
nonisolated vertices, and at most \(12\eta\) arrows have an endpoint
outside \(K\).

\textbf{Proof.} Only vertices in \(A\) can supply nonzero indirect
source rows, and only vertices in \(B\) can supply nonzero recipient
columns. These are the usual camping obstructions at degree zero or one.
In particular \(I\le qr\le m^2/4\).

At most \(d\) arrows end outside \(B\), so at least \(z=\max(0,q-b-d)\)
vertices of \(A\setminus B\) have indegree zero. They emit at least
\(2z\) arrows. Let \(x\) be the total excess indegree above two inside
\(B\). At most \(d-x\) arrows end outside \(B\), and vertices of
indegree at least three have total indegree at most \(3x\), since
\(e\le3(e-2)\) for \(e\ge3\). Thus at most \(3d\) of the \(2z\) arrows
can avoid indegree-two recipients. At least \(z-\lfloor3d/2\rfloor\)
distinct indegree-two recipients have an indegree-zero predecessor.

Each such recipient has a zero indirect column. An indirect source
cannot reach the indegree-zero predecessor and must approach through the
other predecessor, where the pursuer can wait. Therefore
\[C\le r-z+\lfloor3d/2\rfloor
\le r-q+b+d+\lfloor3d/2\rfloor
=b+\lfloor a/2\rfloor+2d.\] The identity uses \(r-q=(a-d)/2\) and the
common parity of \(a,d\).

For the row bound, use the original arrow directions. At least \(r-b-a\)
members of \(B\setminus A\) are sinks. Their incoming arrows force
outdegree-two predecessors except for at most \(3a\) arrows, by the same
surplus count. An outdegree-two source with a sink successor has zero
indirect row: the pursuer camps at the other successor. Repeating the
count gives \(R\le b+2a+\lfloor d/2\rfloor\). This is not an appeal to
symmetry under reversing the pursuit game. The displayed minimum bounds
and \(I\le RC\) follow.

Put \(S=a+d\). Since \(q+r=m-S/2\) and \(I\le qr\),
arithmetic--geometric mean gives \(S\le2m-4\sqrt I=2\eta\). The row and
column estimates also imply \[2\sqrt I\le R+C\le2b+\tfrac52S,\] so
\(b\ge\sqrt I-5S/4\ge m/2-3\eta\). Its upper bound is immediate.

Every nonisolated vertex outside \(A\cup B\) contributes at least one to
its incoming plus outgoing degree. Outgoing degrees outside \(A\) total
at most \(a\), and incoming degrees outside \(B\) total at most \(d\).
There are therefore at most \(S\) such vertices. The nonisolated support
is at most \[q+r-b+S=m-b+S/2
\le m-\sqrt I+\tfrac74S\le m/2+4\eta.\] Finally, the number of arrows
with tail outside \(A\cap B\) is at most \(2(q-b)+a\), and the analogous
incoming count is at most \(2(r-b)+d\). Counting an arrow twice can only
increase this estimate. The number of arrows outside the induced core is
thus at most \[2(q+r-2b)+S=2m-4b\le2m-4\sqrt I+5S\le12\eta.\] This
proves all assertions. \(\square\)

For example, if \(m=2k\) and \(I\ge k(k-3)-\Delta\) with
\(\Delta=O(k)\), then \(\eta=O(1)\). Thus the support has at most
\(k+O(1)\) vertices and all but \(O(1)\) arrows lie inside a core of at
most \(k\) vertices. The constants above are explicit but are not
asserted to be optimal. Theorems 56 and 57 retain their sharper
conclusions in their stated ranges. In particular, this concentration
estimate does not supply the game restriction \(I(K)\ge I(G)\) that is
proved under the more specific hypotheses of Theorem 57.

The argument was independently reviewed and checked against both game
solvers on all 4,165 labelled loopless digraphs through four vertices
and 181 additional graphs. The ordinary proof, finite tests, and exact
integer comparisons are recorded in
\texttt{research/\allowbreak{}routing-\allowbreak{}defect-\allowbreak{}phase2-\allowbreak{}support.\allowbreak{}md} and its two
verification receipts. This theorem is not presently formalized in Lean.

\subsubsection{A bounded host for every remaining small
budget}\label{a-bounded-host-for-every-remaining-small-budget}

The refined outside-vertex counts in the small sparse section make a
stronger finite reduction possible. This is a reduction of the whole
extremal problem, not a claim that an arbitrary induced core preserves a
strategy.

\textbf{Corollary 60.1.} For \(12\le m\le21\), define \(B_m\) by \[
(B_{12},B_{13},\ldots,B_{21})
=(12,13,16,18,15,18,15,18,11,15).
\] Then, for every admissible host order, \[M(n,m)=M(\min\{n,B_m\},m).\]
An explicit graph attaining the value on \(B_m\) vertices supplies an
explicit attainer on every larger host by adding isolated vertices.

\textbf{Proof.} The saved lower scores at these ten budgets are,
respectively, \[9,11,10,11,19,19,30,30,53,56.\] They are attained on
respectively \(7,8,7,7,8,8,9,11,10,11\) vertices. The 19-arrow witness
is the disjoint union of the nine-vertex, 18-arrow graph and a single
directed arrow. Every listed support is at most \(B_m\).

Apply the refined necessary integer-profile bounds to a graph scoring
strictly more than the listed lower score. Complete enumeration of the
integer parameters, independently implemented in two ways, bounds its
nonisolated support by \[12,13,16,18,15,18,15,18,13,15.\] This
enumerates degree counts, not graphs. Every actual graph has one of the
enumerated tuples; a surviving tuple need not be realizable.

Only the 20-arrow bound needs a further observation. Its sole possible
13-vertex tuple has a nine-vertex core in which every vertex has
indegree and outdegree two, together with two one-arrow sources and two
one-arrow sinks. The eligible source--recipient rectangle has at most 56
pairs. An arrow from an outside source into the core kills a whole
indirect recipient column: a core source can reach at most the other
predecessor, where the pursuer waits. An arrow from the core to an
outside sink kills the corresponding source row. Each such row or column
contains at least six eligible pairs, leaving at most 50 guarantees. If
neither interaction occurs, the core has 18 arrows and the outside
arrows form a separate zero-count component. Then \(I\le M_E(18)\le44\).
Neither case can beat 53, so every improvement has support at most
twelve at budget twenty. The twenty-arrow host-compression proposition
in the small sparse section then replaces its only surviving
twelve-vertex form by an eleven-vertex graph of the same score, using a
fresh sink. This gives \(B_{20}=11\); it does not claim that every
improving graph already has support eleven.

For \(n>B_m\), take any \(n\)-vertex graph. If its score is no greater
than the saved lower score, the saved witness already matches or
improves it on at most \(B_m\) vertices. Otherwise delete only its
isolated vertices, using the proved support bound, using the compression
at budget twenty when needed, and pad to \(B_m\) if necessary. This
proves \(M(n,m)\le M(B_m,m)\); adding isolates proves the reverse
inequality. For \(n\le B_m\) the assertion is an identity. \(\square\)

The exact arithmetic, literal lower witnesses and independent replays
are indexed in \texttt{research/\allowbreak{}small-\allowbreak{}case-\allowbreak{}results.\allowbreak{}json}. Together with
the complete joint coverage from order eleven, this confines every
remaining unknown numerical value to a host of order at most eighteen.
The separate recipe audit is still required: a numerical reduction alone
would not prove that an infinite known-value domain has explicit
constructions.

\section{Sparse attaining
constructions}\label{sparse-attaining-constructions}

\subsection{Sparse constructions with a linear
deficit}\label{sparse-constructions-with-a-linear-deficit}

\subsubsection{A safe branching
principle}\label{a-safe-branching-principle}

\textbf{Lemma 10.} Suppose distinct vertices have closed move sets
intersecting in at most one vertex. If a vertex \(r\) has two distinct
outgoing neighbors that each guarantee delivery to \(t\), then \(r\)
guarantees delivery to \(t\).

\textbf{Proof.} For every pursuer position \(c\ne r\), one of the two
successors lies outside \(C(c)\). Move there. It survives every
intervening pursuer response, after which its guaranteed strategy
applies. If the chosen vertex is \(t\), delivery is already terminal.
This is a finite concatenation of actual strategies. \(\square\)

\textbf{Theorem 11.} Let \(a\ge3\), \(h=2a+1\), \(n\ge2h+1\), and
\(\gcd(n,h)=1\). The digraph on \(\mathbb Z/n\mathbb Z\) with outgoing
steps \(1,a,h\) has every pair guaranteed. It has \(3n\) arcs and
\[I(G)=n(n-4).\]

\textbf{Proof.} The positive differences in \(\{0,1,a,2a+1\}\) are
\[1,\ a-1,\ a,\ a+1,\ 2a,\ 2a+1,\] all distinct. Since \(n>2h\), no
positive difference coincides modulo \(n\) with a negative one. Distinct
closed move sets therefore intersect in at most one vertex.

Translate the recipient to zero. Vertices \(0,-1,-a,-h\) are already
delivered or have direct arcs to it. Vertex \(-a-1\) is guaranteed via
its neighbors \(-a,-1\). Thus all four vertices
\[Q(z)=\{z,z-1,z-a,z-a-1\}\] are guaranteed for \(z=0\). If \(Q(z)\) is
guaranteed, apply Lemma 10 in the following order:

{\def\LTcaptype{none} 
\Needspace{8\baselineskip}
\begin{longtable}[]{@{}ll@{}}
\toprule\noalign{}
New source & Two guaranteed successors \\
\midrule\noalign{}
\endhead
\bottomrule\noalign{}
\endlastfoot
\(z-h\) & \(z-a-1,\ z\) \\
\(z-h-1\) & \(z-h,\ z-1\) \\
\(z-h-a\) & \(z-h,\ z-a\) \\
\(z-h-a-1\) & \(z-h-a,\ z-a-1\) \\
\end{longtable}
}

This proves \(Q(z-h)\). Iterating at most \(n-1\) times covers every
vertex because \(h\) is coprime to \(n\). The induction is finite; it
supplies terminating strategies, not merely perpetual evasion. Finally
subtract the \(n\) diagonal pairs and \(3n\) arcs from \(n^2\).
\(\square\)

\begin{figure}
\centering
\pandocbounded{\includegraphics[keepaspectratio,alt={A guaranteed quartet produces the next quartet seven positions earlier in the step-1,3,7 construction.}]{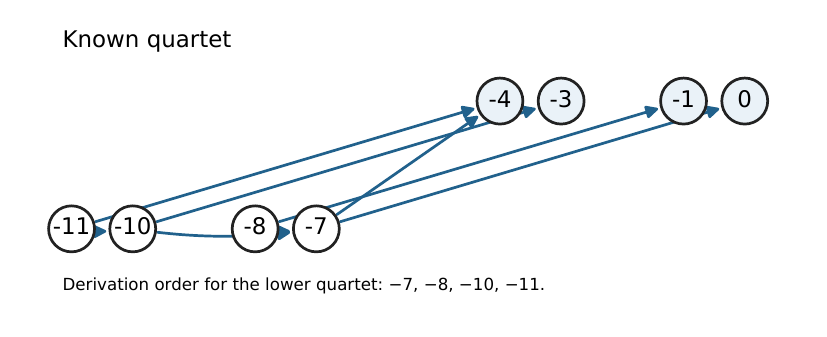}}
\caption{A guaranteed quartet produces the next quartet seven positions
earlier in the step-1,3,7 construction.}
\end{figure}

The choice \(a=3\) proves the theorem for every \(n\ge15\) not divisible
by seven. Adding an isolated vertex gives a linear-deficit construction
for the remaining large sizes as well. The next elementary lemma removes
that padding and improves the constant.

\textbf{Lemma 12.} For every \(n\ge25\), there is an odd integer \(h\)
with \(7\le h<n/2\) and \(\gcd(n,h)=1\).

The proof and its explicit choices are included in the construction
appendix below. Taking \(a=(h-1)/2\) gives, for every \(n\ge25\),
\[n(n-4)\le M_V(n)\le n(n-3),\qquad 3n\le n^2-M_V(n)\le4n.\] These
estimates motivate the sharper endpoint result: Theorem 40 subsequently
proves equality in the upper bound for every \(n\ge12\).

\subsection{Local strategies on cyclic
graphs}\label{local-strategies-on-cyclic-graphs}

The following lemmas separate the finite strategy certificates from the
ordinary simulation and composition arguments used by the sparse
families.

\textbf{Lemma (local simulation).} Represent actual delivery-game states
by states of a local alternating arena, preserving the player to move. A
local winning strategy applies to the actual game if:

\begin{enumerate}
\def\labelenumi{\arabic{enumi}.}
\tightlist
\item
  Its selected messenger moves lift to legal actual moves with the
  prescribed represented successors.
\item
  Every actual pursuer response has a permitted local image, and
  nonterminal collisions are preserved.
\item
  Local delivery success means actual receipt; local safe handoff means
  a noncollision messenger turn after the pursuer response, satisfying
  the next module's source requirements.
\end{enumerate}

Its messenger-move bound is preserved. Encountering the actual recipient
inside a handoff module instead ends play immediately, before collision.

\textbf{Proof.} The first two conditions represent every actual play
following the strategy by a legal local play; the third identifies its
success. Termination follows from the certified nonnegative ranks. Both
conventions used below suffice: strict descent at each nonterminal
half-move, or strict descent at messenger moves and nonincrease at
pursuer responses. The latter strictly decreases rank each round. Rank
zero requires the stated success; nonterminal collision states cannot be
winning. In particular, a handoff is not complete until the pursuer
reply has been survived. \(\square\)

\textbf{Lemma (window projection).} Represent a pursuer exactly in
columns \([-R,R]\) of a ring of circumference \(q\), and otherwise by
\(\Omega\). Confine the messenger to \([-M,M]\), where \(M<R\). Suppose
actual pursuer jumps have absolute column displacement at most
\(J\ge1\), and the local arena includes their represented moves, sends
exits to \(\Omega\), and allows \(\Omega\) to wait or enter any existing
vertex in either boundary's \(J\) layers. If \[q\ge2R+J+1,\] every
actual pursuer move projects legally and no collision is concealed.
Selected messenger moves must separately lift legally.

\textbf{Proof.} Represented columns are distinct. A wrapped jump between
their opposite ends would need displacement at least \(q-2R>J\). Outside
moves project to waiting at \(\Omega\) or entry into a permitted
boundary layer; other moves are literal or exits. Every messenger
position lies in the exact pursuer window, so collisions are detected.
\(\square\)

For a single periodically repeated defect, include all changed outgoing
tails, deleted vertices and exceptional source or recipient positions in
its relative support \([a,b]\). A window meets only anchors in
\([-R-b,R-a]\). Hence \[q\ge2R+b-a+1\] allows at most one periodic copy
to affect the window. Certifying those anchors and the unchanged arena
supplies a module at every center, subject to the move inclusions above.
Several nearby defects instead require a proved covering list of local
crops.

\textbf{Theorem (cyclic strategy composition).} Fix a recipient and an
orientation of its \(q\) columns. Delivery-entry columns form a
consecutive block of length \(1\le L<q\). From every eligible state
outside this block, suppose a simulated transfer uses at most \(\tau\)
messenger moves and ends at a safe handoff advancing by some step in
\(\{1,\ldots,d\}\), with \(1\le d\le L\). From every eligible entry
state, suppose simulated delivery uses at most \(D\) moves. Modules
cover all noncollision pursuer starts, and transfers preserve messenger
eligibility. Then delivery takes at most \[\tau(q-L)+D\] messenger
moves.

\textbf{Proof.} Outside the entry block, use the forward distance to its
first column as potential; it is at most \(q-L\). A transfer either
decreases it or enters the block, since a step of size at most \(L\)
cannot skip the whole block. Thus at most \(q-L\) phases suffice. Each
handoff preserves turn, eligibility and a pursuer position covered by
the next module. Recenter between modules, then deliver. Earlier receipt
only shortens play, and summing the local bounds gives the estimate.
\(\square\)

Unit transfers therefore need at most \(q-1\) phases for one entry
column and \(q-2\) for two adjacent entries. Progress of one or two
likewise reaches two adjacent entries within \(q-2\) phases; it need not
reach one specified column. Eligibility can exclude particular messenger
types or vertices, as in the odd construction.

For optional arrows, restrict the messenger to compulsory moves and give
the local pursuer every optional move: every intermediate graph then
satisfies the simulation inclusions. Omitting an external defect is
sound only if the selected messenger moves remain legal and all actual
pursuer responses remain represented. The following chapters check these
construction-specific facts and retain every finite certificate premise.

\subsection{An infinite family attaining the sparse
bounds}\label{an-infinite-family-attaining-the-sparse-bounds}

The five-type graph used in the dense construction also attains the
sparse vertex and arc-budget bounds. This requires a different strategy:
its messenger may need to travel around a substantial part of the ring.
The proof combines fixed local certificates with a simulation valid at
every sufficiently large circumference, followed by six finite base
cases.

For clarity, define the sparse graph \(H_k\) again. Its vertices are
\[ (i,j)\in (\mathbb Z/k\mathbb Z)\times\{0,1,2,3,4\}.\] The first
outgoing arc follows the cycle
\[(i,0)\to(i,1)\to(i,2)\to(i,3)\to(i,4)\to(i+1,0).\] The second outgoing
arc is specified by \[
\begin{aligned}
(i,0)&\to(i+1,0), & (i,1)&\to(i-1,2),\\
(i,2)&\to(i-1,4), & (i,3)&\to(i-1,3),\\
(i,4)&\to(i-1,1).
\end{aligned}
\] Both arc maps are permutations, are loopless, and have distinct
images at every vertex when \(k\ge4\). Hence \(H_k\) has \(5k\) vertices
and \(10k\) arcs, with indegree and outdegree two everywhere. Waiting is
an additional legal action, as throughout the paper.

\begin{figure}
\centering
\pandocbounded{\includegraphics[keepaspectratio,alt={The five station types in consecutive columns. Black arrows follow the long cycle and blue arrows give the second permutation. Identifying columns modulo k closes the strip into the graph H\_k.}]{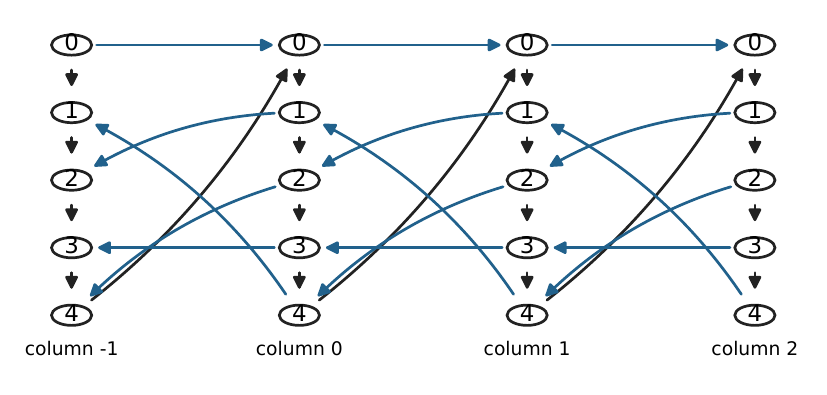}}
\caption{The five station types in consecutive columns. Black arrows
follow the long cycle and blue arrows give the second permutation.
Identifying columns modulo k closes the strip into the graph H\_k.}
\end{figure}

\textbf{Theorem 21 (uniform sparse attainment).} For every integer
\(k\ge4\), every source--recipient pair in \(H_k\) is guaranteed.
Consequently \[
M_V(5k)=5k(5k-3),\qquad
M_E(10k)=5k(5k-3).
\] For \(k\ge10\), a sufficient delivery time is \(8k+16\) messenger
moves. The theorem is computer-assisted: its finite ingredients are the
two local objectives specified below and the six rings \(4\le k\le9\).

\subsubsection{The stronger local arena}\label{the-stronger-local-arena}

Unwrap the column coordinate to the integers, retaining the same two arc
rules. Confine the messenger to columns \(-3,\ldots,3\), disallowing
moves out of this interval. Represent the pursuer by its exact vertex in
columns \(-4,\ldots,4\), together with an exterior state \(\Omega\).
Thus there are 35 messenger positions and 46 pursuer positions.

Within its nine-column window the pursuer has all its usual moves and
its wait; a move leaving the window goes to \(\Omega\). From \(\Omega\)
it may wait or enter any vertex of either boundary column \(-4\) or
\(4\). This grants more freedom than a particular physical pursuer
outside the window. The messenger and \(\Omega\) can never coincide.

The objectives are safe handoff in column zero and immediate delivery to
\((0,t)\), with the distinct stopping rules of local simulation. The
certificates use strict half-move rank descent, seeding handoff only at
noncollision messenger turns and receipt before collision.

The checked finite bounds are:

{\def\LTcaptype{none} 
\Needspace{10\baselineskip}
\begin{longtable}[]{@{}
  >{\raggedright\arraybackslash}p{(\linewidth - 6\tabcolsep) * \real{0.2500}}
  >{\raggedleft\arraybackslash}p{(\linewidth - 6\tabcolsep) * \real{0.2500}}
  >{\raggedleft\arraybackslash}p{(\linewidth - 6\tabcolsep) * \real{0.2500}}
  >{\raggedleft\arraybackslash}p{(\linewidth - 6\tabcolsep) * \real{0.2500}}@{}}
\toprule\noalign{}
\begin{minipage}[b]{\linewidth}\raggedright
Local objective
\end{minipage} & \begin{minipage}[b]{\linewidth}\raggedleft
Initial column
\end{minipage} & \begin{minipage}[b]{\linewidth}\raggedleft
Initial noncollision states
\end{minipage} & \begin{minipage}[b]{\linewidth}\raggedleft
Maximum messenger moves
\end{minipage} \\
\midrule\noalign{}
\endhead
\bottomrule\noalign{}
\endlastfoot
Safe transfer to column zero & 1 & 225 & 8 \\
Delivery to \((0,0)\) & 2 & 225 & 24 \\
Delivery to \((0,1)\) & 2 & 225 & 16 \\
Delivery to \((0,2)\) & 2 & 225 & 16 \\
Delivery to \((0,3)\) & 2 & 225 & 17 \\
Delivery to \((0,4)\) & 2 & 225 & 18 \\
\end{longtable}
}

Each row covers all five source types and all 46 pursuer states other
than the source itself. In particular, it includes an initially exterior
pursuer. Only the bounds 8 and 24 are needed below; optimality of these
local times is not claimed as a separate result.

\subsubsection{Projection to a ring}\label{projection-to-a-ring}

Apply the window-projection lemma with \(M=3\), \(R=4\), \(J=1\) and
circumference \(k\ge10=2R+J+1\). Both arrow rules lift literally, so
every selected messenger move is legal. The local-simulation lemma
therefore transfers both objectives to every translated window. At
\(k=9\) the boundary columns have wraparound adjacencies, which is why
smaller rings retain separate finite certificates.

\subsubsection{Composition and
termination}\label{composition-and-termination}

\textbf{Proof of Theorem 21.} Relative to a recipient in column zero,
transfers move one column backward and delivery starts in column two.
All source types are eligible. Cyclic strategy composition with
\(L=d=1\), \(\tau=8\) and \(D=24\) gives the bound \[8(k-1)+24=8k+16\]
messenger moves.

For \(4\le k\le9\), the explicit finite-ring certificates verify every
source and pursuer position for each of the five recipient types in
column zero. Translating columns gives every recipient. These six
checks, together with the uniform argument for \(k\ge10\), establish
universality for all \(k\ge4\).

Finally, a universal graph with \(n=5k\) vertices and \(2n\) arcs has
\[I(H_k)=n^2-n-2n=n(n-3).\] The vertex upper bound and the even
arc-budget upper bound both equal this number. Hence both asserted
extremal equalities follow. \(\square\)

\textbf{Certificate reproduction and scope.} The file
\texttt{research/\allowbreak{}five-\allowbreak{}type-\allowbreak{}certificates.\allowbreak{}json} contains the fixed local
arena, ranks, and chosen messenger moves, generated and directly checked
by \texttt{src/\allowbreak{}five\_\allowbreak{}type\_\allowbreak{}certificate.\allowbreak{}py}. The independent
implementation \texttt{src/\allowbreak{}five\_\allowbreak{}type\_\allowbreak{}local\_\allowbreak{}verify.\allowbreak{}py} constructs the
explicit alternating-state graph from the arc rules and checks every
winning rank condition, as well as closure of the losing region. Its
full ranks and selected moves are saved in
\texttt{research/\allowbreak{}five-\allowbreak{}type-\allowbreak{}independent-\allowbreak{}local-\allowbreak{}certificates.\allowbreak{}json},
together with the implementation hash. Both implementations agree on the
six local bounds above. The six finite bases are stored in
\texttt{research/\allowbreak{}five-\allowbreak{}type-\allowbreak{}small-\allowbreak{}ring-\allowbreak{}certificates.\allowbreak{}json}, including
complete arc lists, messenger-move layers, and explicit
alternating-state rank certificates for the five representative
recipient types. The script
\texttt{src/\allowbreak{}five\_\allowbreak{}type\_\allowbreak{}small\_\allowbreak{}certificates.\allowbreak{}py} reproduces those checks
using both game solvers and verifies every certificate transition. The
ordinary projection and composition proof explains why these finite
arrays certify infinitely many rings. This theorem is not currently
claimed as a Lean formalization.

\textbf{Corollary 22 (sharp arc-budget leading term).} As
\(m\to\infty\), \[M_E(m)=\frac{m^2}{4}-O(m).\] More explicitly, for
every \(m\ge40\),
\[\frac{m^2}{4}-6m\le M_E(m)\le\left\lfloor\frac m2\right\rfloor^2.\]

\textbf{Proof.} Set \(k=\lfloor m/10\rfloor\ge4\) and \(r=m-10k\), so
\(0\le r\le9\). Take \(H_k\) together with \(r\) disjoint one-way arcs.
Each added component has no indirect guaranteed pair, and additivity
over components preserves \(I(H_k)\). This graph has exactly \(m\) arcs
and \[
I=25k^2-15k
=\frac{(m-r)^2}{4}-\frac{3(m-r)}2
\ge\frac{m^2}{4}-6m.
\] The upper bound follows because there are at most
\(\lfloor m/2\rfloor\) possible sources and at most that many possible
recipients of indirect guaranteed pairs. Thus the quadratic coefficient
\(1/4\) is sharp for arbitrary arc counts, not only the progression
\(m=10k\). \(\square\)

Theorem 40 proves exact vertex maxima at every order from twelve and
exact even arc maxima from 24 arcs; Corollary 45.1 proves exact odd arc
maxima from 25 arcs. The padding argument above is retained as a direct
proof of the leading term. The remaining small exceptions are stated in
the current result summary.

\subsection{Sparse endpoint attainment at every order from
twelve}\label{sparse-endpoint-attainment-at-every-order-from-twelve}

A second periodic construction admits controlled removal of vertices.
Unlike a finite modification checked at one circumference, the removal
is certified in a fixed local game that applies at every translated
position. Separating the removals then provides all residue classes of
the vertex count.

\textbf{Theorem 40 (with finite certificates).} For every integer
\(n\ge12\), there is an explicit two-in/two-out-regular universally
winning digraph on \(n\) vertices. Consequently
\[M_V(n)=M_E(2n)=M(n,2n)=n(n-3).\] Every source--recipient pair in the
construction can be guaranteed within \(9\lceil n/6\rceil+18\) messenger
moves.

\textbf{Proof.} First construct a uniform family for \(n\ge85\). Begin
with \(k\) columns, indexed modulo \(k\), containing vertices \((i,j)\)
of types \(0\le j\le5\). One permutation follows the Hamilton cycle
\[0_i\to1_i\to2_i\to3_i\to4_i\to5_i\to0_{i+1}.\] The other permutation
has arcs \[
\begin{array}{lll}
0_i\to5_i,&1_i\to4_{i+1},&2_i\to1_{i+1},\\
3_i\to0_{i-1},&4_i\to2_{i+1},&5_i\to3_{i-1}.
\end{array}
\] Choose columns \(S\) whose consecutive cyclic gaps are at least
three. At each \(h\in S\), remove \(4_h\) and smooth its path separately
in the two permutations: \[
3_h\to4_h\to5_h\quad\longmapsto\quad3_h\to5_h,
\] \[
1_{h-1}\to4_h\to2_{h+1}\quad\longmapsto\quad1_{h-1}\to2_{h+1}.
\] All remaining arcs stay unchanged. The resulting graph has \(6k-|S|\)
vertices. Smoothing preserves each permutation's bijectivity. There are
no loops or coincident outgoing arcs: the modified first-permutation arc
leaves type three for type five, while its other arc ends at type zero;
the modified second-permutation arc leaves type one for type two two
columns ahead, while its first-permutation arc ends at type two in the
same column. Thus the graph remains simple and two-in/two-out regular.

\begin{figure}
\centering
\pandocbounded{\includegraphics[keepaspectratio,alt={Remove one type-four vertex and reconnect each color separately. All undisplayed arcs stay unchanged. Spaced removals give every residue class of the total number of vertices; the local certificates prove that delivery remains possible.}]{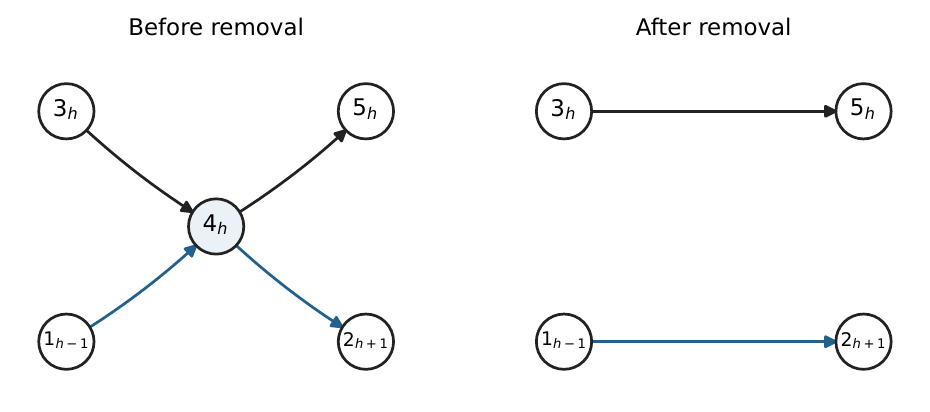}}
\caption{Remove one type-four vertex and reconnect each color
separately. All undisplayed arcs stay unchanged. Spaced removals give
every residue class of the total number of vertices; the local
certificates prove that delivery remains possible.}
\end{figure}

We establish a uniform strategy when \(k\ge15\). Around an arbitrary
reference column zero, restrict the messenger to columns \([-4,4]\) and
represent the pursuer in columns \([-6,6]\), together with a state
\(\Omega\) for all outside positions. Retain both players' waits and
graph arcs. Messenger moves leaving its window are disallowed; pursuer
moves leaving its window lead to \(\Omega\). From \(\Omega\) allow
waiting and entry at any existing vertex in columns \(-6,-5,5,6\).

The removals inside the represented pursuer window form a subset of
\(\{-6,\ldots,6\}\) with consecutive gaps at least three. There are
exactly 189 such subsets. For every one of these patterns, the following
finite local claims have complete rank certificates:

{\def\LTcaptype{none} 
\Needspace{6\baselineskip}
\begin{longtable}[]{@{}
  >{\raggedright\arraybackslash}p{(\linewidth - 4\tabcolsep) * \real{0.3333}}
  >{\raggedleft\arraybackslash}p{(\linewidth - 4\tabcolsep) * \real{0.3333}}
  >{\raggedleft\arraybackslash}p{(\linewidth - 4\tabcolsep) * \real{0.3333}}@{}}
\toprule\noalign{}
\begin{minipage}[b]{\linewidth}\raggedright
Objective
\end{minipage} & \begin{minipage}[b]{\linewidth}\raggedleft
Initial messenger column
\end{minipage} & \begin{minipage}[b]{\linewidth}\raggedleft
Maximum messenger moves
\end{minipage} \\
\midrule\noalign{}
\endhead
\bottomrule\noalign{}
\endlastfoot
Reach column zero safely at the start of a messenger turn & \(-1\) &
\(9\) \\
Deliver to any specified existing vertex in column zero & \(-2\) &
\(27\) \\
\end{longtable}
}

Each claim covers every existing initial messenger type and every
noncolliding initial pursuer position, including \(\Omega\).
Safe-transfer goals are seeded only at messenger turns, after the
preceding pursuer response. Delivery goals instead use immediate receipt
before collision. There are 1,287 certificates: seven objectives for
each of the 189 patterns, with delivery to \(4_0\) omitted in the 36
patterns where that vertex is removed. An independent explicit
alternating-state implementation verifies all 560,885 required initial
states and all 10,271,740 states of these finite games, including
winning rank descent and losing-region closure.

The window-projection lemma applies with \(M=4\), \(R=6\), \(J=2\) and
\(k\ge15\). It remains to check that selecting a crop neither removes an
actual pursuer move nor supplies an illegal messenger move.

A removal outside the represented window does not change any
represented-to-represented arc. The only change to an outward arc occurs
when removing \(4_7\) redirects \(1_6\to4_7\) to \(1_6\to2_8\); both
heads represent \(\Omega\). Outside-to-inside changes are already
included among its permitted boundary entries. At the messenger
boundary, any analogous changed outgoing arc still leaves the allowed
messenger window. The local model therefore depends only on the removals
inside the pursuer window. Their consecutive gaps are at least three, so
it always has one of the 189 certified patterns.

All existing messenger types are covered by each certified crop.
Transfers advance one column, and delivery starts two columns before the
recipient. Cyclic strategy composition with \(L=d=1\), \(\tau=9\) and
\(D=27\) therefore gives \[9(k-1)+27=9k+18.\] This proves universal
winning for every permitted set \(S\).

Finally, given \(n\ge85\), put \(k=\lceil n/6\rceil\) and \(r=6k-n\).
Then \(k\ge15\) and \(0\le r\le5\). Choose \[S=\{0,3,\ldots,3(r-1)\}\]
when \(r>0\), and take \(S\) empty when \(r=0\). The displayed gaps are
three, and the wraparound gap is at least \(15-12=3\). The graph has
exactly \(n\) vertices and \(2n\) arcs. Universal winning gives
\(I=n(n-1)-2n=n(n-3)\), and the vertex and even arc-budget upper bounds
give the three stated equalities for \(n\ge85\).

The complete primary certificates and their direct checker are
\texttt{research/\allowbreak{}sparse-\allowbreak{}composition-\allowbreak{}multi-\allowbreak{}certificates.\allowbreak{}json.\allowbreak{}gz} and
\texttt{src/\allowbreak{}sparse\_\allowbreak{}composition\_\allowbreak{}multi\_\allowbreak{}certificate.\allowbreak{}py}. The
independently written verifier uses direct arc formulas and an explicit
alternating graph: \texttt{src/\allowbreak{}six\_\allowbreak{}type\_\allowbreak{}independent\_\allowbreak{}spacing.\allowbreak{}py}, with
full ranks in
\texttt{research/\allowbreak{}six-\allowbreak{}type-\allowbreak{}independent-\allowbreak{}spacing-\allowbreak{}3-\allowbreak{}certificates.\allowbreak{}json.\allowbreak{}gz}
and a summary in
\texttt{research/\allowbreak{}six-\allowbreak{}type-\allowbreak{}independent-\allowbreak{}spacing-\allowbreak{}3-\allowbreak{}summary.\allowbreak{}json}. The
ordinary argument uses the two certified local claims. For the
prescribed removals giving every order \(n\ge85\),
\texttt{SixTypeConstruction.\allowbreak{}lean} now instantiates the graph, canonical
crop coverage, local certificates, projection, composition and exact
maximum. This does not assert Lean coverage of every permitted removal
pattern or of all the finite witnesses below 85.

For completeness, the finite range \(12\le n\le84\) is supplied by 73
explicit graphs. Each graph has \(2n\) distinct nonloop arcs and both
incoming and outgoing degrees two. Two independent exact game algorithms
agree, and every actual-game rank obligation is checked for every
recipient and every allowed initial pursuer. The finite witnesses also
satisfy the displayed delivery-time bound. Their arc lists are in
\texttt{research/\allowbreak{}sparse-\allowbreak{}finite-\allowbreak{}bridge-\allowbreak{}witnesses.\allowbreak{}json}, the complete
ranks in \texttt{research/\allowbreak{}sparse-\allowbreak{}finite-\allowbreak{}bridge-\allowbreak{}certificates.\allowbreak{}jsonl.\allowbreak{}gz},
and the verification receipt in
\texttt{research/\allowbreak{}sparse-\allowbreak{}finite-\allowbreak{}bridge-\allowbreak{}verification.\allowbreak{}json}. The unified
generator \texttt{src/\allowbreak{}sparse\_\allowbreak{}finite\_\allowbreak{}bridge\_\allowbreak{}construct.\allowbreak{}py} returns the
corresponding finite witness below 85 and the uniform construction above
it. This proves the remaining finite cases of the theorem. \(\square\)

Together with the complete small-order tables, this settles the vertex
maximum at every order. The separate Theorem 45 addresses the odd
arc-budget bound. \texttt{OddEndpointConstruction.\allowbreak{}lean} formalizes its
uniform family from support order98; \texttt{OddSparseConstruction.\allowbreak{}lean}
and \texttt{OddExplicitAttainer.\allowbreak{}lean} now join every finite support
from13. The verification section distinguishes these maximum and recipe
statements from stronger component properties.

\subsection{The odd endpoint and completion of its degree
defect}\label{the-odd-endpoint-and-completion-of-its-degree-defect}

The odd arc-budget equality theorem prescribes a particularly small
defect: on \(n=k+1\) vertices there is one vertex \(u\) of outdegree
one, a different vertex \(v\) of indegree one, all other indicated
degrees are two, and the arc \(u\to v\) is absent. The candidate
extremal value is \[M_E(2n-1)=(n-2)(n-3).\] The equality theorem reduces
attainment to guaranteeing every missing pair whose source differs from
\(u\) and whose recipient differs from \(v\). It does not make arbitrary
deletion from, or completion to, a two-regular graph preserve the
relevant guarantees.

\subsubsection{Two necessary structural
properties}\label{two-necessary-structural-properties}

\textbf{Lemma.} Every odd-extremal graph with the preceding degree
pattern is strongly connected. Moreover, for each vertex \(c\), the
graph with \(c\) deleted has a directed path from every
\(s\notin\{u,c\}\) to every \(t\notin\{v,c\}\).

\textbf{Proof.} Each such ordered pair with \(s\ne t\) is either a
direct arc or a guaranteed missing pair. Against a pursuer initially at
\(c\) and waiting there, any successful delivery path avoids \(c\)
internally. Since its endpoints differ from \(c\), it gives the asserted
path in the deleted graph.

To prove strong connectivity without deleting a vertex, all vertices
different from \(u\) already reach all vertices different from \(v\).
The unique predecessor \(b\) of \(v\) is different from \(u\), because
\(u\to v\) is absent. Thus all vertices other than \(u\) reach \(v\) by
first reaching \(b\). The unique successor \(a\) of \(u\) differs from
\(u\), and \(u\) reaches every vertex by first going to \(a\).
\(\square\)

These are necessary conditions, not substitutes for the full game: the
pursuer can move, so ordinary connectivity and stationary separators do
not establish an attaining strategy.

\subsubsection{Completing an odd extremizer can
fail}\label{completing-an-odd-extremizer-can-fail}

\textbf{Proposition (finite certified example).} There is a graph \(G\)
on 19 vertices with 37 arcs attaining \(I(G)=272=M_E(37)\) such that
adding the unique degree-completing arc yields a two-in/two-out graph
which is not universally winning. Its completed indirect count is 276,
whereas a universally winning two-regular graph on 19 vertices would
have 304.

The example has vertex set \(\{0,\ldots,18\}\) and the following
outgoing neighbors:

{\def\LTcaptype{none} 
\Needspace{14\baselineskip}
\begin{longtable}[]{@{}rlrl@{}}
\toprule\noalign{}
Vertex & Outgoing neighbors & Vertex & Outgoing neighbors \\
\midrule\noalign{}
\endhead
\bottomrule\noalign{}
\endlastfoot
0 & 5, 12 & 10 & 2, 17 \\
1 & 15, 10 & 11 & 1, 6 \\
2 & 7, 0 & 12 & 3, 8 \\
3 & 1, 18 & 13 & 4, 5 \\
4 & 15 & 14 & 4, 12 \\
5 & 11, 9 & 15 & 0, 16 \\
6 & 8, 13 & 16 & 7, 9 \\
7 & 13, 17 & 17 & 14, 11 \\
8 & 10, 16 & 18 & 2, 6 \\
9 & 3, 14 & & \\
\end{longtable}
}

Its exceptional vertices are \(u=4\) and \(v=18\). The finite
actual-game certificate verifies every guarantee required by odd
equality and the count 272. The upper bound follows from Theorem 7 with
\(k=18\).

There is also a short ordinary obstruction after adding \(4\to18\). Fix
recipient zero and start the messenger at 1, the pursuer at 4. Define \[
S_4=\{1,3,5,6,7,8,9,11,12,13,14,16,17\},\qquad
S_{13}=\{6,11,12,14\}.
\] The following four classes of messenger-turn states form a pursuer
trap: \[
(s,4)\ (s\in S_4),\qquad (10,18),\qquad (17,6),\qquad
(s,13)\ (s\in S_{13}).
\] They contain the initial state and exclude every state where the
messenger has reached zero. The pursuer acts as follows, in terms of the
messenger's new vertex \(x\):

{\def\LTcaptype{none} 
\Needspace{8\baselineskip}
\begin{longtable}[]{@{}
  >{\raggedleft\arraybackslash}p{(\linewidth - 2\tabcolsep) * \real{0.5000}}
  >{\raggedright\arraybackslash}p{(\linewidth - 2\tabcolsep) * \real{0.5000}}@{}}
\toprule\noalign{}
\begin{minipage}[b]{\linewidth}\raggedleft
Current pursuer
\end{minipage} & \begin{minipage}[b]{\linewidth}\raggedright
Response
\end{minipage} \\
\midrule\noalign{}
\endhead
\bottomrule\noalign{}
\endlastfoot
4 & If \(x=10\), move to 18. If \(x=15\) or \(18\), capture at \(x\).
Otherwise wait at 4. \\
18 & If \(x=2\), capture there. If \(x=17\), move to 6. If \(x=10\),
wait. \\
6 & If \(x=11\) or \(14\), move to 13. If \(x=17\), wait. \\
13 & If \(x\in\{1,3,4,8\}\), move to 4, capturing if \(x=4\). Otherwise
wait at 13. \\
\end{longtable}
}

An immediate collision with the current pursuer already loses. Direct
inspection of the outgoing-neighbor table shows that every other
possible move either permits the stated capture or returns to one of the
four trap classes. Waiting is included. Since infinite play without
receipt loses, the messenger cannot deliver to zero. This explicitly
identifies how the new arc helps the pursuer: the transition from cop
position 4 to 18 after the messenger moves from 1 to 10 was unavailable
before completion.

\texttt{src/\allowbreak{}odd\_\allowbreak{}endpoint\_\allowbreak{}completion.\allowbreak{}py} contains the fixed graph and
runs both independent game solvers, checking every winning rank and
every losing closure obligation.
\texttt{research/\allowbreak{}odd-\allowbreak{}endpoint-\allowbreak{}completion-\allowbreak{}certificate.\allowbreak{}json} stores the
complete finite certificates for the original and completed graphs. The
ordinary trap proves nonuniversality without computation; the original
graph's attainment and the completed count 276 use the verified finite
certificates. This example is not presently embedded in Lean. The
subsequent construction settles the odd endpoint at every sufficiently
large budget.

\subsubsection{A useful degree-preserving search
operation}\label{a-useful-degree-preserving-search-operation}

Given a two-in/two-out graph, split one vertex \(w\) into vertices
\(v,u\). Assign one old incoming and one old outgoing arc to each new
vertex, then add \(v\to u\). This adds one vertex and one arc, giving
exactly the odd equality degree pattern: \(v\) has indegree one and
outdegree two, while \(u\) has indegree two and outdegree one. The
forbidden opposite arc \(u\to v\) is absent because the original graph
was loopless.

This operation supplies a structured family to investigate, but it is
not a preservation theorem. A split of the 20-vertex five-type graph
attains \(M_E(41)=342\), while all twenty canonical splits of its
30-vertex member fail to attain the corresponding odd bound. Likewise,
changing the two horizontal lane speeds to three produces attaining
single-arc deletions at 40 and 50 vertices but failures at neighboring
and larger orders. These counterchecks prevent promoting isolated
success to a uniform construction. A different six-type family supplies
the following uniform positive result.

\subsubsection{Every sufficiently large odd
budget}\label{every-sufficiently-large-odd-budget}

\textbf{Theorem 45 (the odd endpoint).} For every integer \(n\ge98\), \[
M_E(2n-1)=(n-2)(n-3).
\] Equivalently, for every odd \(m\ge195\), \[
M_E(m)=\frac{(m-3)(m-5)}4.
\] An explicit attaining graph has \(n\) vertices. Every pair required
for attainment has a delivery strategy using at most
\(34\lceil n/6\rceil+27\) messenger moves. The construction and its
composition proof are ordinary finite arguments; the bounded local
strategy lemmas are supplied by independently verified finite
certificates.

\textbf{Proof.} Put \(q=\lceil n/6\rceil\) and \(r=6q-n\), so
\(0\le r\le5\). The sufficient circumference for each residue is given
below. We first describe a graph on the surviving vertices of
\((\mathbb Z/q\mathbb Z)\times\{0,\ldots,5\}\). Its two initial
permutations are \[
\alpha(i,j)=
\begin{cases}(i,j+1)&j<5,\\(i+1,0)&j=5,\end{cases}
\qquad
\beta(i,j)=(i+e_j,\sigma(j)).
\] Use the following four choices of \((\sigma,e)\), whose entries are
in type order \(0,1,\ldots,5\):

{\def\LTcaptype{none} 
\Needspace{8\baselineskip}
\begin{longtable}[]{@{}lll@{}}
\toprule\noalign{}
Pattern & \(\sigma\) & \(e\) \\
\midrule\noalign{}
\endhead
\bottomrule\noalign{}
\endlastfoot
\(A\) & \((0,2,5,1,4,3)\) & \((-1,-1,1,1,-1,1)\) \\
\(B\) & \((0,5,1,4,3,2)\) & \((-1,1,1,-1,-1,1)\) \\
\(C\) & \((3,5,0,2,4,1)\) & \((-1,1,-1,-1,-1,0)\) \\
\(D\) & \((4,1,5,2,0,3)\) & \((-1,-1,1,1,-1,1)\) \\
\end{longtable}
}

For a set \(H\) of deleted vertices, smooth each permutation separately:
from a surviving vertex, follow that permutation until first reaching
another surviving vertex, and use the resulting arc. Finally delete one
specified arc of the smoothed graph. The complete choices are:

{\def\LTcaptype{none} 
\Needspace{10\baselineskip}
\begin{longtable}[]{@{}rclcr@{}}
\toprule\noalign{}
\(r\) & Pattern & Tail \(u\) of deleted arc & Which arc & Direction
\(d\) \\
\midrule\noalign{}
\endhead
\bottomrule\noalign{}
\endlastfoot
0 & \(A\) & \((0,0)\) & \(\alpha\) & \(1\) \\
1 & \(B\) & \((1,4)\) & \(\beta\) & \(1\) \\
2 & \(C\) & \((0,5)\) & \(\beta\) & \(-1\) \\
3 & \(D\) & \((0,1)\) & \(\alpha\) & \(1\) \\
4 & \(D\) & \((0,1)\) & \(\alpha\) & \(1\) \\
5 & \(C\) & \((0,5)\) & \(\beta\) & \(-1\) \\
\end{longtable}
}

The corresponding deleted vertex sets are:

\begin{itemize}
\tightlist
\item
  \(r=0\): \(H=\varnothing\).
\item
  \(r=1\): \(H=\{(0,3)\}\).
\item
  \(r=2\): \(H=\{(0,0),(0,1)\}\).
\item
  \(r=3\): \(H=\{(0,2),(0,3),(0,4)\}\).
\item
  \(r=4\): \(H=\{(-1,4),(0,2),(0,3),(0,4)\}\).
\item
  \(r=5\): \(H=\{(0,0),(0,1),(1,0),(1,1),(1,5)\}\).
\end{itemize}

Write \(v\) for the head of the final deleted arc. Every type cycle of
each displayed \(\sigma\) has a nonzero total displacement \(e\) with
\(|e|\) no greater than its length \(\ell\). Its cyclic lift has length
\(\ell q/\gcd(q,e)\ge q\), whereas \(\alpha\) is a cycle of length
\(6q\). At most five vertices are removed, so no permutation cycle
disappears. The smoothing operation therefore preserves the permutation
property on surviving vertices. Inspection of these six fixed patterns
shows that the two smoothed images at each surviving vertex are
distinct, are not the vertex itself, and change the column by at most
two. The changes to outgoing arcs are confined to tails in columns
\([-2,2]\); removed vertices lie in columns \([-1,1]\). These are finite
checks on the displayed tables, independent of \(q\). The sufficient
circumferences are at least 15, so reduction modulo \(q\) introduces
neither a coincidence nor a loop. Thus the final graph has \(n=6q-r\)
vertices, \(2n-1\) arcs, precisely one outdegree-one vertex \(u\),
precisely one indegree-one vertex \(v\), and all other corresponding
degrees two. Moreover \(u\ne v\) and \(u\to v\) is absent.

Use the retained local certificates with messenger radius \(M\), pursuer
radius \(R\), and an outside state \(\Omega\). Both players may wait;
messenger moves leaving its window are forbidden. A pursuer leaving its
window goes to \(\Omega\), which permits waiting and entry into its
\(J\) boundary layers. Let \([a,b]\) contain every changed outgoing
tail, deleted vertex and degree-exception vertex. The exact parameters
are:

{\def\LTcaptype{none} 
\Needspace{10\baselineskip}
\begin{longtable}[]{@{}
  >{\raggedleft\arraybackslash}p{(\linewidth - 14\tabcolsep) * \real{0.1200}}
  >{\raggedleft\arraybackslash}p{(\linewidth - 14\tabcolsep) * \real{0.1200}}
  >{\raggedleft\arraybackslash}p{(\linewidth - 14\tabcolsep) * \real{0.1200}}
  >{\raggedleft\arraybackslash}p{(\linewidth - 14\tabcolsep) * \real{0.1200}}
  >{\centering\arraybackslash}p{(\linewidth - 14\tabcolsep) * \real{0.1600}}
  >{\raggedleft\arraybackslash}p{(\linewidth - 14\tabcolsep) * \real{0.1200}}
  >{\raggedleft\arraybackslash}p{(\linewidth - 14\tabcolsep) * \real{0.1200}}
  >{\raggedleft\arraybackslash}p{(\linewidth - 14\tabcolsep) * \real{0.1200}}@{}}
\toprule\noalign{}
\begin{minipage}[b]{\linewidth}\raggedleft
\(r\)
\end{minipage} & \begin{minipage}[b]{\linewidth}\raggedleft
\(M\)
\end{minipage} & \begin{minipage}[b]{\linewidth}\raggedleft
\(R\)
\end{minipage} & \begin{minipage}[b]{\linewidth}\raggedleft
\(J\)
\end{minipage} & \begin{minipage}[b]{\linewidth}\centering
\([a,b]\)
\end{minipage} & \begin{minipage}[b]{\linewidth}\raggedleft
Sufficient \(q\)
\end{minipage} & \begin{minipage}[b]{\linewidth}\raggedleft
Transfer
\end{minipage} & \begin{minipage}[b]{\linewidth}\raggedleft
Delivery
\end{minipage} \\
\midrule\noalign{}
\endhead
\bottomrule\noalign{}
\endlastfoot
0 & 4 & 6 & 2 & \([0,0]\) & 15 & 14 & 26 \\
1 & 4 & 6 & 1 & \([-1,1]\) & 15 & 9 & 22 \\
2 & 4 & 6 & 2 & \([-1,1]\) & 15 & 17 & 26 \\
3 & 5 & 7 & 2 & \([-1,1]\) & 17 & 34 & 23 \\
4 & 5 & 7 & 2 & \([-1,1]\) & 17 & 28 & 27 \\
5 & 5 & 7 & 2 & \([-1,2]\) & 18 & 9 & 26 \\
\end{longtable}
}

The table's support is obtained without a search: trace each removed
vertex backwards in each permutation to its first surviving predecessor,
and include the final deleted arc and its two endpoints. Outside those
tails, outgoing arcs are unchanged. All actual column jumps are at most
\(J\); in residue zero the allowed entry depth two is conservative.

For every represented defect translation, safe transfer goes from column
\(-d\) to column zero on a noncolliding messenger turn, excluding \(u\)
at both endpoints. Delivery goes from column \(-2d\) to every surviving
recipient of column zero other than \(v\), again excluding initial
source \(u\). In residue five, transfer may instead end in column zero
or \(-1\), and delivery starts in column 2 or 3. Receipt is terminal
before collision. The table lists the maximum certified messenger moves.

The frozen independent archives contain all 847 local objectives: 97 for
residue zero and 750 for the nonempty patches. The direct replay checks
all 470,026 required initial states, 8,417,696 alternating states, and
24,787,015 legal-transition obligations, including losing closure. For
residue zero the anchors are \(-6,\ldots,6\); for the others they are
\(-10,\ldots,10\), always with an additional unchanged arena. These are
finite local lemmas, not tests extrapolated over ring lengths.

The window-projection lemma and its single-defect condition require
\[q\ge2R+J+1,\qquad q\ge2R+b-a+1.\] The relevant anchor interval
\([-R-b,R-a]\) lies within the certified anchors in every table row. On
represented outgoing tails the selected one-patch arena agrees with the
smoothed arrow formulas; other incoming effects are covered by the
permitted exterior entries. If no patch is relevant, the unchanged arena
applies. Thus local simulation holds at every center. The displayed
circumferences satisfy both inequalities and cover every \(n\ge98\); the
last excluded residue representative is \(n=97\).

Fix a source other than \(u\) and a recipient other than \(v\).
Eligibility means that the messenger avoids \(u\) at module starts and
handoffs; every certified transfer preserves it. For \(r\le4\),
transfers advance one column in the displayed direction \(d\), and
delivery starts in the single column two steps before the recipient.
Cyclic strategy composition uses at most \(q-1\) transfers. For \(r=5\),
transfers advance one or two columns leftward, and the delivery entries
are the two adjacent columns 2 and 3 relative to the recipient. The same
theorem, with maximum progress two and entry-block length two, uses at
most \(q-2\) transfers. In particular, the uniform bound \(34q+27\)
stated above remains valid.

Every missing pair with source different from \(u\) and recipient
different from \(v\) is therefore guaranteed. Their number is
\((n-1)^2-(n-2)-(2n-3)=(n-2)(n-3)\). The odd arc-budget upper bound
proves equality. \(\square\)

When \(r=0\), the same argument gives the sharper bound \(14q+12\) for
all \(q\ge15\), retaining the original congruence-family theorem. The
generator is \texttt{src/\allowbreak{}odd\_\allowbreak{}endpoint\_\allowbreak{}all\_\allowbreak{}residues.\allowbreak{}py}; the unchanged
proof archives are
\texttt{research/\allowbreak{}six-\allowbreak{}type-\allowbreak{}independent-\allowbreak{}odd-\allowbreak{}certificates.\allowbreak{}json.\allowbreak{}gz} and
\texttt{research/\allowbreak{}six-\allowbreak{}type-\allowbreak{}independent-\allowbreak{}odd-\allowbreak{}residues-\allowbreak{}certificates.\allowbreak{}json.\allowbreak{}gz}.
\texttt{src/\allowbreak{}odd\_\allowbreak{}projection\_\allowbreak{}threshold\_\allowbreak{}check.\allowbreak{}py\ -\allowbreak{}\/\allowbreak{}-\allowbreak{}check} replays
their full ranks, checks exact defect support and the projection
inequalities, and compares its result with
\texttt{research/\allowbreak{}odd-\allowbreak{}projection-\allowbreak{}threshold-\allowbreak{}audit.\allowbreak{}json}. The original
generation receipts retain their historical threshold 145 and source
hashes; no frozen certificate was regenerated or repinned. The original
primary certificates are retained as additional evidence.
\texttt{OddEndpointConstruction.\allowbreak{}lean} now instantiates the prescribed
family for every \(n\ge98\), including the 847 local objectives,
actual-game projection, exact counts and attaining maximum. The stronger
delivery bounds retain their explicitly stated domains.

\textbf{Corollary 45.1 (the odd endpoint from 25 arcs).} For every
\(n\ge13\), \[
M(n,2n-1)=M_E(2n-1)=(n-2)(n-3).
\] Consequently the exact odd arc-only extremum is
\[M_E(m)=\frac{(m-3)(m-5)}4\qquad(m\ge25\text{ odd}).\]

\textbf{Proof.} Theorem 45 covers \(n\ge98\). For every one of the
remaining 85 integers \(13\le n\le97\), the supplementary data give a
fixed simple graph with \(2n-1\) arcs, one outdegree-one vertex \(u\),
one different indegree-one vertex \(v\), all other corresponding degrees
two, and \(u\to v\) absent. The displayed six-type residue families
supply most of these graphs. Deleting one arc from independently
verified finite two-regular examples supplies the remaining orders from
19 upward. Six further fixed witnesses cover orders 13 through 18.

For every recipient other than \(v\), the finite certificates specify an
integer rank for every explicit messenger-turn and pursuer-turn state.
Terminal receipt has rank zero, nonterminal collision has rank \(-1\),
and every other nonnegative messenger rank has a legal successor of
smaller nonnegative rank. Every legal successor of a nonnegative pursuer
rank has smaller nonnegative rank. The direct checker verifies these
conditions against the original graph, together with the entire
complementary losing closure. Every initial messenger state whose source
differs from \(u\) and whose pursuer position differs from the source
has nonnegative rank. Strict rank descent therefore forces finite
delivery from each required state against every pursuer policy. This
proves all required guarantees for each listed graph.

The retained full archive, which also covers the now unnecessary orders
98 through 144, contains 10,236 recipient objectives and 215,967,810
explicit alternating states. The same degree count as in Theorem 45 then
gives \((n-2)(n-3)\) guaranteed missing pairs. The arc-only upper bound
applies to all support sizes, and hence also to the fixed support size
\(n\). This proves both equalities. \(\square\)

The 126 graphs at orders 19 through 144 are in
\texttt{research/\allowbreak{}odd-\allowbreak{}endpoint-\allowbreak{}finite-\allowbreak{}bridge-\allowbreak{}graphs.\allowbreak{}json}. Their complete
ranks are in the streamable archive
\texttt{research/\allowbreak{}odd-\allowbreak{}endpoint-\allowbreak{}finite-\allowbreak{}bridge-\allowbreak{}certificates.\allowbreak{}jsonl.\allowbreak{}gz}, with
source hashes and a full numerical summary in
\texttt{research/\allowbreak{}odd-\allowbreak{}endpoint-\allowbreak{}finite-\allowbreak{}bridge-\allowbreak{}verification.\allowbreak{}json}. The
producer \texttt{src/\allowbreak{}odd\_\allowbreak{}endpoint\_\allowbreak{}finite\_\allowbreak{}certify.\allowbreak{}cpp} uses an
explicit alternating-game attractor. The separate direct verifier
\texttt{src/\allowbreak{}odd\_\allowbreak{}endpoint\_\allowbreak{}finite\_\allowbreak{}verify.\allowbreak{}py} checks the original legal
moves; it does not rerun or trust the producer's queue logic. Its
accelerated array checks have an equivalent standard-library scalar
mode. Scalar replay of all objectives at orders 19 through 24 also
passes, and the earlier two independent game algorithms agree with
selected saved certificates at orders 19, 28 and 144. The complete
direct rank checks, rather than these additional sample comparisons,
prove the finite bridge. The six smaller witnesses are fully checked by
both original game algorithms and the direct transition verifier in
\texttt{src/\allowbreak{}odd\_\allowbreak{}endpoint\_\allowbreak{}small\_\allowbreak{}verify.\allowbreak{}py}; their complete proof data
and summary are
\texttt{research/\allowbreak{}odd-\allowbreak{}endpoint-\allowbreak{}small-\allowbreak{}certificates.\allowbreak{}json.\allowbreak{}gz} and
\texttt{research/\allowbreak{}odd-\allowbreak{}endpoint-\allowbreak{}small-\allowbreak{}verification.\allowbreak{}json}. The six
witnesses at orders 13 through 18 are now imported in
\texttt{SmallOddConstruction.\allowbreak{}lean}.
\texttt{OddEndpointFiniteCorollaries.\allowbreak{}lean} also covers orders 88, 89,
90, 94, 95 and 96 using the uniform local catalogue. All73 remaining
support orders are now joined in \texttt{OddSparseConstruction.\allowbreak{}lean} and
\texttt{OddExplicitAttainer.\allowbreak{}lean}, giving the maximum and a prescribed
graph at every support from13, with isolate padding to every larger
host. Original recipient-specific bounds remain attached to their
certificate statements.

\textbf{Proposition (small odd exceptions; exhaustive finite check).}
For \(6\le k\le11\), the general odd upper bound is not attained:
\[M_E(2k+1)\le k(k-3)+1.\] In particular, \(M_E(23)\le89\).

\textbf{Proof.} Equality in the upper bound \(k(k-3)+2\) would, by the
equality characterization, have precisely \(n=k+1\) nonisolated
vertices, one outdegree-one vertex \(u\), one distinct indegree-one
vertex \(v\), all other corresponding degrees two, and \(u\to v\)
absent. Add this missing arc. The resulting two-in/two-out digraph is a
bipartite two-regular graph when its sources and recipients are
represented in separate parts. Each of its bipartite cycles has even
length, so alternating edge colors decompose it into two permutation
factors \(A,B\). Choose \(A\) to contain \(u\to v\) and label
\(u=0,v=1\).

Up to relabeling fixing these two vertices, the \(A\)-cycle through zero
is \((0,1,\ldots,\ell-1)\) for one \(\ell\ge2\). Its other cycles are
represented by an integer partition of \(n-\ell\) into parts at least
two. For each such choice, enumerating all permutations \(B\) with
\(B(i)\ne i,A(i)\) covers every candidate graph; remove \(0\to1\) before
applying any game test.

Two necessary obstructions safely prune this enumeration. Once the
outgoing neighborhood of an eligible source \(s\) is complete, if some
\(c\ne s\) has \(N^+(s)\subseteq N^+[c]\) and there is an eligible
missing recipient for \(s\), the pursuer can wait at \(c\) and capture
the first nonwaiting messenger move. Once the incoming neighborhood of
an eligible recipient \(t\) is complete, if \(N^-(t)\subseteq N^+[c]\)
and there is an eligible source outside \(N^-(t)\cup\{t,c\}\), waiting
at \(c\) until the messenger first enters \(N^-(t)\) likewise permits
capture. In a partial assignment the cop's outgoing row may be
incomplete: coverage already present persists under later additions. The
source or recipient neighborhood used in the inclusion is always
complete. The code checks explicitly that the required eligible missing
partner exists.

All remaining completed graphs are evaluated by the exact finite
reachability game, with receipt priority and finite termination.
Enumeration finishes at every \(7\le n\le12\) and finds no attaining
graph. Its counts are:

{\def\LTcaptype{none} 
\Needspace{10\baselineskip}
\begin{longtable}[]{@{}
  >{\raggedleft\arraybackslash}p{(\linewidth - 6\tabcolsep) * \real{0.2500}}
  >{\raggedleft\arraybackslash}p{(\linewidth - 6\tabcolsep) * \real{0.2500}}
  >{\raggedleft\arraybackslash}p{(\linewidth - 6\tabcolsep) * \real{0.2500}}
  >{\raggedleft\arraybackslash}p{(\linewidth - 6\tabcolsep) * \real{0.2500}}@{}}
\toprule\noalign{}
\begin{minipage}[b]{\linewidth}\raggedleft
\(n\)
\end{minipage} & \begin{minipage}[b]{\linewidth}\raggedleft
Distinguished cycle partitions
\end{minipage} & \begin{minipage}[b]{\linewidth}\raggedleft
Partial nodes
\end{minipage} & \begin{minipage}[b]{\linewidth}\raggedleft
Completed graphs surviving local obstructions
\end{minipage} \\
\midrule\noalign{}
\endhead
\bottomrule\noalign{}
\endlastfoot
7 & 7 & 756 & 0 \\
8 & 11 & 6,338 & 21 \\
9 & 15 & 59,375 & 2,092 \\
10 & 22 & 800,507 & 64,992 \\
11 & 30 & 10,079,292 & 767,210 \\
12 & 42 & 164,276,910 & 16,902,091 \\
\end{longtable}
}

These finite exclusions and the equality characterization prove the
claimed strict bounds. \(\square\)

The complete enumerator is \texttt{src/\allowbreak{}odd\_\allowbreak{}endpoint\_\allowbreak{}exhaustive.\allowbreak{}cpp};
its finished receipts are
\texttt{research/\allowbreak{}odd-\allowbreak{}endpoint-\allowbreak{}small-\allowbreak{}exhaustive.\allowbreak{}json} and
\texttt{research/\allowbreak{}odd-\allowbreak{}endpoint-\allowbreak{}n12-\allowbreak{}exhaustive.\allowbreak{}json}. A known attaining
13-vertex graph passes every partial pruning test as a positive control.
The independent review
\texttt{research/\allowbreak{}odd-\allowbreak{}endpoint-\allowbreak{}twelve-\allowbreak{}independent-\allowbreak{}audit.\allowbreak{}md} checks the
permutation coverage, both pruning arguments, and the final game
semantics. It is an independent source and proof audit, not a separately
implemented repetition of the enumeration. This small-exception
exclusion remains outside Lean; the preceding odd endpoint attainment
has the partial formal coverage stated above. The proposition does not
determine the exact small odd maxima.

The general support reduction in Theorem 57 and the exact finite
exclusions in Theorem 53 supersede the earlier bespoke 23-arc reductions
and intermediate bounds. Their original arguments and receipts are
preserved in
\texttt{research/\allowbreak{}history/\allowbreak{}manuscript-\allowbreak{}superseded-\allowbreak{}arguments-\allowbreak{}2026-\allowbreak{}09-\allowbreak{}13.\allowbreak{}md}.
The following eleven-vertex core lemma underlies the retained
small-order exclusions.

\textbf{Lemma (eleven-vertex cores).} Every graph on 11 vertices with 22
or 23 arcs satisfies \(I(G)\le81\).

\textbf{Proof.} A vertex of outdegree at most one leaves at most ten
contributing source rows, each of size at most eight, giving \(I\le80\).
A low-indegree vertex gives the same column bound. We may therefore
assume minimum indegree and outdegree two. At 22 arcs all degrees equal
two. At 23 arcs exactly one vertex has outdegree three and exactly one
has indegree three.

If a source camping obstruction occurs, at least eight missing pairs
fail at 22 arcs, and at least seven fail at 23 arcs. If a recipient
camping obstruction occurs at \(t\) with pursuer position \(c\), every
source outside \(N^-(t)\cup\{t,c\}\) fails. This gives at least seven
failed pairs at 22 arcs and at least six at 23 arcs. The total numbers
of missing pairs are 88 and 87, respectively, so every graph with one of
these obstructions has \(I\le81\).

It remains to enumerate the graphs with neither obstruction. The 22-arc
case uses the two-permutation decomposition already described. For 23
arcs, the bipartite graph still has a perfect matching: a set \(S\) of
source vertices emits at least \(2|S|\) edges, whereas a set of possible
heads \(T\) can receive at most \(2|T|+1\). Thus \(|N(S)|\ge|S|\), which
implies a perfect matching by the augmenting-path criterion. Its removal
leaves residual degree one at every source except the unique
degree-three source, where two choices remain; recipient capacities are
analogous.

Fix the matching permutation to each cycle partition with no singleton.
Its cycle rotations and permutations of equal-length cycles allow the
degree-three source to be placed at one representative vertex for each
distinct cycle length. Enumerate all positions of the degree-three
recipient, then every residual row subject to the prescribed column
capacities. This covers every graph in question. Only complete source
rows or recipient columns are used in camping prunes; a partial cop row
suffices because any coverage already present persists under later
additions.

The complete 22-arc enumeration has 120,302 graphs surviving these
prunes, with maximum \(I=66\). The complete 23-arc enumeration has
792,080 such graphs, with maximum \(I=74\). Both use the exact finite
game recurrence and were freshly recompiled and replayed with all counts
and witnesses matching. Thus the graphs surviving the prunes also have
\(I\le81\), proving the lemma. \(\square\)

The source and mathematical audit, including the matching argument,
symmetry reduction and rollback checks, is
\texttt{research/\allowbreak{}odd-\allowbreak{}endpoint-\allowbreak{}eleven-\allowbreak{}enumeration-\allowbreak{}audit.\allowbreak{}md}. The
enumerators are \texttt{src/\allowbreak{}small\_\allowbreak{}exceptional\_\allowbreak{}regular\_\allowbreak{}max.\allowbreak{}cpp} and
\texttt{src/\allowbreak{}small\_\allowbreak{}exceptional\_\allowbreak{}one\_\allowbreak{}excess.\allowbreak{}cpp}; matching replay
receipts and hashes are indexed by
\texttt{research/\allowbreak{}odd-\allowbreak{}endpoint-\allowbreak{}eleven-\allowbreak{}enumeration-\allowbreak{}audit.\allowbreak{}json}. The saved
maximizing witnesses additionally pass both original game solvers and
all direct rank checks. This is independent review and re-execution of
the same enumeration sources, rather than a second independently
implemented enumeration.

Theorem 53 combines the exact value \(M(11,23)=78\), the shared
source/sink and positive-support bounds, and the three finite
exceptional-family exclusions to prove \(M_E(23)=83\). Its verified
twelve-vertex witness also proves \(M(n,23)=83\) for every \(n\ge12\)
after adding isolates.

\subsection{A packed sparse construction with optional
arcs}\label{a-packed-sparse-construction-with-optional-arcs}

Deleting a short, specified cluster of vertices from the six-type family
allows an explicit sparse construction at every order at least 115. Near
each deletion, only seven optional arrows must be suppressed. The local
strategies allow both of two adjacent columns as entries for final
delivery.

\textbf{Theorem 52 (with finite local certificates).} Let \(n\ge115\),
put \(k=\lceil n/6\rceil\) and \(r=6k-n\in\{0,\ldots,5\}\). For every
integer \[2n\le m\le20k-11r,\] there is an explicit universally winning
digraph with \(n\) vertices and \(m\) arcs. Consequently
\[M(n,m)=n(n-1)-m.\] The construction has \(2n\) compulsory arcs and
exactly \(8k-9r\) optional arcs; \textbf{every subset} of the optional
arcs works. Its range includes every \(2n\le m\le3n\), and every
delivery takes at most \[9\lceil n/6\rceil+21\] messenger moves.

\textbf{Proof.} We specify the graph, prove the exact arc count,
describe the finite strategy certificates and their complete scope, and
lift those strategies to every order in the theorem.

\subsubsection{The base graph and the seven
exclusions}\label{the-base-graph-and-the-seven-exclusions}

Begin with the six-type graph on
\((\mathbb Z/k\mathbb Z)\times\{0,\ldots,5\}\). Write \(j_i\) for
\((i,j)\). The first permutation is the Hamilton cycle
\[0_i\to1_i\to2_i\to3_i\to4_i\to5_i\to0_{i+1},\] and the second
permutation has arrows
\[0_i\to5_i,\quad1_i\to4_{i+1},\quad2_i\to1_{i+1},\quad
3_i\to0_{i-1},\quad4_i\to2_{i+1},\quad5_i\to3_{i-1}.\] Delete \(4_h\) in
the canonical columns \[S=\{0,3,\ldots,3(r-1)\},\] with
\(S=\varnothing\) when \(r=0\). Reconnect each permutation through each
deleted vertex. The changed arrows are
\[3_h\to5_h,\qquad1_{h-1}\to2_{h+1}.\] There are exactly \(n=6k-r\)
vertices and \(2n\) base arcs. Both indegree and outdegree remain two.
The displayed formulas, \(k\ge20\), and the gaps of at least three
between deleted columns show that the two outgoing arrows at each
surviving vertex are distinct and nonloop. In particular the smoothed
type 1 arrow has head of type 2 two columns ahead, distinct from its
type 2 head in the same column.

The available optional arrow types \((i,j)\to(i+d,t)\) are \[
(j,d,t)\in\{(1,-3,5),(3,-3,5),(4,-3,0),(4,3,1),
(2,3,0),(3,3,0),(3,-3,0),(5,-2,0)\}.
\] Omit arrows whose source was deleted. In addition, for each
\(h\in S\), omit precisely the following seven arrows, where a quadruple
\((\delta,j,d,t)\) denotes \((h+\delta,j)\to(h+\delta+d,t)\): \[
\begin{split}
\mathcal B=\{&(0,1,-3,5),(0,3,-3,5),\\
             &(1,1,-3,5),(1,3,-3,5),(1,3,3,0),\\
             &(2,3,-3,5),(2,4,-3,0)\}.
\end{split}
\] All column indices in this definition are modulo \(k\). These
optional arcs are distinct nonloop arcs outside the base graph, and
their target types are 0, 1, or 5, so their targets always survive.

Originally there are \(8k\) optional arrows. Each deletion removes the
two arrows with its absent type 4 source and the seven further arrows in
\(\mathcal B\). None of the seven has the deleted source. The
source-column sets \(\{h,h+1,h+2\}\) for different holes are disjoint,
including across the cyclic gap, so these losses do not overlap. Thus
the optional count is exactly \(8k-9r\). Moreover \[8k-9r-n=2k-8r\ge0\]
because \(k\ge20\) and \(r\le5\). The full envelope has
\[2n+(8k-9r)=20k-11r\] arcs and contains every budget through \(3n\).

\subsubsection{The finite local strategy
statements}\label{the-finite-local-strategy-statements}

Unwrap a window around a prospective destination column, now called
zero. Retain messenger columns \([-5,5]\) and cop columns \([-8,8]\).
Messenger exits are disallowed. A cop exit enters an outside state
\(\Omega\); from \(\Omega\) the cop may wait or reenter any existing
vertex in columns \(-8,-7,-6,6,7,8\).

For a local deletion crop \(H\subseteq[-8,8]\), construct the smoothed
base arrows literally. The messenger uses only these base arrows,
together with waiting. The cop receives all surviving optional arrows,
except the seven relative exclusions for each hole in \(H\), as well as
its base arrows and waiting. The finite certificates prove the following
statements for every crop that can arise from the canonical
construction:

\begin{itemize}
\tightlist
\item
  From every existing vertex in column \(-1\), the messenger can reach
  column zero within nine messenger moves and hand off a messenger-turn
  state in which it is distinct from the cop.
\item
  From every existing vertex in column \(-3\), the messenger can deliver
  to any specified existing vertex in column zero within 39 messenger
  moves.
\item
  From every existing vertex in column \(-2\), the corresponding
  delivery bound is 38 messenger moves.
\end{itemize}

Each assertion quantifies over every cop position other than the source,
including \(\Omega\). Receipt is terminal before capture; safe transfer
instead requires surviving the cop reply before the handoff.

Here the finite evidence consists of complete ranks and messenger
policies, not just successful simulations. There are 101 crops and 696
objectives: 101 safe-transfer objectives and 595 specified-target
delivery objectives. Every delivery objective checks both entry columns.
Ranks use the round convention of local simulation: the selected
messenger move decreases rank, every cop reply is safe and
nonincreasing, and zero messenger ranks have the required terminal
condition. These checked obligations prove the stated local bounds.

\subsubsection{\texorpdfstring{Why those crops cover every
\(k\ge20\)}{Why those crops cover every k\textbackslash ge20}}\label{why-those-crops-cover-every-kge20}

Enumerate \(k=20,\ldots,29\), \(r=0,\ldots,5\), and every choice of
center. Cropping the canonical hole sets to \([-8,8]\) gives exactly 101
distinct sets. For any \(k\ge29\), the gap from the last hole back to
the first is at least \[k-3(r-1)\ge k-12\ge17.\] A 17-column crop has
span 16, so it cannot contain holes from both sides of this gap. Its
holes form a contiguous portion of one arithmetic progression of step
three, containing at most five terms. Every such cropped progression
already occurs in the \(k=29\) enumeration, by taking the appropriate
number of holes and translating the center. No new local pattern can
therefore occur at larger \(k\).

This is a statement about the specified canonical cluster. It does not
assert that arbitrary deletion sets with minimum spacing three admit
these local strategies.

\subsubsection{Projection of every actual cop
move}\label{projection-of-every-actual-cop-move}

Smoothed base arrows have column jumps at most two and optional arrows
at most three. Thus the window-projection lemma applies with \(M=5\),
\(R=8\), \(J=3\) and \(k\ge20\). The preceding crop-coverage argument
supplies a certified arena at every center. Its remaining move
inclusions require the following check on omitted holes.

Omitting holes outside the crop is conservative. A hole immediately to
the right may change a represented boundary arrow from an outside head
to another outside head; both project to \(\Omega\). No other omitted
hole changes a represented-to-represented base arrow. An omitted hole to
the left can suppress optional arrows through its relative exclusions;
the cropped arena instead permits these extra cop moves. Messenger moves
inside \([-5,5]\) are modeled literally. Consequently the projection
neither grants the messenger an illegal move nor removes an actual cop
response, and it preserves collisions and both terminal conditions.

The certificates grant the cop every permitted optional arrow at once,
while the messenger uses only compulsory arrows. They therefore remain
valid after selecting any subset of the optional arrows.

\subsubsection{Global delivery and the time
bound}\label{global-delivery-and-the-time-bound}

For a recipient in column \(t\), the delivery entries \(t-3,t-2\) form a
consecutive block of length two. Every existing source type is covered,
and safe transfer advances one column. Cyclic strategy composition with
\(L=2\), \(d=1\), \(\tau=9\) and \(D=39\) bounds the total messenger
moves by \[9(k-2)+39=9k+21.\] Hence every missing ordered pair is
guaranteed. The universal missing-pair upper bound gives
\(M(n,m)=n(n-1)-m\), as claimed. \(\square\)

\textbf{Certificate boundary and reproducibility.} The complete finite
ranks and policies are saved in
\texttt{research/\allowbreak{}dense-\allowbreak{}phase2-\allowbreak{}sparse-\allowbreak{}two-\allowbreak{}entry-\allowbreak{}certificates.\allowbreak{}json.\allowbreak{}gz}.
The checker
\texttt{src/\allowbreak{}dense\_\allowbreak{}phase2\_\allowbreak{}sparse\_\allowbreak{}two\_\allowbreak{}entries.\allowbreak{}py\ -\allowbreak{}\/\allowbreak{}-\allowbreak{}check}
reconstructs every arena, checks exact coverage of the 101 crops and 696
objectives, and directly verifies every supplied rank, policy, and
required entry state. The reconstruction spells out all six smoothed
source-type maps; the independent generator uses explicit messenger-turn
and cop-turn states with an AND/OR priority queue, rather than the
earlier collapsed-round solver. The summary is
\texttt{research/\allowbreak{}dense-\allowbreak{}phase2-\allowbreak{}sparse-\allowbreak{}two-\allowbreak{}entry-\allowbreak{}check.\allowbreak{}json}. Literal
global graph and exact optional-count checks are recorded in
\texttt{research/\allowbreak{}dense-\allowbreak{}phase2-\allowbreak{}sparse-\allowbreak{}global-\allowbreak{}map.\allowbreak{}json}.

The smaller messenger window \([-4,4]\) does not establish all the new
column \(-3\) entry statements. Enlarging it to \([-5,5]\) is part of
this proof, and the failed smaller-window checks are retained in the
research log. The earlier broad exclusion buffers and their projection
proofs remain research history; they are unnecessary for the theorem
above. Negative certificate ranks make no assertion of global
impossibility or optimal delivery time. The graph construction,
canonical crop coverage, projection and global composition are also
instantiated in \texttt{LowBandConstruction.\allowbreak{}lean}, with all 101 crops
and 696 local objectives. Its exact optional-prefix attainers cover
\(n\ge115\) through the displayed endpoint \(20k-11r\). The finite
supplement below retains its separate proof boundary.

\subsubsection{Finite cyclic supplement below order
115}\label{finite-cyclic-supplement-below-order-115}

The same adjacency recipe gives a further computer-assisted exact
interval at orders \[
\mathcal S=\{60\}\cup[64,66]\cup[68,72]\cup[75,78]
\cup[80,90]\cup[92,114],
\] where the brackets here denote all integers in the displayed range.
For every \(n\in\mathcal S\), put \(k=\lceil n/6\rceil\) and \(r=6k-n\).
Every budget \[2n\le m\le20k-11r\] has the explicit graph above, and
every subset of its \(8k-9r\) optional arrows is universally winning.
These 47 finite envelopes have delivery bounds at most 72 messenger
moves. In particular, at \(n=114\) the full interval is 228 through 380,
with bound 66, and it includes the low-budget interval 229 through 341.

\textbf{Proof and finite evidence.} Retain the full cyclic vertex set
for both players; no window or exterior-state abstraction is used. Give
the messenger only the two compulsory permutations and waiting, while
giving the cop all compulsory and optional arrows and waiting. For each
recipient, the saved policy assigns to every nonterminal noncollision
messenger state a compulsory move. It either reaches the recipient
immediately, or every cop response is noncollision and strictly
decreases the nonnegative messenger rank. Induction on that rank proves
finite delivery. The policy remains legal in every graph between the
compulsory graph and its full optional envelope, and all possible cop
replies remain covered.

When there are no holes, simultaneous one-column rotation is an
automorphism of both the compulsory and optional arrow sets, so six
representative targets suffice. This identity is checked on both
complete arrow sets. Otherwise every literal recipient is certified. All
source/cop positions are checked for every representative target. The
optional count is the same displayed count as above; the finite records
separately check distinctness, absence of loops, and the complete
compulsory and optional lists. Taking the appropriate lexicographic
prefix gives each intermediate integer budget.

The 47 finite canonical certificates are contained in the single bundle
\texttt{research/\allowbreak{}sparse-\allowbreak{}interface-\allowbreak{}third-\allowbreak{}certificates.\allowbreak{}zip}, indexed by
\texttt{research/\allowbreak{}sparse-\allowbreak{}interface-\allowbreak{}third-\allowbreak{}proof-\allowbreak{}registry.\allowbreak{}json}. The
producer reconstructs and smooths the two permutation cycles. The
checker \texttt{src/\allowbreak{}sparse-\allowbreak{}interface-\allowbreak{}third-\allowbreak{}independent.\allowbreak{}py} imports no
project graph builder or solver: it instead reconstructs each
source-type arrow literally and directly checks the rank after the
messenger's chosen move against every actual cop reply. Its receipt is
\texttt{research/\allowbreak{}sparse-\allowbreak{}interface-\allowbreak{}third-\allowbreak{}independent.\allowbreak{}json}. Together
these are a finite proof of the displayed domain, not an extrapolation
to smaller or larger uncertified rings. Theorem 52 still supplies the
ordinary infinite continuation. Neither this supplement nor its finite
graph generation is presently formalized in Lean.

\textbf{The order 109 capacity repair.} At \(n=109\) (\(k=19,r=5\)),
additionally permit \[3_1\to0_4,\qquad3_2\to5_{18}.\] Both were among
the seven suppressed relative arrows near the hole at column zero. These
two fixed additions extend the optional count from 107 to 109 and
therefore the exact budget interval from 218--325 to 218--327, reaching
\(3n\). The same full cyclic rank argument was independently replayed
for every one of the 109 recipients, retaining the messenger's
compulsory arrows and giving the cop the entire enlarged optional set.
Its maximum rank is 72. The bundle contains this 48th certificate under
the variant name \texttt{repair109}; it is a specific finite extension,
not a claim that these two arrows may be restored at arbitrary orders.
Together with the preceding supplement, this gives every
\(2n\le m\le3n\) for all \(108\le n\le114\).

\subsection{A consecutive interval above the two-regular
endpoint}\label{a-consecutive-interval-above-the-two-regular-endpoint}

The five-type construction can tolerate a controlled set of additional
arcs. The proof gives the pursuer every additional move while
restricting the messenger to the original moves. This verifies all
intermediate graphs at once, without assuming that winning is monotone
under adding arcs.

\textbf{Theorem 36 (with finite local certificates).} Let \(n=5k\) with
\(k\ge10\). For every integer \[2n\le m\le\frac{11n}{5},\] there is an
explicit graph with \(n\) vertices and \(m\) arcs satisfying
\[M(n,m)=n(n-1)-m.\] Every pair can be guaranteed within \(8k+18\)
messenger moves.

\textbf{Proof.} Use the five-type graph \(H_k\) of Theorem 21 on
vertices \((i,j)\), where \(i\) is modulo \(k\) and \(0\le j\le4\). Its
two outgoing arcs are \[
\begin{aligned}
\alpha(i,j)&=(i,j+1) &&(j<4),&\alpha(i,4)&=(i+1,0),\\
\beta(i,0)&=(i+1,0),&\beta(i,1)&=(i-1,2),\\
\beta(i,2)&=(i-1,4),&\beta(i,3)&=(i-1,3),&\beta(i,4)&=(i-1,1).
\end{aligned}
\] There are \(2n\) arcs. Add any chosen subset of the \(k\) further
arcs \[ (i,0)\longrightarrow(i-1,4).\] They are distinct from the
original arcs and from each other, so every integer budget in the stated
interval is obtained.

We verify the stronger asymmetric game in which the messenger may use
only \(\alpha,\beta\) and waiting, while the pursuer may also use all
\(k\) additional arcs. Construct the same local arena as in Theorem 21:
messenger columns \(-3,\ldots,3\), pursuer columns \(-4,\ldots,4\), and
a pursuer state \(\Omega\) representing the outside. The messenger is
constrained to its local columns. The pursuer may leave into \(\Omega\),
wait there, and re-enter at any vertex in either boundary column \(-4\)
or \(4\). Within its columns it has the two original arcs and the
additional arc above whenever its type is zero.

Two independently implemented finite attractor computations verify the
following objectives against every legal initial pursuer state,
including \(\Omega\):

{\def\LTcaptype{none} 
\Needspace{10\baselineskip}
\begin{longtable}[]{@{}
  >{\raggedright\arraybackslash}p{(\linewidth - 4\tabcolsep) * \real{0.3333}}
  >{\raggedleft\arraybackslash}p{(\linewidth - 4\tabcolsep) * \real{0.3333}}
  >{\raggedleft\arraybackslash}p{(\linewidth - 4\tabcolsep) * \real{0.3333}}@{}}
\toprule\noalign{}
\begin{minipage}[b]{\linewidth}\raggedright
Objective
\end{minipage} & \begin{minipage}[b]{\linewidth}\raggedleft
Initial messenger column
\end{minipage} & \begin{minipage}[b]{\linewidth}\raggedleft
Maximum messenger moves
\end{minipage} \\
\midrule\noalign{}
\endhead
\bottomrule\noalign{}
\endlastfoot
Reach column zero safely at the beginning of a messenger turn & 1 & 8 \\
Deliver to \((0,0)\) & 2 & 26 \\
Deliver to \((0,1)\) & 2 & 19 \\
Deliver to \((0,2)\) & 2 & 17 \\
Deliver to \((0,3)\) & 2 & 17 \\
Deliver to \((0,4)\) & 2 & 18 \\
\end{longtable}
}

Each row covers all five messenger types and all 45 legal pursuer
positions per type: 225 required initial states. The complete ranks
verify every selected move and every pursuer response, with the distinct
handoff and receipt rules of local simulation.

The window-projection lemma applies with \(M=3\), \(R=4\), \(J=1\) and
\(k\ge10\). The messenger uses only original arrows and the local
pursuer receives every optional arrow, so local simulation applies to
every chosen subset. Cyclic strategy composition, with backward unit
transfers and delivery entry column two, gives at most
\[8(k-1)+26=8k+18.\] Thus every pair is guaranteed, and all \(n(n-1)-m\)
missing pairs are counted. The missing-pair upper bound proves equality.
\(\square\)

The primary generator and full ranks are
\texttt{src/\allowbreak{}five\_\allowbreak{}type\_\allowbreak{}sparse\_\allowbreak{}interval.\allowbreak{}py} and
\texttt{research/\allowbreak{}five-\allowbreak{}type-\allowbreak{}sparse-\allowbreak{}interval-\allowbreak{}certificates.\allowbreak{}json}. The
independent explicit alternating-state implementation is
\texttt{src/\allowbreak{}five\_\allowbreak{}type\_\allowbreak{}augmented\_\allowbreak{}verify.\allowbreak{}py}, with its full certificate
in \texttt{research/\allowbreak{}five-\allowbreak{}type-\allowbreak{}augmented-\allowbreak{}independent-\allowbreak{}certificates.\allowbreak{}json}.
The ring projection and concatenation are ordinary proofs; the bounded
local facts depend on these finite certificates.
\texttt{FiveShortRing10.\allowbreak{}lean} and \texttt{FiveShortRing11.\allowbreak{}lean} now
instantiate the optional envelopes at orders 50 and 55. The full uniform
family and its stated timing bound retain their ordinary proof boundary.
The theorem makes no claim for the small rings \(k=4,\ldots,9\) or for
budgets beyond \(11n/5\).

\textbf{Theorem 37 (with finite local certificates).} Let \(n=5k\) with
\(k\ge16\). At every integer budget \(2n\le m\le3n\), there is an
explicit graph attaining \[M(n,m)=n(n-1)-m.\] Every pair can be
guaranteed within \(8k+16\) messenger moves.

\textbf{Proof.} Keep the original five-type graph, and add any subset of
the \(n\) arcs \[ (i,j)\longrightarrow(i+3,0),\qquad 0\le j\le4.\] They
are pairwise distinct and are not original arcs. It suffices to win with
the messenger restricted to the original graph and the pursuer allowed
all the additional arcs.

Use messenger columns \([-3,3]\), pursuer columns \([-6,6]\), and the
outside state \(\Omega\). Give the represented pursuer its original and
added moves; moves leaving the interval enter \(\Omega\). From
\(\Omega\) permit waiting and re-entry at every vertex in the six
columns \(-6,-5,-4,4,5,6\). Independent synchronous and explicit
alternating-state certificates verify safe transfer from column one to
zero in at most eight moves, and delivery from column two to the five
recipient types at zero in at most \(24,16,16,17,18\) moves. Each
objective covers all 325 legal initial states: five messenger types and
65 pursuer positions per type. The transfer ends on a safe messenger
turn, and delivery has immediate terminal priority, exactly as in
Theorem 36.

Now the window-projection parameters are \(M=3\), \(R=6\), \(J=3\),
requiring \(k\ge2R+J+1=16\). This is the gap condition, not merely
distinctness of represented columns. The original messenger moves lift
literally and all optional pursuer moves are covered. Cyclic strategy
composition with backward unit transfers, delivery entry column two,
\(\tau=8\) and \(D=24\) gives \(8(k-1)+24=8k+16\). Local simulation
applies to every optional subset, so all pairs are guaranteed at each of
the \(n+1\) integer budgets. Counting missing pairs proves equality.
\(\square\)

The full primary ranks are
\texttt{research/\allowbreak{}five-\allowbreak{}type-\allowbreak{}full-\allowbreak{}sparse-\allowbreak{}interval-\allowbreak{}certificates.\allowbreak{}json},
generated by \texttt{src/\allowbreak{}five\_\allowbreak{}type\_\allowbreak{}sparse\_\allowbreak{}interval.\allowbreak{}py}. The
independent 4620-state implementation and certificate are
\texttt{src/\allowbreak{}five\_\allowbreak{}type\_\allowbreak{}long\_\allowbreak{}augmented\_\allowbreak{}verify.\allowbreak{}py} and
\texttt{research/\allowbreak{}five-\allowbreak{}type-\allowbreak{}long-\allowbreak{}augmented-\allowbreak{}independent-\allowbreak{}certificates.\allowbreak{}json}.
The statement is a uniform theorem with a finite local proof ingredient,
not an extrapolation from tested circumferences.
\texttt{FiveFullRing16.\allowbreak{}lean} through \texttt{FiveFullRing22.\allowbreak{}lean}
instantiate the finite envelopes at orders 80, 85, 90, 95, 100, 105 and
110. \texttt{FiveFinite65.\allowbreak{}lean}, \texttt{FiveFinite70.\allowbreak{}lean} and
\texttt{FiveFinite75.\allowbreak{}lean} separately supply three smaller full
envelopes. These selected instances do not formalize this entire uniform
family or its stated timing bound. For orders divisible by five with
\(n\ge80\), Theorems 34 and 37 combine into a consecutive explicit
interval starting at \(2n\) and reaching order \(n^{4/3}\).

\subsection{Explicit attainment at every budget in a sparse
interval}\label{explicit-attainment-at-every-budget-in-a-sparse-interval}

The circulant construction gives an exact joint extremum, not only a
lower bound on the vertex maximum: for every \(n\ge25\), Theorem 11 and
Lemma 12 already show that \(M(n,3n)=n(n-4)\). The next result fills
consecutive budgets above \(3n\). Its proof controls the extra pursuer
moves; arbitrary additions of arcs would not justify that conclusion.

\textbf{Theorem 34.} For every \(n\ge25\), \[M(n,3n)=n(n-4).\] Moreover,
if \(d\ge4\) and \[n>d^3+2d,\] then at every integer budget
\(3n\le m\le dn\) there is an explicit graph attaining
\[M(n,m)=n(n-1)-m.\] In particular, this holds throughout
\(3n\le m\le4n\) for every \(n\ge73\). More generally the explicit
interval reaches arc counts of order \(n^{4/3}\).

\textbf{Proof.} Say that a set \(A\subseteq\mathbb Z/n\mathbb Z\) has
distinct nonzero differences when
\[a-b=a'-b'\ne0,\qquad a,b,a',b'\in A,\] implies \((a,b)=(a',b')\). For
such a set, two distinct translates \(x+A\) and \(y+A\) intersect in at
most one point: two intersections would give two representations of the
nonzero difference \(y-x\).

Choose \(h\) as in Lemma 12 and put \(a=(h-1)/2\). The set
\[A_0=\{0,1,a,h\}\] has distinct nonzero differences by the proof of
Theorem 11. The base graph \(B\) has steps \(1,a,h\). For any larger set
\(A\) with distinct nonzero differences, let \(H_A\) have the nonzero
steps in \(A\); we will keep \(0\in A\) throughout. Every graph \(G\)
satisfying \[B\subseteq G\subseteq H_A\] has distinct closed move sets
intersecting in at most one vertex, because \(C_G(x)\subseteq x+A\).

The proof of Theorem 11 therefore works inside every such \(G\). To make
this preservation explicit, fix the recipient zero. Its direct
predecessors \(-1,-a,-h\) and zero are guaranteed, and \(-a-1\) is
guaranteed using the two successors \(-a,-1\). If
\[Q(z)=\{z,z-1,z-a,z-a-1\}\] is guaranteed, the following four
applications of Lemma 10 establish \(Q(z-h)\):

{\def\LTcaptype{none} 
\Needspace{8\baselineskip}
\begin{longtable}[]{@{}ll@{}}
\toprule\noalign{}
New source & Two guaranteed successors \\
\midrule\noalign{}
\endhead
\bottomrule\noalign{}
\endlastfoot
\(z-h\) & \(z-a-1,\ z\) \\
\(z-h-1\) & \(z-h,\ z-1\) \\
\(z-h-a\) & \(z-h,\ z-a\) \\
\(z-h-a-1\) & \(z-h-a,\ z-a-1\) \\
\end{longtable}
}

All displayed moves belong to \(B\). Coprimality of \(h\) and \(n\)
makes finitely many repetitions cover all sources. The same argument
with all coordinates shifted applies to every recipient, even when \(G\)
itself is not translation invariant. Thus every intermediate graph is
universally winning.

It remains to enlarge \(A_0\). Suppose \(|A|=k\) and its nonzero
differences are distinct. To adjoin \(x\), exclude the following
residues:

\begin{enumerate}
\def\labelenumi{\arabic{enumi}.}
\tightlist
\item
  The \(k\) members of \(A\) itself.
\item
  Residues \(x=a+b-c\) with \(a,b,c\in A\) and \(b\ne c\). These exclude
  every equality between an old nonzero difference and either sign of a
  new difference. There are at most \(k^2(k-1)\) such residues.
\item
  Solutions of \(2x=a+b\) with \(a,b\in A\). These exclude an equality
  \(x-a=b-x\) between the two signs of new differences. There are
  \(k(k+1)/2\) unordered pairs with repetition, and each congruence has
  at most two solutions, hence at most \(k(k+1)\) excluded residues.
\end{enumerate}

Within each sign the new differences are already distinct. The three
exclusions are therefore sufficient, and together rule out at most
\[k+k^2(k-1)+k(k+1)=k^3+2k\] residues. Whenever \(n>k^3+2k\), at least
one admissible residue remains; choose the smallest representative.
Under the theorem's hypothesis this process continues from size four to
size \(d+1\).

The resulting \(H_A\) has \(dn\) arcs and contains the \(3n\) arcs of
\(B\). For a prescribed integer \(m\) in the interval, add exactly
\(m-3n\) of the remaining arcs in lexicographic order. The intermediate
graph is universally winning by the preservation argument, so all
\(n(n-1)-m\) missing pairs are guaranteed. This matches the elementary
upper bound.

For \(d=4\) the hypothesis is \(n>72\). Taking the largest admissible
\(d\) gives \(d\) of order \(n^{1/3}\), and therefore an upper interval
endpoint of order \(n^{4/3}\). The endpoint \(m=3n\) only needs the
original \(n\ge25\) construction, with no extension. \(\square\)

This is a deterministic construction at each stated integer budget. The
underlying delivery argument may require many moves; it is not a
two-move theorem. Its ordinary proof is independent of the finite tests
in \texttt{src/\allowbreak{}sidon\_\allowbreak{}main\_\allowbreak{}interval.\allowbreak{}py}, and no new Lean coverage is
asserted. It provides an initial consecutive part of Raychev's
full-density interval. The stronger constructions and exact-budget
argument in Theorems 38--48 subsequently cover the entire interval at
every sufficiently large order.

\subsection{An explicit algebraic family at intermediate
density}\label{an-explicit-algebraic-family-at-intermediate-density}

The envelope in this construction uses the classical finite-field
parabola Sidon set, the distinct-difference ingredient of the
Erdős--Turán construction (Erdős and Turán 1941); see also Eberhard and
Manners (Eberhard and Manners 2023, Construction 1). That set and its
intersection property are established ingredients. The proof below
supplies the four-generator delivery core and its guarantee throughout
the envelope.

\textbf{Theorem 35.} Let \(p\ge5\) be a prime. For every integer
\[4p^2\le m\le p^2(p-1),\] there is an explicit graph attaining \[
M(p^2,m)=p^2(p^2-1)-m.
\] In particular, the interval extends to \(m=n^{3/2}-n\) at orders
\(n=p^2\). For primes \(p\ge11\), Theorem 34 extends its lower endpoint
to \(3p^2\).

\textbf{Proof.} Put \(A=\{(t,t^2):t\in\mathbb F_p\}\), and let \(H_A\)
be the directed graph with all nonzero moves in \(A\). Let \(B_0\) be
its four-move subgraph with parameters \[t\in\{1,-1,-2,-3\}.\] These are
distinct and nonzero when \(p\ge5\). We prove universal delivery in
every graph \(G\) satisfying \(B_0\subseteq G\subseteq H_A\), including
graphs without translation symmetry.

Every closed outgoing neighborhood in \(G\) is a subset of the
corresponding translate of \(A\). Distinct translates of \(A\) intersect
in at most one point. Indeed, a nonzero difference \[
(t,t^2)-(s,s^2)=(t-s,(t-s)(t+s))
\] determines \(t,s\) uniquely: if the first coordinate is nonzero, it
determines their difference and sum, and division by two is valid; if
the first coordinate is zero, the entire difference is zero. Thus every
nonzero difference has at most one ordered representation.

This intersection property gives the following safe-branching rule. Fix
a recipient. If two distinct outgoing neighbors of a vertex are already
universally winning positions for that recipient, the vertex is
universally winning as well. Against an initial cop at another vertex,
at most one of the two successors lies in the cop's closed outgoing
neighborhood. Move to the other, survive the cop's response, and use its
winning strategy. The rule concerns the actual game, including the
intervening pursuer move.

We first use recipient \(0\) and only moves of the uniform core \(B_0\).
Put \[u=(1,1),\qquad v=(-2,4),\qquad h=(-3,9)=u+2v,\] and define
\[Q(z)=\{z,z-u,z-v,z-u-v\}.\] The vertices \(0,-u,-v\) are initially
winning, by immediate receipt or a direct final move. The remaining
vertex \(-u-v\) has winning successors \(-u,-v\), so all of \(Q(0)\) is
winning.

If every vertex of \(Q(z)\) is winning, so is every vertex of
\(Q(z-h)\). Apply the safe-branching rule in the following order; each
row lists the new vertex and two already winning successors: \[
\begin{array}{c|cc}
z-h & z-u-v & z\\
z-h-u & z-h & z-u\\
z-h-v & z-h & z-v\\
z-h-u-v & z-h-v & z-u-v
\end{array}
\] The corresponding pairs of moves are \((v,h),(u,h),(v,h),(u,h)\),
respectively. They are distinct allowed moves.

Repeated propagation makes the entire line
\(L=\{\lambda h:\lambda\in\mathbb F_p\}\) winning: the nonzero vector
\(h\) has additive order \(p\), so integer multiples cover all the
displayed scalars.

Use the fourth core move \(w=(-1,1)\). We have \[w-v=(1,-3)=-h/3\in L,\]
while \(v\notin L\) because \(\det(h,v)=6\ne0\). If every vertex of a
coset \(L+z\) is winning, then every vertex of \(L+z-v\) has two
distinct winning successors, obtained by adding \(v\) and \(w\). The
safe-branching rule therefore fills that entire coset as well. Iterating
gives all cosets \(L-jv\), which are the \(p\) distinct cosets because
\(v\notin L\). Thus every vertex is winning.

The same proof shifted by the recipient uses only moves of \(B_0\) and
the global closed-neighborhood intersection bound. It applies to every
recipient even if the extra arcs of \(G\) are not translation invariant.
Hence every intermediate graph is universally winning.

The core has \(4p^2\) arcs, and \(H_A\) has \(p^2(p-1)\) arcs. Add the
remaining arcs one at a time in any fixed order until the prescribed
integer budget \(m\) is reached. All \(p^2(p^2-1)-m\) missing nonloop
pairs are guaranteed, attaining the elementary upper bound. Finally,
when \(p\ge11\) we have \(p^2\ge73\), so Theorem 34 covers the adjacent
interval \(3p^2\le m\le4p^2\). \(\square\)

The proof is finite and constructive: initializing a quartet takes at
most two messenger moves, and each quartet propagation increases a bound
by at most four. At most \(p-1\) such propagations and \(p-1\) further
coset steps therefore give the uniform, unoptimized bound \(5p-3\)
messenger moves.

The prime hypothesis is intentional. In an extension field, repeated
addition only covers scalars from its prime subfield, so this four-move
propagation does not justify the same conclusion for arbitrary prime
powers. The construction is checked at \(p=5,7,11\) by
\texttt{src/\allowbreak{}parabola\_\allowbreak{}cayley.\allowbreak{}py}; it verifies the full-parabola
intersection condition, the four-move core's complete safe-branching
closure, and the core/upper-graph inclusions for all 827 consecutive
integer budgets in these examples. At \(p=5\), both full game solvers
and their rank certificates independently verify every recipient. These
finite checks supplement the ordinary proof and are recorded in
\texttt{research/\allowbreak{}parabola-\allowbreak{}cayley-\allowbreak{}checks.\allowbreak{}json}. No new Lean formalization
is claimed.

\section{Protected routing and the main
interval}\label{protected-routing-and-the-main-interval}

\subsection{Protecting the moves used by the
strategy}\label{protecting-the-moves-used-by-the-strategy}

The distinct-difference condition in Theorem 34 is stronger than the
delivery strategy needs. Its branching steps only compare pairs of the
three original outgoing moves. Preserving those comparisons permits a
much larger set of additional moves.

\textbf{Theorem 38.} For every \(n\ge25\) and every integer
\[3n\le m\le n\left(\left\lceil\frac n7\right\rceil-1\right),\] there is
an explicit graph attaining \[M(n,m)=n(n-1)-m.\] The strategy uses only
three outgoing moves per vertex; the remaining arcs may be added
individually up to the stated budget.

\textbf{Proof.} Choose \(h\) using Lemma 12, put \(a=(h-1)/2\), and set
\[S=\{1,a,h\},\qquad A_0=\{0,1,a,h\}\subseteq\mathbb Z/n\mathbb Z.\] The
proof of Theorem 11 shows that \(A_0\) has distinct nonzero differences.
In particular, each of the six ordered differences
\[D=\{u-v:u,v\in S,\ u\ne v\}=\{\pm(a-1),\pm(a+1),\pm2a\}\] has just its
original ordered representation in \(A_0\). We only preserve uniqueness
for these six differences.

Starting from \(A=A_0\), adjoin any residue outside
\[A\cup(A+D),\qquad A+D=\{x+d:x\in A,\ d\in D\}.\] No new ordered pair
involving the adjoined residue has a difference in \(D\), because
\(D=-D\). Thus the six protected differences retain their original
unique representations. When \(|A|=k\), at most \(7k\) residues are
excluded. If \(k<\lceil n/7\rceil\), then \(7k<n\), so the extension
remains possible. Choosing the smallest available representative
constructs a set of size at least \(\lceil n/7\rceil\). For
\(25\le n\le28\), the initial set already has that size.

Let \(B\) be the uniform three-move graph with steps in \(S\), and let
\(H_A\) have all nonzero steps in \(A\). Consider any graph
\[B\subseteq G\subseteq H_A.\] For a source vertex \(r\) and distinct
core moves \(u,v\in S\), no different cop vertex \(c\) can have both
\(r+u\) and \(r+v\) in its closed outgoing neighborhood. Indeed, that
neighborhood is contained in \(c+A\). If both successors belonged to it,
their corresponding elements \(x,y\in A\) would satisfy
\(x-y=u-v\in D\). Protected uniqueness forces \((x,y)=(u,v)\) and
therefore \(c=r\), a contradiction.

Consequently, whenever two successors obtained using core moves are
universally winning for a fixed recipient, their source is universally
winning: at least one of those successors is inaccessible to the initial
cop and survives its response. This argument does not assert that
arbitrary pairs of actual outgoing successors have the same property.

The entire proof of Theorem 11 uses exactly these protected pairs. For
recipient zero, the direct predecessors \(-1,-a,-h\) are winning, and
\(-a-1\) uses the core moves \(1,a\). If \[Q(z)=\{z,z-1,z-a,z-a-1\}\] is
winning, the four steps establishing \(Q(z-h)\) use the move pairs
\[(a,h),\quad(1,h),\quad(a,h),\quad(1,h),\] respectively, in the order
given in Theorem 34. Hence that quartet propagates in \(G\). Since
\(\gcd(h,n)=1\), its successive translates cover every source. Shifting
the proof by the recipient covers every target, regardless of whether
the added arcs are translation invariant.

Thus every graph between \(B\) and \(H_A\) is universally winning. The
first has \(3n\) arcs and the latter has at least
\(n(\lceil n/7\rceil-1)\). Add the extra arcs in lexicographic order
until the prescribed integer budget is reached. All \(n(n-1)-m\) missing
pairs are guaranteed, meeting the elementary upper bound. \(\square\)

The improvement from order \(n^{4/3}\) to order \(n^2\) comes from
preserving only the distinctions actually used by the strategy. Other
differences may repeat many times. The graph may also have larger
closed-neighborhood intersections; the proof protects its particular
pairs of core successors instead of imposing a global intersection
bound.

For example, at \(n=29\) the construction gives core moves \(1,3,7\) and
the envelope \(A=\{0,1,3,7,8\}\). Its protected differences are
\(\pm2,\pm4,\pm6\), each still uniquely represented. The unprotected
difference \(-1\) occurs twice, as \(0-1=7-8\). Thus this envelope
genuinely fails the full distinct-difference condition, yet every
intermediate graph with \(87\le m\le116\) is universally winning.

\textbf{Corollary 39.} There is an absolute \(N\) such that, for every
\(n\ge N\) and every integer \[3n\le m\le n(n-3),\] we have
\[M(n,m)=n(n-1)-m.\] For every sufficiently large multiple of five
\(n\), the same formula holds at every integer budget
\(2n\le m\le n(n-3)\).

\textbf{Proof.} The upper endpoint in Theorem 38 satisfies
\[n(\lceil n/7\rceil-1)\ge n^2/7-n.\] Since
\(\sqrt{n/\log n}\to\infty\), this endpoint eventually exceeds
\(\lceil10n^{3/2}\sqrt{\log n}\rceil\), the lower endpoint in Corollary
26. The two integer intervals overlap. Theorem 38 covers the lower one
and Corollary 26 covers all remaining budgets through \(n(n-3)\),
including its dense end via Theorem 20, now a consequence of the direct
dense envelopes in Theorems 48 and 58. Taking a common sufficiently
large threshold establishes the first statement.

For multiples of five, Theorem 37 supplies every budget from \(2n\)
through \(3n\) once \(n\ge80\). Combining the intervals proves the
second statement. \(\square\)

The later Theorems 44 and 59 replace this eventual overlap by a fully
explicit construction for every \(n\ge115\), including the sparse range.
Theorem 38 retains its smaller-order domain, arbitrary-subset guarantee,
and protected-difference strategy. The probabilistic overlap in
Corollary 39 is included only as a consequence of its underlying
constructions.

The generator \texttt{src/\allowbreak{}main\_\allowbreak{}interval\_\allowbreak{}resume.\allowbreak{}py} checks protected
differences at eight orders from 25 through 1000. It checks every
protected source-pair/cop triple at the seven listed orders through 200,
and verifies 24 selected recipient games with both independent solvers
and the explicit rank-certificate checker. Those tests include
intermediate graphs lacking translation symmetry. The receipt
\texttt{research/\allowbreak{}main-\allowbreak{}interval-\allowbreak{}resume-\allowbreak{}checks.\allowbreak{}json} records their exact
scope; the parameterized proof above is independent of these finite
checks.

\subsection{Protected routing sets at every exact
budget}\label{protected-routing-sets-at-every-exact-budget}

The random construction can be strengthened by reserving a small set of
routing vertices. Only arrows incident with that set are used to certify
delivery. All arrows outside the set can then be chosen freely,
including adversarially, without invalidating any certified route.

\textbf{Theorem 54 (protected routing).} Let \(n\ge4\), let \(H\) be a
specified set of \(b\le n\) vertices, and put \[
L=n(n-1),\qquad Q=b(2n-b-1).
\] Choose an integer \(0<a<Q\) and write \(p=a/Q\). If \[
(b-3)p^2(1-p)>4\log n+\log(8/3), \tag{R}
\] there is a deterministic algorithm, polynomial in \(n\), choosing
exactly \(a\) arrows with at least one endpoint in \(H\) such that
\textbf{every assignment of the arrows between vertices outside \(H\)}
guarantees every source--recipient pair within two messenger moves. In
particular, at every integer budget \[
a\le m\le a+(n-b)(n-b-1),
\] an explicit algorithm attains \(M(n,m)=L-m\).

Separately, for every real \(0<p<1\), without condition (R) or an
integer exact-budget requirement, independently including each of the
\(Q\) incident arrows with probability \(p\) fails to give such a robust
routing set with probability at most \[
n(n-1)^2\bigl(1-p^2(1-p)\bigr)^{b-3}.
\] Thus the theorem includes both the deterministic exact-budget and the
random-construction conclusions.

\textbf{Proof.} Fix \(s\ne t\) and an initial pursuer position
\(c\ne s\). Each \(u\in H\setminus\{s,t,c\}\) certifies a safe two-move
route when \[
s\to u,\qquad u\to t\quad\hbox{are present},\qquad c\to u
\quad\hbox{is absent}.
\] The three variables are distinct, even when \(c=t\). All are incident
with \(H\). For different candidates the sets of three variables are
disjoint: an equality between an incoming-to-candidate variable and an
outgoing-from-candidate variable would force a candidate to equal one of
\(s,t,c\), which is excluded.

These witness tests are specified adjacency-extension patterns, of the
kind studied for existentially closed tournaments (Bonato et al. 2010).
The probabilistic and conditional-expectation methods used here are
standard (Alon and Spencer 2016). We give the details because the
protected routing set, exact incident-arrow count, and arbitrary
exterior arrows are part of the assertion.

There are at least \(b-3\) candidates, each with probability
\(p^2(1-p)\) of success. There are exactly \(n(n-1)^2\) legal triples.
The union bound proves the random assertion. If \(Z\) counts bad
triples, it also gives \[
\mathbb E Z\le n(n-1)^2e^{-(b-3)p^2(1-p)}
<\frac{3}{8n}
\] under (R).

Let \(X\) be the number of chosen incident arrows. Then
\(X\sim\operatorname{Bin}(Q,a/Q)\) and \(a\) is a mode. Since
\(\operatorname{Var}(X)\le Q/4\), Chebyshev's inequality gives
\(\Pr(|X-a|<\sqrt Q)\ge3/4\). Here \(Q\le n(n-1)\), so this interval
contains at most \(2n-1\) integers. Modal maximality implies \[
\Pr(X=a)\ge\frac{3}{4(2n-1)}\ge\frac{3}{8n}.
\] Therefore \(\mathbb E[Z\mid X=a]<1\).

For completeness, the conditional-expectation algorithm remains exact
when restricted to these \(Q\) variables. After a partial assignment,
let \(\ell\) be the number of undecided incident arrows and \(d\) the
number of additional present arrows required. The total number of
completions is \(C=\binom{\ell}{d}\). For a fixed legal triple, each
candidate contributes a factor \((1+x)^q\) if its witness pattern is
already contradicted; otherwise its factor is \((1+x)^q-x^j\), where
\(q\) variables remain undecided and \(j\) of them require presence.
Multiply these factors, multiply by \((1+x)^{q_0}\) for unused
variables, and extract the coefficient of \(x^d\). This counts exactly
the completions in which the triple is bad. Summing over triples gives
an integer \(B\) with initial \(B<C\).

For either decision on the next arrow, child counts satisfy
\(B=B_0+B_1\) and \(C=C_0+C_1\). Hence some feasible child retains
\(B_i<C_i\). At the end \(C=1\) and \(B=0\). There are \(O(n^2)\)
decisions and \(O(n^3)\) triples per decision, with polynomial-size
polynomial products and \(O(n^2+\log n)\)-bit counts. Straightforward
exact integer arithmetic is therefore polynomial in \(n\).

For every legal triple a selected candidate differs from the pursuer and
both endpoints. Moving to it is safe, and the pursuer cannot reach it on
its response. The messenger then reaches the recipient; receipt precedes
collision. No outside-\(H\) arrow appears in these conditions, so the
same strategies work for every assignment of those arrows. There are
\((n-b)(n-b-1)\) of them. Choosing any prescribed number proves the
entire displayed exact-budget interval. \(\square\)

Setting \(H=V\) gives \(Q=L\) and the former exact-budget density
theorem at every original order and budget, with its polynomial
algorithm and two-move guarantee. A smaller routing set adds genuine
robustness and extends the low-density range; it does not merely combine
old intervals.

\subsubsection{An explicit interval from oriented routing
pairs}\label{an-explicit-interval-from-oriented-routing-pairs}

\protect\phantomsection\label{the-scale-of-the-first-universal-two-move-budget}
The separate two-move sparsity question and its bounds are in the
\href{https://angelraychev.com/cop-prisoner/delivery-time/\#the-scale-of-the-first-universal-two-move-budget}{delivery-time
companion}. The exact-budget routing constructions remain here.

\textbf{Orientation specialization.} Fix a routing set \(H\) of size
\(b\), and let \(Q=b(2n-b-1)\). If \[n(n-1)^2(7/8)^{b-3}<1,\] there is a
deterministic polynomial construction using exactly \(Q/2\) incident
arrows that preserves universal two-move delivery under every choice of
exterior arrows. Thus all budgets from \(Q/2\) to \(L-Q/2\) work. A
sharper sufficient condition is
\[n(n-1)(n-2)(7/8)^{b-3}+n(n-1)(3/4)^{b-2}<1.\]

To prove this, independently orient each unordered pair touching \(H\),
with its two directions equally likely. The incident count is exactly
\(Q/2\) for every outcome. For distinct \(s,t,c\), a candidate
\(u\in H\setminus\{s,t,c\}\) is a safe intermediary when \(s\to u\),
\(u\to t\), and \(u\to c\), with probability \(1/8\). The last
orientation excludes \(c\to u\). Different candidates use independent
orientation variables. When \(c=t\), only two variables are needed and
the probability is \(1/4\). Counting the two kinds of triples gives the
sharper displayed failure bound; the coarser one also follows. Exterior
arrows occur in none of these tests.

For exact derandomization, fix orientations in sequence. For a given
triple and candidate with \(j\) free relevant variables, there are
\(2^j\) bad completions if a fixed variable already prevents the
witness, and \(2^j-1\) otherwise. Multiply over candidates, whose
variable sets are disjoint, and account for irrelevant free variables.
Summing these integer counts over triples gives the total number of bad
completions, with multiplicity. Initially it is smaller than the number
of all completions. At each decision choose a branch preserving that
strict inequality; the two child counts sum to their parent counts. At
the end there is one completion and zero failures. There are
polynomially many operations on integers of polynomial bit length. This
specialization replaces coefficient extraction by products, while the
full density and sampling assertions of Theorem 54 remain available.

\textbf{Corollary 55.} Set \[
b_n=\lceil24\log n+4\rceil,\qquad
A_n=\frac{b_n(2n-b_n-1)}2.
\] Whenever \(b_n\le n\), every integer budget \(A_n\le m\le L-A_n\)
admits universal two-move delivery by a deterministic polynomial
algorithm, with arbitrary exterior-arrow subsets as in Theorem 54. For
every \(n\ge4096\), this extends to every integer budget \[
A_n\le m\le n(n-3).
\]

\textbf{Proof of the constructions.} Use the orientation specialization.
Since \((7/8)^x\le e^{-x/8}\) for \(x\ge0\),
\[n(n-1)^2(7/8)^{b_n-3}<n^3 e^{-(b_n-3)/8}
\le e^{-1/8}<1.\] It gives the exact robust interval \([A_n,L-A_n]\).

For \(n\ge4096\), \[
b_n\le24\log n+5\le32\log n+12\le(n-1)/13.
\] At \(4096\), \(\log4096=12\log2<9\) makes the left side of the second
inequality less than \(300\), while its right side equals \(315\). The
difference increases thereafter because \(1/13-32/n>0\). The elementary
inequality \(\log2<3/4\) follows, for example, from
\(1+3/4+(3/4)^2/2>2\).

Consequently \(A_n\le nb_n\le n(n-1)/13\). The existing dense
constructions supply every missing budget from \(2n\) through
\(n\lceil(n-1)/13\rceil\) for \(n\ge32\), in two moves (Theorems 48 and
58). Their interval overlaps \([A_n,L-A_n]\), proving the claimed
extension through \(L-2n=n(n-3)\). \(\square\)

\textbf{Verification boundary.} The proofs above are ordinary
mathematical arguments. The restricted exact-count implementation and
its independent finite checks are recorded in
research/routing-unification-checks.json; the independent proof review
is research/routing-independent-audit.md. No large instance of the
polynomial conditional-expectation algorithm is claimed.
\texttt{ProtectedRouting.\allowbreak{}lean} verifies preservation of safe routing
witnesses under arbitrary exterior-arrow changes and their implication
for actual terminating delivery. The probability estimates and counting
algorithm remain ordinary proofs.

\subsection{Composition through incomparable hub
sets}\label{composition-through-incomparable-hub-sets}

The construction below reduces a delivery problem on an arbitrarily
large outside set to a family of subsets of a small core. The subsets
specify which hubs an outside vertex cannot reach. Incomparability makes
one outside vertex's inaccessible hubs available to another vertex. A
two-part core handles the remaining configurations, including recipients
inside the core. This produces an entire robust interval of exact
budgets.

The antichain condition is the classical Sperner-family condition
(Sperner 1928); equal-size distinct subsets are automatically
incomparable. Its role here is to separate source and pursuer access to
the hubs. The cyclic restrictions and two-part core below supply the
additional conditions needed for every recipient and pursuer position.

\textbf{Theorem 59 (antichain composition).} Let \(k\ge5\). Take two
sets of hubs \(P_i,L_i\), indexed modulo \(k\). Their missing arrows are
\[P_i\to L_i,L_{i-1},\qquad L_i\to P_{i+4},P_{i+6};\] all other distinct
hub pairs are arrows. For each of \(q\) outside vertices \(x\), choose a
set \(K_x\) of hubs subject to these conditions:

\begin{itemize}
\tightlist
\item
  \(|K_x|\ge2\), and \(K_x\) lies entirely in one of the two hub parts;
\item
  a set in the \(P\) part contains no two indices differing by \(\pm2\);
  a set in the \(L\) part contains no adjacent indices;
\item
  a two-element code contains no pair at either cyclic difference
  \(\pm1\) or \(\pm2\);
\item
  the sets \(K_x\) form an antichain: no one contains another.
\end{itemize}

Include every hub-to-outside arrow, and include \(x\to u\) for precisely
the hubs \(u\notin K_x\). Every subset of the arrows between distinct
outside vertices may then be included. Every resulting graph has
universal delivery within two messenger moves.

Writing \(n=2k+q\) and \(W=\sum_x|K_x|\), this gives every exact budget
\[A\le m\le A+q(q-1),\qquad A=2k(2k-3)+4kq-W.\] At the upper endpoint
the number of missing arrows is \(4k+W\). The assertion holds for every
optional subset, not merely a specially chosen subset of each size.

\textbf{Proof.} Let \(F\) denote the missing graph inside the hubs. For
a missing pair \((s,t)\) and a pursuer at \(c\ne s\), we seek a hub
\(u\) with arrows \(s\to u\to t\), with \(u\ne c\) and no arrow
\(c\to u\). The messenger can move to \(u\), survive every pursuer
response, and then reach \(t\). The case \(c=t\) is included.

First suppose all three vertices are hubs. For a missing hub pair, the
set \(D=F^+(s)\cup F^-(t)\) contains two points and two lines. If
\(s=P_i\), the line indices are \(i,i-1\) and the point indices are
either \(i,i+1\) or \(i-1,i\). A point pursuer has its entire missing
row in \(D\) only when it is \(s\); a line pursuer's point pair has
difference two and cannot fit inside the adjacent point pair. If
\(s=L_i\), the point indices are \(i+4,i+6\) and the line indices are
either \(i,i-2\) or \(i,i+2\). The same argument with the two parts
exchanged applies. Thus \(F^+(c)\) contains a hub outside \(D\), which
is the desired intermediary.

Now consider the other nondirect configurations.

\begin{enumerate}
\def\labelenumi{\arabic{enumi}.}
\item
  \textbf{Both endpoints outside.} If the pursuer is outside, choose
  \(u\in K_c\setminus K_s\), which exists by incomparability. It is
  reachable from \(s\), inaccessible from \(c\), and reaches \(t\). If
  the pursuer is a hub, its missing row consists of an adjacent line
  pair or a point pair of difference two. The restrictions on \(K_s\)
  prevent that entire row from lying in \(K_s\). A remaining member
  works.
\item
  \textbf{Source outside, recipient a hub.} Missingness means
  \(t\in K_s\). For an outside pursuer whose code lies in the same part,
  choose \(u\in K_c\setminus K_s\); that part is a clique, so
  \(u\to t\). For an outside pursuer using the other part, \(K_c\) is
  disjoint from \(K_s\). Its code cannot fit inside \(F^-(t)\): a larger
  code has at least three members, and a pair code has neither of that
  row's possible differences, one and two. Thus a member reaches \(t\).
  For a hub pursuer whose missing row lies in the part containing \(t\),
  the code restriction leaves a member outside \(K_s\), and the clique
  supplies the last arrow. For the other hub pursuer part, both missing
  neighbors are reachable from \(s\). Their adjacent or difference-two
  pair cannot equal the opposite kind of pair \(F^-(t)\), so at least
  one reaches \(t\).
\item
  \textbf{Both endpoints hubs, pursuer outside.} The set \(D\) above
  contains only two members in either part, with cyclic difference one
  or two. A larger code cannot fit in that part, and a pair code
  excludes both differences. Hence \(K_c\setminus D\) contains an
  intermediary.
\end{enumerate}

A hub source reaches an outside recipient directly. All intermediaries
used in the proof are hubs, so optional arrows between outside vertices
cannot invalidate a witness. Counting the compulsory and optional arrows
gives the stated interval. \(\square\)

\begin{figure}
\centering
\includegraphics[width=0.9\linewidth,height=\textheight,keepaspectratio,alt={Two incomparable missing-hub sets leave a hub that the source can reach but the pursuer cannot. Every hub reaches the outside recipient, completing a safe two-move delivery. The theorem separately handles recipients and pursuers inside the core.}]{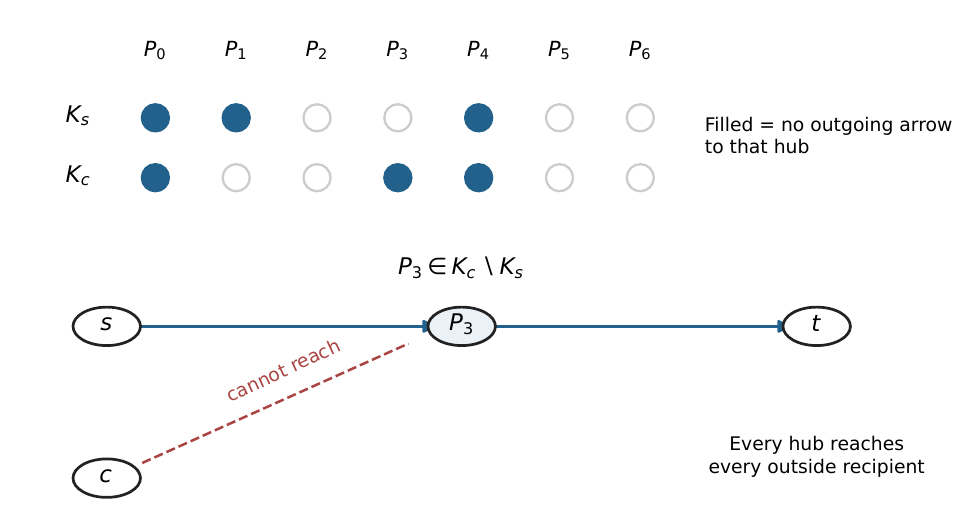}
\caption{Two incomparable missing-hub sets leave a hub that the source
can reach but the pursuer cannot. Every hub reaches the outside
recipient, completing a safe two-move delivery. The theorem separately
handles recipients and pursuers inside the core.}
\end{figure}

\textbf{Corollary (separated pair codes).} For \(k\ge6\), there are
\(k(k-5)\) permitted pair codes. For \(0\le q\le k(k-5)\), put
\(n=2k+q\). Every budget
\[ B\le m\le n(n-3),\qquad B=(4k-2)n-2k(2k+1) \] has an explicit
attainer with two-move delivery and every optional subset.

\textbf{Proof.} Each hub index has \(k-5\) partners after excluding
itself and differences \(\pm1,\pm2\). The two parts supply \(k(k-5)\)
pairs. Distinct pairs form an antichain. List the point pairs, then line
pairs, lexicographically, and take the first \(q\). Theorem 59 with
\(W=2q\) gives base \(B\) and terminal missing count \(4k+2q=2n\). Take
the desired optional prefix in ordered-endpoint order. The original
codes of size at least three, and the \(k=5,q=0\) core, remain
unchanged. \(\square\)

\subsubsection{The main interval above three arrows per
vertex}\label{the-main-interval-above-three-arrows-per-vertex}

\textbf{Theorem (pair-code overlap).} Every \(n\ge79\) and
\(3n\le m\le n(n-3)\) has an explicit attaining graph.

\textbf{Proof.} Choose the least \(k\ge6\) with \(n\le k(k-3)\). Then
\(2k\le n\) and the pair corollary supplies \(B\le m\le n(n-3)\).
Theorem 38 supplies \(3n\le m\le T=n(\lceil n/7\rceil-1)\), with
\(T\ge5n\). The five-offset circle with admissible \(h\) reaches at
least \(n(n-h-1)/2\), starting at \(5n\). The following direct parameter
rules make the intervals meet.

For \(79\le n\le108\), take \(h=7\) unless \(7\mid n\), then \(h=11\).
The \(h=7\) lower cap exceeds \(B\) on \(80\ldots88\) (\(k=11\)) and
\(89\ldots108\) (\(k=12\)): check each left endpoint and the increasing
quadratic difference. At 79 its exact cap is \(2844>B=2812\). The
\(h=11\) cases are \(84,91,98,105\); their lower cap \(n(n-12)/2\)
exceeds \(B\), first by 2 at 84 and by \(17/2\) at 91, then
increasingly.

For \(109\le n\le827\), choose the first coprime \(h\) in
\(7,11,15,19\). Failure requires \(7\cdot11\cdot19\) and either 3 or 5
to divide \(n\), forcing \(n\ge4389\). All circle hypotheses hold. The
cap is at least \(n(n-20)/2\). Here \(k\ge13\), and minimality gives
\(n\ge N=k^2-5k+5\). The cap minus \(B\) increases from \(N\); twice its
value there, with \(u=k-13\), is
\[k^4-18k^3+67k^2-6k-55=u^4+34u^3+379u^2+1398u+205>0.\]

For \(n\ge828\), Theorem 38 alone reaches \(B\). Subtracting \(B\) from
\(n^2/7-n\) gives \(f=n^2/7-(4k-1)n+4k^2+2k\). At \(k=31,n=828\) and
\(k=32,n=869\) it is respectively \(18/7\) and \(11740/7\), and
increases through both ranges. For \(k\ge33\), minimality gives
\(n\ge k^2-5k+5\ge28k-7\), so \(n^2/7-n\ge(4k-2)n>B\). All choices use
arithmetic and optional prefixes, with no graph search. \(\square\)

\texttt{src/\allowbreak{}main\_\allowbreak{}interval\_\allowbreak{}pair\_\allowbreak{}codes.\allowbreak{}py} implements the direct rules;
\texttt{src/\allowbreak{}main\_\allowbreak{}interval\_\allowbreak{}pair\_\allowbreak{}codes\_\allowbreak{}check.\allowbreak{}py} independently
reconstructs witnesses and coverage. Its receipt is
\texttt{research/\allowbreak{}main-\allowbreak{}interval-\allowbreak{}pair-\allowbreak{}code-\allowbreak{}checks.\allowbreak{}json}. The separated
pair-code construction and its overlap from order 79 are also
instantiated in \texttt{PairCodeConstruction.\allowbreak{}lean} and
\texttt{PairCodeOverlap.\allowbreak{}lean}. The retained larger and mixed codes below
keep their separate proof boundaries. The former classification-overlap
proof is preserved in
\texttt{research/\allowbreak{}history/\allowbreak{}before-\allowbreak{}slimming-\allowbreak{}core-\allowbreak{}code-\allowbreak{}composition.\allowbreak{}md}. The
following independent stronger construction results remain part of the
paper.

\subsubsection{Retained larger and mixed
codes}\label{retained-larger-and-mixed-codes}

\textbf{Equal-weight codes at every core parity.} For a cycle of length
\(t\ge3\), write \(C(t,0)=1\), put \(C(t,j)=0\) if \(j<0\) or \(2j>t\),
and otherwise put \[C(t,j)=\frac{t}{t-j}\binom{t-j}{j}.\] This counts
its independent \(j\)-sets: for \(j\ge1\), distinguish a selected vertex
and distribute the unselected vertices among the \(j\) nonempty gaps.
Dividing by the \(j\) choices of distinguished selected vertex gives
\(\frac{t}{j}\binom{t-j-1}{j-1}\), the same expression.

For every \(k\ge5\) and \(3\le r\le\lfloor k/2\rfloor\), the total
number of permitted size-\(r\) codes in the two hub parts is
\[Q(k,r)=C(k,r)+
\begin{cases}
C(k,r),&k\text{ odd},\\
\displaystyle\sum_{j=0}^{r}C(k/2,j)C(k/2,r-j),&k\text{ even}.
\end{cases}\] Indeed, the forbidden graph in the line part is a
\(k\)-cycle. For odd \(k\), multiplication by two identifies this cycle
with the forbidden difference-two graph in the point part. For even
\(k\), the even and odd indices instead form two disjoint cycles of
length \(k/2\). Choosing \(j\) indices in one and \(r-j\) in the other
gives the displayed convolution. The two parts have disjoint hub labels,
so their counts add. All distinct size-\(r\) codes are incomparable.

Consequently every \(q\le Q(k,r)\) supplies the explicit interval in
Theorem 59, with \[A=(4k-r)n-2k(2k+3-r),\qquad L-A-q(q-1)=2n+(r-2)q.\]
List independent sets by taking increasing indices with at least one
unselected index between consecutive selections and across the cyclic
boundary. For even cores, list the two component choices in increasing
order of \(j\) and then in lexicographic order within each component.
List line codes first, then point codes, and take the first \(q\). This
gives a concrete rule from \((k,r,q)\) without pursuit-game search. It
includes all previous odd-core ranges. The implementation and
independent counting and asymmetric witness checks are
\texttt{src/\allowbreak{}core\_\allowbreak{}code\_\allowbreak{}all\_\allowbreak{}parities.\allowbreak{}py} and
\texttt{src/\allowbreak{}core\_\allowbreak{}code\_\allowbreak{}all\_\allowbreak{}parities\_\allowbreak{}check.\allowbreak{}py}; the original odd-core
implementation and its earlier receipts remain available unchanged.

For example, \(Q(10,3)=100\), giving \(37n-400\le m\le n(n-1)-(3n-20)\)
for \(20\le n\le120\). The fixed mixed list at order 105 contains 85
incomparable \(k=10\) codes of total weight 268, giving every budget
\(3472\ldots10612\). At order 116 the saved \(k=11\) list has two
size-five and 92 size-four codes, weight 378, giving
\(4176\ldots12918\): the two larger codes exclude at most ten of the 110
available size-four codes. The literal lists and independent witnesses
remain in \texttt{research/\allowbreak{}core-\allowbreak{}code-\allowbreak{}third-\allowbreak{}mixed-\allowbreak{}105.\allowbreak{}json} and
\texttt{research/\allowbreak{}routing-\allowbreak{}defect-\allowbreak{}phase2-\allowbreak{}mixed-\allowbreak{}code-\allowbreak{}checks.\allowbreak{}json}. These
families retain two-move delivery and every optional subset. Their
stronger two-move sparsity consequences remain in the
\href{https://angelraychev.com/cop-prisoner/delivery-time/\#antichain-two-move-sparsity}{delivery-time
companion}.

\subsubsection{Retained wider protected
circles}\label{retained-wider-protected-circles}

Here we use the protected differences of Theorem 38 without changing its
strategy. Let \(h\) be an allowed odd step, \(a=(h-1)/2\), and
\[A_0=\{0,1,a,2a+1\},\qquad
D=\{\pm(a-1),\pm(a+1),\pm2a\}.\] Put \(p=3a-1\) and
\[J=\{jp+t:0\le j<\lfloor n/p\rfloor,\ 0\le t<a-1\}.\] No difference
between two members of \(J\) belongs to \(D\): within a run it has
magnitude at most \(a-2\), while between runs, including the cyclic gap,
the shortest gap is at least \(2a+1\). For any translate \(J+v\), remove
all members of \(A_0\cup(A_0+D)\) and then adjoin \(A_0\). Every
protected difference still has exactly its original representation. Thus
the result is an allowed envelope in Theorem 38. The excluded set has at
most 19 elements, by expanding the four translates of the six
differences. Averaging over \(v\) gives an envelope of size at least
\[4+\left\lceil\frac{|J|(n-19)}{n}\right\rceil.\] One may then apply the
greedy extension from Theorem 38. In particular the old guarantee
\(|A|\ge\lceil n/7\rceil\) is preserved. Searching the \(n\)
translations is explicit and takes polynomial time.

The literal protected sets at orders \(100\ldots111\) remain in
\texttt{research/\allowbreak{}core-\allowbreak{}code-\allowbreak{}third-\allowbreak{}circle-\allowbreak{}envelopes.\allowbreak{}json}, supplying every
budget \(3n\ldots n(|A|-1)\). At order 117, the explicit envelope
\(\{0,1,20\}\cup[41,59]\cup[100,115]\) with \(h=41\) has 38 members and
exactly one representation of each protected difference, retaining
\(351\ldots4329\). Its independent check and the order-116 mixed codes
remain in
\texttt{research/\allowbreak{}routing-\allowbreak{}defect-\allowbreak{}phase2-\allowbreak{}mixed-\allowbreak{}code-\allowbreak{}checks.\allowbreak{}json}. Original
implementations and receipts remain unchanged.
\texttt{RobustTwoMove.\allowbreak{}lean} proves the semantic implication from
supplied witnesses, not the combinatorial construction or the complete
classification.

\subsection{Routing through a layer of
guides}\label{routing-through-a-layer-of-guides}

A small hub-and-guide interface can serve arbitrarily many additional
vertices. The following attachment lemma separates its delivery proof
from the antichain argument that enters it. This gives the three- and
five-move constructions below without repeating their pursuer cases.

\textbf{Lemma (guide attachment).} Let \(H\) be a set of \(h\) hubs with
\(a\) internal arrows, and let \(D\) be \(g\) disjoint guide vertices.
Guide \(j\) reaches precisely the hubs outside its omitted set \(K_j\),
and no other guide. Put \(W_K=\sum_j|K_j|\). Include all arrows from
\(R\subseteq H\) to the guides and no other hub-to-guide arrows. Choose
\(S\subseteq U\subseteq H\).

Attach \(b\) bulk vertices. Every hub in \(R\) reaches every bulk
vertex; each bulk vertex \(x\) reaches \(S\) and its assigned guide set
\(J_x\). These are all compulsory arrows. The distinct sets \(J_x\) have
size at least two and form an antichain. Put \(W_J=\sum_x|J_x|\). The
optional arrows are all guide-to-bulk arrows, all distinct bulk-to-bulk
arrows, the arrows from a fixed \(C\subseteq H\setminus R\) to every
bulk vertex, and the arrows from each bulk vertex to \(U\setminus S\).
Assume:

\begin{enumerate}
\def\labelenumi{\arabic{enumi}.}
\tightlist
\item
  Every hub in \(R\) has a missing hub successor in \(S\); the successor
  must differ from the hub itself.
\item
  From every hub or guide, delivery takes at most \(d\) moves using only
  hubs as nonterminal intermediaries and compulsory final arrows. This
  interface strategy remains valid against a bulk pursuer allowed every
  hub in \(U\), every guide, and every bulk vertex.
\end{enumerate}

Then every optional subset gives universal delivery within \(d+1\)
moves. Its exact-budget interval is \[
 A\le m\le A+b\bigl(g+b-1+|C|+|U\setminus S|\bigr),
\] where \[
 A=a+(|R|+h)g-W_K+(|R|+|S|)b+W_J.
\] The interface hypothesis is finite: a recipient outside \(H\) is
represented by a terminal outlet reached from \(R\), and a bulk pursuer
by its permitted hub set \(U\) and unrestricted guide access. Hub moves
do not distinguish other bulk vertices. Receipt at the actual recipient
is always terminal, even if that vertex is occupied by the pursuer.

\textbf{Proof.} Let the messenger start at a bulk vertex \(x\). If the
pursuer is a hub in \(R\), the first hypothesis supplies a safe
compulsory move into \(S\). A hub outside \(R\) reaches no guide, so any
assigned guide is safe against it. Against a guide pursuer choose an
assigned guide other than the pursuer, possible because \(|J_x|\ge2\).
Against another bulk pursuer \(c\), choose a guide in
\(J_x\setminus J_c\); incomparability makes this difference nonempty,
and no optional arrow enters a guide from a bulk vertex. Each chosen
vertex lies in the interface and outside the pursuer's closed
neighborhood in the full envelope. After the cop reply the interface
strategy therefore applies. If the chosen vertex is the recipient,
receipt has already ended play. Interface sources need no entry move.
This proves the time bound against the full envelope and hence for every
optional subset.

Count the \(a\) hub arrows, \(|R|(g+b)\) outlet arrows, \(hg-W_K\)
guide-to-hub arrows, and \(|S|b+W_J\) bulk arrows to obtain \(A\). The
four disjoint optional classes have sizes \(gb\), \(b(b-1)\), \(|C|b\),
and \(|U\setminus S|b\). Taking an ordered-endpoint prefix attains every
stated integer budget. \(\square\)

\textbf{Theorem (three-move layered construction).} Choose integers
\(2\le r\le g\le14\) and \(0\le b\le\binom gr\), and put \(n=14+g+b\).
There is a direct construction attaining \(M(n,m)=n(n-1)-m\) for every
\[
 A\le m\le A+b(g+b-1),\qquad A=154+25g+(22+r)b.
\] It guarantees delivery within three messenger moves. Every subset of
the specified optional arrows is permitted. For pairs of guides,
\(r=2\), the compulsory budget is \(A=24n-182+g\).

\textbf{Construction and proof.} Take fourteen hubs \(P_i,L_i\), indexed
modulo seven, with missing arrows
\[P_i\to L_i,L_{i-1},\qquad L_i\to P_{i+4},P_{i+6}.\] All other distinct
hub pairs are arrows, giving \(a=154\). Use the first \(g\) omitted
guide codes in the ordered list \[
 \{L_i,L_{i+2},L_{i+4}\}\quad(0\le i<7),\qquad
 \{P_i,P_{i+1},P_{i+4}\}\quad(0\le i<7).
\] Set \(R=H\), \(C=\varnothing\), and
\[S=U=\{P_0,P_1,P_4,P_5,L_0,L_2,L_4,L_6\}.\] Assign the bulk vertices
the first \(b\) distinct \(r\)-subsets of the guides in lexicographic
order. The guide attachment lemma prescribes all arrows. In particular,
only guide-to-bulk and distinct bulk-to-bulk arrows are optional; take
their first \(m-A\) ordered endpoints for the exact budget.

The displayed guide codes satisfy Theorem 59. Thus hub and guide sources
deliver within two moves against hub and guide pursuers, using only hub
intermediaries. To check a bulk pursuer, put
\[T=H\setminus U=\{P_2,P_3,P_6,L_1,L_3,L_5\}.\] For a missing hub pair
\((s,t)\), its danger set \(F^+(s)\cup F^-(t)\) has only two points and
two lines, so it cannot contain \(T\). A hub in the difference is a safe
intermediary. A hub source reaches any outside recipient directly. For a
guide source, the opposite part from its omitted code contains three
reachable safe hubs in \(T\). All reach an outside recipient. A missing
hub recipient lies in the code's own part, so at most two fail to reach
it. This proves the two-move interface hypothesis. Finally, \(T\)'s
point triple contains no difference-two pair and its line triple
contains no adjacent pair, so \(S\) meets every hub missing row. The
lemma applies with \(d=2\), \(W_K=3g\) and \(W_J=rb\), giving exactly
the stated count and three-move bound. Every missing pair delivers, so
the missing-pair upper bound is attained. \(\square\)

For example, \((g,b,r)=(10,36,2)\) gives every budget 1268 through 2888
at order 60; \((11,55,2)\) gives 1749 through 5324 at order 80. The
larger-\(r\) formulas remain valid when more bulk vertices are needed.
The direct implementation is \texttt{src/\allowbreak{}main\_\allowbreak{}interval\_\allowbreak{}layered.\allowbreak{}py};
its checks replay the two-move interface and safe-entry obligations
against the full optional envelope. They supplement the ordinary
all-parameter proof.

\subsubsection{A smaller core with antipodal guide
codes}\label{a-smaller-core-with-antipodal-guide-codes}

\textbf{Theorem (twelve-hub guide construction).} Choose
\(2\le r\le g\le8\), \(0\le b\le\binom gr\), and put \(n=12+g+b\). Every
exact budget \[
 A\le m\le A+b(g+b-1),\qquad A=106+22g+(19+r)b
\] has a direct attaining construction with delivery within three
messenger moves. Every optional subset is valid.

\textbf{Construction and proof.} Use the same hub missing-arrow formula
modulo six, giving \(h=12\) and \(a=108\). Use the first \(g\) guide
codes in \[
 \{L_0,L_2,L_4\},\quad\{L_1,L_3,L_5\},\quad
 \{P_i,P_{i+3}\}\ (0\le i<3),\quad
 \{L_i,L_{i+3}\}\ (0\le i<3).
\] Set \(R=H\), \(C=\varnothing\), and
\[S=U=H\setminus T,\qquad T=\{P_0,P_3,L_0,L_2,L_4\}.\] Again assign the
first \(b\) lexicographic \(r\)-subsets of guides and use the attachment
lemma's adjacency and optional-prefix rule.

All guide codes satisfy generalized Theorem 59, including its separated
pair condition. Its hub-intermediary proof supplies the interface
against hub and guide pursuers. Against a bulk pursuer, a missing hub
pair has a safe intermediary in \(T\): its line triple cannot fit in the
danger set's two lines. For a point-code guide source, all three safe
lines are reachable, and at most two fail to reach a missing hub
recipient, which is a point. For a line-code guide, the two safe points
are reachable. A missing hub recipient is a line, whose missing
predecessors are adjacent points; they cannot contain the antipodal pair
\(\{P_0,P_3\}\). Either choice of opposite part also reaches every
outside recipient. Hub sources reach outside recipients directly. This
proves the two-move interface.

The safe point pair has difference three and the safe line triple has no
adjacent pair, so \(S\) meets every hub missing row. The attachment
lemma now uses \(d=2\), \(|S|=7\), \(W_K=2g+2\) and \(W_J=rb\), giving
\(108+24g-(2g+2)+(19+r)b=A\) and the claimed optional count. \(\square\)

For pair guide codes the base is \(21n-146+g\). At
\((n,g,b,r)=(40,7,21,2)\) this gives 701 through 1268. Taking
\(g=8,r=3\) gives \[
 22n-158\le m\le22n-158+(n-20)(n-13),\qquad20\le n\le76.
\] At orders 50 and 60 the intervals are respectively 942 through 2052
and 1162 through 3042. These are all-parameter consequences of the
ordinary proof, not extrapolations from numerical games.

\texttt{lean/\allowbreak{}GuideAttachment.\allowbreak{}lean} formalizes antichain safe entry
followed by a supplied interface strategy. \texttt{TwelveHubGraph.\allowbreak{}lean},
\texttt{TwelveHubCounts.\allowbreak{}lean} and \texttt{TwelveHubCodes.\allowbreak{}lean} also
instantiate the twelve-hub equal-weight family, using numeric mask order
for its codes. Other interfaces and code variants retain their separate
ordinary and finite proof boundaries.

\subsection{Circle constructions with three-way
escape}\label{circle-constructions-with-three-way-escape}

The following construction gives the messenger several alternative
moves. No pursuer at a different vertex can cover all alternatives. A
finite propagation argument then makes every source winning for every
recipient. The construction tolerates every subset of its optional
arrows.

\textbf{Theorem (half-density circles and a wraparound gap).} Let
\(h\ge7\), \(h\equiv3\pmod4\), \(2h+2\le n\), and \(\gcd(h,n)=1\). In
\(\mathbb Z/n\mathbb Z\), take the compulsory offsets
\[C=\{1,h-3,h-2,h,h+2\}.\] One may use either of the following closed
offset sets, whose representatives are in \(\{0,\ldots,n-1\}\):

\begin{enumerate}
\def\labelenumi{\arabic{enumi}.}
\tightlist
\item
  If \(4\mid n\), take
  \[A=\{x:x\equiv0,1\pmod4\}\cup\{h\}\setminus\{2h-1\}.\]
\item
  At any order satisfying the hypotheses, take
  \[A=\{0\le x\le n-h:x\equiv0,1\pmod4\}\cup\{h\}
  \setminus\{2h-1\}.\]
\end{enumerate}

In both formulas the deletion applies to the whole preceding union.
Include all arrows \(x\to x+c\) for \(c\in C\), and an arbitrary subset
of the arrows \(x\to x+a\) for \(a\in A\setminus(C\cup\{0\})\). Every
resulting graph has universal delivery and hence attains
\[M(n,m)=n(n-1)-m\qquad\text{for every }5n\le m\le n(|A|-1).\] For an
exact budget, include the first \(m-5n\) optional arrows in increasing
lexicographic order of their ordered endpoints.

The first version has \(|A|=n/2\). For the second version, put
\(L=n-h\), write \(L+1=4q+r\) with \(0\le r<4\), and put
\[e=\mathbf1_{2h-1\le L}.\] Then \(|A|=2q+\min(r,2)+1-e\). These are
explicit adjacency and counting rules, with no graph search.

\textbf{Proof.} All five compulsory offsets are distinct and nonzero.
They belong to either envelope: \(h+2\le n-h\), and none equals the
deleted offset. The offset zero is allowed for waiting but is not
counted as an arrow. The cardinalities follow by counting the residue
classes and the one exceptional offset \(h\).

We first prove a local escape property. Relative to a messenger at zero,
the following sets of alternative moves cannot all be reached in one
move (including waiting) by any pursuer at \(q\ne0\):
\[\{h,1\},\qquad\{h,h-3,h-2\},\qquad\{h,h-2,h+2\}.\]

For a full modulo-four envelope, suppose the pursuer reaches \(h\). A
pursuer congruent to one modulo four cannot do so. A pursuer congruent
to zero must be zero, because \(h\) is the only envelope element
congruent to three. A pursuer congruent to two reaches no zero-class
point and at most one one-class point. A pursuer congruent to three
reaches no one-class point. Each displayed triple contains \(h\)
together with either a zero-class and a one-class point, or two distinct
one-class points. This proves their escape property. For the pair
\(\{h,1\}\), the only remaining possible pursuer is \(q=1-h\). Reaching
\(h\) would require the deleted offset \(2h-1\), a contradiction.

The second envelope has a gap of \(h-1\) consecutive absent offsets at
the end of its representative interval. Every displayed alternative set
has diameter at most \(h-1\). If, after subtracting a pursuer position,
its representatives wrapped around zero, their extreme representatives
would differ by at least \(n-h+1\). They cannot both lie in \([0,n-h]\).
Thus any putative covering uses a common integer lift of the entire
alternative set. The same residue argument applies to that lift. Its
only exceptional pursuer is the integer zero, which represents the
excluded modular pursuer zero. The pair argument again requires a
deleted offset. This proves the local escape property for every order in
the second version.

Fix recipient zero. Call a vertex winning when it guarantees eventual
delivery from every distinct pursuer position. The recipient and all
direct compulsory predecessors are winning. Whenever the compulsory
successors of a vertex contain one of the displayed alternative sets of
already winning vertices, that vertex is winning too: choose a successor
outside the pursuer's full closed envelope. After the pursuer response,
the positions remain distinct and the winning successor strategy
applies. Each use of this rule increases a finite delivery rank, so the
conclusion is delivery in finite time, not merely indefinite evasion.

Let \(P(z)=\{z-2,z-1,z,z+2\}\). The pattern \(P(-h)\) is winning: three
of its vertices are direct compulsory predecessors; its remaining vertex
\(-h-1\) uses offsets \(h,1\) to reach predecessors \(-1,-h\). Given
\(P(z)\), make the following four additions in order:

\begin{itemize}
\tightlist
\item
  \(z-h+2\) uses \(h,h-3,h-2\), reaching \(z+2,z-1,z\);
\item
  \(z-h\) uses \(h,h-2,h+2\), reaching \(z,z-2,z+2\);
\item
  \(z-h-1\) uses \(h,1\), reaching \(z-1,z-h\);
\item
  \(z-h-2\) uses \(h,1\), reaching \(z-2,z-h-1\).
\end{itemize}

Every successor belongs to the original pattern or an earlier addition.
The resulting set contains \(P(z-h)\). Since \(\gcd(h,n)=1\), repeated
shifts by \(-h\) visit every residue. Modular coincidences only identify
vertices already known to be winning; they do not invalidate this finite
induction. Translating this proof gives every recipient. All moves of
the messenger are compulsory, and the pursuer has throughout been
allowed the entire envelope. Consequently the assertion holds for every
optional subset, including nonsymmetric ones. The arc count and the
missing-pair upper bound give the stated value of \(M\). \(\square\)

\textbf{Parameters at every order divisible by four.} The full-envelope
version applies to every \(n\ge16\) divisible by four. If \(8\mid n\),
choose \(h=n/2-1\). Otherwise write \(n=4k\) with odd \(k\ge5\): choose
\(h=k+2\) when \(k\equiv1\pmod4\), and \(h=k+4\) when
\(k\equiv3\pmod4\). The relevant gcd divides eight or sixteen and is
odd, so is one. The remaining inequalities follow directly. Thus the
complete direct interval in that subfamily is \[5n\le m\le n(n/2-1).\]

For other orders, take the smallest admissible \(h\), or choose among
the finitely many arithmetic candidates to maximize the explicit
envelope count. This parameter calculation is separate from graph or
strategy search. In particular \(h=7\) applies whenever \(n\ge16\) and
\(7\nmid n\), regardless of the parity of \(n\). The second family
therefore bridges the old divisibility restriction while losing only a
bounded number of offsets when \(h\) is fixed.

\textbf{Relation to the earlier eight-move construction.} The previous
construction had the additional compulsory offsets \(5,h-7,h+1\) and
removed both \(2h-5\) and \(2h-1\), under the stronger assumption
\(h\ge11\). Its compulsory graph contains the present one and its
envelope is contained in the present envelope. Thus the present theorem
preserves every old optional-subset guarantee while lowering the
compulsory budget from \(8n\) to \(5n\) and restoring one allowed
offset. In particular the earlier interval \(8n\le m\le n(n/2-2)\) at
\(4\mid n\), \(n\ge24\), remains valid.

The direct implementation is
\texttt{src/\allowbreak{}main\_\allowbreak{}interval\_\allowbreak{}half\_\allowbreak{}density.\allowbreak{}py}. Finite structural replays
are recorded separately. \texttt{CircleConstruction.\allowbreak{}lean} and
\texttt{CircleGuardedConstruction.\allowbreak{}lean} also instantiate the full and
guarded five-offset constructions with exact counts and optional
prefixes; \texttt{CircleParameters.\allowbreak{}lean} supplies the displayed
divisible-by-four domain.

\subsection{The main interval above three arrows per
vertex}\label{the-main-interval-above-three-arrows-per-vertex-1}

The circle and guide constructions reduce the order needed for the broad
main interval. The distinction between its first linear segment and its
remaining budgets is now explicit.

\textbf{Theorem (expanded main interval).} There are explicit attaining
constructions with \[M(n,m)=n(n-1)-m\] throughout both of the following
domains:

\begin{enumerate}
\def\labelenumi{\arabic{enumi}.}
\tightlist
\item
  Every \(n\ge43\) and \(3n\le m\le n(n-3)\).
\item
  Every \(n\ge104\) and \(2n\le m\le n(n-3)\).
\end{enumerate}

The thresholds are sufficient, not claimed necessary. Every construction
comes with a concrete adjacency rule and arithmetic parameter selection;
no new graph search or derandomization is required for a given
\((n,m)\). The individual families retain their optional-subset and
delivery-time guarantees. The combined statement imposes no uniform
two-move deadline.

\textbf{Proof.} For \(n\ge106\), the previous complete main-interval
theorem already applies. It remains to cover the finite range
\(43\le n\le105\). Here we combine the following proved intervals,
taking only parameters satisfying their stated arithmetic hypotheses:

\begin{itemize}
\tightlist
\item
  The protected three-jump circle, from \(3n\) through
  \(n(\lceil n/7\rceil-1)\), reaches at least \(5n\).
\item
  The new five-jump circle reaches \(n(|A|-1)\), using a full
  modulo-four envelope or its guarded version as appropriate.
\item
  The twelve- and fourteen-hub guide constructions provide their
  displayed exact intervals, with \(2\le r\le g\) and enough distinct
  guide subsets.
\item
  The equal-weight antichain construction of Theorem 59 provides its
  explicit interval at each admissible core size and code weight.
\item
  The earlier dense envelopes finish at \(n(n-3)\), starting no later
  than \(n(n-1)-n\max\{3,\lceil(n-1)/13\rceil\}\).
\end{itemize}

The finite parameter certificate
\texttt{research/\allowbreak{}main-\allowbreak{}interval-\allowbreak{}expanded-\allowbreak{}checks.\allowbreak{}json} lists a covering
chain at every order \(43,\ldots,105\). Each chain starts at \(3n\), or
at \(2n\) for orders 104 and 105, and adjacent integer intervals meet
without a missing budget. All parameters are integers checked against
the preceding formulas. The two smaller starting budgets use the already
certified canonical sparse envelopes: order 104 covers 208 through 316,
and order 105 covers 210 through 327. Both reach \(3n\).

This is a finite certificate of parameters and inequalities for proved
families, rather than a collection of unverified game outcomes. The
selector enumerates only the displayed arithmetic parameters, takes a
covering family in a fixed order, and selects the prescribed number of
optional arrows. The direct constructor
\texttt{src/\allowbreak{}main\_\allowbreak{}interval\_\allowbreak{}expanded.\allowbreak{}py} implements this rule. Its
separate check constructs graphs at the ends and middle of each selected
interval and checks their exact arrow counts. The original constructor
is used unchanged above order 105. This proves both assertions.
\(\square\)

For example, at order 60 the new circle covers 300 through 1740 and the
twelve-hub guide family covers 1162 through 3042. The preceding sparse,
protected-circle and dense constructions cover the two ends of the main
interval. Their overlap fills every budget from 120 through 3420.

\begin{figure}
\centering
\includegraphics[width=1\linewidth,height=\textheight,keepaspectratio,alt={At order 60, the five-jump circle and guide layer bridge the previously uncovered budgets 841 through 1655. The horizontal bars are exact construction intervals; their overlap provides every integer budget.}]{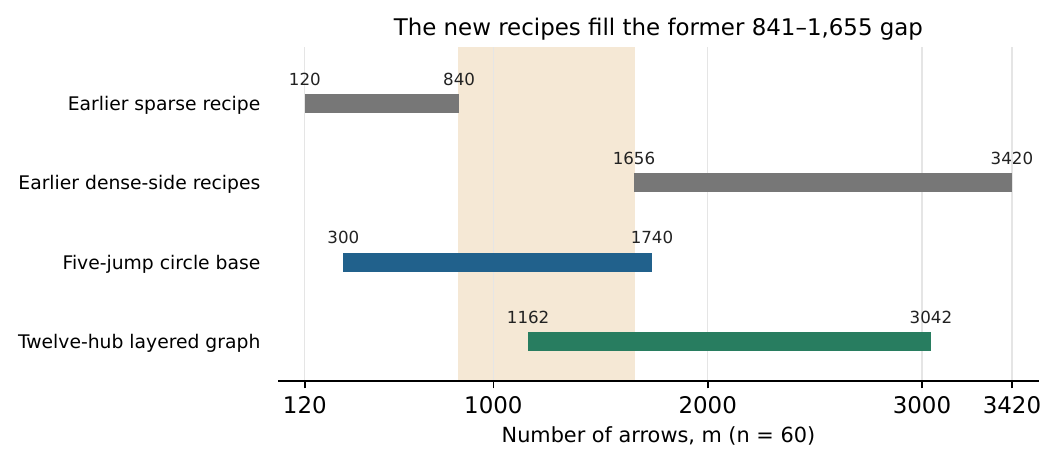}
\caption{At order 60, the five-jump circle and guide layer bridge the
previously uncovered budgets 841 through 1655. The horizontal bars are
exact construction intervals; their overlap provides every integer
budget.}
\end{figure}

Combining the second domain with the previously proved sparse and dense
formulas strengthens Corollary 47 to every \(n\ge104\) and \(m\ge22\).
Budgets 12 through 21 remain a separate small-support problem. At orders
below 104, the remaining gap inventory distinguishes \(2n<m<3n\) from
the uncovered main budgets at orders at most 42. Neither the new domain
nor a successful finite construction implies attainment at an unlisted
pair.

The ordinary circle and guide proofs have independent local and
actual-game rank checks. \texttt{RobustLayers.\allowbreak{}lean} proves that supplied
strictly decreasing base/envelope source ranks imply actual delivery,
and that safe entry into a routing region adds one move. Concrete circle
and guide-code modules now also discharge selected construction
premises, counts and overlaps. \texttt{MainUpperInterval.\allowbreak{}lean} proves
the full interval above \(3n\) from order25. The later
\texttt{MainIntervalClassification.\allowbreak{}lean} joins the entire main interval
from order12, including this theorem's domain. Stronger optional-subset
and timing guarantees remain attached to their individual constructions;
the remaining full-classification obligations are stated in the
verification section.

\subsection{Exact-budget consequences of protected
routing}\label{exact-budget-consequences-of-protected-routing}

The full-vertex random and deterministic constructions are
specializations of the protected-routing theorem. Their common counting
proof is Theorem 54.

\textbf{Theorem 43.} Put \(L=n(n-1)\) and \(p=m/L\). Under the
hypotheses of Theorem 25, there is a deterministic algorithm, polynomial
in \(n\), constructing an \(n\)-vertex, exactly \(m\)-arc graph in which
every pair is guaranteed within two messenger moves.

\textbf{Proof.} Set \(H=V\), \(b=n\), and \(a=m\) in Theorem 54. The
condition \((n-3)p^2(1-p)>5\log n\) implies the required condition
because \(\log n>\log(8/3)\) for \(n\ge4\). The restricted completion
counts become exactly the full-vertex counts. \(\square\)

The algorithm's polynomial bound does not imply a small practical
runtime. \texttt{src/\allowbreak{}exact\_\allowbreak{}budget\_\allowbreak{}derandomize.\allowbreak{}py} retains the earlier
full-vertex implementation and small exhaustive checks;
\texttt{src/\allowbreak{}routing\_\allowbreak{}unification.\allowbreak{}py} generalizes those counts to a
specified routing set. Neither receipt claims execution on a large
asymptotic instance.

\subsubsection{An explicit threshold for the entire main
interval}\label{an-explicit-threshold-for-the-entire-main-interval}

\textbf{Theorem 44.} For every integer \(n\ge106\) and every integer
\[2n\le m\le n(n-3),\] an explicit polynomial-time construction produces
a graph with exactly \(n\) vertices and \(m\) arcs in which every pair
is guaranteed. Consequently \(M(n,m)=n(n-1)-m\) throughout this
interval.

\textbf{Proof.} First let \(106\le n\le114\). The finite cyclic
supplement to Theorem 52 covers every budget \(2n\le m\le3n\), using the
two specified restored optional arrows at order 109. The literal
protected-circle envelopes at these nine orders have maximum budget at
least \(37n-400\). With \(k=10\) and code weight three, Theorem 59
supplies 100 codes, more than the required \(q=n-20\), and its interval
is \[37n-400\le m\le n(n-1)-(3n-20).\] Its upper endpoint overlaps the
two-move dense bridge \(n(n-1)-3n\le m\le n(n-3)\). Thus the sparse,
circle, code and dense intervals cover every budget at these nine
orders. The explicit circle lists are recorded in
\texttt{research/\allowbreak{}core-\allowbreak{}code-\allowbreak{}third-\allowbreak{}circle-\allowbreak{}envelopes.\allowbreak{}json} for 106 through
111 and in \texttt{research/\allowbreak{}core-\allowbreak{}code-\allowbreak{}circle-\allowbreak{}overlap-\allowbreak{}screen.\allowbreak{}json} for
112 through 114. Their six protected differences and all integer
overlaps are independently checked in
\texttt{research/\allowbreak{}sparse-\allowbreak{}interface-\allowbreak{}third-\allowbreak{}cross-\allowbreak{}lane-\allowbreak{}review.\allowbreak{}json}. The
executable adjacency recipe is \texttt{src/\allowbreak{}main\_\allowbreak{}interval\_\allowbreak{}composed.\allowbreak{}py}.
Its finite extension is checked against the sparse certificate graphs,
the circle conditions and the full code/dense witness conditions in
\texttt{research/\allowbreak{}main-\allowbreak{}interval-\allowbreak{}composed-\allowbreak{}third-\allowbreak{}checks.\allowbreak{}json}; the earlier
outputs at the listed regression orders are unchanged.

For \(n\ge115\), Theorem 52 covers \(2n\le m\le3n\). The separated-pair
overlap theorem following Theorem 59 covers \(3n\le m\le n(n-3)\) at
every \(n\ge79\), hence at each remaining order. Together these two
intervals cover every stated budget. Each construction lists its base
arrows and optional arrows; include the prescribed number of optional
arrows in lexicographic order. All relevant lists can be constructed in
polynomial time. The missing-pair bound gives the reverse inequality.
\(\square\)

The sufficient order threshold 106 is not asserted to be optimal. This
proof is entirely by explicit combinatorial constructions and finite
local strategy certificates. At \(n\ge115\), the sparse part retains its
\(9\lceil n/6\rceil+21\) move bound; its finite 106--114 supplement uses
at most 72 moves. These are bounds for the sparse part; there is no
uniform two-move assertion for the entire interval. The
protected-routing theorem remains independently useful because it allows
a specified routing set and a variable incident density, with an
explicit random success probability and exact-budget derandomization.

\subsection{Six-offset circles and finite changes of the compulsory
moves}\label{six-offset-circles-and-finite-changes-of-the-compulsory-moves}

The five-offset circle uses a protected pair of escape moves. Its
envelope therefore cannot exceed half density: two translates of that
envelope have at most one common point. Adding a sixth compulsory move
replaces the pair by a triple. The same four-position pattern still
propagates, but the permitted graph may now contain substantially more
arrows.

\textbf{Theorem (three escape triples).} Work in
\(\mathbb Z/n\mathbb Z\) and choose \(h\) with \(\gcd(h,n)=1\). Require
the offsets \[C=\{1,h-3,h-2,h,h+1,h+2\}\] to be six distinct nonzero
residues. Let \(A\) contain \(C\cup\{0\}\). Suppose that, for every
nonzero residue \(q\), none of
\[S_0=\{h,1,h+1\},\qquad S_1=\{h,h-3,h-2\},\qquad
S_2=\{h,h-2,h+2\}\] is contained in \(q+A\). Include every arrow whose
difference belongs to \(C\) and any subset of the other arrows whose
nonzero difference belongs to \(A\). Every resulting graph has universal
delivery. Consequently this is an explicit attaining construction at
every exact budget \[6n\le m\le n(|A|-1),\qquad M(n,m)=n(n-1)-m.\]

\textbf{Proof.} Fix recipient zero and call a source winning when
delivery is guaranteed from every distinct pursuer position. A winning
set of vertices supports an additional source whenever all its
successors from one of the three displayed triples are already winning.
The escape condition selects one successor outside the pursuer's full
closed envelope. After every pursuer response, the messenger remains
uncaptured and uses the successor strategy. Every such addition has a
finite strictly decreasing strategy rank when followed during play.

Put \(P(z)=\{z-2,z-1,z,z+2\}\). Every member of \(P(-h)\) is a direct
compulsory predecessor of zero; the new offset \(h+1\) supplies its
previously missing member. From a winning \(P(z)\), add the following
four sources in order:

\begin{itemize}
\tightlist
\item
  \(z-h+2\) uses \(S_1\), reaching \(z+2,z-1,z\);
\item
  \(z-h\) uses \(S_2\), reaching \(z,z-2,z+2\);
\item
  \(z-h-1\) uses \(S_0\), reaching \(z-1,z-h,z\);
\item
  \(z-h-2\) uses \(S_0\), reaching \(z-2,z-h-1,z-1\).
\end{itemize}

This yields \(P(z-h)\). Coprimality makes its successive translates
cover every vertex. All additions use previously established finite
ranks, even when some modular vertices coincide. Translating the
argument establishes every recipient. The messenger uses only compulsory
arrows and the pursuer may use the full envelope throughout; therefore
every optional subset is covered, including nonsymmetric subsets. Add
the first \(m-6n\) optional arrows in lexicographic order of ordered
endpoints to obtain the exact budget. All missing pairs are winning,
proving equality in the missing-pair upper bound. \(\square\)

\textbf{Two periodic recipes.} The escape hypotheses admit direct
residue rules.

\begin{enumerate}
\def\labelenumi{\arabic{enumi}.}
\item
  If \(4\mid n\), \(h\ge7\), \(h\equiv3\pmod4\), \(2h+2\le n\), and
  \(\gcd(h,n)=1\), take \[A=\{x:x\equiv0,1\pmod4\}\cup\{h\}.\] There are
  no deleted offsets. This gives \(|A|=n/2+1\) and the whole interval
  \(6n\le m\le n^2/2\). The parameter rule from the five-offset theorem
  provides such an \(h\) at every \(4\mid n\), \(n\ge16\).

  To check escape, a pursuer of residue one cannot reach \(h\). A
  pursuer of residue zero reaches \(h\) only when it is zero. A pursuer
  of residue two reaches no zero-class point and at most one one-class
  point; a pursuer of residue three reaches no one-class point. Every
  required triple contains \(h\) and either a zero-class and a one-class
  point or two distinct one-class points. None can be covered.
\item
  If \(5\mid n\), \(h\ge8\), \(h\equiv3\pmod5\), \(2h+2\le n\), and
  \(\gcd(h,n)=1\), take
  \[A=\bigl(\{x:x\equiv0,1,4\pmod5\}\cup\{h\}\bigr)
      \setminus\{h+3,h-4,2h-1\}.\] The three deleted offsets are
  distinct, lie in the periodic part, and avoid \(C\cup\{0,h\}\). Thus
  \(|A|=3n/5-2\), giving \[6n\le m\le n(3n/5-3).\] For every odd
  multiple of five with \(n\ge25\), one may simply take \(h=8\). The
  rule also applies to all other displayed arithmetic parameters, with
  no graph search.

  For the escape check, the periodic residue set is three consecutive
  residues modulo five, whereas the residues of each required triple are
  not three consecutive residues. A covering translate must therefore
  use the exceptional envelope point \(h\). For \(S_0\) this leaves
  possible pursuers \(q=0,1-h,1\); for \(S_1\) they are \(0,-3,-2\); and
  for \(S_2\) they are \(0,-2,2\). Ignore the forbidden pursuer zero.
  For \(S_0\), \(q=1\) would require \(h-1\), whose residue is absent,
  and \(q=1-h\) requires the deleted \(2h-1\). For \(S_1\), \(q=-2\)
  requires the absent-residue offset \(h-1\) and \(q=-3\) requires the
  deleted \(h+3\). For \(S_2\), \(q=-2\) requires the absent-residue
  offset \(h+4\), and \(q=2\) requires the deleted \(h-4\). These
  exhaust the possible covers. \(\square\)
\end{enumerate}

\textbf{Finite circle envelopes.} The literal offsets in
\texttt{research/\allowbreak{}six-\allowbreak{}circle-\allowbreak{}recipes.\allowbreak{}json} supply the theorem's escape
premises at every order from 16 through 42. They are finite construction
data; constructing a graph reads those data and performs no
optimization. The direct checker separately verifies every translated
triple against every nonzero pursuer position.

At some small orders, the fixed three-triple test stops short of the
next construction. One can instead enlarge the envelope and make some
newly available offsets compulsory. The messenger can then use the new
arrows to counter the pursuer's new moves. These finitely many promoted
envelopes carry a different certificate: a literal nonnegative
source-rank vector \(r\) for recipient zero. It has \(r(0)=0\), and
every other source either has a direct compulsory arrow to zero or,
against each distinct pursuer, has a compulsory successor of smaller
rank outside the pursuer's closed envelope. The usual finite-rank
induction proves actual delivery. For recipient \(t\), use \(r(s-t)\);
the checker also replays every target explicitly. This applies to all
optional subsets and does not infer intermediate-budget success merely
from successful endpoints.

For example, at order 24 the six-offset envelope is
\[A_0=\{0,1,4,5,7,8,9,12,13,16,17,20,21\},\] covering
\(144\le m\le288\). Add offset three to the envelope and take every
nonzero offset in the enlarged envelope except offset one as compulsory.
The resulting \(288\le
m\le312\) family has the rank vector \[\begin{split}
(&0,3,4,1,1,3,6,1,1,3,6,1,1,5,2,1,1,1,4,1,1,1,4,7).
\end{split}\] The only optional arrows in this second envelope have
difference one. Thus every budget between 144 and 312 is covered by a
direct finite recipe.

In the same way, the six-offset and promoted envelopes together provide
every \(6n\le m\le U_n\) for orders 19 through 27, where \[
(U_{19},U_{20},U_{21},U_{22},U_{23},U_{24},U_{25},U_{26},U_{27})
=(228,240,252,264,299,312,350,364,378).
\] The supplied source ranks, rather than the three-triple theorem,
justify the promoted portions. This is a finite certificate boundary,
not an assertion that an optimization result proves an infinite family.

The direct constructor is \texttt{src/\allowbreak{}main\_\allowbreak{}interval\_\allowbreak{}six\_\allowbreak{}circle.\allowbreak{}py}.
Its check replays the escape clauses or supplied ranks, the four
propagation identities, exact arc counts, and every recipient for the
promoted rows. \texttt{research/\allowbreak{}six-\allowbreak{}circle-\allowbreak{}checks.\allowbreak{}json} records the
executable checks. \texttt{CircularSourceRanks.\allowbreak{}lean} and
\texttt{SixCircleFinite.\allowbreak{}lean} now instantiate the retained finite
source-rank rows and their exact optional-prefix intervals. The general
three-escape-triple theorem retains its ordinary parameterized proof
boundary.

\subsection{A smaller incoming and outgoing hub
interface}\label{a-smaller-incoming-and-outgoing-hub-interface}

The guide attachment lemma also permits fewer hub arrows at both ends.
Only its fixed interface proof changes: it now needs four moves,
followed by the same single safe-entry move from a bulk vertex.

\textbf{Theorem (reduced twelve-hub interface).} Let \(2\le g\le7\).
Assign the \(b\) bulk vertices any distinct guide subsets
\(J_x\subseteq\{0,\ldots,g-1\}\) which have size at least two and form
an antichain. Put \(n=12+g+b\) and \(W=\sum_x|J_x|\). There is an
explicit attaining construction at every budget \[
 A\le m\le A+b(g+b+8),\qquad A=106+18g+12b+W.
\] Every subset of the specified optional arrows is permitted, and
delivery takes at most five messenger moves. In particular, distinct
size-\(r\) codes give \(A=106+18g+(12+r)b\) whenever \(2\le r\le g\le7\)
and \(0\le b\le\binom gr\).

\textbf{Construction.} Label the hubs \(P_i=i\) and \(L_i=6+i\), with
indices modulo six. Their missing arrows are
\[P_i\to L_i,L_{i-1},\qquad L_i\to P_{i+4},P_{i+6};\] all other distinct
hub pairs are arrows. Use the first \(g\) guide codes \[
 \{L_0,L_2,L_4\},\quad\{L_1,L_3,L_5\},\quad
 \{P_0,P_3\},\quad\{P_1,P_4\},\quad\{P_2,P_5\},\quad
 \{L_0,L_3\},\quad\{L_1,L_4\}.
\] The parameters of the guide attachment lemma are \[
\begin{split}
 R&=\{P_1,P_2,P_3,P_4,L_0,L_1,L_3,L_4\},\\
 S&=\{P_1,P_4,L_1,L_3\},\qquad T=\{P_0,P_2,P_3\},\\
 U&=H\setminus T,\qquad C=H\setminus R.
\end{split}
\] Thus precisely \(R\) reaches every guide compulsorily, each guide
reaches the hubs outside its code, and bulk vertex \(x\) reaches
\(S\cup J_x\). The optional arrows are all guide-to-bulk and distinct
bulk-to-bulk arrows, all arrows from \(H\setminus R\) to bulk vertices,
and all arrows from a bulk vertex to \(H\setminus(T\cup S)\). In
particular a bulk pursuer may reach nine hubs, whereas the messenger
needs only four.

For equal-size codes, take the first \(b\) lexicographic \(r\)-subsets.
For an exact budget, take the first \(m-A\) optional arrows in
ordered-endpoint order. These are direct finite listing rules with no
game search.

\textbf{Proof.} We verify the attachment lemma's interface hypothesis
with \(d=4\), allowing a bulk pursuer every guide and every hub outside
\(T\). A hub move is safe against hub \(c\) exactly when it belongs to
the missing row \(F^+(c)\), against a guide when it belongs to that
guide's code, and against a bulk pursuer whenever it belongs to \(T\).
All nonterminal intermediaries below are hubs.

For a hub recipient, Theorem 59 supplies the two-move proof against hub
and guide pursuers; removing hub-to-outside arrows cannot affect its hub
witnesses. Against a bulk pursuer, a missing hub pair's danger set has
only two points, so it misses a member of \(T\). For a point-code guide
source, a missing recipient is a point; at least one of \(T\)'s three
points lies outside its two-element code and reaches the recipient
through the point clique. For a line-code guide, all three safe points
are reachable, and at most two fail to reach a line recipient. Thus all
interface sources reach every hub recipient within two moves.

For an outside recipient, the hubs in \(R\) reach it directly. The
following table specifies the further source-rank layers and their
reachable, already winning hub sets, using the numeric hub labels.

{\def\LTcaptype{none} 
\Needspace{10\baselineskip}
\begin{longtable}[]{@{}
  >{\raggedright\arraybackslash}p{(\linewidth - 4\tabcolsep) * \real{0.3333}}
  >{\raggedleft\arraybackslash}p{(\linewidth - 4\tabcolsep) * \real{0.3333}}
  >{\raggedright\arraybackslash}p{(\linewidth - 4\tabcolsep) * \real{0.3333}}@{}}
\toprule\noalign{}
\begin{minipage}[b]{\linewidth}\raggedright
New sources
\end{minipage} & \begin{minipage}[b]{\linewidth}\raggedleft
Delivery bound
\end{minipage} & \begin{minipage}[b]{\linewidth}\raggedright
Reachable already winning hubs
\end{minipage} \\
\midrule\noalign{}
\endhead
\bottomrule\noalign{}
\endlastfoot
\(P_0\) & 2 & \(\{1,2,3,4,7,9,10\}\) \\
\(P_5\) & 2 & \(\{1,2,3,4,6,7,9\}\) \\
\(L_2\) & 3 & \(\{1,3,4,5,6,7,9,10\}\) \\
\(L_5\) & 3 & \(\{0,1,2,4,6,7,9,10\}\) \\
Point-code guide with code \(K\) & 3 &
\((H\setminus\{L_2,L_5\})\setminus K\) \\
Remaining guide with code \(K\) & 4 & \(H\setminus K\) \\
\end{longtable}
}

Each displayed set meets the missing row of every hub pursuer other than
the source, every other guide code, and \(T\). These are exactly the
safe-intermediary tests, with a source's own pursuer position excluded.
We verify the tests explicitly. In the first two rows, the four point
members meet every difference-two pair; the listed lines meet every
adjacent pair except the source's own missing row. Both rows meet the
first guide triple in \(L_4\) or \(L_0\), respectively, the second in
\(L_1,L_3\), all three point codes, and the two retained line-pair
codes. They meet \(T\) in \(P_2,P_3\).

In the next two rows the only omitted lines form the antipodal pair
\(\{L_2,L_5\}\), so every adjacent line pair is met. Their point members
meet every difference-two pair except the source's own row. Directly
reading the seven codes shows that both sets meet each; they meet \(T\)
in \(P_3\) and in \(P_0,P_2\), respectively. For a point-code guide,
removing its antipodal point pair leaves a point in every difference-two
pair. The four retained lines meet every adjacent pair and every line
code; the other point codes are disjoint from its own. At least one
member of \(T\) remains. Finally \(H\setminus K\) for a remaining guide
meets every hub missing row because \(K\) contains neither an adjacent
line pair nor a difference-two point pair. It meets every other guide
code by incomparability, and contains all of \(T\), since \(K\) is a
line code. This verifies the table.

Every safe hub has smaller delivery rank, independently of the cop's
reply. Rank induction proves the four-move interface, including against
the allowed bulk pursuer. Optional hub-to-bulk arrows do not change any
missing hub row. Finally \(S\) meets exactly the eight hub missing rows
indexed by \(R\); the four exceptions are \(P_0,P_5,L_2,L_5\). This is
the attachment lemma's remaining hypothesis.

Apply that lemma with \(h=12\), \(a=108\), \(|R|=8\), \(|S|=4\),
\(W_K=2g+2\), \(W_J=W\), \(|C|=4\), and \(|U\setminus S|=5\). It gives
the compulsory count \[108+20g-(2g+2)+12b+W=A,\] the optional count
\(b(g+b+8)\) and the five-move bound. Every missing pair therefore
delivers, attaining the missing-pair upper bound throughout the
interval. \(\square\)

For pair codes the base is \(14n-62+4g\). At \((n,g,b)=(25,5,8)\) the
interval is \(308\le m\le476\); at \((27,5,10)\) it is
\(336\le m\le566\). These cover, respectively, the earlier outstanding
budgets 326 through 339 and 352 through 375. The optional-subset
guarantee and five-move bound preserve the earlier family's stronger
three-move claim as a separate result.

The eighth code \(\{L_2,L_5\}\) is deliberately absent: it misses the
first two safe-intermediary sets. This is an obstruction to that
interface, not to every eight-guide construction or to the game. The
direct constructor is \texttt{src/\allowbreak{}main\_\allowbreak{}interval\_\allowbreak{}reduced\_\allowbreak{}hubs.\allowbreak{}py}. Its
checks replay the fixed hub ranks, safe-entry premises, full-capacity
equal-weight families, a mixed antichain and the two targeted intervals.
\texttt{ReducedTwelveHubGraph.\allowbreak{}lean} and
\texttt{ReducedTwelveHubCounts.\allowbreak{}lean} now instantiate the reduced
twelve-hub equal-weight family, with numeric-mask code selection, exact
counts and optional-prefix maxima. Other code variants and supplementary
timing claims retain their stated boundaries.

\subsection{The complete main interval from order
twelve}\label{the-complete-main-interval-from-order-twelve}

\textbf{Theorem (main interval).} For every integer \(n\ge12\) and every
\(2n\le m\le n(n-3)\), \[M(n,m)=n(n-1)-m.\] There is an explicit
attaining graph at every pair. Its construction uses arithmetic families
and fixed finite arrow lists; it performs no game search or
derandomization. Finite lists cover a bounded bridge; the circle and
separated-pair constructions cover all budgets above \(3n\) from order
79. The finite-values section treats smaller orders and the other budget
ranges.

\subsubsection{Stages with changing compulsory
arrows}\label{stages-with-changing-compulsory-arrows}

Let \(B\subseteq E\) be two graphs on the same vertices. For each
recipient \(t\) suppose that every nonterminal, noncollision position
\((s,c)\) has a legal compulsory messenger move \(u\), and either
\(u=t\), or \(u\ne c\) and every legal cop reply \(c'\) in \(E\)
satisfies \[u\ne c',\qquad 0\le R_t(u,c')<R_t(s,c).\] Here \(R_t\) is a
nonnegative integer rank on the noncollision positions. Receipt takes
precedence over collision at \(t\), as throughout this paper. Induction
on this rank proves delivery in every graph \(B\subseteq G\subseteq E\).
The messenger uses only \(B\), while the cop is allowed all of \(E\).

One may now replace \(B\) between successive budget intervals. In
particular, if \(G_0\subseteq\cdots\subseteq G_q\) and each
\(G_i/G_{i+1}\) pair has such a certificate, every budget from
\(|E(G_0)|\) through \(|E(G_q)|\) is supplied. Earlier additions are
available to the messenger in subsequent stages. No strategy is required
to survive all final cop arrows using only the original sparse base.
Every optional subset within each certified stage works; no
arbitrary-subset assertion across different stages is inferred.

The same observation works when discovering graphs by deleting arrows
from a dense example. Read the resulting list backwards. Discovery
direction is irrelevant to the proof: every saved certificate is for a
literal lower base and upper envelope, with all intervening budgets
obtained by prefixes.

\subsubsection{Finite bridge and explicit
recipe}\label{finite-bridge-and-explicit-recipe}

The construction data in \texttt{research/\allowbreak{}sparse-\allowbreak{}staged-\allowbreak{}records.\allowbreak{}json}
specify a base arrow list, an ordered list of additional arrows, and
certified prefix intervals. At budget \(m\), take exactly the first
\(m-|E(B_0)|\) arrows in that list. A stored stage containing \(m\)
supplies the delivery proof. The selected base need not be the same
graph as a previously known construction at an overlapping endpoint;
each interval has its own complete certificate.

The finite certificates check both positions, waits, receipt before
collision, and strict descent after every legal cop reply. Only their
stated intervals are used: a stored arrow list can extend beyond its
certified stages, and those extra prefixes are not inferred to work.
Combine these domains.

\begin{enumerate}
\def\labelenumi{\arabic{enumi}.}
\tightlist
\item
  For \(12\le n\le78\), use the finite stages, protected and
  five/six-offset circles, layered and reduced hub interfaces, and
  separated pair codes. The existing literal dense envelopes at orders
  13 and 14 supply the remaining budgets \(129\) and
  \(147,148,149,151\), respectively. The independent integer-interval
  check covers this entire finite range. No general larger-code capacity
  theorem or dense connector is needed by this cover.
\item
  For \(79\le n\le114\), the pair-code overlap theorem covers every
  budget from \(3n\). The retained sparse stages and certified sparse
  envelopes cover \(2n\) through \(3n\), including the repaired
  order-109 envelope.
\item
  For \(n\ge115\), Theorem 52 covers \(2n\) through \(3n\), and the
  pair-code overlap theorem covers the remainder through \(n(n-3)\).
\end{enumerate}

\texttt{src/\allowbreak{}main\_\allowbreak{}interval\_\allowbreak{}complete.\allowbreak{}py} selects these intervals and
constructs their adjacency.
\texttt{src/\allowbreak{}main\_\allowbreak{}interval\_\allowbreak{}pair\_\allowbreak{}codes\_\allowbreak{}check.\allowbreak{}py} independently
reconstructs the finite interval union from formulas and saved stage
endpoints, and checks actual constructor outputs. The current
source-pinned receipts are
\texttt{research/\allowbreak{}main-\allowbreak{}interval-\allowbreak{}pair-\allowbreak{}code-\allowbreak{}checks.\allowbreak{}json} and
\texttt{research/\allowbreak{}main-\allowbreak{}interval-\allowbreak{}pair-\allowbreak{}complete-\allowbreak{}checks.\allowbreak{}json}. Earlier
constructors, certificate files, stronger guarantees, and historical
receipts remain intact.

Since every resulting missing nonloop pair is guaranteed, each graph has
\(I(G)=n(n-1)-m\), matching the elementary missing-pair upper bound.
This proves the theorem. The finite bridge is computer-assisted; the
circle and pair-code overlap have ordinary proofs. The sparse infinite
family retains its stated finite local certificate ingredient.

\subsubsection{Formal verification and possible
simplifications}\label{formal-verification-and-possible-simplifications}

\texttt{lean/\allowbreak{}StagedEnvelopes.\allowbreak{}lean} proves the generic asymmetric
rank-envelope implication for every intermediate graph and arbitrary
legal cop play. Concrete stage lists are instantiated at orders12
through18 in \texttt{MainTwelve.\allowbreak{}lean} and
\texttt{MainSmallClassification.\allowbreak{}lean}. \texttt{LowBandConstruction.\allowbreak{}lean}
supplies the uniform lower interval from115, and
\texttt{MainUpperInterval.\allowbreak{}lean} the interval above \(3n\) from25. The
remaining finite envelopes are now joined in
\texttt{MainIntervalClassification.\allowbreak{}lean} and
\texttt{MainExplicitAttainer.\allowbreak{}lean}: the entire main interval, with a
prescribed graph, is formal at every order from12. The complete
classification still has the separate finite upper obligations stated in
the verification section.

A single repeatable small attachment rule remains a possible
simplification, not a missing step in this finite-bridge theorem. Two
fixed attachments transported whole budget families from order 25 to 26
and 26 to 27. Reusing only the newest vertex's ports subsequently
creates a recipient whose predecessors one cop can cover. That
obstruction rules out the proposed repeated tip attachment; it does not
obstruct the explicit constructions above. The private extension
research is retained separately from the main theorem.

\section{Dense constructions and the final
staircase}\label{dense-constructions-and-the-final-staircase}

\subsection{The dense boundary}\label{the-dense-boundary}

\textbf{Theorem 13.} For \(n\ge2\),
\[M(n,n(n-2))=\left\lfloor\frac{n-2}{2}\right\rfloor,\] and \(M(n,m)=0\)
for every admissible \(m>n(n-2)\).

\textbf{Proof.} If \((s,t)\) is indirect guaranteed, row \(s\) misses
\(s\to t\). By Lemma 1, every other row \(c\) must miss an arc into
\(N^+(s)\). Thus every row misses at least one nonloop arc, and
\(m\le n(n-2)\).

Now suppose \(m=n(n-2)\). If \(I=0\), the upper bound is immediate.
Otherwise the preceding argument shows that every row misses exactly one
arc; write \(f(v)\) for its endpoint. Necessarily \(f(v)\ne v\). A
possible indirect pair is \((s,f(s))\). Lemma 1 says
\[ f^{-1}(s)=\varnothing,\qquad f^{-1}(f(s))=\{s\}. \tag{2}\] These
conditions suffice: for initial pursuer \(c\ne s\), choose the
intermediate \(f(c)\). It is different from \(s,f(s)\), is inaccessible
on the pursuer's next move, and has an arc to \(f(s)\). This is a
two-move guarantee.

Pairs satisfying (2) have disjoint source and target vertices. If there
are \(j\) pairs, the remaining vertices contain a cycle of the loopless
functional graph of \(f\) and therefore number at least two. Hence
\(2j\le n-2\). For attainment use pairs \(s_i\to t_i\) in the
missing-arc graph, send every \(t_i\) to \(a\), put
\(a\leftrightarrow b\), and send any remaining vertex to \(a\). Take all
nonloop arcs except these missing arcs. \(\square\)

\begin{figure}
\centering
\pandocbounded{\includegraphics[keepaspectratio,alt={The missing arcs at the dense boundary form a functional graph. Each displayed two-vertex branch contributes one guaranteed indirect pair.}]{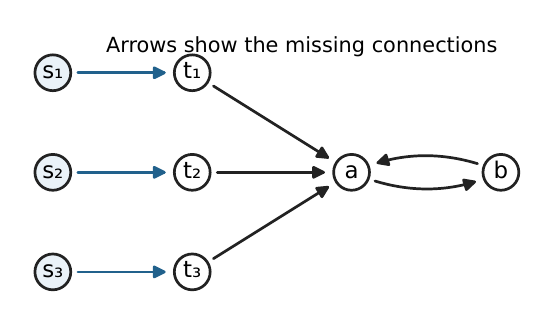}}
\caption{The missing arcs at the dense boundary form a functional graph.
Each displayed two-vertex branch contributes one guaranteed indirect
pair.}
\end{figure}

\textbf{Theorem 14.} A noncomplete graph with every pair guaranteed has
at most \(n(n-3)\) arcs.

\textbf{Proof.} Such a graph has an indirect pair. No row can be
complete: its vertex would block every indirect source other than itself
and have no indirect target of its own. If a row \(c\) missed only
\(c\to t\), then \(N^+(t)\subseteq C(c)\). Vertex \(t\) also has an
indirect recipient, contradicting Lemma 1. Every row therefore misses at
least two arcs. \(\square\)

The next section gives a uniform construction attaining this density on
an infinite set of orders. Its proof is independent of the sparse-game
certificate for the complementary graph proved later.

\subsubsection{A periodic family with three-move
delivery}\label{a-periodic-family-with-three-move-delivery}

Let \(n=5k\), \(k\ge4\), and write vertices as \(5i+j\), with \(i\)
modulo \(k\) and \(0\le j\le4\). Put \(\alpha(v)=v+1\) and \[
\begin{aligned}
\beta(5i)&=5(i+1),&\beta(5i+1)&=5(i-1)+2,\\
\beta(5i+2)&=5(i-1)+4,&\beta(5i+3)&=5(i-1)+3,\\
\beta(5i+4)&=5(i-1)+1.
\end{aligned}
\] Let \(H\) have arcs \(v\to\alpha(v),v\to\beta(v)\), and let \(G\) be
its loopless directed complement. Both maps are permutations, without
fixed points or coincident images at a vertex; thus \(G\) has \(n(n-3)\)
arcs.

\textbf{Theorem 15.} Every pair in \(G\) is guaranteed within three
messenger moves. Consequently \(M(5k,5k(5k-3))=10k\) for every
\(k\ge4\).

\textbf{Proof.} For a missing pair \((s,t)\) and initial pursuer
\(c\ne s\), a safe two-move intermediate is any member of
\[\{\alpha(c),\beta(c)\}\setminus F(s,t),\qquad
F(s,t)=\{s,t,\alpha(s),\beta(s),\alpha^{-1}(t),\beta^{-1}(t)\}.\] Its
exclusion from \(F\) ensures the two required arcs in \(G\), while its
membership in the displayed pair makes it inaccessible to the pursuer's
next move.

Translate the source by a multiple of five so that \(0\le s\le4\). A
failure requires \(c\in F-1\), since \(\alpha(c)\in F\). The complete
local calculation is:

{\def\LTcaptype{none} 
\Needspace{14\baselineskip}
\begin{longtable}[]{@{}rrll@{}}
\toprule\noalign{}
\(s\) & \(t\) & \(F(s,t)\) & \((c,\beta(c))\), with \(c\in F-1\),
\(c\ne s\) \\
\midrule\noalign{}
\endhead
\bottomrule\noalign{}
\endlastfoot
0 & 1 & \(\{0,1,5,9\}\) & \((-1,-9),(4,-4),(8,3)\) \\
0 & 5 & \(\{0,1,4,5\}\) & \((-1,-9),(3,-2),(4,-4)\) \\
1 & 2 & \(\{-3,1,2,6\}\) & \((-4,-8),(0,5),(5,10)\) \\
1 & \(-3\) & \(\{-4,-3,1,2\}\) & \((-5,0),(-4,-8),(0,5)\) \\
2 & 3 & \(\{-1,2,3,8\}\) & \((-2,-7),(1,-3),(7,4)\) \\
2 & \(-1\) & \(\{-2,-1,2,3\}\) & \((-3,-6),(-2,-7),(1,-3)\) \\
3 & 4 & \(\{-2,3,4,7\}\) & \((-3,-6),(2,-1),(6,2)\) \\
3 & \(-2\) & \(\{-3,-2,3,4\}\) & \((-4,-8),(-3,-6),(2,-1)\) \\
4 & 5 & \(\{-4,0,4,5\}\) & \((-5,0),(-1,-9),(3,-2)\) \\
4 & \(-4\) & \(\{-5,-4,4,5\}\) & \((-6,-14),(-5,0),(3,-2)\) \\
\end{longtable}
}

The only candidate with \(\beta(c)\in F\) is \((s,t,c)=(4,5,-5)\). Every
difference being tested has absolute value at most 19, so no extra
equality appears modulo any \(5k\ge20\). In that exceptional state move
from 4 to 0. This is legal and safe: the pursuer at \(-5\) has no arc to
0. Afterwards the pair \((0,5)\) is covered by the nonexceptional second
row for every surviving pursuer position, and delivery takes at most two
more moves. \(\square\)

\subsection{Direct missing graphs unify the dense
constructions}\label{direct-missing-graphs-unify-the-dense-constructions}

The following direct constructions replace the high-girth existence
arguments for the numerical dense results. Its proof examines the
possible safe intermediaries directly, so no girth hypothesis is needed.
The parameterized families have ordinary proofs; a few small orders use
explicitly listed finite constructions and independently checked
witnesses. All strategies here take at most two messenger moves. Write
\(L=n(n-1)\) and \(h=L-m\).

For a missing-arc graph \(F\), a missing pair \((s,t)\) has a safe
intermediary against pursuer \(c\ne s\) precisely when
\[N_F^+(c)\setminus\bigl(N_F^+(s)\cup N_F^-(t)\bigr)\ne\varnothing.\] We
call the union in this expression the danger set. It contains both
endpoints. A vertex outside it supplies both required messenger arcs,
while being a missing outgoing neighbor of the pursuer makes it safe.

There is also a robust version. Suppose \(B\subseteq F\) are missing
graphs, and for every arc \((s,t)\) of \(F\) and every \(c\ne s\),
\[N_B^+(c)\not\subseteq N_F^+(s)\cup N_F^-(t). \tag{R}\] Then every
intermediate missing graph \(B\subseteq F'\subseteq F\) guarantees all
pairs in two moves: retain the witness from \(B\), while the danger set
only shrinks. In particular, optional arcs can be chosen individually
and in arbitrary order.

\textbf{Theorem 58 (a unified dense construction).} There are
deterministic polynomial constructions with the following properties.

\begin{enumerate}
\def\labelenumi{\arabic{enumi}.}
\tightlist
\item
  At every \(n\ge10\), a graph with \(n(n-3)\) arcs has all \(2n\)
  missing pairs guaranteed within two moves. Thus \(M(n,n(n-3))=2n\).
\item
  At every \(n\ge16\), one compulsory set of \(2n\) missing arcs and a
  list of \(n\) optional missing arcs satisfy (R). Hence every integer
  \(2n\le h\le3n\) attains \(M(n,L-h)=h\), with two-move delivery.
\item
  At every \(n\ge20\), the entire final dense band is exact:
  \[M\bigl(n,n(n-3)+r\bigr)
  =2n-r-\left\lceil\frac r2\right\rceil-1,
  \qquad1\le r\le n.\] Every counted pair in an attaining construction
  delivers within two moves. The first step \(r=1\) already holds at
  every \(n\ge16\).
\item
  For \(n=10q+s\ge30\), \(0\le s\le9\), put \[H_5(n)=10q^2+5q+s^2+s.\] A
  compulsory set of \(2n\) missing arcs extends to an envelope with
  \(H_5(n)\) missing arcs satisfying (R). For \(n=4q+s\ge44\),
  \(0\le s\le3\), the same assertion holds with
  \[H_4(n)=4q^2-26q+s^2+s.\] Thus every integer \(2n\le h\le H(n)\) has
  \(M(n,L-h)=h\), where \(H(n)=H_5(n)\) for \(30\le n<44\) and
  \(H(n)=\max\{H_5(n),H_4(n)\}\) for \(n\ge44\). Every optional subset
  works, and every missing pair delivers within two moves.
\end{enumerate}

These statements determine every missing budget \(n\le h\le H(n)\) at
every \(n\ge30\), and every \(n\le h\le3n\) at every \(n\ge20\). The
value is the displayed staircase for \(h<2n\) and \(h\) for \(h\ge2n\).
The new cap strictly exceeds the cap of Theorem 48 throughout their
common domain \(n\ge30\), preserving its time and arbitrary-subset
guarantees; part 2 also covers all of its smaller orders. From \(n=52\)
onward, \(H(n)=H_4(n)=n^2/4-O(n)\).

\subsubsection{The two-type core}\label{the-two-type-core}

Use vertices \(P_i,L_i\), with \(i\) modulo \(k\), and the compulsory
arcs \[
P_i\longrightarrow L_i,L_{i-1},\qquad
L_i\longrightarrow P_{i+4},P_{i+6}. \tag{B}
\] There are exactly two incoming and two outgoing arcs at every vertex.
The two types are disjoint, so no loop occurs. For \(k\ge5\), all pairs
in the directed complement are guaranteed in two moves.

To see this, translate the source index to zero. A pursuer \(P_i\) is
blocked precisely when the line part of the danger set contains both
\(L_i,L_{i-1}\). A pursuer \(L_i\) is blocked precisely when its point
part contains both \(P_{i+4},P_{i+6}\). Thus line indices are tested for
adjacency, and point indices for a difference of two. The four possible
missing pairs have the following point and line danger sets:

\begin{itemize}
\tightlist
\item
  For \(P_0\to L_0\): point indices \(\{0,1\}\), line indices
  \(\{0,-1\}\).
\item
  For \(P_0\to L_{-1}\): point indices \(\{-1,0\}\), line indices
  \(\{0,-1\}\).
\item
  For \(L_0\to P_4\): point indices \(\{4,6\}\), line indices
  \(\{0,-2\}\).
\item
  For \(L_0\to P_6\): point indices \(\{4,6\}\), line indices
  \(\{0,2\}\).
\end{itemize}

For \(k\ge5\), the only blocked pursuer in each case is the source
itself, which is not a legal initial position. This gives every even
order \(n=2k\ge10\).

For odd order \(n=2k+1\ge13\), introduce a vertex \(w\). Replace the
arcs \(P_0\to L_0\) and \(P_3\to L_3\) by
\[P_0\to w\to L_0,\qquad P_3\to w\to L_3. \tag{O}
\] Each degree is still two. We verify the two-move property directly.
For an old missing pair and a pursuer other than \(P_0,P_3,w\), its old
witnesses are unchanged: the danger set only gains the new vertex \(w\)
and may lose old vertices. The pursuer at \(w\) has missing neighbors
\(L_0,L_3\), which cannot both occur in any of the four old line danger
sets, whose only nonzero differences are one or two.

A pursuer at \(P_a\), \(a\in\{0,3\}\), can use \(w\) unless the source
is the other modified point or the recipient is \(L_0\) or \(L_3\). In
those cases \(L_{a-1}\) is safe. Explicitly, for source \(P_b\) with
\(a\ne b\) its remaining old recipient is \(L_{b-1}\) and the line
danger set is \(\{L_{b-1}\}\). For old recipient \(L_b\), the remaining
source is \(P_{b+1}\) and the line danger set is \(\{L_b,L_{b+1}\}\).
None contains \(L_{a-1}\) when \(k\ge6\) and \(a,b\in\{0,3\}\).

For a new pair \(P_b\to w\), the danger set has point indices
\(\{0,3\}\) and the single line index \(b-1\), in addition to \(w\).
Neither an adjacent line pair nor a point pair of difference two fits;
the other modified point pursuer also retains its different line
neighbor. For a new pair \(w\to L_b\), the line indices are \(\{0,3\}\)
and the point index is \(b+1\), in addition to \(w\). Again all legal
pursuers have a safe neighbor.

The remaining order \(n=11\) has the following explicit missing rows:

\[\begin{aligned}
0&:\{4,8\},&1&:\{9,10\},&2&:\{8,9\},&3&:\{7,10\},\\
4&:\{0,1\},&5&:\{2,6\},&6&:\{3,5\},&7&:\{3,4\},\\
8&:\{2,7\},&9&:\{0,5\},&10&:\{1,6\}.
\end{aligned}\]

Every incoming and outgoing degree is two. For each source, its two
recipients in the displayed order give the following two danger sets:

\[\begin{array}{c|c|c}
s&D(s,t_1)&D(s,t_2)\\\hline
0&\{0,4,7,8\}&\{0,2,4,8\}\\
1&\{1,2,9,10\}&\{1,3,9,10\}\\
2&\{0,2,8,9\}&\{1,2,8,9\}\\
3&\{3,7,8,10\}&\{1,3,7,10\}\\
4&\{0,1,4,9\}&\{0,1,4,10\}\\
5&\{2,5,6,8\}&\{2,5,6,10\}\\
6&\{3,5,6,7\}&\{3,5,6,9\}\\
7&\{3,4,6,7\}&\{0,3,4,7\}\\
8&\{2,5,7,8\}&\{2,3,7,8\}\\
9&\{0,4,5,9\}&\{0,5,6,9\}\\
10&\{1,4,6,10\}&\{1,5,6,10\}
\end{array}\]

Inspection of the eleven rows shows that each danger set contains only
its source's row. Every legal pursuer therefore has a safe missing
neighbor. This finite set-containment argument proves the eleven-vertex
case and completes part 1. The accompanying certificate also lists all
220 pair/pursuer witnesses explicitly.

\subsubsection{An envelope with one optional arc per
vertex}\label{an-envelope-with-one-optional-arc-per-vertex}

For even order \(n=2k\), retain (B) as \(B\) and add the optional arcs
\[P_i\longrightarrow P_{i+6},\qquad
L_i\longrightarrow L_{i+4}. \tag{E}
\] For \(k\ge11\), condition (R) holds for the full envelope \(F\).
There are six types of missing pairs. Their point and line danger
indices are respectively:

\begin{itemize}
\tightlist
\item
  \(P_0\to L_0\): \(\{6,0,1\}\) and \(\{0,-1,-4\}\);
\item
  \(P_0\to L_{-1}\): \(\{6,-1,0\}\) and \(\{0,-1,-5\}\);
\item
  \(P_0\to P_6\): \(\{6,0\}\) and \(\{0,-1,2\}\);
\item
  \(L_0\to P_4\): \(\{4,6,-2\}\) and \(\{4,0,-2\}\);
\item
  \(L_0\to P_6\): \(\{4,6,0\}\) and \(\{4,2,0\}\);
\item
  \(L_0\to L_4\): \(\{4,6,5\}\) and \(\{4,0\}\).
\end{itemize}

Test line adjacency and point differences of two as before. The only
blocked compulsory pair belongs to the source. All index differences
tested have absolute value at most eight; \(k\ge11\) introduces no
additional modular adjacency or difference of two. Thus (R) holds.

For odd order \(n=2k+1\), with \(k\ge16\) or \(k=14\), use a different
placement of the new vertex \(w\): replace \(P_0\to L_0\) and
\(P_7\to L_7\) by paths through \(w\). Retain every optional arc in (E),
and add the single optional arc \(w\to P_2\). The compulsory graph has
\(2n\) arcs and the full envelope \(3n\) arcs. We check (R).

For an old pair, every unmodified old pursuer retains its old safe
witness because the only new danger vertex is \(w\). The compulsory
neighbors of the pursuer at \(w\) are \(L_0,L_7\). No translated line
danger set in the six-case list contains both: their integer differences
have absolute value at most six, while seven is not congruent to any
such nonzero difference modulo the stated \(k\).

The remaining old pursuers are \(P_a\), \(a\in\{0,7\}\), with compulsory
neighbors \(w,L_{a-1}\). If \(w\) is in the danger set and the source
differs from \(P_a\), only the following cases arise:

\begin{itemize}
\tightlist
\item
  The source is the other modified point \(P_b\). Its retained old
  targets are \(L_{b-1}\) and \(P_{b+6}\). The respective line danger
  sets are \(\{b-1,b-5\}\) and \(\{b-1,b,b+2\}\).
\item
  The recipient is \(L_b\), \(b\in\{0,7\}\). The remaining old sources
  are \(P_{b+1}\) and \(L_{b-4}\). Their line danger sets are
  \(\{b+1,b,b-4\}\) and \(\{b,b-4\}\).
\item
  The recipient is \(P_2\). Its old sources are
  \(L_{-2},L_{-4},P_{-4}\), and their line danger sets are contained in
  \(\{2,0,-2,-4,-5\}\).
\end{itemize}

None contains \(a-1\) in its applicable case, for \(k\ge16\) or
\(k=14\). Thus the remaining compulsory neighbor is safe.

It remains to check the new missing pairs. In the following lists, the
danger set also contains \(w\):

\begin{itemize}
\tightlist
\item
  For \(P_b\to w\), \(b\in\{0,7\}\), the point indices are
  \(\{0,7,b+6\}\) and the only line index is \(b-1\).
\item
  For \(w\to L_b\), the point indices are \(\{2,b+1\}\) and the line
  indices are \(\{0,7,b-4\}\).
\item
  For \(w\to P_2\), the point indices are \(\{2,-4\}\) and the line
  indices are \(\{0,7,-2,-4\}\).
\end{itemize}

The point sets contain no pair of difference two and the line sets no
adjacent pair at the stated moduli. For example, the largest relevant
point difference is 13, so it cannot become two when \(k\ge16\); modulo
14 it becomes \(-1\), which is also harmless. The line differences have
absolute value at most eleven, so they introduce no modular adjacency. A
modified point pursuer, when legal, retains a line neighbor \(L_{-1}\)
or \(L_6\) outside its danger set; the pursuer at \(w\) is either the
source or sees only one of its two line neighbors in that set. This
proves (R) at the stated odd orders.

Four explicit repairs fill the remaining odd orders. On \(2k+1\)
vertices, subdivide \(P_0\to L_0\) and \(P_j\to L_j\) through \(w\) as
above, retain the optional arcs (E), and add \(w\to P_z\). Make the
following substitutions in the optional list; every index is modulo
\(k\):

\begin{itemize}
\tightlist
\item
  At \(n=23\), use \((k,j,z)=(11,5,0)\), replacing \(L_1\to L_5\) by
  \(L_1\to L_7\) and \(L_7\to L_0\) by \(L_7\to L_2\).
\item
  At \(n=25\), use \((k,j,z)=(12,5,2)\), replacing \(L_1\to L_5\) by
  \(L_1\to P_{10}\).
\item
  At \(n=27\), use \((k,j,z)=(13,6,2)\), replacing \(L_2\to L_6\) by
  \(L_2\to P_{12}\).
\item
  At \(n=31\), use \((k,j,z)=(15,5,0)\), replacing \(L_1\to L_5\) by
  \(L_1\to P_{10}\).
\end{itemize}

In each case the compulsory set has \(2n\) arcs and the optional list
has \(n\) distinct additional arcs. These are finite constructions:
\texttt{research/\allowbreak{}dense-\allowbreak{}phase2-\allowbreak{}certificates.\allowbreak{}json} supplies, for every
full-envelope pair and every legal pursuer, a compulsory neighbor
outside its danger set. The direct checker
\texttt{src/\allowbreak{}dense\_\allowbreak{}phase2\_\allowbreak{}certificates.\allowbreak{}py} verifies each listed witness
using literal legal-move sets; an independent check is recorded in
\texttt{research/\allowbreak{}dense-\allowbreak{}phase2-\allowbreak{}independent-\allowbreak{}root-\allowbreak{}check.\allowbreak{}json}. The
certificate proves (R), so it covers every optional subset rather than
just the chosen total counts. Together with the parameterized families,
this establishes the bridge at every \(n\ge22\).

\paragraph{Six smaller envelopes}\label{six-smaller-envelopes}

For \(n=16,\ldots,21\), use the compulsory core (B), or its odd
modification (O). Identify \(P_i=i\), \(L_i=k+i\), and, in odd order,
\(w=2k\), where \(k=\lfloor n/2\rfloor\). The following ordered pairs
are optional missing arcs:

\[\begin{aligned}
O_{16}=\{&(0,12),(1,13),(2,13),(3,15),(4,15),(5,9),(6,10),\\
&(7,10),(8,1),(9,2),(10,3),(11,4),(12,5),(13,6),\\
&(14,7),(15,0)\}
\end{aligned}\]

\[\begin{aligned}
O_{17}=\{&(0,7),(0,11),(1,13),(2,14),(3,14),(4,16),(5,9),\\
&(6,8),(7,10),(8,3),(10,3),(12,5),(13,6),(14,10),\\
&(15,6),(16,1),(16,13)\}
\end{aligned}\]

\[\begin{aligned}
O_{18}=\{&(0,14),(1,14),(2,15),(4,16),(5,17),(6,9),(7,10),\\
&(8,11),(8,13),(10,1),(11,2),(12,3),(12,4),(13,5),\\
&(14,6),(14,9),(15,7),(17,8)\}
\end{aligned}\]

\[\begin{aligned}
O_{19}=\{&(0,3),(0,12),(1,2),(2,15),(3,15),(4,5),(5,9),\\
&(6,10),(6,18),(7,11),(8,11),(8,13),(12,18),(13,4),\\
&(14,6),(15,10),(16,9),(17,8),(17,13)\}
\end{aligned}\]

\[\begin{aligned}
O_{20}=\{&(0,13),(0,16),(1,14),(2,15),(3,18),(4,10),(4,18),\\
&(5,10),(9,12),(9,15),(10,9),(11,17),(12,2),(13,3),\\
&(13,17),(14,5),(15,6),(17,11),(18,14),(19,8)\}
\end{aligned}\]

\[\begin{aligned}
O_{21}=\{&(0,3),(1,14),(1,16),(2,18),(3,18),(6,11),(7,11),\\
&(8,12),(8,14),(9,12),(11,1),(11,15),(12,2),(13,20),\\
&(15,6),(15,19),(17,6),(17,7),(18,8),(18,10),(18,14)\}
\end{aligned}\]

Each list contains exactly \(n\) distinct arcs outside the compulsory
set of \(2n\) arcs. For every missing pair in the full envelope and
legal initial pursuer, the archived certificate explicitly supplies
\(u\in N_B^+(c)\) with \((s,u),(u,t)\notin F\) and \(u\notin\{s,t,c\}\).
These are precisely the obligations in (R). The archives are
\texttt{research/\allowbreak{}routing-\allowbreak{}defect-\allowbreak{}phase2-\allowbreak{}dense-\allowbreak{}tiny-\allowbreak{}checks.\allowbreak{}json} for
orders \(16,17,18\) and
\texttt{research/\allowbreak{}routing-\allowbreak{}defect-\allowbreak{}phase2-\allowbreak{}dense-\allowbreak{}small-\allowbreak{}checks.\allowbreak{}json} for
\(19,20,21\); at order 18 the archived full envelope is larger, so the
displayed subset inherits its witnesses.

Separate literal readback programs, importing neither the producer nor
its graph helpers, check every stored witness and its complete
pair/pursuer coverage. Their receipts are
\texttt{research/\allowbreak{}routing-\allowbreak{}defect-\allowbreak{}phase2-\allowbreak{}dense-\allowbreak{}tiny-\allowbreak{}independent.\allowbreak{}json} and
\texttt{research/\allowbreak{}routing-\allowbreak{}defect-\allowbreak{}phase2-\allowbreak{}dense-\allowbreak{}small-\allowbreak{}independent.\allowbreak{}json}. In
addition, \texttt{SmallDenseEnvelopes.\allowbreak{}lean} checks all six row lists,
the base and envelope arrow counts, and every required safe witness by
kernel reduction. The robust-envelope semantic theorem then proves
actual two-move delivery for every intermediate graph, against arbitrary
legal pursuer choices. This verifies every optional subset, not just a
single graph at \(h=3n\). Together with the preceding families, it
completes part 2 for every \(n\ge16\).

\subsubsection{Quadratic envelopes from a shared danger-set
calculation}\label{quadratic-envelopes-from-a-shared-danger-set-calculation}

The interval can be enlarged while keeping just two compulsory missing
neighbors per pursuer. In a core with indices modulo \(k\), use
\[P_i\to L_i,L_{i-1},\qquad L_i\to P_{i+\alpha},P_{i+\beta},\] and
optional arcs \(P_i\to P_{i+d}\) for \(d\in D\) and \(L_i\to L_{i+e}\)
for \(e\in E\). Translate the source index to zero. The six
full-envelope danger sets have the following point and line index parts:

\[\begin{array}{c|c|c}
\text{pair}&\text{point indices}&\text{line indices}\\\hline
P_0\to L_0&D\cup\{0,1\}&\{0,-1\}\cup(-E)\\
P_0\to L_{-1}&D\cup\{-1,0\}&\{0,-1\}\cup(-1-E)\\
P_0\to P_d&D\cup(d-D)&\{0,-1,d-\alpha,d-\beta\}\\
L_0\to P_\alpha&\{\alpha,\beta\}\cup(\alpha-D)&E\cup\{0,\alpha-\beta\}\\
L_0\to P_\beta&\{\alpha,\beta\}\cup(\beta-D)&E\cup\{0,\beta-\alpha\}\\
L_0\to L_e&\{\alpha,\beta,e,e+1\}&E\cup(e-E)
\end{array} \tag{T}\]

Condition (R) asks whether the line part contains an adjacent pair, or
the point part contains a compulsory pair of difference
\(\beta-\alpha\), belonging to a pursuer other than the source. This
single calculation will establish both caps in part 4.

\paragraph{A residue-five family at every order from
30}\label{a-residue-five-family-at-every-order-from-30}

Write \(n=10q+s\), where \(q\ge3\) and \(0\le s\le9\), and put \(k=5q\).
Take \((\alpha,\beta)=(4,6)\) and
\[D=\{10,15,\ldots,k-5\},\qquad E=\{5,10,\ldots,k-5\}.\] Add vertices
\(x_0,\ldots,x_{s-1}\) with compulsory missing rows
\[N_B^+(x_a)=\{P_a,L_{a+3}\}.\] There are no missing arcs from the core
to these vertices. Permit every directed nonloop arc between added
vertices as an additional optional missing arc.

First consider only old sources and pursuers. The multiples-of-five sets
in (T) contain neither adjacent indices nor indices of difference two.
In the four compulsory-pair rows of (T), a compulsory cop row can be
contained only when it is the source's own row. For \(P_0\to P_d\), the
line set \(\{0,-1,d-4,d-6\}\) acquires an additional adjacent pair only
at \(d=2,4,5,7\); values \(3\) and \(6\) merely duplicate the source
row. Only \(5\) is a multiple of five, and it is excluded from \(D\).
For \(L_0\to L_e\), the apparent extra difference-two pair at \(e=5\) is
just the existing pair \(\{4,6\}\). Thus old legal pursuers always have
a witness. These modular statements hold because \(5\mid k\) and
\(k\ge15\).

An added cop has point/line index difference three. For an optional
\(P\to P\) pair, the point danger has the source's residue, while its
line residues differ by \(0,1,4\); difference three is absent. For an
optional \(L\to L\) pair, the line danger has one residue and the point
residues differ from it by \(0,1,4\), again excluding the row. The four
compulsory cases in (T) give the same exclusion except possibly
\((P_{i+1},L_{i-1})\) or \((P_{i+4},L_{i+2})\). Their actual
line-minus-point difference is \(-2\), which cannot equal three modulo
\(k\ge15\). Hence added cops also retain witnesses.

For a new source \(x_a\) targeting \(P_a\), the old danger parts are
\[\{P_a\}\cup\{P_{a-d}:d\in D\},\qquad
\{L_{a+3},L_{a-4},L_{a-6}\}.\] The first lies in one residue class. The
second has pairwise differences \(7,9,2\), none congruent to \(\pm1\)
modulo \(k\ge15\). Neither old row type is blocked. An added row from
\(x_b\) would require \(b-a\) to be a multiple of five and one of
\(0,-7,-9\), forcing \(b=a\), the excluded source. For recipient
\(L_{a+3}\), the old point danger is \(\{P_a,P_{a+3},P_{a+4}\}\), whose
differences \(1,3,4\) cannot be \(\pm2\); the old line danger has one
residue class. An added cop would require \(b-a\in\{0,3,4\}\) to be a
multiple of five, again forcing \(b=a\).

Finally, an optional pair \(x_a\to x_b\) has old danger exactly
\(\{P_a,L_{a+3}\}\). It contains neither an old row nor a distinct new
row. Incoming arrows from added vertices only add new vertices to other
danger sets, whereas every compulsory row lies in the old core. They
cannot spoil any witness. This proves (R).

There are \(4k+2s=2n\) compulsory missing arcs, \(k((q-2)+(q-1))\)
optional core arcs, and \(s(s-1)\) optional added arcs. Their total is
exactly \(H_5(n)=10q^2+5q+s^2+s\).

\paragraph{A parity family with quarter-density
cap}\label{a-parity-family-with-quarter-density-cap}

Write \(n=4q+s\), \(q\ge11\), \(0\le s\le3\), and put \(k=2q\). Take
\((\alpha,\beta)=(4,7)\). Let \(D\) consist of all nonzero even residues
except \(\pm2,\pm4,\pm6,\pm8\), and let \(E\) consist of all nonzero
even residues except \(\pm2,\pm4,8,10,12\). For each \(0\le a<s\) add a
vertex with compulsory row \[N_B^+(x_a)=\{P_a,P_{a+5}\},\] and again
permit all nonloop missing arcs between added vertices, with none from
old vertices to new ones.

Both forbidden row differences in (T) are odd: one for a point cop,
three for a line cop. Since \(D,E\) contain only even indices, testing
the exceptional even/odd pairs in each row of (T) gives the following
possible offsets for a blocked old cop other than the source:

\begin{itemize}
\tightlist
\item
  \(P_0\to L_0\): \(4,-2\) in \(D\), or \(2\) in \(E\);
\item
  \(P_0\to L_{-1}\): \(2,-4\) in \(D\), or \(-2\) in \(E\);
\item
  \(P_0\to P_d\): \(d=2,8\);
\item
  \(L_0\to P_4\): \(-6\) in \(D\), or \(-2,-4\) in \(E\);
\item
  \(L_0\to P_7\): \(6\) in \(D\), or \(2,4\) in \(E\);
\item
  \(L_0\to L_e\): \(e=0,10\).
\end{itemize}

All these offsets are excluded. Coincidences that merely reproduce the
source row do not block another cop: adjacent pairs and pairs of
difference three have unique source indices when \(k\ge22\).

An added cop has a point pair of difference five. In the same six
point-danger sets, such a pair would require respectively offsets
\(6,-4\) in \(D\); \(4,-6\) in \(D\); no possibility in the all-even
third case; \(2,-8\) in \(D\); \(-2,8\) in \(D\); or \(-2,8,2,12\) in
\(E\). Every exception is excluded, so added cops are safe against old
pairs.

For a new source \(x_a\) targeting \(P_a\), its old point danger is
\[\{P_a,P_{a+5}\}\cup\{P_{a-d}:d\in D\},\] and its old line danger is
\(\{L_{a-4},L_{a-7}\}\). The latter is not adjacent. A difference-three
point pair would require \(d=-2\) or \(-8\), excluded from \(D\). For
recipient \(P_{a+5}\), the point danger is
\[\{P_a,P_{a+5}\}\cup\{P_{a+5-d}:d\in D\};\] the same test requires the
excluded offsets \(2,8\). Its line danger again has difference three. A
distinct added row could be contained only if \(b=a\) or \(b=a\pm5\)
modulo \(k\): one parity occurs just once in the point danger. The first
is the source and the others are impossible for \(a,b\in\{0,1,2\}\) and
\(k\ge22\).

For a new-to-new pair, the old danger is the source's difference-five
point pair, containing no old row and no distinct added row. As before,
new danger vertices belong to no compulsory row. This proves (R).

The eight exclusions in \(D\) and seven in \(E\) are distinct because
\(k\ge22\), so \(|D|=q-9\) and \(|E|=q-8\). There are \(2n\) compulsory
missing arcs and \(k(2q-17)+s(s-1)\) optional ones, totaling
\(H_4(n)=4q^2-26q+s^2+s\). Part 4 follows by selecting any desired
number of optional arcs.

\paragraph{Comparison and proof
boundary}\label{comparison-and-proof-boundary}

The improvement retains the old small-order domains. For \(n\ge30\),
\[H_5(n)-3n=5q(2q-5)+(s-1)^2-1\ge14.\] Also \[\begin{aligned}
13H_5(n)-n(n+11)
&=30q^2-45q-20qs+12s^2+2s\\
&\ge\frac{260(q-1)^2-261}{12}>0.
\end{aligned}\] Since \(\lceil(n-1)/13\rceil\le(n+11)/13\), this
strictly exceeds the cap \(n\max\{3,\lceil(n-1)/13\rceil\}\) of Theorem
48. For \(H_4\), the difference from \(3n\) is at least \(65\) and
\[13H_4(n)-n(n+11)=36q^2-382q-8qs+12s^2+2s.\] At \(q=11\) its values are
\(154,80,30,4\), and each increases with \(q\). Thus this family also
strictly exceeds the old cap at every \(n\ge44\).

The twenty residue classes modulo \(20\) show \(H_4(n)>H_5(n)\) for
every \(n\ge52\). Indeed it holds at \(52,\ldots,71\), and if
\(n=20t+u\), the difference increases on replacing \(n\) by \(n+20\) by
\[120t-80+40\bigl(\lfloor u/4\rfloor-\lfloor u/10\rfloor\bigr)>0.\] For
the old orders \(n=19,\ldots,29\), that cap is exactly \(3n\), already
covered by part 2. Thus the present envelopes subsume its entire
numerical statement, preserving the time, algorithm, and
arbitrary-subset conclusions at every old order. No earlier domain is
lost.

The two cap formulas have ordinary proofs above. The independent review
\texttt{research/\allowbreak{}dense-\allowbreak{}phase2-\allowbreak{}uniform-\allowbreak{}independent-\allowbreak{}audit.\allowbreak{}md} checks their
modular exceptions and counts. The primary receipts
\texttt{research/\allowbreak{}dense-\allowbreak{}phase2-\allowbreak{}uniform-\allowbreak{}envelope.\allowbreak{}json} and
\texttt{research/\allowbreak{}dense-\allowbreak{}phase2-\allowbreak{}quarter-\allowbreak{}envelope.\allowbreak{}json} additionally check
bounded full-envelope games; those samples are not used to infer the
formulas.

\subsubsection{The first forced loss at every order from
16}\label{the-first-forced-loss-at-every-order-from-16}

Use an even core of order \(2k\ge12\), or its odd modification (O) of
order \(2k+1\ge13\). Add four vertices \(c,x,y,z\) and missing arcs
\[c\to x,\quad x\to y,z,\quad
y\to a,b,\quad z\to d,e. \tag{A}
\] In the even case choose \((a,b,d,e)=(P_0,P_3,P_1,P_4)\). In the odd
case choose \((a,b,d,e)=(P_1,P_4,P_2,P_5)\). The four attachments are
distinct; each branch pair has index difference three, hence neither one
nor two modulo \(k\ge6\).

For an old missing pair, old pursuers retain their old safe neighbors,
since the danger set only gains \(y\) or \(z\). A pursuer at \(c\) uses
\(x\). A pursuer at \(x\) can use \(y\) or \(z\): no old recipient
receives both branches. The pursuers at \(y,z\) have point pairs of
difference three, which cannot fit into the old point danger pairs of
difference one or two. In the odd core there is also the point danger
set \(\{P_0,P_3\}\) for a pair ending at \(w\); neither chosen branch is
that pair. The other new old-core pairs have at most one point in their
danger set.

For a new pair with source \(y\) or \(z\) and recipient an attachment
\(P_t\), ordinary point pursuers have adjacent line neighbors while the
old incoming line pair of \(P_t\) has difference two. Ordinary line
pursuers have point-neighbor difference two while the branch pair has
difference three. In the odd core, the modified point pursuers use
\(w\), and the pursuer at \(w\) has line neighbors of difference three,
which cannot both lie in the incoming pair of difference two. The
opposite branch is disjoint, \(x\) can use that opposite branch, and
\(c\) can use \(x\). Thus these four pairs win. The pair \(c\to x\) also
wins: no missing neighbor enters \(c\), and every other pursuer has two
missing neighbors, at least one outside \(\{c,x\}\).

Exactly the pairs \((x,y),(x,z)\) fail. A pursuer at \(c\) has closed
move set \(V\setminus\{x\}\) in the delivery graph and blocks every
nonterminal departure from \(x\). All other initial pursuers have safe
neighbors for those pairs. The construction has \(2n-1\) missing arcs,
so \(2n-3\) pairs count. This proves the first-step assertion for all
\(n\ge16\).

\subsubsection{Functional roots and the rest of the
staircase}\label{functional-roots-and-the-rest-of-the-staircase}

Use the spare functional vertex when the order is odd. On \(r\ge2\)
vertices put \(j=\lfloor(r-2)/2\rfloor\), retain the missing chains
\(s_i\to t_i\to a\) and the cycle \(a\leftrightarrow b\), and, for odd
\(r\), put the remaining vertex \(z\to a\). The set of vertices with no
missing incoming arc is
\[R=\{s_1,\ldots,s_j\}\quad\hbox{together with $z$ when $r$ is odd}.\]
Consequently \(|R|=\lfloor(r-1)/2\rfloor\).

Add \(k\) vertices, each with a distinct two-element missing-neighbor
set in \(R\). This is possible whenever
\[k\le\binom{\lfloor(r-1)/2\rfloor}{2}. \tag{F}
\] It guarantees exactly \(j+2k\) counted pairs. Against a functional
pursuer, its unique missing neighbor is outside \(R\) and outside the
added vertices, so it is safe for every added source and its two
recipients. Against another added pursuer, choose a member of its
two-element set outside the source's set. Any two distinct roots are
joined in the delivery graph. The original pairs \((s_i,t_i)\) also
survive: an added pursuer has a root neighbor different from \(s_i\),
and that root has a direct delivery arc to \(t_i\). This argument
includes the extra root \(z\). Theorem 27 gives the matching upper
bound.

Now fix \(n\ge20\) and \(1\le r\le n\), and put \(k=n-r\). The first
step is already proved. For \(r\ge2\) and \(k\ge10\), join the
functional \(r\)-vertex graph to the direct \(k\)-vertex core; Lemma 28
preserves \(j+2k\) guarantees. If \(1\le k\le9\) remains, then
\(r\ge11\), so (F) holds because
\[\binom{\lfloor(r-1)/2\rfloor}{2}\ge\binom52=10\ge k.\] Finally \(k=0\)
is the functional endpoint in Theorem 13. Since
\[j+2k=2n-r-\lceil r/2\rceil-1,\] this proves part 3. More precisely,
the same formula holds whenever \(r\ge2\) and either the core condition
or (F) applies, including smaller orders than 20; the endpoint \(r=n\)
holds at every \(n\ge2\).

All recipes consist of modular arc lists and elementary subset choices,
and therefore take polynomial time. The receipt
\texttt{research/\allowbreak{}dense-\allowbreak{}unified-\allowbreak{}core-\allowbreak{}checks.\allowbreak{}json} records bounded
complete two-move and robustness checks from
\texttt{src/\allowbreak{}dense\_\allowbreak{}unified\_\allowbreak{}core.\allowbreak{}py}. The parameterized proofs above do
not depend on those finite tests. The explicitly isolated eleven-vertex
core and small optional envelopes have the finite proofs and
certificates described above. The direct even and odd cores are also
instantiated in \texttt{DenseTwoType.\allowbreak{}lean} and
\texttt{DenseOddCore.\allowbreak{}lean}, alongside the finite core and envelopes in
\texttt{ElevenVertexDense.\allowbreak{}lean} and \texttt{SmallDenseEnvelopes.\allowbreak{}lean}.
\texttt{DenseSharpClassification.\allowbreak{}lean} formalizes all dense extrema and
explicit attainers from order eleven. Alternative attachment families
and stronger guarantees retain their separately stated proof boundaries.

\subsection{Functional interfaces in the final dense
band}\label{functional-interfaces-in-the-final-dense-band}

All graphs in this section are missing-arrow graphs; the delivery graph
is their loopless directed complement. Put \(n=r+k\) and
\(j=\lfloor(r-2)/2\rfloor\). The functional graph on \(r\ge2\) vertices
has \(a\leftrightarrow b\), chains \(s_i\to t_i\to a\) for
\(1\le i\le j\), and a spare \(z\to a\) when \(r\) is odd. Its roots
\(R\) have missing indegree zero and \(|R|=\lfloor(r-1)/2\rfloor\);
exactly its \(j\) pairs \((s_i,t_i)\) are guaranteed.

\textbf{Theorem (root-colored attachment).} The exact value
\[M(n,n(n-3)+r)=j+2k=2n-r-\lceil r/2\rceil-1\] has a direct attaining
construction in the following domains:

\begin{itemize}
\tightlist
\item
  \(k=0\); or \(r\ge5\) and \(k=1\) or even; or \(r\ge7\), with two-move
  delivery;
\item
  \(r\ge5\) and arbitrary \(k\), with three-move delivery.
\end{itemize}

\textbf{Construction and proof.} For \(k=0\) use the functional graph.
For \(k=1\) give the new vertex two distinct roots as its missing
neighbors. For \(k\ge2\), put a missing cycle on the new vertices, with
successor \(\sigma\), and add \(x\to\gamma(x)\) for a proper root
coloring \(\gamma\). Two alternating roots suffice for even \(k\); an
odd cycle uses a third root at its last vertex. Each new row has two
missing neighbors.

Fix a new missing pair, from \(x\) to either \(\gamma(x)\) or
\(\sigma(x)\). An old pursuer \(c\) has a unique missing neighbor in the
old nonroots, which is a safe intermediary to either recipient. A new
pursuer \(c\ne x\) with different root color uses \(\gamma(c)\). If the
colors agree, use \(\sigma(c)\): injectivity gives
\(\sigma(c)\ne\sigma(x)\), and proper coloring ensures its delivery
arrow to \(\gamma(x)\). It also has an arrow to \(\sigma(x)\), since a
missing arrow there would force \(\sigma(c)=x\), contradicting the equal
colors. For an old counted pair \((s_i,t_i)\), a new pursuer uses
\(\sigma(c)\); the source reaches every new vertex and each new vertex
reaches the old nonroot \(t_i\). The \(k=1\) case instead uses a root
different from the functional source. These intermediaries prove every
claimed two-move guarantee.

Only odd \(k\ge3\) with two roots remains. Use an alternating cycle on
\(k-1\) vertices and a singleton \(p\) whose missing neighbors are the
two roots. A cycle-source pair against pursuer \(p\) uses the other
root; an old counted pair uses a root different from its source. Thus
all these pairs retain their two-move guarantees. From \(p\) to a root,
an old pursuer again supplies its unique nonroot missing neighbor.
Against a cycle pursuer \(c\), move safely to \(\sigma(c)\) and either
deliver directly or invoke its already proved two-move strategy. This
takes at most three moves.

The graph has \(r+2k\) missing arrows and \(j+2k\) guarantees, matching
Theorem 27. The existing fixed-graph certificates check the same
constructions; the displayed argument proves their unbounded domains.
\(\square\)

\textbf{Corollary (the two-move band).} For every \(n\ge16\) and
\(1\le r\le n\), the staircase value has an attaining construction with
two-move delivery.

\textbf{Proof.} The first step uses the four-vertex attachment of
Theorem 58. For \(2\le r\le6\), the remaining core has \(k=n-r\ge10\)
vertices, so functional-block addition to the direct core applies. For
\(r\ge7\) use the preceding colored attachment. The case \(k=0\) is
functional. \(\square\)

\subsubsection{Exact addition across missing-arrow
blocks}\label{exact-addition-across-missing-arrow-blocks}

\textbf{Lemma (arbitrary-duration missing-block addition).} Let \(F\) be
a disjoint union of loopless directed graphs \(F_1,\ldots,F_b\), every
vertex of which has positive outgoing degree. Let \(G_i\) be the
loopless complement of \(F_i\) on its own vertex set, and let \(G\) be
the loopless complement of \(F\). Then \[I(G)=\sum_{i=1}^b I(G_i).\]
More precisely, a missing pair inside block \(i\) is guaranteed in \(G\)
if and only if it is guaranteed in \(G_i\). If its original guarantee
uses at most \(D\ge1\) messenger moves, at most \(\max\{2,D+1\}\) moves
suffice in \(G\).

\textbf{Proof.} Every arrow between distinct blocks is present in \(G\),
so all missing pairs lie within one block. Fix such a pair \((s,t)\) in
block \(i\).

Suppose first it is guaranteed in \(G_i\). If the pursuer starts outside
block \(i\), choose one of its missing neighbors \(w\) in its own block.
The source has a delivery arrow to \(w\), the pursuer cannot reach \(w\)
in its reply, and \(w\) has a delivery arrow to \(t\). Two moves
suffice. Otherwise follow the original within-block strategy while the
pursuer remains in block \(i\). If the pursuer first leaves that block
after a nonterminal messenger move, apply the same two-move escape
argument from the messenger's present vertex. A play in which the
pursuer stays inside follows the original winning strategy. If the
strategy has a bound \(D\), any departure occurs after at most \(D-1\)
messenger moves, giving at most \(D+1\) moves in total. All waits are
included in these strategies. The chosen \(w\) is not the pursuer's
current vertex, since \(F\) has no loops, so waiting cannot invalidate
its safety.

Conversely, suppose the pair is not guaranteed in \(G_i\). Choose a
permitted initial pursuer position and its avoiding strategy in that
finite reachability game. Keep the pursuer inside the block and follow
that strategy. If the messenger also stays inside, it cannot guarantee
receipt. If it first leaves, its new vertex is outside the recipient's
block, so receipt has not occurred, and the pursuer immediately captures
it using the complete cross-block arrow. Thus an outside excursion
cannot turn an internal losing pair into a winning one. This proves the
equivalence and summing over missing pairs proves the identity.
\(\square\)

The two-move version of the original block lemma retains its sharper
two-move time conclusion. The present lemma removes that hypothesis for
the count identity; it does not assert that every copied strategy
continues to use two moves.

\textbf{Proposition (a nine-vertex core).} There is an explicit graph
with \(n=9\), \(m=54\) and \(I=18\). Label vertices \(0,\ldots,8\) and
take all nonloop delivery arrows except the following missing-neighbor
rows: \[
\begin{array}{c|ccccccccc}
v&0&1&2&3&4&5&6&7&8\\\hline
N_F^+(v)&\{1,3\}&\{4,7\}&\{5,6\}&\{0,2\}&\{6,8\}&
\{2,8\}&\{0,4\}&\{1,5\}&\{3,7\}
\end{array}
\] Every missing pair is guaranteed within three messenger moves. Four
of the eighteen pairs need more than two moves under the exact two-move
criterion. The literal graph is a permitted finite exception: it is
verified by two independent actual-game solvers and a full replay of all
alternating-state rank certificates, in
\texttt{research/\allowbreak{}dense-\allowbreak{}third-\allowbreak{}core-\allowbreak{}compose-\allowbreak{}certificates.\allowbreak{}json}. The
missing-pair upper bound proves \(M(9,54)=18\).

\subsubsection{One additional missing arrow supplies every first
step}\label{one-additional-missing-arrow-supplies-every-first-step}

\textbf{Lemma (enlarged direct cores).} At every order \(N\ge9\) there
is an explicit missing graph with \(2N+1\) arrows whose complement
guarantees all of them. Three messenger moves suffice for \(N\ge10\);
the retained nine-core uses four.

\textbf{Proof.} For even \(N=2q\ge10\), use the two-type core (B) of
Theorem 58 and add \(P_0\to P_1\). Only pairs leaving \(P_0\), entering
\(P_1\), and the new pair can lose an old safe intermediary. For the new
pair the point and line danger indices are \(\{0,1\}\) and
\(\{0,-1,-3,-5\}\). For \(q\ge7\), the latter contains only the source's
adjacent pair \(\{0,-1\}\). The affected old pairs have just one
additional blocked state: \[(P_0,L_{-1},L_{-5}),\] where the entries are
source, recipient and pursuer. Move safely to \(P_{-1}\), whose missing
pair to \(L_{-1}\) still has its universal two-move guarantee. This
proves the assertion for \(q\ge7\).

For odd \(N=2q+1\ge13\), use (O) and add the same \(P_0\to P_1\). For
\(q\ge8\) the only additional blocked states and safe first moves are
\[(P_0,L_{-1},L_{-5})\longmapsto P_{-1},\qquad
(P_0,w,L_{-3})\longmapsto P_3.\] Both resulting pairs retain their
universal two-move guarantees. Indeed the new pair has line danger
indices \(\{-1,-3,-5\}\), point indices \(\{0,1\}\), and vertex \(w\).
No ordinary cop row fits; the modified row \(\{w,L_2\}\) fits only when
\(q=7\). The old pair to \(w\) has point danger \(\{0,1,3\}\) and the
old pair to \(L_{-1}\) has point danger \(\{-1,0,1\}\), giving exactly
the two states displayed. Pairs entering \(P_1\) introduce no others for
\(q\ge6\).

The small moduli have the following complete exceptional-state lists;
each safe first move leads to a universally two-move pair:

{\def\LTcaptype{none} 
\Needspace{8\baselineskip}
\begin{longtable}[]{@{}
  >{\raggedright\arraybackslash}p{(\linewidth - 4\tabcolsep) * \real{0.3333}}
  >{\raggedright\arraybackslash}p{(\linewidth - 4\tabcolsep) * \real{0.3333}}
  >{\raggedright\arraybackslash}p{(\linewidth - 4\tabcolsep) * \real{0.3333}}@{}}
\toprule\noalign{}
\begin{minipage}[b]{\linewidth}\raggedright
Core
\end{minipage} & \begin{minipage}[b]{\linewidth}\raggedright
Source, recipient, pursuer
\end{minipage} & \begin{minipage}[b]{\linewidth}\raggedright
Safe first move
\end{minipage} \\
\midrule\noalign{}
\endhead
\bottomrule\noalign{}
\endlastfoot
even \(q=5\) & \((P_0,L_4,L_0)\); \((L_2,P_1,L_4)\) & \(P_4\);
\(P_0\) \\
even \(q=6\) & \((P_0,P_1,P_1)\); \((P_0,L_5,L_1)\) & \(L_1\);
\(P_5\) \\
odd \(q=6\) & \((P_0,L_5,L_1)\); \((P_0,w,L_3)\) & \(P_5\); \(P_3\) \\
odd \(q=7\) & \((P_0,P_1,P_3)\); \((P_0,L_6,L_2)\); \((P_0,w,L_4)\) &
\(L_2\); \(P_6\); \(P_3\) \\
\end{longtable}
}

For the eleven-core listed in Theorem 58, add missing arrow \(7\to0\).
The five affected danger sets, for pairs
\((7,0),(7,3),(7,4),(4,0),(9,0)\), are respectively
\[\{0,3,4,7,9\},\ \{0,3,4,6,7\},\ \{0,3,4,7\},\
\{0,1,4,7,9\},\ \{0,4,5,7,9\}.\] Each contains only its source's missing
row. Thus every missing pair still wins in two moves. Finally, the
nine-core above with \(2\to4\) added is the retained finite all-delivery
witness with 19 missing arrows in
\texttt{research/\allowbreak{}small-\allowbreak{}dense-\allowbreak{}anchor-\allowbreak{}certificates.\allowbreak{}json}.

The modular danger-set proof establishes the unbounded families.
\texttt{src/\allowbreak{}dense\_\allowbreak{}functional\_\allowbreak{}transfer\_\allowbreak{}check.\allowbreak{}py} independently checks
their exceptional states and safe moves, the five literal eleven-core
sets, and both exact game algorithms on the fixed small graphs. The
receipt \texttt{research/\allowbreak{}dense-\allowbreak{}functional-\allowbreak{}transfer-\allowbreak{}audit.\allowbreak{}json} records
the checks; no new Lean instantiation is asserted. \(\square\)

\textbf{Corollary (the first dense step).} For every \(n\ge11\),
\[M(n,n(n-3)+1)=2n-3.\]

\textbf{Proof.} Adjoin a disjoint missing two-cycle to the enlarged core
of order \(n-2\). Missing-block addition preserves all \(2n-3\) core
guarantees; there are \(2n-1\) missing arrows, meeting Theorem 27. Four
moves suffice for \(n\ge12\) by the lemma; the retained full certificate
gives the same bound at \(n=11\). The earlier two-move construction
remains available for \(n\ge16\). The retained finite attachments at
\(n=14,15\) additionally retain their three-move bounds; their full
ranks are in
\texttt{research/\allowbreak{}dense-\allowbreak{}third-\allowbreak{}core-\allowbreak{}compose-\allowbreak{}certificates.\allowbreak{}json}.
\(\square\)

\subsubsection{Equality in the second dense
step}\label{equality-in-the-second-dense-step}

\textbf{Proposition (the second step reduces exactly to a smaller dense
core).} Let \(n\ge5\) and put \(k=n-2\). Then \[M(n,n(n-3)+2)=2n-4
\quad\Longleftrightarrow\quad M(k,k(k-3))=2k.\] In particular,
\[M(n,n(n-3)+2)\le
\max\{M(k,k(k-3)),\,2n-5\}.\] An attaining graph at the left-hand
endpoint must have a missing graph consisting of a two-cycle disjoint
from a \(k\)-vertex missing core whose outdegree is exactly two
everywhere and whose complement guarantees all \(2k\) missing pairs. No
missing-indegree-two condition or two-move assumption is imposed on this
core.

\textbf{Proof.} Suppose the first maximum is attained at \(I=2n-4\).
There are \(h=2n-2\) missing arrows, so the failed-pair count is
\(\Delta=2\). The proof of Theorem 27 applies because \(I>0\). In its
notation, every missing outdegree is positive, the number \(\ell\) of
vertices with missing outdegree one satisfies \(\ell\ge2\), and
\(\Delta\ge\lceil\ell/2\rceil+1\). Therefore \(\ell=2\). The total
degree sum then forces all other missing outdegrees to be two.

Write \(S\) for the two degree-one vertices and \(f\) for their unique
missing-neighbor map. If \(b=|f(S)\setminus S|\ge1\), the stronger
inequality in that proof gives
\[\Delta\ge\left\lceil\frac{\ell+3b}{2}\right\rceil\ge3,\] which is
impossible. Thus \(f(S)\subseteq S\), and looplessness forces the
two-cycle on \(S\). Its two missing pairs already fail, as both vertices
lie in \(f(S)\).

Every missing pair with source outside \(S\) must consequently be
guaranteed. If such a source \(v\) has a missing arrow to \(a\in S\),
let \(b\) be the other cycle vertex. The pursuer at \(b\) has closed
delivery neighborhood \(V\setminus\{a\}\). The messenger at \(v\) cannot
move to \(a\), so every legal move, including waiting, is intercepted
unless it is immediate receipt; a missing recipient cannot be reached
immediately. Hence the entire missing source row of \(v\) fails, a
contradiction. All missing arrows from outside \(S\) therefore stay
outside \(S\).

The missing graph is the disjoint union of the two-cycle and a
\(k\)-vertex outdegree-two graph. Exact missing-block addition gives
\(I=I(G_k)\), because the complement of the two-cycle has no indirect
guarantees. This yields the right-hand equality.

Conversely, a graph attaining \(M(k,k(k-3))=2k\) guarantees every
missing pair. Its missing graph has positive outgoing degree at every
vertex; otherwise the corresponding universal pursuer neighborhood would
force \(I=0\). Adjoin a disjoint missing two-cycle. Exact addition
preserves \(2k=2n-4\) guarantees, matching the staircase bound. The
maximum inequality follows because any larger integer than \(2n-5\) must
equal the staircase upper bound \(2n-4\), and the equivalence then
forces the smaller-core maximum to have that value. \(\square\)

The already completed finite joint tables give \(M(5,10)=3\) and
\(M(6,18)=6\). Hence two of the gaps at that stage had stricter upper
bounds than the eventual staircase: \[M(7,30)\le9,\qquad M(8,42)\le11.\]
These immediate bounds inherit the existing exhaustive finite table
proofs. The complete small-core checks below strengthen them and also
exclude second-step equality at orders nine and ten. Their completeness,
rather than a failed bounded construction search, supplies the
exclusion.

\subsubsection{A safely reached root and the staircase from order
eleven}\label{a-safely-reached-root-and-the-staircase-from-order-eleven}

A rooted missing core consists of vertices \(K\) and one additional root
\(\rho\), with missing rows
\(F(v)\subseteq(K\setminus\{v\})\cup\{\rho\}\) for \(v\in K\). In its
restricted delivery game, surviving pursuer positions are in \(K\).
Receipt at the actual recipient \(\rho\) remains immediate. For a
recipient in \(K\), a messenger move to \(\rho\) is useful only when
\(\rho\in F(c)\) for the current pursuer \(c\): otherwise the pursuer
can capture there. After a safe arrival the root delivers directly to
the recipient on the next move. Restricting surviving pursuer positions
does not remove this capture at an unsafe root arrival.

\textbf{Lemma (designation and exact functional replacement).}
Designating any vertex of an all-delivery missing core as \(\rho\) and
forgetting its own row yields an all-delivery rooted interface. More
generally, let every core row contain a missing neighbor in \(K\).
Attach any loopless functional missing graph \(H\) at \(\rho\), with no
missing arrow in \(H\) entering \(\rho\), and no other cross missing
arrows. Then \[I(G)=I_{\rm root}(K,\rho)+I(\overline H),\] where the
rooted count includes missing pairs whose recipient is \(\rho\).

\textbf{Proof.} Designation follows the old strategy until receipt or
safe arrival at \(\rho\); in the latter case the strengthened root row
delivers next. Receipt when \(\rho\) is the actual recipient is
immediate, even if the pursuer is there.

For a core pair, an \(H\) pursuer's unique missing neighbor lies in
\(H\setminus\{\rho\}\). It is safely reachable from any core position
and directly delivers to every core recipient and to \(\rho\). Thus a
rooted winning strategy extends whenever the pursuer leaves \(K\). A
losing pursuer can remain in \(K\), capturing every messenger excursion
into \(H\setminus\{\rho\}\); visits to \(\rho\) obey precisely the
rooted safe-arrival rule. Hence the core guarantees are exactly
preserved.

For an \(H\) pair the missing recipient differs from \(\rho\). Against a
core pursuer choose its missing neighbor in \(K\), which is safely
reachable from \(H\) and directly delivers to that recipient. This also
handles a pursuer leaving \(H\) during an internal winning strategy.
Conversely, a losing pursuer can remain in \(H\) and capture every
excursion into \(K\). Thus the \(H\) guarantees are also exactly
preserved. All waits and terminal receipt are included. \(\square\)

For a functional missing map \(f\), its pair \((s,f(s))\) is guaranteed
exactly when \(s\) has missing indegree zero and \(f(s)\) has missing
indegree one. Otherwise a predecessor of \(s\), or a second predecessor
of \(f(s)\), camps and blocks every nonterminal move. If both conditions
hold, \(f(c)\) is a safe intermediary against every \(c\ne s\). This
proves the functional count used by the attachment lemma directly.

\textbf{Proposition (rooted cores).} An all-delivery rooted interface
with two missing neighbors per core vertex exists at every core order
\(k\ge7\). At \(k=7\), with root 7, its rows are
\[\{2,7\},\ \{4,5\},\ \{0,1\},\ \{4,6\},\
\{0,3\},\ \{6,7\},\ \{2,5\}.\] Its rooted guarantee has the retained
finite actual-game proof in
\texttt{research/\allowbreak{}small-\allowbreak{}dense-\allowbreak{}rooted-\allowbreak{}certificates.\allowbreak{}json} and
\texttt{research/\allowbreak{}small-\allowbreak{}dense-\allowbreak{}rooted-\allowbreak{}deletion-\allowbreak{}certificates.\allowbreak{}json}. At
\(k=8\), designate vertex 8 of the nine-core as root; at \(k\ge9\)
designate a vertex of the direct core on \(k+1\) vertices. Every row has
an internal missing neighbor, so exact functional replacement applies.
In particular, the seven-core with \(\rho\to a\leftrightarrow b\) gives
the retained small value \(M(10,73)=14\). The audit checks both positive
and negative transfer masks on fixed interfaces; the lemma itself is an
ordinary proof.

\textbf{Corollary (the complete staircase from order eleven).} For every
\(n\ge11\) and \(1\le r\le n\),
\[M(n,n(n-3)+r)=2n-r-\lceil r/2\rceil-1.\] Four messenger moves suffice;
at every \(n\ge16\) the earlier uniform two-move construction remains
available.

\textbf{Proof.} The first step was proved above. For \(r=2\) adjoin a
disjoint missing two-cycle to the ordinary core on \(n-2\ge9\) vertices.
For \(r=3,4\), use the rooted core on \(k=n-r\ge7\) vertices and attach
the functional \(r\)-vertex graph, choosing one of its roots as
\(\rho\). The exact replacement identity gives
\(2k+\lfloor(r-2)/2\rfloor\). For \(r\ge5\) use the root-colored cycle
and singleton construction, including \(k=0,1\). Every count meets
Theorem 27. The fixed rooted seven-core and nine-core certificates
retain their four-move attached bounds; larger designated cores have
two-move strategies before replacement and need at most three
afterwards. \(\square\)

\subsubsection{Strict small values in the second through fourth
steps}\label{strict-small-values-in-the-second-through-fourth-steps}

\textbf{Theorem.} The following are exact values, with literal attaining
graphs: \[
\begin{array}{c|rrrr}
 &n=7&n=8&n=9&n=10\\\hline
r=2&8&10&12&15\\
r=3&7&9&11&14\\
r=4&6&8&10&12
\end{array}
\qquad(m=n(n-3)+r).
\] In particular the uniform staircase threshold eleven is sharp:
\(M(10,74)=12<13\), the staircase value at that cell.

The ordinary seven- and eight-vertex endpoint cores needed below admit
complete small normal forms. In an all-delivery two-out missing core,
the two members of any row have disjoint missing rows and no missing
arrow between them. Otherwise source or predecessor camping loses a
missing pair. Normalize \(F(0)=\{1,2\}\). If either neighbor points back
to 0, the other rows can be labeled \(F(1)=\{0,3\}\) and
\(F(2)=\{4,5\}\); otherwise they can be labeled \(F(1)=\{3,4\}\) and
\(F(2)=\{5,6\}\). These exhaust the possibilities, without imposing any
missing-indegree condition. Complete enumeration of all remaining rows
has no survivors at either order. The independent chronological replay
uses 227 and 4,307 prefixes respectively, recorded in
\texttt{research/\allowbreak{}small-\allowbreak{}dense-\allowbreak{}endpoint-\allowbreak{}independent-\allowbreak{}checks.\allowbreak{}json}. This is
the finite endpoint exclusion used here, not an assumption about
two-move delivery strategies.

Here are the finite and ordinary components of the proof. Write
\(q(v)=d_F^+(v)\), \(\Delta=2n-r-I\), \(S=\{v:q(v)=1\}\), \(\ell=|S|\),
and \(e=\sum_v\max\{q(v)-2,0\}\). Positive \(I\) implies \(q(v)>0\)
everywhere; hence \(\ell=r+e\). From the degree-one proof,
\(\Delta\ge\lceil\ell/2\rceil+1\). Put \(U=f(S)\), \(a=|U\cap S|\),
\(b=|U\setminus S|\), and let \(\rho_0\) count vertices in
\(S\setminus U\) whose image has exactly one preimage in \(S\). Its
stronger bounds are \[
\Delta\ge\ell-\rho_0+2b,\qquad
\Delta\ge
\begin{cases}
\max\{a+2b,\ell-a+b\},&b>0,\\
\max\{a,\ell-a+2\},&b=0.
\end{cases}
\] Every row with source in \(U\) is killed. An outside row pointing
into \(U\) is killed as well by the corresponding degree-one pursuer.

At \(r=4\), equality in the staircase would require \(\Delta=3\),
forcing \(\ell=4,e=0,b=0,a=3,\rho_0=1\). The functional graph on \(S\)
is therefore \(\rho\to t_0\to a\leftrightarrow b\). Its three killed
rows exhaust the losses, and all other rows avoid \(U\). Equality is
equivalent to an all-delivery safe-root two-out core on \(n-4\)
vertices.

At \(r=3\), equality again means \(\Delta=3\). For \(\ell=4,e=1\) the
same functional path leaves a rooted core on \(n-4\) vertices with one
degree-three row and all other degrees two. For \(\ell=3,e=0\), the case
\(b\ge2\) is excluded. The case \(b=1\) would force \(a=1,\rho_0\ge2\):
the two vertices outside \(U\) would need distinct unique images, but
the vertex of \(S\cap U\) must also map to the outside image, a
contradiction. Thus \(b=0\), and \(S\) is a three-cycle or a two-cycle
with one incoming tail. All three pairs fail, leaving a rooted two-out
core on \(n-3\) vertices (with unused root in the three-cycle case).

Complete safe-root enumeration excludes two-out cores of orders 3
through 6 and one-degree-three cores of orders 3 through 6. The source
test is \(F(c)\subseteq F(s)\cup\{s\}\) for distinct core vertices: it
kills row \(s\). The predecessor test requires a cop row contained in
the core whose neighbors' completed rows have a common missing
recipient. A cop row using the root is not pruned by this latter test.
There is no missing-indegree restriction. Two-out rows using the root
have normalized initial rows \(F(0)=\{1,\rho\}\) and either
\(F(1)=\{0,2\}\) or \(F(1)=\{2,3\}\). Rows without root use the two
ordinary depth-two forms. A unique degree-three row is normalized to
\(\{1,2,\rho\}\) or \(\{1,2,3\}\). Every other loopless row is
enumerated.

The two-out families have 4, 14, 72, 1,168 prefix nodes at orders 3, 4,
5, 6; only the last has static leaves, six in total. The
one-degree-three families have 4, 74, 1,442, 29,852 nodes, with twelve
static leaves at order 6 only. All eighteen leaves fail in the full
actual game after functional attachment; both winning and losing ranks
are checked. These complete finite exclusions prove the strict third-
and fourth-step upper bounds in the table. The seven-core attachment
supplies its exceptional \(r=3\) entry at order 10.

For \(r=2\), the earlier equivalence and complete ordinary
seven/eight-core exclusions first rule out \(\Delta=2\) at orders 7
through 10. At order 10 this gives the claimed upper bound 15. At orders
7, 8, 9 suppose \(\Delta=3\). Then \(\ell=2,3,4\). The \(\ell=4\) case
leaves the four-vertex functional path and a rooted core of order
\(n-4\) with excess 2: one degree-four row or two degree-three rows. All
six profiles at orders 3, 4, 5 have zero static leaves in complete
all-row enumeration. Their prefix counts are 1, 3, 8, 69, 396, 4,211.
The \(\ell=3\) case leaves a rooted core of order \(n-3\) with one
degree-three row, already excluded.

When \(\ell=2\), all other degrees are two. A closed two-cycle forces
every outside row to avoid it, since an additional killed row would give
four losses. Exact missing-block addition then bounds \(I\) by the
ordinary core values at orders 5, 6, 7, all below \(2n-5\). Otherwise
the only shape with at most three losses is \(c\to d\to x\), with
degrees 1, 1, 2. The rows at \(d,x\) exhaust the losses. Put \(Y=F(x)\).
The degree-one pursuer at \(d\) kills every missing column into \(Y\);
\(d\notin Y\), or the additional pair \((c,d)\) would fail. Every
remaining good row avoids \(\{d,x\}\cup Y\), except \((c,d)\). If
\(c\in Y\), removing these four special vertices leaves an ordinary
universal two-out core of order \(n-4\le5\), impossible by the small
tables. Otherwise it leaves a rooted two-out core of order
\(n-5=2,3,4\), root \(c\). The order-two case is immediately
source-camped; the other two are among the complete exclusions above.
This proves \(\Delta\ge4\) and the strict upper \(2n-6\) at orders 7, 8,
9.

For the lower bounds at \(r=4\), take a missing two-cycle
\(0\leftrightarrow1\) and two free roots \(2\to0\), \(3\to0\). Attach
the alternating colored cycle and, when necessary, singleton on vertices
4 onward, using roots 2, 3. All four functional pairs fail, while the
previous routing proof guarantees all \(2(n-4)\) new pairs in at most
three moves. This attains \(2n-8\). At \(r=3\), orders 7, 8, 9, add
missing arrow \(4\to6\) to those respective graphs. At \(r=2\), add also
\((5,6)\), \((5,7)\), or \((4,8)\) respectively. For order 10, \(r=2\),
add missing arrow \(1\to3\) to the seven-core three-vertex attachment.
These finite modifications have full game certificates, without any
assertion that arbitrary additions preserve their guarantees.

The independent finite review covers every normalized case and degree
profile above, plus all eighteen negative attached graphs, in
\texttt{research/\allowbreak{}small-\allowbreak{}case-\allowbreak{}rooted-\allowbreak{}independent-\allowbreak{}review.\allowbreak{}json}. The literal
positive graphs and ranks are in
\texttt{research/\allowbreak{}small-\allowbreak{}dense-\allowbreak{}two-\allowbreak{}root-\allowbreak{}certificates.\allowbreak{}json} and
\texttt{research/\allowbreak{}small-\allowbreak{}dense-\allowbreak{}second-\allowbreak{}step-\allowbreak{}certificates.\allowbreak{}json}. The
ordinary reduction and finite game checks are distinct from the generic
Lean results retained elsewhere in the paper.

\subsubsection{The remaining first
steps}\label{the-remaining-first-steps}

\textbf{Theorem.} The first dense-step values at orders seven through
ten are
\[M(7,29)=8,\quad M(8,41)=11,\quad M(9,55)=14,\quad M(10,71)=16.\]
Together with the preceding results, this determines every dense cell in
the original unresolved list. Together with the existing complete tables
through order 6 and the constructions for \(r\ge5\), every admissible
cell of the final dense band is exact and has an explicit attaining
graph. The eventual staircase is exact throughout the final band at
every order at least eleven.

\textbf{Proof of the upper bounds.} The seven-vertex value comes from
the independently replayed complete finite table. At the other orders
put \(h=2n-1\) and \(\Delta=h-I\). The staircase bound gives
\(\Delta\ge2\). As above, \(\ell=1+e\) and
\(\Delta\ge\lceil\ell/2\rceil+1\).

Suppose first that \(\Delta=2\). If \(\ell=2\), the degree-one image
count forces a disjoint missing two-cycle and an ordinary universal core
on \(k=n-2\) vertices with \(2k+1\) missing arrows. At \(n=8\), the
six-vertex maximum excludes this; at \(n=9\), the exact value
\(M(7,27)=10\) excludes it. At \(n=10\), it would be a universal
\((8,39)\) core. Such a core has no missing degree zero or one, so one
degree is three and the other seven are two. Normalize the unique
degree-three row to \(F(0)=\{1,2,3\}\). Complete enumeration of every
other loopless degree-two row, using the two camping predicates, has no
survivors. The primary search visits 3,322 prefixes; the independent
chronological replay visits 189,190 prefixes.

If \(\ell=1\), all other degrees are two. Write the unique degree-one
row as \(c\to x\) and put \(Y=F(x)\). The two pairs from \(x\) exhaust
the losses. Every missing column into \(Y\) fails under the pursuer at
\(c\), and every other nonfunctional row pointing to \(x\) would be
killed. Thus all other good rows avoid \(Y\cup\{x\}\). If \(c\in Y\),
removing \(c,x\) and the other member of \(Y\) leaves an ordinary
universal two-out core of order \(n-3=5,6,7\), excluded by the
small-order results. Otherwise removing \(\{c,x\}\cup Y\) leaves a
safely rooted two-out core of order \(n-4=4,5,6\), root \(c\), among the
complete exclusions above. Therefore \(\Delta\ge3\) at orders 8, 9, 10,
giving respective upper bounds 12, 14, 16.

At order 8 it remains to exclude \(\Delta=3\). Here \(\ell=1,2,3,4\).
For \(\ell=4\), the four-vertex functional path leaves a rooted
four-core with total excess 3. Its degree profiles are \((5,2,2,2)\),
\((4,3,2,2)\), and \((3,3,3,2)\). Complete all-row enumeration gives no
static leaves; degree 5 is an empty row domain. For \(\ell=3\), the
closed functional triple leaves a rooted five-core with excess 2, namely
the one-degree-four or two-degree-three profiles already excluded above.

For \(\ell=2\), the other degrees are two except for one degree-three
row. A closed two-cycle would reduce the score to a six-vertex core and
is excluded by \(M_V(6)=7\). Otherwise only \(c\to d\to x\) with
\(q(c)=q(d)=1,q(x)=2\) is compatible with three losses; a higher degree
at \(x\), two distinct outside images, or a common outside image gives
at least four. The earlier reduction with \(Y=F(x)\) leaves either an
ordinary four-core with all degrees at least two, impossible by its
known maximum, or a rooted three-core with all degrees two and at most
one degree-three row. Both rooted profiles were excluded.

For \(\ell=1\), normalize \(F(0)=\{1\}\). Up to relabeling,
\(F(1)=\{0,2\}\) or \(F(1)=\{2,3\}\), and all other rows have degree
two. Complete enumeration of all such rows has no completion with at
most three distinct forced losses. The two searches use different row
orders: the primary combined count is 578 prefixes; the independent
counts are 12,160 and 9,934. To avoid double counting, each source
stores the set of its marked missing recipients. A source-camping
witness marks its entire row. A predecessor-camping witness at \(c\)
marks the shared missing recipient column, except for initial source
\(c\), where that cop placement is excluded. The size of the union, not
the sum of witness counts, is the pruning bound. Every pruned prefix
therefore forces at least four different lost pairs in every completion.
This completes the upper bound 11 at order 8.

\textbf{Literal lower bounds.} The following missing rows, listed in
source order, give the four attaining graphs: \[
\begin{array}{c|l}
n&F(0),F(1),\ldots,F(n-1)\\\hline
7&1;\ 0;\ 01;\ 0;\ 256;\ 346;\ 23\\
8&1;\ 0;\ 0;\ 0;\ 256;\ 367;\ 247;\ 34\\
9&123;\ 245;\ 06;\ 16;\ 05;\ 46;\ 7;\ 8;\ 7\\
10&27;\ 345;\ 014;\ 46;\ 03;\ 67;\ 25;\ 8;\ 9;\ 8
\end{array}
\] Each string denotes the set of its digit labels; all labels are
single digits in this table. The order-nine graph is the three-vertex
attachment of a rooted six-core with two degree-three rows; its root is
6. The literal full graphs have respective exact counts 8, 11, 14, 16
and delivery bounds 3, 3, 5, 4 from both game solvers and full rank
replay. These finite bounds retain the four-move staircase guarantee
from order 11 and the uniform two-move guarantee from 16.

The constructor is \texttt{src/\allowbreak{}small\_\allowbreak{}dense\_\allowbreak{}first\_\allowbreak{}step.\allowbreak{}py}, with fixed
game records in
\texttt{research/\allowbreak{}small-\allowbreak{}dense-\allowbreak{}first-\allowbreak{}step-\allowbreak{}certificates.\allowbreak{}json}. The
independent upper-bound replay is
\texttt{research/\allowbreak{}small-\allowbreak{}case-\allowbreak{}first-\allowbreak{}step-\allowbreak{}independent-\allowbreak{}review.\allowbreak{}json}; the
literal game replay is
\texttt{research/\allowbreak{}small-\allowbreak{}case-\allowbreak{}first-\allowbreak{}step-\allowbreak{}literal-\allowbreak{}review.\allowbreak{}json}. These are
ordinary reductions and finite certificate proofs, with no claim that
the new small classification has been formalized in Lean. \(\square\)

\subsection{Missing-graph criteria and the dense
endpoint}\label{missing-graph-criteria-and-the-dense-endpoint}

The following criterion applies to arbitrary missing graphs. The
numerical constructions are now supplied by the stronger direct Theorem
58.

\textbf{Lemma 18 (high-girth complement).} Let \(F\) be an orientation
of a finite simple undirected graph of girth at least six. Suppose every
vertex has outdegree at least two in \(F\). Let \(G\) contain exactly
the nonloop arcs absent from \(F\). Then every source--recipient pair in
\(G\) is guaranteed within two messenger moves.

\textbf{Proof.} Diagonal and direct pairs need no argument. For a
missing pair \((s,t)\) of \(G\), the arc \(s\to t\) belongs to \(F\).
Set \[
B=N_F^+(s)\cup N_F^-(t).
\] Both \(s\) and \(t\) belong to \(B\). In the underlying undirected
graph, the vertices of \(B\) form an induced double star: its centers
are \(s,t\), the other outneighbors of \(s\) are leaves at \(s\), and
the other predecessors of \(t\) are leaves at \(t\). The two sets of
leaves are disjoint, since an overlap would make a triangle. An
additional edge inside \(B\) would make a triangle or quadrilateral.
Thus the distance between any two vertices in this induced tree is at
most three.

Fix an initial pursuer vertex \(c\ne s\). If \(c\notin B\), two
neighbors of \(c\) in \(B\) would form an undirected cycle of length at
most five, so \(c\) has at most one outgoing \(F\)-neighbor in \(B\).
The same bound holds when \(c\in B\setminus\{s\}\): a leaf at \(s\) has
its only incident edge in the tree directed toward it; a leaf at \(t\)
has only its edge toward \(t\); and all tree edges incident with \(t\)
point into \(t\).

Since \(c\) has at least two outgoing \(F\)-neighbors, choose
\[u\in N_F^+(c)\setminus B.\] Then \(u\ne s,t\), and \(s\to u\) and
\(u\to t\) are arcs of \(G\). Also \(u\ne c\) and \(c\to u\) is absent
from \(G\). The messenger moves to \(u\), survives every pursuer
response, and delivers to \(t\) on its next turn. \(\square\)

The lemma only requires a lower bound on outgoing degrees. In
particular, edges can later be added to the underlying graph and
oriented arbitrarily, provided its girth remains at least six.

\textbf{Theorem 19 (the densest noncomplete all-delivery networks).} At
every \(n\ge10\), \[M(n,n(n-3))=2n.\] An attaining graph guarantees
every pair within two messenger moves.

\textbf{Proof.} Use the direct two-type missing graph, with its
odd-order subdivision when needed, and the explicit eleven-vertex core,
in Theorem 58. It has exactly \(2n\) missing arcs and a two-move
guarantee for every pair. The missing-pair upper bound proves equality.
\(\square\)

The density \(n(n-3)\) is maximal for a noncomplete all-delivery graph
by Theorem 14. The earlier high-girth construction at orders \(n\ge52\)
is no longer needed for this numerical theorem. Lemma 18 remains useful
as a criterion for arbitrary missing graphs; its structural hypothesis
is stronger than the condition needed by the direct construction.

\textbf{Theorem 20 (a retained dense parameterization).} Fix integers
\(n\ge52\) and \(5\le D\le n-1\), and put \[
B_D=D^4-2D^3+2D^2+1,\qquad
L_D(n)=\left\lfloor\frac{Dn-(D-4)B_D}{2}\right\rfloor.
\] Every integer \(2n\le h\le L_D(n)\) satisfies \[M(n,n(n-1)-h)=h,\]
with a deterministic two-move attaining construction. An empty interval
makes no assertion.

\textbf{Proof.} A nonempty interval forces \(n\ge B_D\ge B_5=426\). For
\(D\ge5\), the displayed polynomial gives \(B_D-1\ge13D/2\). Hence \[
L_D(n)\le\frac{Dn}{2}
\le n\left\lceil\frac{n-1}{13}\right\rceil.
\] The robust bridge of Theorem 58 and the protected envelope of Theorem
48 cover every missing budget from \(2n\) to the latter endpoint. Their
constructions have two-move delivery. \(\square\)

This parameterization is a consequence of the stronger dense envelope,
retained for references to its explicit \(D\)-dependent bound. It
introduces no additional graph search, girth-growth certificate, or
separate construction family.

\subsection{Attainment across densities and losses near
completeness}\label{attainment-across-densities-and-losses-near-completeness}

Angel Raychev proposed that, for every sufficiently large \(n\), the
missing-pair bound should be attained at \textbf{every} integer budget
\[2n\le m\le n(n-3).\] Theorem 44 proves this conjecture at every
\(n\ge115\) using the explicit constructions of Theorems 38, 52, 58, and
59. This section gives the separate exact-budget probability argument
and the dense counting and composition principles. Raychev also proposed
separating existence from explicit constructions and studying how
domination destroys missing pairs as the graph approaches completeness.
These distinctions motivate the statements below; their proofs were
developed in the September 2026 follow-up investigation.

Write \(L=n(n-1)\) and \(h=L-m\). Thus \(h\) is the number of missing
nonloop arcs, and always \(I(G)\le h\).

\subsubsection{An exact-budget existence
theorem}\label{an-exact-budget-existence-theorem}

\textbf{Theorem 25.} Let \(n\ge4\), let \(0<m<L\) be an integer, and put
\(p=m/L\). If \[
(n-3)p^2(1-p)>5\log n,
\] then some graph with exactly \(n\) vertices and \(m\) arcs guarantees
every source--recipient pair within two messenger moves. In particular,
\[M(n,m)=L-m.\]

\textbf{Proof.} Apply Theorem 54 with the routing set equal to the
entire vertex set: \(H=V\), \(b=n\), and \(a=m\). Since \(n\ge4\),
\[5\log n>4\log n+\log(8/3),\] the displayed hypothesis implies its
condition. Theorem 54 supplies an exactly \(m\)-arc graph with a safe
two-move intermediary for every legal source--recipient--pursuer triple,
proving the assertion. \(\square\)

This is an exact integer-budget existence theorem, rather than a
statement about an expected number of arcs. It also gives a finite
randomized search: sample graphs with exactly \(m\) arcs and check the
safe-intermediate condition. No practical running-time bound near the
threshold is asserted.

\textbf{Corollary 26.} There is an absolute \(N\) such that, for every
\(n\ge N\) and every integer budget \[
\left\lceil10n^{3/2}\sqrt{\log n}\right\rceil\le m\le n(n-3),
\] we have \(M(n,m)=n(n-1)-m\). Attainment can always be achieved with
delivery in at most two moves.

\textbf{Proof.} Let \(p=m/L\) and \(h=L-m\). For \(n\ge6\), we have
\(n-3\ge n/2\). If \(p\le1/2\), the lower bound on \(m\) gives
\(p\ge10\sqrt{\log n/n}\) and therefore \[
(n-3)p^2(1-p)\ge25\log n>5\log n.
\] If \(p\ge1/2\) and \(h\ge100n\log n\), then
\(1-p=h/L\ge100\log n/n\), giving \[
(n-3)p^2(1-p)\ge\frac{100}{8}\log n>5\log n.
\] Theorem 25 covers both cases.

For the remaining dense budgets use Theorem 20 with
\[D=\left\lfloor(n/2)^{1/4}\right\rfloor.\] For all sufficiently large
\(n\) its hypotheses hold, and its upper endpoint satisfies
\[L_D(n)\ge n(D/4+1)-1.\] Since \(D\) grows as a positive constant times
\(n^{1/4}\), this endpoint eventually exceeds \(100n\log n\). Theorem 20
supplies every integer missing-arc count from \(2n\) to \(L_D(n)\), and
the probabilistic range supplies the rest. The intervals overlap, so no
admissible integer budget is omitted. Both proofs give two-move
delivery. \(\square\)

The absolute threshold is not optimized or asserted to be moderate. In
particular, this corollary does not establish the full interval at
\(n=52\). The general sparse interval between \(2n\) and order
\(n^{3/2}\sqrt{\log n}\) remains outside its conclusion.

\subsubsection{A growing loss in the final dense
band}\label{a-growing-loss-in-the-final-dense-band}

Let \(F\) be the loopless directed complement of \(G\), and write
\(q(v)=d_F^+(v)\). The number \[\Delta=h-I(G)\] counts missing pairs
that fail to guarantee delivery.

\textbf{Theorem 27.} For \(n\ge2\) and every admissible integer
\(n\le h<2n\), \[
I(G)\le h-\left\lceil\frac{2n-h}{2}\right\rceil-1.
\] Equivalently, writing \(r=2n-h\), we have \[
M\bigl(n,n(n-3)+r\bigr)\le2n-r-\lceil r/2\rceil-1
\qquad(1\le r\le n).
\]

\textbf{Proof.} If \(I(G)=0\), the conclusion holds: the right-hand side
is increasing with \(h\) and at \(h=n\) equals
\(\lfloor(n-2)/2\rfloor\ge0\). Assume henceforth that \(I(G)>0\).

Every \(q(v)\) is positive. Otherwise a vertex \(c\) with \(q(c)=0\) can
intercept every nonterminal departure from any other source, while its
own source row has no missing pair. Put \[
S=\{c:q(c)=1\},\qquad \ell=|S|.
\] Since \(h\ge\ell+2(n-\ell)\), we have \(\ell\ge2n-h\ge1\). For
\(c\in S\) let \(f(c)\) be its unique missing neighbor; thus
\(f(c)\ne c\). Set
\[U=f(S),\qquad a=|U\cap S|,\qquad b=|U\setminus S|.\]

Every source \(x\in U\) loses all of its \(q(x)\) missing pairs. To see
this, choose \(c\in S\) with \(f(c)=x\). The pursuer's closed outgoing
neighborhood in \(G\) is \(V\setminus\{x\}\). It can wait at \(c\) while
the messenger waits at \(x\) and intercept every nonterminal departure.
An indirect recipient cannot be reached by a direct departure.

A second obstruction applies when two distinct vertices \(c,c'\in S\)
have the same image \(x\). The missing pair \((c,x)\) loses against a
pursuer starting at \(c'\), whose closed outgoing neighborhood contains
every possible nonterminal intermediate vertex.

Let \(R\) be the set of \(c\in S\setminus U\) for which \(f(c)\) has
exactly one preimage in \(S\), and put \(\rho=|R|\). Every missing edge
with tail in \(S\setminus R\) loses by one of the two obstructions. The
killed rows in \(U\setminus S\) have at least two missing edges each,
with tails outside \(S\). These losses are disjoint, so
\[\Delta\ge\ell-\rho+2b.\]

If \(b\ge1\), then \(\rho\le\ell-a\) because
\(R\subseteq S\setminus U\), and \(\rho\le a+b\) because its members
have distinct images in \(U\). Consequently \[
\Delta\ge\max\{a+2b,\ell-a+b\}
\ge\left\lceil\frac{\ell+3b}{2}\right\rceil
\ge\lceil\ell/2\rceil+1.
\] If \(b=0\), the function \(f\) maps \(S\) into itself without fixed
points. Its finite functional graph has a directed cycle of length at
least two. Every cycle vertex belongs to \(U\) and already has a
preimage in \(U\), so it cannot also be the image of a member of \(R\):
that member is outside \(U\) and would give a second preimage. Thus
\(\rho\le a-2\), as well as \(\rho\le\ell-a\). It follows that \[
\Delta\ge\max\{a,\ell-a+2\}\ge\lceil\ell/2\rceil+1.
\] Using \(\ell\ge2n-h\) proves the bound. \(\square\)

Thus the forced loss increases as \(2,2,3,3,4,4,\ldots\) as the arc
budget moves one, two, three, four, five, six, and further steps above
\(n(n-3)\). At the other endpoint, \(h=n\), the bound agrees with
Theorem 13.

\subsubsection{Combining two-move
constructions}\label{combining-two-move-constructions}

\textbf{Lemma 28.} Suppose a missing-edge graph \(F\) is a disjoint
union of graphs \(F_i\), each with positive outgoing degree at every
vertex. Let \(G_i\) be the complement on its own block. If an indirect
pair in \(G_i\) can be guaranteed within two messenger moves, then it
remains guaranteed in the complement \(G\) of the entire \(F\).

\textbf{Proof.} Let the source and recipient lie in block \(i\). If the
initial pursuer also lies there, the original two-move strategy chooses
an intermediate vertex outside its closed outgoing neighborhood. The
same is true in \(G\). The pursuer may respond by moving outside the
block, but this cannot prevent receipt on the messenger's next move.

If the pursuer starts at \(c\) in another block, choose any outgoing
\(F\)-neighbor \(u\) of \(c\) in that block. The arcs \(s\to u\) and
\(u\to t\) are present in \(G\) because they cross blocks, whereas
\(c\to u\) is absent. Also \(u\ne c,s,t\). This gives the required safe
two-move delivery. \(\square\)

The lemma concerns two-move strategies. No preservation theorem for
arbitrary longer strategies under this join is assumed.

\textbf{Theorem 29.} Let \(r\ge2\) and \(k=n-r\). If \(k\ge10\), then
\[M\bigl(n,n(n-3)+r\bigr)=2n-r-\lceil r/2\rceil-1.\] An explicit
polynomial construction attains this value with every counted pair
delivered within two messenger moves.

\textbf{Proof.} Use the functional missing graph on \(r\) vertices, with
\(j=\lfloor(r-2)/2\rfloor\) guaranteed pairs. On the other \(k\)
vertices use the direct missing core in Theorem 58, which has \(2k\)
guaranteed missing pairs. Both blocks have positive missing outdegree,
so Lemma 28 preserves all \(j+2k\) guarantees in the complement of their
disjoint union. There are \(r+2k=2n-r\) missing arcs, and
\[j+2k=2n-r-\lceil r/2\rceil-1.\] Theorem 27 gives optimality.
\(\square\)

This includes the former \(k\ge52\) range and the separate 26-vertex
block, without any finite girth-growth dependency. In missing-budget
notation the range is \(n+10\le h\le2n-2\).

\subsubsection{A restriction on permutation
constructions}\label{a-restriction-on-permutation-constructions}

\textbf{Proposition 30.} Suppose the two outgoing neighbors of each
vertex are \(\alpha(v)\) and \(\beta(v)\), where \(\alpha,\beta\) are
commuting permutations, both images differ from \(v\), and
\(\alpha(v)\ne\beta(v)\). Then \(I(G)\le n\). In fact all indirect
guarantees are decided by a two-move criterion.

\textbf{Proof.} Fix a recipient \(t\) and put
\(s_t=\alpha^{-1}\beta^{-1}(t)\). Commutation gives \[
N^+(s_t)=\{\beta^{-1}(t),\alpha^{-1}(t)\}=N^-(t).
\] The predecessor camping lemma excludes every indirect source except
possibly \(s_t\): a pursuer starting at \(s_t\) can intercept the
messenger when it first reaches a predecessor of \(t\). Thus each
recipient contributes at most one pair.

When \(s_t\ne t\), the pair \((s_t,t)\) is missing because neither
predecessor of a loopless recipient equals that recipient. This pair is
guaranteed exactly when no \(c\ne s_t\) has \(N^-(t)\subseteq C(c)\).
Necessity is the same camping lemma. For sufficiency, from each legal
initial pursuer position choose \(u\in N^-(t)\setminus C(c)\). The
messenger moves from \(s_t\) to this safe predecessor and then delivers.
\(\square\)

The bound is sharp: on \(\mathbb Z/7\mathbb Z\), take arcs \(x\to x+1\)
and \(x\to x+3\). For recipient zero the predecessors are \(6,4\), and
the only closed outgoing neighborhood containing both is that of \(3\).
The criterion gives \((3,0)\); translation gives exactly the seven pairs
\((s,s+4)\). Since \(n(n-3)>n\) for \(n\ge5\), commuting permutations
cannot supply a universal vertex-extremal graph in that range.
Noncommutation is necessary, not asserted sufficient.

\subsubsection{An extension that fills the small dense
offsets}\label{an-extension-that-fills-the-small-dense-offsets}

\textbf{Proposition 31.} Let \(r\ge2\), \(k\ge1\), and \(n=r+k\). If
\[k\le\binom{\lfloor(r-1)/2\rfloor}{2},\] then a deterministic
polynomial construction attains
\[M\bigl(n,n(n-3)+r\bigr)=2n-r-\lceil r/2\rceil-1,\] with a two-move
strategy for every counted pair.

\textbf{Proof.} The functional construction has
\(j=\lfloor(r-2)/2\rfloor\) guaranteed pairs. Its vertices with no
missing incoming arc are the \(j\) branch sources, together with the
spare vertex when \(r\) is odd. There are exactly
\(q=\lfloor(r-1)/2\rfloor\) such roots. Give each added vertex a
different two-element missing-neighbor set among these roots. The
functional-root argument in Theorem 58 proves all \(2k\) new missing
pairs and preserves the original \(j\) pairs in two moves. The identity
\(j+2k=2n-r-\lceil r/2\rceil-1\) and Theorem 27 prove equality.
\(\square\)

This strengthens the former bound
\(k\le\binom{\lfloor(r-2)/2\rfloor}{2}\) when \(r\) is odd and preserves
it when \(r\) is even. In particular, all the earlier \(n\ge75\),
\(1\le n-r\le51\) cases remain covered.

Theorem 58 supplies the common direct construction and complete proof.
The counting bound and binary two-move block-composition principle are
also formalized in \texttt{DenseStaircase.\allowbreak{}lean} and
\texttt{ComplementarySum.\allowbreak{}lean}. The commuting-permutation restriction
and supplementary construction claims retain their separate proof
boundaries.

\subsection{Exact final dense band}\label{exact-final-dense-band}

\textbf{Theorem 32.} For every \(n\ge16\),
\[M\bigl(n,n(n-3)+1\bigr)=2n-3.\] A deterministic attaining graph
guarantees all but two of its missing pairs within two messenger moves.

\textbf{Proof.} Use the four-vertex attachment to a direct two-type core
in Theorem 58. The graph has \(2n-1\) missing arcs. Exactly the two
missing pairs leaving the attachment's vertex \(x\) are blocked by a
pursuer at \(c\); every other missing pair has a safe two-move
intermediary. Theorem 27 gives the matching upper bound. \(\square\)

\textbf{Corollary 33 (the complete final dense band).} For every
\(n\ge20\) and every integer \(1\le r\le n\), \[
\boxed{M\bigl(n,n(n-3)+r\bigr)
=2n-r-\left\lceil\frac r2\right\rceil-1.}
\] Equivalently, for every integer \(n\le h\le2n-1\),
\[M\bigl(n,n(n-1)-h\bigr)
=h-\left\lceil\frac{2n-h}{2}\right\rceil-1.\] There is a deterministic
polynomial attaining construction at every such budget, with a two-move
strategy for every counted pair.

\textbf{Proof.} Theorem 27 gives the upper bound. The first step is
Theorem 32. For \(r\ge2\), put \(k=n-r\). If \(k\ge10\), Theorem 29
applies. Every remaining positive \(k\) is at most nine, so \(r\ge11\)
and \[\binom{\lfloor(r-1)/2\rfloor}{2}\ge\binom52=10\ge k.\] Proposition
31 applies. Finally \(k=0\) is the functional endpoint of Theorem 13.
\(\square\)

Together with Theorem 19 at \(m=n(n-3)\) and Theorem 13 above
\(m=n(n-2)\), this determines the joint maximum for every \(n\ge20\) and
\(m\ge n(n-3)\). At the boundary the value is \(2n\); the next step has
value \(2n-3\); the staircase then reaches \(\lfloor(n-2)/2\rfloor\)
before becoming zero. The proof does not assert that 20 is the smallest
possible whole-band threshold.

The core and attachment formulas and their complete ordinary proofs are
given together in Theorem 58. The eleven-vertex core is listed
explicitly and independently checked; its two-move guarantee is also
formalized in Lean. \texttt{DenseFirstStep.\allowbreak{}lean} also proves the
first-step maximum and explicit attainment from order 16, and
\texttt{DenseSharpClassification.\allowbreak{}lean} supplies the complete dense
values and explicit attainers from order eleven. Supplementary
guarantees retain their separately stated proof boundaries.

\begin{figure}
\centering
\pandocbounded{\includegraphics[keepaspectratio,alt={The complete final dense band at n=100. The solid curve is the proved exact maximum at every integer budget; the dashed line counts all missing pairs and does not account for forced failures.}]{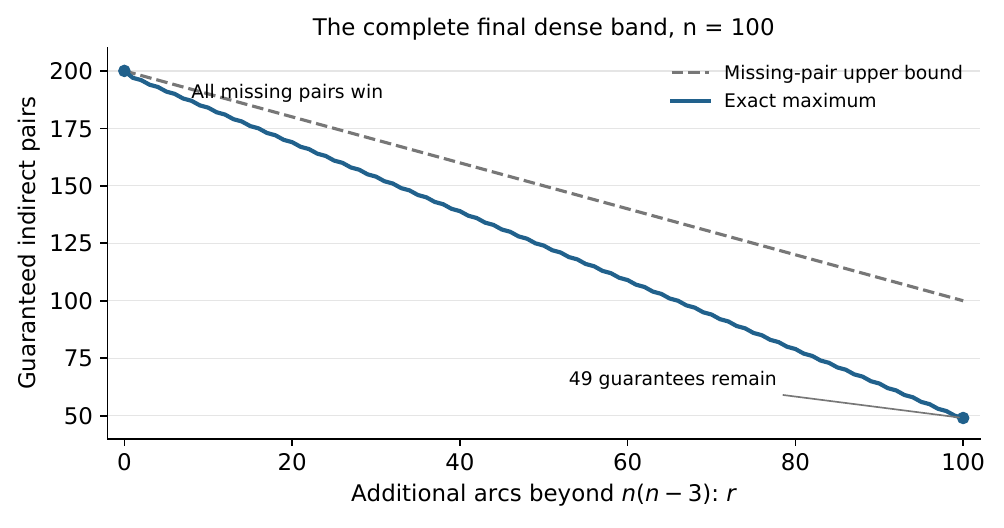}}
\caption{The complete final dense band at n=100. The solid curve is the
proved exact maximum at every integer budget; the dashed line counts all
missing pairs and does not account for forced failures.}
\end{figure}

\subsection{A quadratic dense interval from three safe
intermediaries}\label{a-quadratic-dense-interval-from-three-safe-intermediaries}

The sparse missing graph need not retain large girth while its density
increases. Three designated outgoing missing arcs per vertex suffice,
provided certain differences between their endpoints stay protected. The
resulting construction covers a quadratic interval of missing-arc counts
by a direct modular rule.

\textbf{Theorem 48.} For every integer \(n\ge19\), put
\[K(n)=\max\left\{3,\left\lceil\frac{n-1}{13}\right\rceil\right\}.\] At
every integer missing-arc count \(3n\le h\le nK(n)\), an explicit graph
attains \[M(n,n(n-1)-h)=h.\] Every source--recipient pair in the
construction is guaranteed within two messenger moves. Arbitrary
intermediate subsets of the specified optional missing arcs are
permitted.

\textbf{Proof.} Work in \(\mathbb Z/n\mathbb Z\), and put
\[S=\{1,3,7\},\qquad D=\{\pm2,\pm4,\pm6\},
\qquad E=\{\pm3,\pm5,\pm9\}.\] We shall construct a set \(A\) containing
\(S\), not containing zero, with the following properties:

\begin{enumerate}
\def\labelenumi{\arabic{enumi}.}
\tightlist
\item
  Each difference in \(D\) has exactly its original ordered
  representation between two members of \(S\), and no other ordered
  representation between members of \(A\).
\item
  No difference between two members of \(A\) belongs to \(E\).
\end{enumerate}

For \(n\ge19\), these properties hold for \(A=S\): its six ordered
nonzero differences are precisely the six distinct members of \(D\), and
none belongs to \(E\).

Whenever \(|A|=k\), adjoin any residue outside
\[\{0\}\cup A\cup(A+D)\cup(A+E).\] The exclusion of \(A+D\) preserves
all protected representations, and the exclusion of \(A+E\) preserves
the second property; both difference sets are symmetric. At most
\(13k+1\) residues are excluded. If \(k<\lceil(n-1)/13\rceil\), this
quantity is strictly less than \(n\), so the process continues. Starting
with three members, it produces \(|A|=K(n)\); choosing the smallest
permissible representative makes the rule deterministic.

Let \(F\) be any loopless digraph satisfying
\[\{(x,x+d):x\in\mathbb Z/n\mathbb Z, d\in S\}
\subseteq E(F)\subseteq
\{(x,x+d):x\in\mathbb Z/n\mathbb Z, d\in A\}.\] Let \(G\) be its
directed complement. We prove all missing pairs of \(G\) have safe
two-move intermediaries, using only the three designated missing
neighbors of each initial pursuer position.

Fix a missing pair \((s,t)\) and initial pursuer position \(c\ne s\).
Translate so \(s=0\), and write \(a=t-s\in A\) and \(z=c-s\ne0\). The
three candidates are the members of \(z+S\). A candidate is suitable if
it lies outside \[A\cup(a-A).\] Indeed, this ensures both messenger arcs
are present in \(G\), while the candidate is a missing outgoing neighbor
of \(c\) and therefore cannot be reached on the pursuer's response.

At most one candidate belongs to \(A\). Otherwise, for distinct
\(u,v\in S\), the difference of \(z+u,z+v\in A\) is the protected
difference \(u-v\). Its unique representation forces \(z+u=u\) and
\(z+v=v\), contrary to \(z\ne0\).

At most two candidates belong to \(a-A\). If the candidates with steps
\(u,v\) belong there, the two elements \(a-z-u,a-z-v\) of \(A\) have
protected difference \(v-u\). Its unique representation is \((v,u)\), so
\[a-z=u+v.\] If the third candidate, with step \(w\), also belonged to
\(a-A\), the same argument applied to \(u,w\) would give \(a-z=u+w\), a
contradiction.

If all three candidates were blocked, the remaining candidate \(z+w\)
would have to belong to \(A\), while the other two are in \(a-A\).
Setting \(b=z+w\in A\), the displayed equality gives \[b-a=w-u-v.\] As
\(w\) runs over \(1,3,7\) and \(u,v\) are the other two steps, this
difference is respectively \(-9,-5,3\). Each belongs to \(E\),
contradicting the second property of \(A\). Hence a suitable candidate
exists.

The suitable candidate differs from \(c\) because \(0\notin S\). It
differs from both endpoints: \(a\in A\) and \(0\in a-A\). Thus the
messenger first moves to it, survives every response, and then reaches
\(t\); receipt precedes any subsequent interception. This argument uses
only the envelope containments, so no translation symmetry of the chosen
intermediate \(F\) is needed.

The compulsory missing graph has \(3n\) arcs and the full envelope has
\(nK(n)\) arcs. Adding any chosen number of optional missing arcs gives
each integer \(h\) in this interval, with all \(h\) missing pairs
guaranteed. The missing-pair upper bound proves equality. \(\square\)

Theorem 58 now subsumes this entire numerical range, including the
orders 19 through 21, while preserving two-move delivery and every
optional subset. The proof above is retained as a distinct sufficient
condition on arbitrary protected difference sets. The main numerical
catalogue therefore no longer counts this as an additional dense
interval.

The executable construction and bounded checks are in
\texttt{src/\allowbreak{}dense\_\allowbreak{}protected\_\allowbreak{}envelope.\allowbreak{}py}; their exact finite scope is
recorded in \texttt{research/\allowbreak{}dense-\allowbreak{}protected-\allowbreak{}envelope-\allowbreak{}checks.\allowbreak{}json}. The
ordinary theorem is independent of those finite pursuit-game checks. No
new Lean coverage is asserted.

\section{Finite extrema}\label{finite-extrema}

\subsection{Completion of the finite sparse
budgets}\label{completion-of-the-finite-sparse-budgets}

The ordinary source/sink and positive-support lemma supplies one common
reduction. We use the exact smaller budgets in increasing order,
followed by the complete finite exclusions listed below. The game and
counted pairs are unchanged: messenger first, waiting allowed, receipt
before capture, and distinct nonadjacent source--recipient pairs.

\textbf{Theorem (finite sparse-budget values).} The budgets twelve
through twenty-one have the following global maxima and least attaining
host orders:

{\def\LTcaptype{none} 
\Needspace{6\baselineskip}
\begin{longtable}[]{@{}lrrrrrrrrrr@{}}
\toprule\noalign{}
\(m\) & 12 & 13 & 14 & 15 & 16 & 17 & 18 & 19 & 20 & 21 \\
\midrule\noalign{}
\endhead
\bottomrule\noalign{}
\endlastfoot
\(M_E(m)\) & 9 & 11 & 13 & 15 & 19 & 28 & 33 & 40 & 53 & 60 \\
Least host \(h\) & 7 & 8 & 8 & 8 & 8 & 10 & 9 & 11 & 10 & 12 \\
\end{longtable}
}

For every \(n\ge h\), \(M(n,m)=M_E(m)\). The exceptional smaller-host
values are \[M(9,17)=24,\quad M(9,19)=35,\quad M(10,19)=38,\]
\[M(10,21)=56,\qquad M(11,21)=58.\]

\subsubsection{The finite domains, in budget
order}\label{the-finite-domains-in-budget-order}

Write \(A_j=M_E(j)\), \(D_m\) for (56.3), and \(B_m\) for the
positive-support bound (56.4). The exact values through eleven arrows
are already established. At each following budget, use only the
preceding exact values.

\textbf{Budgets twelve through sixteen.} In a strict counterexample to
the displayed value, a degree-one source or sink is impossible because
\(A_{m-1}\) is smaller. A source or sink of other degree at least three
has score at most \[A_{m-3}+\lfloor(m-3)/2\rfloor=9,11,13,15,17\] at the
five respective budgets, again at or below the displayed target. Thus
every nonisolated source or sink has other degree exactly two. Every
nonisolated vertex has total degree at least two, so the support has
order at most \(m\). This directly justifies the comparison hosts 12 and
13 for their complete exclusions: a smaller support is padded with
isolates. No final least-attaining-host value is used to prove either
comparison.

For budgets fourteen through sixteen, the deleted degree-two row or
column must be live, because otherwise deletion gives the smaller value
\(A_{m-2}\). A source-to-sink arrow would contradict its required two
nonsink successors or two positive-indegree predecessors. These
observations also hold at twelve and thirteen, although their full-host
exclusions already cover the permitted graphs.

For completeness, the simultaneous live capacities used by the retained
enumerators follow from the same lemma. Let \(S\) be the total outgoing
degree of all indegree-zero vertices and \(T\) the total incoming degree
of all sinks. If \(R,C\) are the numbers of live rows and columns, then
\[R\le\lfloor(m-T)/2\rfloor,\qquad C\le\lfloor(m-S)/2\rfloor.\] With
active sets of sizes \(q,r\) and overlap \(b\), the count for \(u,v\)
actual live rows and columns gives \[
I\le\max_{\substack{0\le u\le\min(q,R)\\0\le v\le\min(r,C)}}
\{uv-(b-q-r+u+v)_+-(2u+2v-m)_+\}. \tag{F.1}
\] Routes through nonlive vertices remain available; this counts pairs,
not an induced game on live vertices.

The complete twelve- and thirteen-arrow exclusions have respectively 65
and 88 outgoing-degree profiles. The fourteen- through sixteen-arrow
packages cover all permitted positive and degree-two-zero profiles on
supports nine through \(m\), together with the reviewed eight-vertex
bounds 13,15,19. The retained seven-vertex bounds at those budgets are
10,11,12. Consequently every support is covered. The literal attainers
give equality, allowing the induction to continue.

\textbf{Budgets seventeen through twenty-one.} Once the preceding exact
values are known, every graph with an actual source or sink is bounded
as follows:

{\def\LTcaptype{none} 
\Needspace{5\baselineskip}
\begin{longtable}[]{@{}lrrrrr@{}}
\toprule\noalign{}
\(m\) & 17 & 18 & 19 & 20 & 21 \\
\midrule\noalign{}
\endhead
\bottomrule\noalign{}
\endlastfoot
\(D_m\) & 22 & 28 & 36 & 42 & 53 \\
\end{longtable}
}

All are strictly below the global target. For positive support, (56.4)
leaves the following finite domains; support cannot exceed \(m\).

{\def\LTcaptype{none} 
\Needspace{9\baselineskip}
\begin{longtable}[]{@{}
  >{\raggedright\arraybackslash}p{(\linewidth - 4\tabcolsep) * \real{0.3333}}
  >{\raggedright\arraybackslash}p{(\linewidth - 4\tabcolsep) * \real{0.3333}}
  >{\raggedright\arraybackslash}p{(\linewidth - 4\tabcolsep) * \real{0.3333}}@{}}
\toprule\noalign{}
\begin{minipage}[b]{\linewidth}\raggedright
\(m\)
\end{minipage} & \begin{minipage}[b]{\linewidth}\raggedright
Supports requiring finite exclusions
\end{minipage} & \begin{minipage}[b]{\linewidth}\raggedright
Larger-support arithmetic bound
\end{minipage} \\
\midrule\noalign{}
\endhead
\bottomrule\noalign{}
\endlastfoot
17 & 9,10 & \(B_m\le28\) for \(N\ge11\) \\
18 & 9,10,11 & \(B_m\le30\) for \(N\ge12\) \\
19 & 9,10,11 & \(B_m\le38\) for \(N\ge12\) \\
20 & 10,11 & \(B_m\le48\) for \(N\ge12\) \\
21 & 10,11 & \(B_m\le60\) at \(N=12\), and \(\le50\) for \(N\ge13\) \\
\end{longtable}
}

At seventeen arrows, the retained eight-vertex value 19 excludes every
smaller host by isolate padding. At eighteen use \(M_V(8)=32\). At
nineteen, twenty and twenty-one, the missing-pair bounds on the
respective smaller supports are 37,52,51. These are independent of the
global values being proved.

The seventeen- and eighteen-arrow joint-degree exclusions therefore give
the claimed maxima and \(M(9,17)=24\). At nineteen arrows the complete
positive-support exclusions give 35,38,40 on orders 9,10,11. The general
source/sink bound 36 suffices at ten and eleven. At nine vertices it
sharpens to 25: degree-one deletion gives \(M(8,18)=21\), degree two
gives \(M(8,17)+6=25\), and degree at least three gives \(A_{16}+5=24\).
An isolate leaves the retained value \(M(8,19)=22\). Thus the fixed-host
claims also cover every graph.

At twenty arrows, the positive ten-vertex possibilities above 53 are the
regular family, the two directional single-transfer families, and eleven
joint-degree multisets comprising double transfers or a transfer in both
directions. The existing regular and directional bounds are 53; the
complete eleven-profile exclusion supplies the same bound for the
remainder. On eleven vertices, exceeding 53 requires \(q=r=9\): there
are two degree-one vertices in each direction and all other degrees two.
The complete two-defect exclusion retains all overlaps of the low sets
and gives 53. This finishes the global bound without the historical
nineteen-arrow upper 53 or the comparison-host argument.

At twenty-one arrows, the ten-vertex exclusions bound the one-excess
family and the twelve remaining joint-degree multisets by 56. All other
positive degree counts have cap at most 54. At eleven vertices, the
one-low-in/one-low-out family has upper 56; the remaining cases above 58
have \((q,r)=(9,10)\) or \((10,9)\) and are included in the complete
eleven-profile exclusion with upper 58. The five \(q=r=9\) profiles in
that inventory already have cap 57. The source/sink bound 53 applies on
both hosts, so their literal values 56 and 58 are exact. The
twelve-vertex literal attains the arithmetic upper 60. No separate
nineteen-arrow upper 49 is needed.

Every lower bound is a saved literal graph with complete winning ranks
and a closed losing region independently checked. Adding isolates
preserves its score. Least attainment at seven and eight follows from
\(M(6,12)=7\) and \((M(7,13),M(7,14),M(7,15),M(7,16))=(10,10,11,12)\).
The displayed smaller-host values handle budgets 17,19,21; at eighteen
use \(M_V(8)=32\), and at twenty a nine-vertex graph has only52 missing
pairs. This proves all assertions. \(\square\)

The next endpoints use exactly the same reduction: at twenty-two arrows,
\(D_{22}=63\), positive supports at least twelve have cap 74, and
supports at most ten have at most 68 missing pairs. Only order eleven
remains. At twenty-three, \(D_{23}=79\), positive supports at least
thirteen have cap 76, and supports at most ten have at most 67 missing
pairs. The twenty-two- and twenty-three-arrow sections complete those
residual families, in that order.

\subsubsection{What a complete finite exclusion
establishes}\label{what-a-complete-finite-exclusion-establishes}

Each producer lists every permitted degree profile. Joint profiles
assign incoming degrees to sorted outgoing degrees, sorting only within
equal-outdegree classes. Zero entries are retained whenever the ordinary
reduction permits them. Adjacency rows are assigned completely. Future
vertices may be permuted only when their prescribed degree pairs and
already fixed incoming prefixes agree; these permutations fix every
completed row. Choosing every possible cardinality and a prefix in each
such class retains a representative of every completion. Earlier
vertices remain individual choices, loops are forbidden, and opposite
arrows are independent. Matching-based families retain their separately
proved complete factorization and normalization.

Pruning subtracts permanent losses as a union. Source camping requires a
different initial cop; predecessor camping retains the source-equals-cop
exception. A column loss applied to an unfinished row uses a full
column, so a future direct arrow cannot be subtracted twice.
Guarded-path tests grant the messenger every individually possible
future arrow and keep only fixed cop arrows, permitting the initial
escape and terminal receipt. Failure in this easier game is permanent.
Every surviving leaf is evaluated in the actual alternating game. Early
exit uses the current count plus every unprocessed eligible pair as an
upper bound. An incumbent rises only after a literal larger score is
evaluated; earlier prunes remain valid under the larger final bound.
Complete traversal therefore excludes every score above that bound.

\subsubsection{Evidence and computational
boundary}\label{evidence-and-computational-boundary}

The following records, under \texttt{research/\allowbreak{}}, specify the exact
completed families, source hashes, literal graphs, and independent
reviews. Counts refer to residual families; inherited finite results
retain their original evidence.

\begin{itemize}
\tightlist
\item
  \textbf{Budgets 12--13.} \texttt{endpoint40\_\allowbreak{}tiny\_\allowbreak{}live\_\allowbreak{}review.\allowbreak{}json}
  supplies complete 65- and 88-profile exclusions, independently counted
  inventories, and independent checks on 6,651 graphs, 44,376 prefixes
  and 662,312 specific camping states. It does not claim a second full
  graph traversal.
\item
  \textbf{Budgets 14--16.} \texttt{finish\_\allowbreak{}budget\_\allowbreak{}manifest.\allowbreak{}json} and
  \texttt{finish-\allowbreak{}budget-\allowbreak{}ordinary-\allowbreak{}review.\allowbreak{}json} supply complete
  prescribed-degree residual families, with 118,547 leaf games
  independently replayed, together with the reviewed eight-vertex
  bounds.
\item
  \textbf{Budgets 17--18.}
  \texttt{finish\_\allowbreak{}eighteen\_\allowbreak{}review\_\allowbreak{}manifest.\allowbreak{}json} and
  \texttt{finish\_\allowbreak{}eighteen\_\allowbreak{}review\_\allowbreak{}seventeen.\allowbreak{}json} supply exact
  joint-profile set comparisons, with 73,472 package leaf games and
  26,662 stronger seventeen-arrow leaf games independently replayed.
\item
  \textbf{Budget 19.}
  \texttt{sparse\_\allowbreak{}finish\_\allowbreak{}nineteen\_\allowbreak{}exact\_\allowbreak{}manifest.\allowbreak{}json} and
  \texttt{finish\_\allowbreak{}nineteen\_\allowbreak{}review\_\allowbreak{}manifest.\allowbreak{}json} supply 153 complete
  joint profiles and 291,301 independently replayed leaf games; the
  literal host values are 35, 38 and 40.
\item
  \textbf{Budget 20.} \texttt{finish-\allowbreak{}sparse-\allowbreak{}review.\allowbreak{}json} supplies
  independent coverage and game review of the new residual families,
  with 19,286 leaf games. Inherited regular and directional-transfer
  bounds are explicitly identified.
\item
  \textbf{Budget 21.} \texttt{sparse\_\allowbreak{}finish\_\allowbreak{}twentyone\_\allowbreak{}manifest.\allowbreak{}json}
  and \texttt{finish-\allowbreak{}twentyone-\allowbreak{}review.\allowbreak{}json} supply complete matching and
  joint-degree families, independent alternating-game leaf replays and
  permanent-loss checks; the literal host values are 56, 58 and 60.
\end{itemize}

The independent game replays use a separately written alternating-state
attractor. For literal attainers, finite winning ranks and all
losing-closure conditions are checked. Except where a record explicitly
supplies a second generator, replay reuses the independently reviewed
producer traversal. The unrestricted upper exclusions at budgets 12
through 21 remain ordinary finite proofs with auditable computational
premises. Exact retained lower graph scores are separately formalized in
\texttt{SmallSparseWitnesses.\allowbreak{}lean}; those witnesses do not establish the
unrestricted upper bounds.

The eight-vertex row is now completely classified: the sparse,
lower-main and upper-main receipts are
\texttt{endpoint40\_\allowbreak{}n8\_\allowbreak{}sparse\_\allowbreak{}receipt.\allowbreak{}json},
\texttt{endpoint40\_\allowbreak{}n8\_\allowbreak{}low\_\allowbreak{}receipt.\allowbreak{}json}, and
\texttt{endpoint40\_\allowbreak{}n8\_\allowbreak{}high\_\allowbreak{}receipt.\allowbreak{}json}, with the separate
high-budget review in \texttt{finish-\allowbreak{}n8-\allowbreak{}review.\allowbreak{}json}; retained formulas
and literal completions cover the intervening and dense ranges. Together
with the other completed nine- and ten-vertex families, these close
every finite exception. The final main-cell values are \(M(9,20)=37\),
\(M(9,21)=39\), \(M(9,22)=40\), \(M(9,23)=44\), and
\((M(10,22),M(10,23),M(10,24),M(10,25))=(60,61,61,65)\). The complete
finite tables identify their matching graph records; the approved
small-case registry pins their separate exclusions and independent
reviews.

\subsection{Exact extrema with at most seven
arcs}\label{exact-extrema-with-at-most-seven-arcs}

\textbf{Proposition 23.} Allowing arbitrarily many vertices, including
isolated vertices, \[
\bigl(M_E(0),M_E(1),\ldots,M_E(6)\bigr)
  =(0,0,0,0,1,1,2).
\]

\textbf{Proof.} Write \(A=\{v:d^+(v)\ge2\}\), \(B=\{v:d^-(v)\ge2\}\),
\(q=|A|\), \(r=|B|\), and \(b=|A\cap B|\). Every counted pair belongs to
\(A\times B\), excludes the diagonal, and excludes arcs. Consequently \[
I(G)\le qr-b-e(A,B),\qquad
e(A,B)\ge\max\{0,2q+2r-m\},\qquad
q,r\le\lfloor m/2\rfloor.                 \tag{23.1}
\] For completeness, the second inequality follows by counting the arcs
leaving \(A\) and those entering \(B\): their respective totals are at
least \(2q\) and \(2r\), and their union has at most \(m\) arcs.

A counted pair \((s,t)\) requires two arcs leaving \(s\) and two
entering \(t\). These four arcs are distinct because \(s\to t\) is
absent. Thus \(I(G)=0\) when \(m\le3\).

We will use two elementary camping observations. If \(d^-(z)=0\),
\(z\to t\), and \(d^-(t)=2\), then no indirect pair ending at \(t\) is
guaranteed. Indeed, let \(p\) be its other predecessor. An indirect
source is neither \(z\) nor \(p\). A pursuer starting at \(p\) can wait
there: the messenger cannot reach \(z\), and cannot enter \(p\) without
being captured before reaching \(t\). Dually, if \(d^+(z)=0\),
\(s\to z\), and \(d^+(s)=2\), then no indirect pair starting at \(s\) is
guaranteed. A pursuer waits at the other successor of \(s\); the
messenger's only alternative departure reaches a sink different from the
indirect recipient.

Suppose first that \(m=4\) or \(m=5\). Then \(q,r\le2\). If either is
zero there is nothing to count, and if both are one the bound is
immediate. When \(q=r=2\), (23.1) gives \(I(G)\le m-4-b\le1\). When
\(q=1\) and \(r=2\), two counted pairs would have the same source and
distinct recipients. The source's two outgoing arcs are then disjoint
from the at least four arcs entering those recipients, which would
require \(m\ge6\). The case \(q=2,r=1\) follows by the same count with
sources and recipients interchanged: at least four arcs leave the two
counted sources, and the recipient's two incoming arcs are disjoint from
them. Hence \(I(G)\le1\).

Now let \(m=6\), so that \(q,r\le3\). If \(q=1\), three counted pairs
would require two outgoing arcs from their common source and six
disjoint arcs entering their three recipients; this is impossible. Thus
\(I(G)\le2\). The same argument applies when \(r=1\); cases with a zero
coordinate are again immediate. If \(q=r=2\), (23.1) gives
\(I(G)\le2-b\le2\).

Consider \(q=2,r=3\). Every arc ends in \(B\), and every vertex of \(B\)
has indegree exactly two. If \(A\subseteq B\), then \(b=2\) and (23.1)
gives \(I(G)\le6-2-4=0\). Otherwise choose \(z\in A\setminus B\). It has
indegree zero and at least two distinct successors in \(B\). The first
camping observation eliminates both corresponding recipient columns. At
most one column remains, containing at most the two sources in \(A\).
Thus \(I(G)\le2\).

For \(q=3,r=2\), every arc starts in \(A\) and every vertex of \(A\) has
outdegree exactly two. If \(B\subseteq A\), the same count gives zero.
Otherwise a vertex \(z\in B\setminus A\) has outdegree zero and at least
two predecessors in \(A\). The second camping observation eliminates
their two source rows, leaving at most one row and two recipients. Again
\(I(G)\le2\).

Finally, if \(q=r=3\), all six arcs go from \(A\) to \(B\), and all
relevant indegrees and outdegrees equal two. If \(A=B\), (23.1) gives
\(I(G)\le9-3-6=0\). Otherwise a vertex of \(A\setminus B\) has indegree
zero and eliminates two recipient columns by the first camping
observation. The remaining column has two distinct predecessors in
\(A\), so at most one of the three sources can contribute an indirect
pair. This completes the upper bounds in all cases.

For the lower bounds, the four-arc diamond \[
E=\{(0,1),(0,2),(1,3),(2,3)\}
\] guarantees the indirect pair \((0,3)\): the messenger chooses the
branch not occupied by the pursuer and delivers on the next move. Adding
a disjoint arc gives the five-arc example. The six-arc graph \[
E=\{(0,1),(0,2),(1,3),(1,4),(2,3),(2,4)\}
\] guarantees both \((0,3)\) and \((0,4)\). If the pursuer starts on one
branch, choose the other; if it starts at either sink, choose either
branch. In each case it cannot reach the selected branch on its next
move, and receipt then occurs on the messenger's second move. Hence
\(I(G)=2\), as required. \(\square\)

These arguments establish the arc-only maximum over every finite vertex
set; no restriction on the number of nonisolated vertices is used.

\textbf{Proposition 24 (with a finite exhaustive computation).} We have
\[
M_E(7)=2.
\]

\textbf{Proof.} First consider a weakly connected digraph with at most
seven arcs and at least seven vertices. Forget the directions but retain
one undirected edge for each arc, so that opposite arcs become parallel
edges. The resulting connected multigraph has at most as many edges as
vertices. It is therefore a tree or has exactly one cycle; a parallel
pair is allowed as a cycle of length two.

In a tree, any reachable indirect source--recipient pair has an internal
vertex separating its endpoints. A pursuer can start and wait there,
preventing delivery. Thus no indirect pair is guaranteed. In the
unicyclic case, a pair with either endpoint outside the unique cycle
likewise has an internal separating vertex, unless its endpoints are
adjacent in the underlying graph. For adjacent endpoints off the cycle,
a directed route exists only if their connecting arc has the required
direction, making the pair direct. Thus every counted pair must have
both endpoints on the cycle.

There are exactly two simple undirected paths along that cycle between
two distinct cycle vertices. If only one is a directed path from the
source to the recipient, it is either a direct arc or has an internal
vertex where the pursuer can wait. If neither is directed, the recipient
is unreachable. Excursions into attached trees cannot bypass an internal
cycle vertex. Consequently a counted pair requires both cycle paths to
be directed from its source to its recipient. This forces the cycle's
unique divergence at that source and unique convergence at that
recipient; every other cycle vertex has one entering and one leaving
cycle edge. A fixed orientation can therefore support at most one such
ordered pair. A cycle of length two consists of opposite arcs and
supports no indirect pair. We have proved that every weakly connected
digraph with at most seven arcs and at least seven vertices has
\(I(G)\le1\).

For a general digraph with seven arcs, distinct weak components
contribute additively to \(I(G)\). Each component with a positive
contribution uses at least four arcs by Proposition 23, so at most one
component contributes. If that component uses at most six arcs,
Proposition 23 bounds its contribution by two. If it uses seven arcs and
has at least seven vertices, the preceding argument bounds its
contribution by one. In the only remaining case it has at most six
vertices. Adding isolated vertices to make six does not change its
indirect guarantee count, and the complete six-vertex enumeration gives
a maximum of two at seven arcs.

This last finite assertion is the computational dependency of the
proposition. It is reproduced by \texttt{src/\allowbreak{}enumerate\_\allowbreak{}small\_\allowbreak{}n6.\allowbreak{}cpp};
the result is recorded with \texttt{complete:\ true} in
\texttt{research/\allowbreak{}solver-\allowbreak{}enumeration-\allowbreak{}n6.\allowbreak{}json}. The enumeration covers
every six-vertex digraph up to the degree-ordering reduction described
in the computational section; the witness is also checked by the
separate explicit-turn solver. The structural reduction above is what
extends that finite calculation to unrestricted vertex counts.

Finally, add one disjoint arc to the six-arc, two-recipient construction
in Proposition 23. This preserves its two guarantees and supplies seven
arcs, proving the matching lower bound. \(\square\)

Theorem 49 below settles the next four budgets, eight through eleven.

\subsection{Exact arc maxima from eight through
eleven}\label{exact-arc-maxima-from-eight-through-eleven}

The six-vertex table supplies useful lower bounds, but does not itself
settle an arc budget when the number of vertices is unrestricted. At
eleven arcs the distinction is essential: seven vertices attain eight
guarantees, exceeding the six-vertex maximum of seven.

\textbf{Theorem 49 (with finite exhaustive certificates).} Over all
finite vertex sets,
\[\bigl(M_E(8),M_E(9),M_E(10),M_E(11)\bigr)=(4,5,6,8).\]

\textbf{Proof of the support reduction.} Indirect guarantees add over
weak components. Every component with a positive contribution has at
least four arcs. At these budgets there can be at most two positive
components. If there are two, each has at most seven arcs; Propositions
23 and 24 give a combined maximum of at most three. This cannot improve
any of the displayed lower bounds.

There is therefore only one relevant positive component. We argue in
increasing order of the budget \(m=8,9,10,11\). If this component has
fewer than \(m\) arcs, the already established smaller-budget maximum
bounds its contribution by the displayed value at \(m\). If it uses all
\(m\) arcs and has at least \(m\) vertices, the tree-or-unicyclic
argument in Proposition 24 gives at most one indirect guarantee. That
argument only needs a connected underlying multigraph with at most as
many edges as vertices, so applies unchanged at the present budgets.

In the remaining case the component has at most \(m-1\) vertices. Pad it
with isolates to exactly \(m-1\) vertices; its arc count and indirect
count do not change. Hence each upper bound reduces to the following one
finite counterexample test:

{\def\LTcaptype{none} 
\Needspace{8\baselineskip}
\begin{longtable}[]{@{}lrr@{}}
\toprule\noalign{}
Arc budget \(m\) & Host order & Excluded threshold for \(I\) \\
\midrule\noalign{}
\endhead
\bottomrule\noalign{}
\endlastfoot
8 & 7 & 5 \\
9 & 8 & 6 \\
10 & 9 & 7 \\
11 & 10 & 9 \\
\end{longtable}
}

\textbf{Complete finite upper checks.} Both enumerators list all
prescribed outdegree profiles with total \(m\) in sorted order and then
all loopless outgoing rows. As in the canonical reduction of Theorem 46,
future vertices with equal prescribed degree and equal incoming pattern
from completed rows can be permuted while leaving every completed row
unchanged. Choosing the new row's neighbors as a prefix of each such
contiguous cell retains a representative of every graph. The second
implementation instead uses descending degrees and suffixes. Both
reduction arguments fix the current vertex and all completed vertices.

The source and predecessor camping conditions of Theorem 46 supply
necessary local rejections. A completed degree-two source that points to
a prescribed sink also contributes zero, by the earlier sink
obstruction. Otherwise a source of degree \(d\ge2\) contributes at most
\(n-1-d\). A target of current indegree \(d\) contributes at most
\(n-1-\max(2,d)\), or zero if its indegree cannot reach two.

One additional numerical bound makes these sparse checks short. Let
\(Q\) contain all sources not yet excluded, and put \(q=|Q|\). Every
member of \(Q\) has prescribed outdegree at least two. If the eventual
number of indegree-at-least-two targets is \(r\), the same
degree-incidence count as (1) gives
\[I\le qr-\max(0,q+r-n)-\max(0,2q+2r-m).\] To bound \(r\) at a partial
assignment, let \(R\) be the number of arcs still to be assigned. A
target currently of indegree \(d\) needs at least \(\max(0,2-d)\) of
these arcs to reach indegree two. Discard targets with too few remaining
possible predecessors, and select as many of the others as possible
using the smallest such costs first, with total cost at most \(R\). The
resulting \(r_{\max}\) is an upper bound: each remaining arc enters only
one target. Maximizing the displayed bound over \(0\le r\le r_{\max}\)
is therefore a safe partial rejection. No sufficiency is assumed for
these degree or camping tests.

Every remaining complete graph is evaluated by the original finite
reachability game. The primary implementation uses synchronous
messenger-position masks; the independent implementation uses an
explicit alternating-state queue attractor. Both start with terminal
receipt states and compute the least winning region, including waiting
and all pursuer responses. They reject every candidate at every one of
the four thresholds. The exact enumeration counts are supplied in the
reproduction receipt; both runs report complete enumeration rather than
termination at their time limit. This proves the four upper bounds.

\textbf{Lower witnesses.} The witnesses at eight, nine and ten arcs are
the graphs in the six-vertex table. The eleven-arc witness has seven
vertices and the following outgoing neighborhoods:

{\def\LTcaptype{none} 
\Needspace{11\baselineskip}
\begin{longtable}[]{@{}ll@{}}
\toprule\noalign{}
Vertex & Outgoing neighbors \\
\midrule\noalign{}
\endhead
\bottomrule\noalign{}
\endlastfoot
0 & 1 \\
1 & 0 \\
2 & 3 \\
3 & 4, 5 \\
4 & 0, 2 \\
5 & 3, 6 \\
6 & 1, 2 \\
\end{longtable}
}

Its guaranteed indirect pairs are
\[ (3,0),(3,1),(3,2),(4,1),(5,0),(5,1),(5,2),(6,0). \] Both independent
game solvers check each of the four lower witnesses, and the saved
explicit ranks satisfy every legal-move obligation. They attain the
upper bounds and complete the proof. \(\square\)

The complete reproduction command is
\texttt{python3\ src/\allowbreak{}small\_\allowbreak{}arcs\_\allowbreak{}resume.\allowbreak{}py}. Its two exhaustive sources
are \texttt{src/\allowbreak{}small\_\allowbreak{}arcs\_\allowbreak{}sparse\_\allowbreak{}primary.\allowbreak{}cpp} and
\texttt{src/\allowbreak{}small\_\allowbreak{}arcs\_\allowbreak{}sparse\_\allowbreak{}independent.\allowbreak{}cpp}; the receipt, source
hashes, complete enumeration summaries and all four lower rank
certificates are in
\texttt{research/\allowbreak{}small-\allowbreak{}arcs-\allowbreak{}resume-\allowbreak{}verification.\allowbreak{}json}. The final full
reproduction took about eighteen seconds. The unrestricted maxima at
budgets 8, 9 and 10 are now formalized in
\texttt{SmallBudgetEight.\allowbreak{}lean}, \texttt{SmallBudgetNine.\allowbreak{}lean} and
\texttt{SmallBudgetTen.\allowbreak{}lean}. At budget 11,
\texttt{SmallSparseWitnesses.\allowbreak{}lean} formalizes the lower witness, while
the unrestricted upper bound retains its ordinary and finite
computer-assisted proof boundary.

\subsection{Small vertex orders}\label{small-vertex-orders}

The smallest order attaining the vertex bound is twelve. This statement
is about attaining the bound, not about the first noncomplete network in
which all pairs are guaranteed: the latter already occurs at seven
vertices.

\textbf{Theorem 46 (with finite exhaustive certificates).} The seven-,
eight-, and eleven-vertex maxima are
\[M_V(7)=21,\qquad M_V(8)=32,\qquad M_V(11)=85.\] In addition,
\(M(11,23)=78\) and \(M(11,24)=79\). The intermediate vertex bounds
proved here at orders nine through eleven are:

{\def\LTcaptype{none} 
\Needspace{7\baselineskip}
\begin{longtable}[]{@{}lrr@{}}
\toprule\noalign{}
\(n\) & Lower bound for \(M_V(n)\) & Upper bound for \(M_V(n)\) \\
\midrule\noalign{}
\endhead
\bottomrule\noalign{}
\endlastfoot
9 & 48 & 53 \\
10 & 62 & 69 \\
11 & 85 & 85 \\
\end{longtable}
}

There is no universally winning two-in/two-out regular digraph on any of
\(7,8,9,10,11\) vertices. Together with the previous exact results
through six vertices and Theorem 40, this implies that twelve is the
smallest order \(n\ge4\) for which \(M_V(n)=n(n-3)\).

\textbf{Proof of the lower bounds.} The finite witnesses at orders seven
through ten have respectively \(21,24,24,28\) arcs and guarantee every
pair. Their complete graph and winning-strategy certificates are
provided in
\texttt{research/\allowbreak{}small-\allowbreak{}exceptional-\allowbreak{}general-\allowbreak{}certificates.\allowbreak{}json}, with the
improved nine-vertex witness in
\texttt{research/\allowbreak{}root-\allowbreak{}small-\allowbreak{}structured-\allowbreak{}9-\allowbreak{}certificate.\allowbreak{}json}. The earlier
eleven-vertex witness with 26 arrows remains in
\texttt{research/\allowbreak{}small-\allowbreak{}exceptional-\allowbreak{}structured-\allowbreak{}11-\allowbreak{}certificate.\allowbreak{}json}. A
new 25-arrow witness improves its count from 84 to 85; its outgoing rows
are \[
(1,7,9),\ (2,5),\ (3,8),\ (5,10),\ (0,6),\ (1,4,6),\
(8,10),\ (1,8),\ (6,7,9),\ (3,4),\ (0,2).
\] Two independent attractor algorithms agree; each stored
alternating-state rank satisfies all messenger-choice and universal
pursuer-response obligations. Thus \(I=n(n-1)-m\) gives
\(21,32,48,62,85\) respectively. The new literal witnesses are
consolidated in \texttt{research/\allowbreak{}small-\allowbreak{}case-\allowbreak{}results.\allowbreak{}json} and
independently replayed by
\texttt{src/\allowbreak{}small\_\allowbreak{}case\_\allowbreak{}certificate\_\allowbreak{}review.\allowbreak{}py}. For concreteness, the
seven-vertex witness is the following graph, whose guarantees require at
most four messenger moves:

{\def\LTcaptype{none} 
\Needspace{11\baselineskip}
\begin{longtable}[]{@{}ll@{}}
\toprule\noalign{}
Vertex & Outgoing neighbors \\
\midrule\noalign{}
\endhead
\bottomrule\noalign{}
\endlastfoot
0 & 1, 3, 6 \\
1 & 0, 4, 6 \\
2 & 0, 5, 6 \\
3 & 1, 2, 5 \\
4 & 0, 2, 3 \\
5 & 1, 3, 4 \\
6 & 2, 4, 5 \\
\end{longtable}
}

\textbf{Seven-vertex upper bound.} If a graph has at most thirteen arcs,
the arc-only bounds give \(I\le20\); it suffices here to use
monotonicity of \(M_E\) (pad by disjoint one-way arcs) and
\(M_E(13)\le20\). A graph with at least twenty-one arcs has at most
\(42-21=21\) missing pairs. It remains to exclude \(I\ge22\) for
\(14\le m\le20\).

We describe the finite reduction and every rejection rule. Write
\(O_s=N^+(s)\), \(P_t=N^-(t)\) and \(C_c=\{c\}\cup N^+(c)\). A counted
pair \((s,t)\) necessarily satisfies

\begin{enumerate}
\def\labelenumi{\arabic{enumi}.}
\tightlist
\item
  \(s\ne t\), \(t\notin O_s\), \(|O_s|\ge2\) and \(|P_t|\ge2\);
\item
  \(O_s\not\subseteq C_c\) for every \(c\ne s\);
\item
  \(P_t\not\subseteq C_c\) for every \(c\ne s\).
\end{enumerate}

The degree conditions are the one-neighbor camping obstructions. For
condition 2, a pursuer waiting at \(c\) captures the messenger after its
first departure from \(s\). For condition 3, it waits at \(c\) and
captures the messenger upon the first entry into a predecessor of \(t\),
before the next messenger move can deliver. Indefinite waiting by the
messenger also fails to deliver. These arguments respect receipt
priority because an indirect source is not itself a predecessor of its
recipient.

Let \(U(G)\) be the number of pairs satisfying these three conditions.
Then \(I(G)\le U(G)\). A source contributes at most four pairs, so
\(I\ge22\) requires at least six sources of outdegree at least two.
Every graph can be relabeled in nondecreasing outdegree order, allowing
only the first row to have degree below two. We enumerate every loopless
outgoing row mask in that order, with total degree between fourteen and
twenty. This includes every possible graph, with harmless repetitions.
At a partial assignment, the following bounds allow safe pruning:

\begin{itemize}
\tightlist
\item
  A completed row of degree \(d\ge2\) contributes at most \(6-d\), or
  zero if a known completed row already violates condition 2. Each
  remaining row contributes at most \(6-\max(2,d_0)\), where \(d_0\) is
  the smallest outdegree still allowed by the ordering.
\item
  A target currently of indegree \(d\) contributes at most
  \(6-\max(2,d)\). Its contribution is zero if the remaining rows cannot
  raise its indegree to two.
\item
  The total arc count must remain extendible to the interval
  \([14,20]\).
\end{itemize}

No omitted edge can repair a rejected completed source: the covering
closed neighborhood only grows. The incoming-degree bounds likewise use
upper estimates on all possible completions. At every surviving complete
graph we count \(U(G)\) directly. The exhaustive enumeration has
\(109,941,961\) recursion nodes and \(4,299,624\) complete leaves before
an optional indegree tie normalization. None has \(U(G)\ge22\). A
separate implementation uses boolean adjacency, descending outdegree
order, and no indegree tie normalization; it completes \(123,280,819\)
nodes and the same \(4,299,624\) leaves, again with no survivor. This
proves \(I\le21\) in the remaining arc range. The lower witness proves
equality.

\textbf{Eight-vertex upper bound.} At most fifteen arcs give
\(I\le M_E(15)\le30\); at least twenty-four arcs leave at most
thirty-two missing pairs. To exclude \(I\ge33\), it therefore suffices
to examine \(16\le m\le23\). Each source contributes at most five pairs,
so at least seven sources must have outdegree at least two. The same
camping conditions and partial bounds used above apply with eight
vertices.

For this slightly larger search we also remove labeling redundancy. Fix
a nondecreasing outdegree sequence. After the outgoing rows numbered
less than \(r\) have been chosen, partition the unassigned vertices into
cells with equal prescribed outdegree and identical adjacency from all
completed rows. Vertices in such a cell may be permuted without changing
any completed row or degree requirement. Therefore, when choosing row
\(r\), its neighbors among the future vertices of each cell may be
required to form a prefix of that cell. This loses no graph: permute
within each cell to put the neighbors first, fixing row \(r\) and all
previous vertices, and continue recursively. The new cells are obtained
by splitting the old cells, so their contiguity is preserved. The
analogous suffix rule and reversed degree order are equally complete.

The primary search completes all \(178\) degree sequences, with
\(946,676,955\) recursion nodes and \(22,465,453\) complete leaves. Only
\(373\) leaves satisfy the necessary pair bound \(U(G)\ge33\). Exact
least-attractor evaluation rejects all of them; the largest actual \(I\)
among these survivors is seventeen. Thus no graph has \(I\ge33\),
proving \(M_V(8)=32\). The ordinary canonicalization and every pruning
rule were independently audited; as an additional small positive check,
every labeled loopless four-vertex graph was canonicalized in both
prefix and suffix conventions and its relabeling and prefix conditions
checked. Independent complete suffix enumeration and the explicit
alternating-state game provide a second finite verification: its
completed intervals \([0,160)\) and \([160,178)\) cover all degree
profiles. They visit \(1,179,975,770\) nodes and \(63,717,034\) complete
leaves, with \(1,115\) static survivors, all rejected. The differing
counts reflect the opposite degree order and suffix convention, which
choose a different redundant collection of representatives. Both
implementations give the same zero-survivor conclusion.

\textbf{Eleven vertices and twenty-three arrows.} We prove
\(M(11,23)=78\). If \(q\) vertices have outdegree at least two, the
other \(11-q\) vertices emit at most \(11-q\) arrows and cannot be
indirect sources. The eligible rows therefore contain at most
\[10q-\bigl(23-(11-q)\bigr)=9q-12\] missing pairs. This is at most 78
for \(q\le10\). The incoming argument is identical: an indirect
recipient requires at least two predecessors. Consequently, a graph with
at least 79 guarantees must have all indegrees and outdegrees at least
two. There is exactly one outdegree-three vertex \(U\) and one
indegree-three vertex \(V\), which may coincide; all other degrees are
two.

The bipartite arrow graph has a perfect matching. Indeed, a source set
\(S\) emits at least \(2|S|\) arrows, whereas fewer than \(|S|\)
recipient columns have total capacity at most \(2(|S|-1)+1\). This
contradicts a failure of Hall's condition. The matching is a
fixed-point-free permutation. Up to simultaneous relabeling there are
fourteen possible cycle partitions of eleven into parts at least two.
Rotating cycles and swapping cycles of equal length places \(U\) at the
first vertex of one cycle of each distinct length. Allowing all eleven
positions for \(V\) retains every relation between the exceptional
vertices and gives 330 patterns in total. The remaining arrows have
outgoing degree one except two at \(U\) and incoming capacity one except
two at \(V\). Enumerating all residual rows, avoiding loops and matching
arrows, covers the entire class. The two residual neighbors at \(U\) are
chosen in increasing order.

During this enumeration, source and predecessor camping identify
permanent losing pairs, retaining the exception when the initial source
equals the camping vertex. No camping class is excluded in advance. A
partial cop row contains only arrows that persist; a predecessor set is
used only when its prescribed capacity is full. The permanent losing
pairs are copied separately at each recursion depth. A row contributes
at most the smaller of its final missing-pair count and the number of
targets not already direct, diagonal, or permanently losing.

The guarded-path bound uses every still-possible messenger arrow and
only fixed cop arrows. It permits the first move to escape the cop's
closed neighborhood and the final move to deliver anywhere; only
interior positions must avoid that neighborhood. Thus a failed
optimistic route certifies a losing pair in every completion. Surviving
graphs receive exact recipient least-attractor evaluation. Solved
recipient counts plus the remaining potentially winning pairs form an
upper bound, so the game calculation may stop once this is at most 78.

The complete ranges of cycle partitions \(1\)--\(6\) and \(7\)--\(14\)
visit 121,790,823 nodes and reach 1,152,430 exact game decisions. Every
decision rejects a score of at least 79. The recorded value \(-1\) is an
early rejection sentinel, not a graph score. The source, complete
receipts, and independent ordinary audit are respectively
\texttt{src/\allowbreak{}small\_\allowbreak{}sparse\_\allowbreak{}one\_\allowbreak{}excess.\allowbreak{}cpp},
\texttt{research/\allowbreak{}small-\allowbreak{}sparse-\allowbreak{}one-\allowbreak{}excess11-\allowbreak{}partitions1-\allowbreak{}6.\allowbreak{}json},
\texttt{research/\allowbreak{}small-\allowbreak{}sparse-\allowbreak{}one-\allowbreak{}excess11-\allowbreak{}partitions7-\allowbreak{}14.\allowbreak{}json}, and
\texttt{research/\allowbreak{}small-\allowbreak{}plateau-\allowbreak{}one-\allowbreak{}excess-\allowbreak{}review.\allowbreak{}md}. The independently
rank-verified literal 78 witness is included in
\texttt{research/\allowbreak{}small-\allowbreak{}case-\allowbreak{}results.\allowbreak{}json}, proving equality.

\textbf{Eleven vertices and twenty-four arrows.} We prove
\(M(11,24)=79\). The certified attaining graph has outgoing rows \[
(1,7,9),\ (2,4),\ (3,8),\ (4,6),\ (0,5),\ (1,6),\
(7,10),\ (1,8),\ (5,7,9),\ (3,10),\ (0,2).
\] Suppose a graph at this budget has \(I\ge80\). If exactly ten source
rows have outdegree at least two, at least 23 arrows start in those
rows, so they contribute at most \(100-23=77\) pairs. At most nine
eligible rows give at most \(9\cdot8=72\). The incoming-degree count is
identical. Thus every indegree and outdegree is at least two. The total
surplus above two is two on each side, so the outgoing degree profiles
are \((4,2,\ldots,2)\) and \((3,3,2,\ldots,2)\), and no degree exceeds
four.

There are 86 missing pairs, permitting at most six failures. Source
camping would make an entire row fail. To lose at most six pairs that
source would need outdegree four. It then uses all outgoing surplus,
while every other cop has a closed neighborhood of size three and cannot
cover its four successors. Hence source camping is impossible.

First suppose there is no predecessor camping either. The preceding
canonical row method, in descending degree order, enumerates both degree
profiles. Incoming feasibility requires enough unassigned tails to reach
indegree two and to escape every completed cop neighborhood covering a
current predecessor set. If no single tail escapes them all, at least
two new predecessors are required; summing these lower bounds over
recipients cannot exceed 24. At a complete graph, each counted pair must
also have a route to a predecessor of its recipient through vertices
outside every stationary cop's closed neighborhood. The final arrow into
the recipient is exempt from this avoidance requirement, so receipt
priority is respected. This gives a necessary pair-count filter before
exact game evaluation.

The complete target-80 run visits 186,865,392 nodes, 2,805,991 complete
camping-free graphs, and 371,060 guarded-path survivors. Its largest
actual count is 79. Every pruned branch has a proved upper below 80. The
frozen receipt is
\texttt{research/\allowbreak{}small-\allowbreak{}loss-\allowbreak{}n11-\allowbreak{}m24-\allowbreak{}camping-\allowbreak{}free-\allowbreak{}80.\allowbreak{}json}.

If predecessor camping occurs at \(t\) with cop \(c\), at least
\(10-d^-(t)-1\) missing pairs fail. Thus \(d^-(t)\) is three or four,
and no other column may be camped. For degree four there are three
cases: \(c\notin N^-(t)\cup\{t\}\), \(c\in N^-(t)\), or \(c=t\). In each
case \(c\) has outdegree four, and its outgoing set is respectively
\(N^-(t)\), \(\{t\}\cup(N^-(t)\setminus\{c\})\), or \(N^-(t)\). All
other degrees are two. For degree three, the exceptional cop-source pair
must be indirect, so \(c\notin N^-(t)\cup\{t\}\). The cop has outdegree
three or four. In the latter case its extra neighbor can be normalized
to one outside vertex, and the second indegree-three vertex has ten
possible labels. In the former case the second outdegree-three vertex
has three symmetry types: \(t\), a predecessor of \(t\), or an outside
vertex. Retaining all ten possible labels of the second indegree-three
vertex gives another 30 configurations. Together with the three
degree-four cases, this is \(3+10+30=43\) configurations.

Exact incoming caps, the prescribed predecessor set at \(t\), and
colored prefix cells preserve all these cases. The colors retain
degrees, recipient status, and predecessor membership. Source camping
and camping at every other column remain forbidden. A complete
enumeration of all 43 configurations visits 113,000,606 nodes, 603,453
leaves, and 66,362 guarded-path survivors; its largest evaluated count
is 73, and no graph reaches 80. This proves \(M(11,24)\le79\), matching
the witness. The source and complete receipt are
\texttt{src/\allowbreak{}small\_\allowbreak{}loss\_\allowbreak{}camped\_\allowbreak{}profiles.\allowbreak{}cpp} and
\texttt{research/\allowbreak{}small-\allowbreak{}loss-\allowbreak{}n11-\allowbreak{}m24-\allowbreak{}camped80.\allowbreak{}json}; the ordinary
completeness audit is
\texttt{research/\allowbreak{}small-\allowbreak{}plateau-\allowbreak{}camped80-\allowbreak{}review.\allowbreak{}md}.

At most 21 arrows give \(I\le71\) by the retained arc bounds. At budgets
22 and 23 the bounds are respectively 79 and 78. We have just proved the
value 79 at budget 24. At least 25 arrows leave at most 85 missing
pairs. The universal 25-arrow witness therefore proves \(M_V(11)=85\).

\textbf{Nonattainment at orders eight through eleven.} By the vertex
equality criterion, an attaining graph would be universally winning and
have indegree and outdegree exactly two everywhere. Any such directed
graph is the union of two arc-disjoint permutation factors. Indeed, its
bipartite tail/head incidence graph is two-regular, so alternating the
edges around each even cycle gives the two perfect matchings. Both
permutations are fixed-point-free. Up to simultaneous relabeling, the
first permutation may be put in canonical cycle form for each partition
of \(n\) having no part one; enumerate every second permutation avoiding
both the diagonal and the first permutation's arcs.

The source and predecessor camping conditions above reject a branch as
soon as a completed two-element neighborhood is covered by a known
closed outgoing neighborhood. For the predecessor condition, at least
\(n-3\ge2\) sources are indirect, so a covering vertex can be avoided as
the initial source. Every remaining complete graph is evaluated by the
finite least-attractor algorithm. The exhaustive counts are:

{\def\LTcaptype{none} 
\Needspace{9\baselineskip}
\begin{longtable}[]{@{}
  >{\raggedright\arraybackslash}p{(\linewidth - 8\tabcolsep) * \real{0.1579}}
  >{\raggedleft\arraybackslash}p{(\linewidth - 8\tabcolsep) * \real{0.2105}}
  >{\raggedleft\arraybackslash}p{(\linewidth - 8\tabcolsep) * \real{0.2105}}
  >{\raggedleft\arraybackslash}p{(\linewidth - 8\tabcolsep) * \real{0.2105}}
  >{\raggedleft\arraybackslash}p{(\linewidth - 8\tabcolsep) * \real{0.2105}}@{}}
\toprule\noalign{}
\begin{minipage}[b]{\linewidth}\raggedright
\(n\)
\end{minipage} & \begin{minipage}[b]{\linewidth}\raggedleft
Cycle partitions
\end{minipage} & \begin{minipage}[b]{\linewidth}\raggedleft
Recursion nodes
\end{minipage} & \begin{minipage}[b]{\linewidth}\raggedleft
Locally admissible complete graphs
\end{minipage} & \begin{minipage}[b]{\linewidth}\raggedleft
Universally winning graphs
\end{minipage} \\
\midrule\noalign{}
\endhead
\bottomrule\noalign{}
\endlastfoot
7 & 4 & 142 & 0 & 0 \\
8 & 7 & 1,446 & 0 & 0 \\
9 & 8 & 10,996 & 535 & 0 \\
10 & 12 & 231,676 & 35,610 & 0 \\
11 & 14 & 2,177,946 & 120,302 & 0 \\
\end{longtable}
}

An independent implementation reverses the permutation assignment order.
Its local tests forbid directed transitive triangles, two distinct
length-two walks with the same ordered endpoints, and duplicate
completed outgoing pairs; these are precisely the same camping
obstructions in adjacency-matrix form. It then uses an explicit
alternating-state AND/OR queue attractor instead of the primary
synchronous bitmask iteration. Both complete enumerations give exactly
the displayed counts. As a positive control, a known universal
twelve-vertex graph passes every incremental pruning test and the game
check of the primary enumeration. Hence equality \(M_V(n)=n(n-3)\) is
impossible at these orders, which lowers the integer upper bound by one.
\(\square\)

The seven-vertex example also gives the first noncomplete universally
winning graph. At orders at most three, the previously established zero
maximum excludes any missing guaranteed pair. At orders four through
six, the density bound for a noncomplete universal graph would require
\(m\le n(n-3)\) and hence \(I\ge2n\), contradicting the exact maxima
\(2,4,7\).

\textbf{Proposition (complete seven-vertex joint table).} Every arrow
budget at order seven is determined:

{\def\LTcaptype{none} 
\Needspace{19\baselineskip}
\begin{longtable}[]{@{}rrrrrr@{}}
\toprule\noalign{}
\(m\) & \(M(7,m)\) & \(m\) & \(M(7,m)\) & \(m\) & \(M(7,m)\) \\
\midrule\noalign{}
\endhead
\bottomrule\noalign{}
\endlastfoot
0 & 0 & 15 & 11 & 30 & 8 \\
1 & 0 & 16 & 12 & 31 & 7 \\
2 & 0 & 17 & 12 & 32 & 6 \\
3 & 0 & 18 & 14 & 33 & 5 \\
4 & 1 & 19 & 15 & 34 & 4 \\
5 & 1 & 20 & 18 & 35 & 2 \\
6 & 2 & 21 & 21 & 36 & 0 \\
7 & 2 & 22 & 14 & 37 & 0 \\
8 & 4 & 23 & 12 & 38 & 0 \\
9 & 5 & 24 & 12 & 39 & 0 \\
10 & 6 & 25 & 11 & 40 & 0 \\
11 & 8 & 26 & 10 & 41 & 0 \\
12 & 9 & 27 & 10 & 42 & 0 \\
13 & 10 & 28 & 9 & & \\
14 & 10 & 29 & 8 & & \\
\end{longtable}
}

\textbf{Proof.} The previously established small-arrow bounds and
witnesses give budgets zero through eleven. The universal 21-arrow graph
gives the central value 21, and the dense results give budgets 33
through 42. The twenty other budgets are certified by complete
independent descending-degree, suffix-cell enumeration, using the
normalization already proved above and an explicit alternating-state
game solver.

For completeness, the additional sparse pruning bound is as follows. If
\(q\) rows and \(r\) columns contain winning indirect pairs, they
include at least \(\max(0,q+r-7)\) diagonal coincidences and at least
\(\max(0,2q+2r-m)\) direct arrows between them. These exclusions are
disjoint, giving \[I\le qr-\max(0,q+r-7)-\max(0,2q+2r-m).\] Partial
assignments overestimate the possible live rows using the permanent
source-camping obstructions. For each potential live column, the cost of
reaching indegree two is \(\max(0,2-d^-_{\rm current})\). The remaining
arrow budget bounds how many of the cheapest columns can become live.
Maximizing the displayed inequality over both dimensions below these
caps is therefore a valid upper bound. Fixed source and column capacity
sums give further necessary bounds. Every surviving complete graph is
checked by the explicit game attractor with receipt priority and
waiting.

The independent runs cover every degree profile at every one of the
twenty budgets, visiting 193,302,277 nodes. None reaches one above the
displayed value. They reproduce all eighteen completed primary
exclusions and also finish the two cases \(m=22,23\); in particular, the
sharp values there are 14 and 12. Literal attaining graphs have complete
winning and losing rank certificates independently replayed by
\texttt{src/\allowbreak{}small\_\allowbreak{}case\_\allowbreak{}certificate\_\allowbreak{}review.\allowbreak{}py}. The upper receipts and
ordinary audit are \texttt{research/\allowbreak{}small-\allowbreak{}case-\allowbreak{}n7-\allowbreak{}independent.\allowbreak{}json} and
\texttt{research/\allowbreak{}small-\allowbreak{}case-\allowbreak{}n7-\allowbreak{}independent-\allowbreak{}review.\allowbreak{}md}; the consolidated
records are in \texttt{research/\allowbreak{}small-\allowbreak{}case-\allowbreak{}results.\allowbreak{}json}. \(\square\)

\textbf{Proposition (further fixed-host values).}
\[M(8,12)=M(9,12)=9,\qquad M(8,13)=11,\qquad M(8,39)=14.\] In
particular, \(M_E(13)\ge11\); no unrestricted equality at budget
thirteen is asserted.

\textbf{Proof.} The seven-vertex 12-arrow witness embeds by isolates at
orders eight and nine. At order eight and budget thirteen, the graph
with rows \[
(5,6,7),\ (6,7),\ (1,5),\ (2,4),\ (0,3),\ (4),\ (1),\ ()
\] has eleven guarantees. Its 1,024 alternating-state rank obligations
are checked directly. Complete fixed-host exclusions of scores 10, 10,
and 12 at \((8,12)\), \((9,12)\), and \((8,13)\) cover respectively 55,
61, and 52 degree profiles. They use the descending-degree, suffix-cell
method and the two-dimensional incidence bound proved above. The source
differs from the independently audited seven-vertex implementation only
in comments, array capacity, and the host-order guard. The ordinary
coverage review and literal rank replay are recorded by
\texttt{src/\allowbreak{}small\_\allowbreak{}case\_\allowbreak{}tiny\_\allowbreak{}review.\allowbreak{}py} in
\texttt{research/\allowbreak{}small-\allowbreak{}case-\allowbreak{}tiny-\allowbreak{}review.\allowbreak{}json}. This is independent
source review, not a second complete negative enumeration at these host
orders.

For \((8,39)\) there are seventeen missing arrows. To exceed fourteen
guarantees, at most two missing pairs could fail. A zero missing row
makes every indirect pair fail. If \(\ell\) missing rows have degree
one, their degree-one loss bound forces \(\ell\le2\). When \(\ell=2\),
the two functional rows must form a two-cycle; the other rows avoid that
cycle and would give a universal six-vertex core with fifteen missing
pairs, contradicting \(M_V(6)=7\). When \(\ell=1\), its unique image row
already accounts for two failures. The remaining rows avoid its forced
bad targets. Restriction gives either an impossible universal
five-vertex core or a four-vertex core with one safe root. The latter
has one of the four degree profiles \((2,2,2,2)\), \((3,2,2,2)\),
\((4,2,2,2)\), or \((3,3,2,2)\), all excluded by the reviewed rooted
certificates. The root is useful only if the current cop cannot capture
there; root arrival is not declared unconditional success.

With no degree-one row, the missing degrees are \((3,2,\ldots,2)\).
Normalize the degree-three source and its neighbors, then enumerate
every remaining missing row. A union of permanent source and recipient
camping losses rejects only prefixes with at least three forced
failures. The complete search has 46,344 prefixes and 3,168 surviving
leaves. An independent orbit reconstruction covers every leaf using 22
fixed representatives; all 22,528 negative-state rank obligations are
checked. Their largest score is thirteen, while discarded graphs have
score at most fourteen. Thus the upper bound is fourteen, met by the
separately rank-verified literal graph in
\texttt{research/\allowbreak{}small-\allowbreak{}case-\allowbreak{}results.\allowbreak{}json}. The complete degree
reduction, corrected root semantics, source audit, and full leaf-orbit
review are in \texttt{research/\allowbreak{}small-\allowbreak{}sparse-\allowbreak{}eight39-\allowbreak{}review.\allowbreak{}md/\allowbreak{}json},
checked by \texttt{src/\allowbreak{}small\_\allowbreak{}sparse\_\allowbreak{}eight39\_\allowbreak{}review.\allowbreak{}py}. This review
does not claim a second traversal of the full search tree or a Lean
proof of its enumeration. \(\square\)

\textbf{Proposition (complete finite intervals).} The missing-pair upper
bound is attained at every budget in each of the following intervals: \[
\begin{array}{c|c}
n&\text{budgets }m\text{ with }M(n,m)=n(n-1)-m\\\hline
8&24\le m\le26\\
9&24\le m\le54\\
10&28\le m\le70\\
11&25\le m\le88
\end{array}
\]

\textbf{Proof.} Saved staged envelopes cover \([24,26]\) at order eight,
\([24,30]\cup[29,53]\) at order nine, \([28,69]\) at order ten, and
\([25,87]\) at order eleven. At each stage the messenger uses only
arrows of a fixed base graph \(B\), while the rank obligations permit
every cop move in the larger envelope \(E\). Hence the same decreasing
ranks prove delivery for every graph \(B\subseteq G\subseteq E\),
including every optional subset and every intermediate arrow budget. The
ranks depend on both positions and are checked at every recipient and
initial cop position. The order-nine stages obtained from the two ends
overlap; no inference from winning endpoints is used. The previously
certified dense endpoints supply \((9,54)\), \((10,70)\), and
\((11,88)\). All resulting counts equal the missing-pair upper bound.
The consolidated literal recipes and certificates are indexed in
\texttt{research/\allowbreak{}small-\allowbreak{}case-\allowbreak{}results.\allowbreak{}json};
\texttt{src/\allowbreak{}small\_\allowbreak{}case\_\allowbreak{}stage\_\allowbreak{}review.\allowbreak{}py} independently replays the
interval obligations. \(\square\)

\textbf{Verification boundary.} These are finite computer-assisted
proofs, with ordinary completeness and pruning arguments. Lean now
verifies every host table through order6 and further specified cells at
orders7 through11; the exact current domains are listed in the
verification section. The complete finite table is not yet a Lean
theorem. The seven-vertex vertex-maximum exclusion needs no game solver
after the static camping reduction; its complete joint table also uses
explicit game evaluation. The lower witnesses have saved rank
certificates checked directly. Every remaining degree-two graph and
every eight-vertex static survivor is evaluated by the exact game
algorithms in the independent exhaustive runs; full per-graph ranks for
those negative enumerations are not stored. The vertex maxima at orders
nine and ten remain open. The eleven-vertex maximum is exact, but this
does not determine every joint budget below 25. Absence of a two-regular
universal graph alone does not determine the vertex maximum.

Primary sources are \texttt{src/\allowbreak{}small\_\allowbreak{}exceptional\_\allowbreak{}n7\_\allowbreak{}exact.\allowbreak{}cpp} and
\texttt{src/\allowbreak{}small\_\allowbreak{}exceptional\_\allowbreak{}regular.\allowbreak{}cpp}; independent sources are
\texttt{src/\allowbreak{}small\_\allowbreak{}exceptional\_\allowbreak{}n7\_\allowbreak{}independent.\allowbreak{}cpp} and
\texttt{src/\allowbreak{}small\_\allowbreak{}exceptional\_\allowbreak{}independent.\allowbreak{}cpp}. The eight-vertex
sources are \texttt{src/\allowbreak{}small\_\allowbreak{}exceptional\_\allowbreak{}canonical.\allowbreak{}cpp} and
\texttt{src/\allowbreak{}small\_\allowbreak{}exceptional\_\allowbreak{}canonical\_\allowbreak{}independent.\allowbreak{}cpp}; the
independent ordinary audit is
\texttt{research/\allowbreak{}small-\allowbreak{}eight-\allowbreak{}canonical-\allowbreak{}independent-\allowbreak{}audit.\allowbreak{}md}. The
reproduction receipt is
\texttt{research/\allowbreak{}small-\allowbreak{}exceptional-\allowbreak{}verification.\allowbreak{}json} and the research
record is \texttt{research/\allowbreak{}small-\allowbreak{}exceptional-\allowbreak{}progress.\allowbreak{}md}.

\subsection{Stronger support bounds and the last small even
endpoint}\label{stronger-support-bounds-and-the-last-small-even-endpoint}

The shared deletion and positive-support argument also preserves the
stronger support results. Complete finite regular-core exclusions then
supply the remaining twenty-two-arrow case. The exact smaller budgets
are proved in the finite sparse-budget section before they are used
here.

\subsubsection{Two stronger support
consequences}\label{two-stronger-support-consequences}

\textbf{Proposition (strict improvement of the even support cutoff).}
For \(k\ge7\), a graph with \(2k\) arrows either has a
two-in/two-out-regular nonisolated support of order \(k\), or
\[I(G)\le U:=k^2-4k+2.\]

\textbf{Proof.} Positive supports are covered by (56.5) and the
missing-pair count. For an actual source or sink of other degree one,
deletion gives \(I\le A_{2k-1}\le k^2-5k+6<U\). Other degree at least
three gives \[I\le A_{2k-3}+k-2\le k^2-6k+10<U.\] For degree two, write
\(H\) for the deleted graph. The general cap is \((k-1)(k-3)=U+1\). A
score above \(U\) must attain this integer cap, so \(I(H)=(k-1)(k-4)\)
and the deleted row or column contributes \(k-1\). The even equality
statement makes \(H\) a regular \((k-1)\)-vertex core plus isolates. A
live deleted row needs two nonsink successors; a live deleted column
needs two positive-indegree predecessors. Both neighbors must therefore
be in that core. The two direct pairs leave at most \(k-3\) possible
additions, a contradiction. \(\square\)

\textbf{Proposition (original-graph support reduction).} For \(k\ge7\),
a graph with \(2k\) arrows and \(I(G)>V:=k^2-5k+8\) has exactly \(k\)
nonisolated vertices, whose degrees need not all be two. Consequently
\[M_E(2k)\le\max\{M(k,2k),k^2-5k+8\}.\]

\textbf{Proof.} Positive supports other than \(k\) are excluded by
(56.5) and the missing-pair count. Source/sink degrees one or at least
three give the same bounds as above, both below \(V\). For other degree
two, \[I(H)\ge I(G)-(k-1)>k^2-6k+9>(k-2)(k-4).\] Theorem 56 at parameter
\(k-1\) forces \(H\) to have a regular \((k-1)\)-vertex core plus
isolates. The deleted row or column must be live, since otherwise
\(I(G)\le I(H)\le(k-1)(k-4)<V\). Its two neighbors therefore lie in the
core, as in the preceding proof. Adjoining the deleted vertex gives
exactly \(k\) nonisolated vertices in the original graph. \(\square\)

At twenty-two arrows, any graph beating74 already fits on eleven
vertices. These implications concern the original graph, not a claim
that arbitrary arrow deletion or contraction preserves delivery. Their
earlier numerical consequences \(M_E(14)\le22\) and \(M_E(16)\le32\)
follow using \(M_V(7)=21,M_V(8)=32\); the finite-budget theorem gives
the sharper exact values 13 and 19.

\subsubsection{Regular cores and their finite
exclusions}\label{regular-cores-and-their-finite-exclusions}

\textbf{Lemma (two-row loss).} In an \(N\)-vertex two-in/two-out-regular
graph, a source camping obstruction forces at least \(2(N-3)\) missing
pairs to fail.

\textbf{Proof.} Suppose \(s\ne c\) and
\(N^+(s)\subseteq\{c\}\cup N^+(c)\). If \(c\notin N^+(s)\), the two
outgoing pairs are equal, and both rows \(s,c\) are zero. Otherwise
write \(N^+(s)=\{c,x\}\). The arrow \(c\to x\) makes \(N^-(x)=\{s,c\}\),
so row \(s\) and column \(x\) are both zero by camping. Their
intersection is the direct pair \((s,x)\), so their \(N-3\) missing
pairs are disjoint. \(\square\)

Without source camping, every predecessor-camping cop is outside the
recipient's predecessor set. Otherwise, if \(N^-(t)=\{c,x\}\) and
\(c\to x\), the row \(c\) is covered by the closed neighborhood of
\(x\), a source obstruction. Two cops camping the same recipient would
have identical outgoing pairs and likewise create source camping. Thus
distinct predecessor-camping pairs have distinct recipients, and two
such pairs destroy at least \(2(N-4)\) missing pairs.

The complete retained finite exclusions for zero or exactly one such
pair give the following bounds. All degrees in this table are two.

{\def\LTcaptype{none} 
\Needspace{7\baselineskip}
\begin{longtable}[]{@{}
  >{\raggedright\arraybackslash}p{(\linewidth - 8\tabcolsep) * \real{0.1579}}
  >{\raggedleft\arraybackslash}p{(\linewidth - 8\tabcolsep) * \real{0.2105}}
  >{\raggedleft\arraybackslash}p{(\linewidth - 8\tabcolsep) * \real{0.2105}}
  >{\raggedleft\arraybackslash}p{(\linewidth - 8\tabcolsep) * \real{0.2105}}
  >{\raggedleft\arraybackslash}p{(\linewidth - 8\tabcolsep) * \real{0.2105}}@{}}
\toprule\noalign{}
\begin{minipage}[b]{\linewidth}\raggedright
\(N\)
\end{minipage} & \begin{minipage}[b]{\linewidth}\raggedleft
Maximum with no camping
\end{minipage} & \begin{minipage}[b]{\linewidth}\raggedleft
Maximum with exactly one predecessor-camping pair
\end{minipage} & \begin{minipage}[b]{\linewidth}\raggedleft
Bound with at least two pairs
\end{minipage} & \begin{minipage}[b]{\linewidth}\raggedleft
Resulting bound
\end{minipage} \\
\midrule\noalign{}
\endhead
\bottomrule\noalign{}
\endlastfoot
9 & 30 & 24 & 44 & 44 \\
10 & 38 & 45 & 58 & 58 \\
11 & 66 & 71 & 74 & 74 \\
\end{longtable}
}

Source camping gives the smaller bound \(N(N-3)-2(N-3)\). This proves,
in particular, the retained eleven-vertex regular-core upper 74. The
exactly-one-pair enumerations evaluate 1,567,20,658,484,471 leaves at
the three orders, without time or target termination. The eleven-vertex
run covers all 14 derangement cycle types of one permutation factor,
with 7,016,903 partial nodes; an independent scalar solver agrees on482
selected leaves, including every successive maximizing witness. The
source and coverage review and replay are retained in
\texttt{research/\allowbreak{}sparse-\allowbreak{}third-\allowbreak{}audit.\allowbreak{}md} and
\texttt{research/\allowbreak{}sparse-\allowbreak{}third-\allowbreak{}enumeration.\allowbreak{}json}. They are complete
source-reviewed exclusions, not second full graph enumerations or Lean
proofs.

A further complete ten-vertex regular enumeration sharpens58 to53. It
covers 12 first-factor cycle types,5,652,419 partial nodes and 289,163
actual-game evaluations. The retained literal attains53. The
factorization, coverage and threshold-dependent camping review is
\texttt{research/\allowbreak{}small-\allowbreak{}plateau-\allowbreak{}regular-\allowbreak{}family-\allowbreak{}review.\allowbreak{}md}. The separate
incoming and outgoing single-transfer families on ten vertices also have
upper 53, with complete receipts
\texttt{research/\allowbreak{}small-\allowbreak{}sparse-\allowbreak{}transfer-\allowbreak{}in10-\allowbreak{}threshold53.\allowbreak{}json} and
\texttt{research/\allowbreak{}small-\allowbreak{}sparse-\allowbreak{}transfer-\allowbreak{}out10-\allowbreak{}threshold53.\allowbreak{}json}. These
three family results are premises of the finite twenty-arrow proof,
independent of its conclusion.

For reference, the more restrictive eligible-camping-free directional
transfer families retain their sharper finite maxima:

{\def\LTcaptype{none} 
\Needspace{6\baselineskip}
\begin{longtable}[]{@{}
  >{\raggedright\arraybackslash}p{(\linewidth - 4\tabcolsep) * \real{0.2727}}
  >{\raggedleft\arraybackslash}p{(\linewidth - 4\tabcolsep) * \real{0.3636}}
  >{\raggedleft\arraybackslash}p{(\linewidth - 4\tabcolsep) * \real{0.3636}}@{}}
\toprule\noalign{}
\begin{minipage}[b]{\linewidth}\raggedright
Order
\end{minipage} & \begin{minipage}[b]{\linewidth}\raggedleft
One outgoing degree1 and one outgoing degree3
\end{minipage} & \begin{minipage}[b]{\linewidth}\raggedleft
One incoming degree1 and one incoming degree3
\end{minipage} \\
\midrule\noalign{}
\endhead
\bottomrule\noalign{}
\endlastfoot
9 & 13 & 29 \\
10 & 32 & 50 \\
\end{longtable}
}

All other corresponding degrees, and every degree in the opposite
direction, are two. Hall's condition follows directly from these degree
patterns: source sets emit at least \(2|S|-1\) arrows against recipient
capacity two, or at least \(2|S|\) against capacity at most \(2|T|+1\).
A perfect matching can therefore be normalized by cycle type; the
high-degree vertex uses one representative of each distinct cycle length
and the low-degree vertex is unrestricted. Incoming and outgoing
families are evaluated separately. The complete scopes and audits remain
in the sparse-third-transfer-* evidence. None of the directional
statements assumes graph-reversal invariance.

\subsubsection{Exact attainment at twenty-two
arrows}\label{exact-attainment-at-twenty-two-arrows}

\textbf{Theorem (the twenty-two-arrow maximum).} For every \(n\ge11\),
\[M(n,22)=M_E(22)=79.\]

\textbf{Proof.} The preceding exact budgets give
\(D_{22}=\max(A_{21},A_{20}+10)=63\). Positive support at least twelve
has bound 74 from (56.4); support at most ten has at most 68 missing
pairs. At support eleven, a coordinate below eleven gives
\[B_{22}(11,10,11)=79,\] with all smaller coordinates bounded as in the
arithmetic lemma. If \(q=r=11\), every degree is two and the retained
regular-core upper 74 applies. This proves the global upper 79 without
either stronger support proposition being a premise.

For attainment, take vertices \(0,\ldots,10\) with outgoing pairs
\[(1,8),\ (2,5),\ (0,10),\ (0,4),\ (5,7),\ (0,6),\
(7,9),\ (2,8),\ (4,9),\ (1,10),\ (3,6).\] Every outdegree is two; vertex
0 has indegree three, vertex 3 indegree one, and all other indegrees
two. The saved full certificate verifies 79 guaranteed indirect pairs.
Both independent solvers agree, and the literal verifier checks every
winning and losing condition at all 2,662 alternating states across the
eleven recipients. The adjacency and ranks are in
\texttt{research/\allowbreak{}sparse-\allowbreak{}third-\allowbreak{}one-\allowbreak{}transfer-\allowbreak{}in-\allowbreak{}certificate.\allowbreak{}json}. Adding
isolates supplies every larger host. \(\square\)

The lower bound uses this literal graph, not completeness of its
discovery search. The upper uses the ordinary kernel and the retained
regular-core exclusions. \texttt{SparseWitness11.\allowbreak{}lean} separately
verifies an attaining eleven-vertex, twenty-two-arrow graph of exact
score 79; the unrestricted twenty-two-arrow upper bound remains outside
Lean.

\subsubsection{Additional outside-vertex
bounds}\label{additional-outside-vertex-bounds}

\textbf{Lemma (outside-vertex refinement).} Put \[\begin{gathered}
A=\{v:d^+(v)\ge2\},\qquad B=\{v:d^-(v)\ge2\},\\
q=|A|,\quad r=|B|,\quad b=|A\cap B|,\quad a=m-2q,\quad d=m-2r.
\end{gathered}\] Let \(Z\) be the nonisolated vertices outside
\(A\cup B\), and write \[z=|Z|,\qquad o=\sum_{v\in Z}d^+(v),\qquad
i=\sum_{v\in Z}d^-(v).\] Then \(0\le o\le a\), \(0\le i\le d\), and
\(\max(o,i)\le z\le o+i\). The nonisolated support is \(q+r-b+z\). If
\(R,C\) count the source rows and recipient columns containing at least
one guaranteed indirect pair, then
\[R\le\min\{q,b+a+\lfloor(d-o-z)/2\rfloor\},\]
\[C\le\min\{r,b+d+\lfloor(a-i-z)/2\rfloor\}.\] Moreover,
\[I\le qr-b-\max\{0,m-\min(a,r-b+o)-\min(d,q-b+i)\},\] and
\[I\le RC-\max\{0,R+C-q-r+b\}-\max\{0,2R+2C-m\}.\]

\textbf{Proof.} Every vertex of \(Z\) has indegree and outdegree at most
one, with positive sum, proving the elementary constraints and support
count. Let \(x=\sum_{v\in B}(d^-(v)-2)\). Exactly \(d-x\) arrows end
outside \(B\), of which exactly \(i\) end in \(Z\). Thus exactly
\(q-b-d+x+i\) vertices of \(A\setminus B\) have indegree zero; they emit
at least twice that many arrows. The \(z-i\) indegree-zero vertices in
\(Z\) each emit one further arrow. At most \(d-x\) of these arrows end
outside \(B\), and recipients of indegree at least three receive at most
\(3x\) arrows in total. Consequently at least \(\max(0,2q-2b-3d+i+z)\)
selected arrows enter indegree-two recipients. Dividing by their
capacity two gives that many distinct zero columns after taking the
ceiling. Such a recipient has an inaccessible predecessor; a cop waits
at its other predecessor. The cop-source exception is a direct pair and
is not counted. Rearranging, using the common parity of \(a,d\), gives
the displayed column bound.

For the row bound, use sinks in the original directions. If
\(y=\sum_{v\in A}(d^+(v)-2)\), exactly \(r-b-a+y+o\) members of
\(B\setminus A\) are sinks, and \(z-o\) further sinks in \(Z\) each
receive one arrow. At most \(a-y+3y\) of these incoming arrows can avoid
outdegree-two predecessors. Thus at least \(\max(0,2r-2b-3a+o+z)\)
arrows come from degree-two sources with a sink successor. Each such
source has a zero indirect row, since the cop waits at its other
successor. Dividing by two and rearranging proves the row bound; no
reversal invariance of the game is used.

Total outdegree outside \(A\) is at most \(\min(a,r-b+o)\), and total
indegree outside \(B\) is at most \(\min(d,q-b+i)\). Subtracting both
from \(m\) bounds the arrows in \(A\times B\) from below. Its \(b\)
diagonal pairs are disjoint from these arrows, giving the first pair
bound. The live row and column sets intersect in at least
\(\max(0,R+C-q-r+b)\) vertices. Their total outgoing and incoming
degrees are at least \(2R\) and \(2C\), forcing at least
\(\max(0,2R+2C-m)\) direct pairs in their rectangle. Removing these
arrows and diagonals gives the second bound. \(\square\)

Dropping the outside-vertex subtractions also gives the simpler bounds
\[R\le\min(q,b+a+\lfloor d/2\rfloor),\qquad
C\le\min(r,b+\lfloor a/2\rfloor+d).\] The full integer constraints imply
that a 20-arrow graph with \(I>50\) has at most 13 nonisolated vertices,
and an 18-arrow graph with \(I>30\) has at most 15. The independent
replay retains respectively 22 and 67 triples \((q,r,b)\), not that many
graphs or degree sequences. Its source and complete arithmetic receipt
are \texttt{research/\allowbreak{}small-\allowbreak{}case-\allowbreak{}support-\allowbreak{}review-\allowbreak{}2026-\allowbreak{}09-\allowbreak{}13.\allowbreak{}py} and the
adjacent JSON. These support bounds do not assert that deleting
exceptional vertices preserves a strategy.

\subsubsection{Consequences for the earlier small-budget
comparisons}\label{consequences-for-the-earlier-small-budget-comparisons}

The exact finite-budget theorem gives \(M_E(18)=33\) and \(M_E(20)=53\).
It therefore preserves the earlier intermediate conclusions
\(M_E(18)\le44\), \(M_E(20)\le58\), and \(53\le M_E(20)\le57\) without
their former separate optimization arguments. The regular and
directional family exclusions used in the final proof were stated
independently above.

\textbf{Proposition (comparison host for twenty arrows).} For every
\(n\ge11\), \[M(n,20)=M(11,20)=53.\] Every twenty-arrow graph has an
explicit eleven-vertex twenty-arrow comparison graph with at least as
many guarantees.

\textbf{Proof.} Pad the retained ten-vertex score 53 literal with one
isolate. The finite-budget theorem bounds every original graph by 53.
Further isolate padding proves the assertion for every larger host.
\(\square\)

The historical stronger task of modifying a hypothetical larger-support
improvement remains documented in
\texttt{research/\allowbreak{}small-\allowbreak{}sparse-\allowbreak{}twenty-\allowbreak{}host-\allowbreak{}compression.\allowbreak{}md} and its
independent review. It is no longer a dependency of the exact
classification. Its frozen evidence and the earlier intermediate
exclusions are preserved unchanged.

\subsection{The exact 23-arc endpoint}\label{the-exact-23-arc-endpoint}

The shared deletion and positive-support bounds reduce the
twenty-three-arrow endpoint to the same three finite families as the odd
support theorem.

\textbf{Theorem 53 (with finite exhaustive checks).}
\[\boxed{M_E(23)=83.}\] Moreover, \(M(n,23)=83\) for every \(n\ge12\).

\textbf{Proof.} The retained twelve-vertex literal has 83 guaranteed
indirect pairs. For the upper bound, the preceding exact budgets give
\[D_{23}=\max\{M_E(22),M_E(21)+10\}=79.\] Thus a graph with \(I\ge84\)
has positive indegree and outdegree after isolates are removed. Support
at most ten has at most 67 missing pairs; support at least thirteen has
at most 76 by (56.4). The retained eleven-vertex value \(M(11,23)=78\)
excludes that order.

At support twelve, \(q,r\le11\), and a coordinate at most ten gives
\[B_{23}(12,10,11)=81.\] Consequently \(q=r=11\). Degree sums force one
degree-one vertex in each direction and every other corresponding degree
two. The two exceptional vertices are either distinct, with their
connecting arrow absent or present, or coincide. Their eligible
missing-pair counts are respectively 90,89,89. This proves the full
finite domain directly; neither Theorem 57 nor a historical
eleven-vertex upper 81 is required for this endpoint.

We exhaust these families using their two-permutation descriptions. For
the first, add the missing exceptional arc \(u\to v\); for the second
add another copy of the existing arc; for the third add a dummy loop at
the common exceptional vertex. The resulting bipartite incidence
multigraph is two-regular and decomposes into two perfect matchings.
This follows directly by alternating the edges of each even cycle,
including a doubled-edge cycle. Regard the matchings as permutations
\(P,Q\), placing the added arc in \(P\).

Label the exceptional source 0 and, when distinct, the exceptional
recipient 1. Up to relabeling, enumerate the cycle of \(P\) containing
the added arc and the unordered cycle lengths of its remaining vertices.
In the common-exception case, \(P\) fixes 0. Every other cycle has
length at least two. Then enumerate every admissible permutation \(Q\),
allowing only the prescribed doubled exceptional arc in the second
family, and delete the added \(P\) arc before any game calculation. This
covers every graph in the three families.

Partial assignments may be discarded when a stationary pursuer covers
the completed outgoing neighbors of a contributing source, or the
completed predecessors of a contributing recipient. A source obstruction
destroys at least eight eligible pairs. A recipient obstruction destroys
at least seven, after omitting a possible source coinciding with the
pursuer. Since there are at most 90 eligible pairs, all discarded graphs
satisfy \(I\le83\). Testing coverage with a partially assigned pursuer
neighborhood is also safe: that neighborhood can only grow.

Every completed graph surviving these obstructions is evaluated by the
exact least reachability fixed point. Evaluation may stop at a leaf once
its remaining possible score falls below 84. The three finished
enumerations give the following counts; none has \(I\ge84\).

{\def\LTcaptype{none} 
\Needspace{7\baselineskip}
\begin{longtable}[]{@{}
  >{\raggedright\arraybackslash}p{(\linewidth - 6\tabcolsep) * \real{0.2000}}
  >{\raggedleft\arraybackslash}p{(\linewidth - 6\tabcolsep) * \real{0.2667}}
  >{\raggedleft\arraybackslash}p{(\linewidth - 6\tabcolsep) * \real{0.2667}}
  >{\raggedleft\arraybackslash}p{(\linewidth - 6\tabcolsep) * \real{0.2667}}@{}}
\toprule\noalign{}
\begin{minipage}[b]{\linewidth}\raggedright
Exceptional configuration
\end{minipage} & \begin{minipage}[b]{\linewidth}\raggedleft
Permutation cases
\end{minipage} & \begin{minipage}[b]{\linewidth}\raggedleft
Partial nodes
\end{minipage} & \begin{minipage}[b]{\linewidth}\raggedleft
Completed graphs surviving camping
\end{minipage} \\
\midrule\noalign{}
\endhead
\bottomrule\noalign{}
\endlastfoot
Distinct, arc absent & 42 & 164,276,910 & 16,902,091 \\
Distinct, arc present & 42 & 13,448,216 & 1,125,588 \\
Common vertex & 14 & 49,814,092 & 3,558,672 \\
\end{longtable}
}

This excludes every proposed counterexample and proves \(I\le83\). The
verified witness gives equality. Adding isolates preserves its
guaranteed indirect pairs, proving the fixed-order assertion for every
\(n\ge12\). \(\square\)

\textbf{Reproducibility and proof boundary.} The new source
\texttt{src/\allowbreak{}odd\_\allowbreak{}endpoint\_\allowbreak{}near23\_\allowbreak{}threshold.\allowbreak{}cpp} is a separately saved,
parameterized copy of the earlier enumerator; it accepts cutoffs at
least 84. The coverage, game calculation, and pruning are unchanged. The
complete receipts are
\texttt{research/\allowbreak{}odd-\allowbreak{}endpoint-\allowbreak{}threshold84-\allowbreak{}absent.\allowbreak{}json},
\texttt{research/\allowbreak{}odd-\allowbreak{}endpoint-\allowbreak{}threshold84-\allowbreak{}present.\allowbreak{}json}, and
\texttt{research/\allowbreak{}odd-\allowbreak{}endpoint-\allowbreak{}threshold84-\allowbreak{}same.\allowbreak{}json}. Their cutoff,
completion flags, exact counts, source hash, and dependencies are
consolidated in \texttt{research/\allowbreak{}odd-\allowbreak{}endpoint-\allowbreak{}threshold84-\allowbreak{}summary.\allowbreak{}json}.
The unchanged dependency on the eleven-vertex enumeration is audited in
\texttt{research/\allowbreak{}odd-\allowbreak{}endpoint-\allowbreak{}eleven-\allowbreak{}enumeration-\allowbreak{}audit.\allowbreak{}md}. The
source/sink and positive-support argument above supplies the current
all-host reduction. Its arithmetic is checked separately by
\texttt{src/\allowbreak{}sparse\_\allowbreak{}upper\_\allowbreak{}kernel\_\allowbreak{}check.\allowbreak{}py}. The earlier specialized
23-arc reduction and the unchanged cutoff-84 pruning are recorded in
\texttt{research/\allowbreak{}odd-\allowbreak{}endpoint-\allowbreak{}threshold84-\allowbreak{}ordinary-\allowbreak{}audit.\allowbreak{}md}. A separate
software audit of the coverage, pruning, and evaluation is
\texttt{research/\allowbreak{}odd-\allowbreak{}endpoint-\allowbreak{}threshold84-\allowbreak{}independent-\allowbreak{}audit.\allowbreak{}md}. It also
compares 37 leaf-game and early-exit cases across all three families
with an explicit alternating-state solver, including full rank and
losing-closure checks; the receipt is
\texttt{research/\allowbreak{}odd-\allowbreak{}endpoint-\allowbreak{}threshold84-\allowbreak{}independent-\allowbreak{}tests.\allowbreak{}json}. The
lower-bound certificate is
\texttt{research/\allowbreak{}odd-\allowbreak{}endpoint-\allowbreak{}m23-\allowbreak{}lower-\allowbreak{}certificate.\allowbreak{}json.\allowbreak{}gz}, checked by
two independently implemented game solvers and every saved rank
obligation. The exhaustive upper-bound computation is reviewed at the
source and mathematical level; it is not a separately implemented second
enumeration or a Lean proof.

\section{Verification and scope
boundary}\label{verification-and-scope-boundary}

\subsection{Exact computation and formal
verification}\label{exact-computation-and-formal-verification}

For a fixed recipient \(t\), let \(X_j(r)\) be the set of initial
pursuer positions from which delivery can be forced in at most \(j\)
messenger moves. Start with \(X_0(t)=V\) and \(X_0(r)=\varnothing\)
otherwise. Set \[
S_j(v)=\{c:C(c)\subseteq X_j(v)\}\quad(v\ne t),\qquad S_j(t)=V,
\] and update \[
X_{j+1}(r)=\left(\bigcup_{v\in C(r)}S_j(v)\right)\setminus\{r\}\quad(r\ne t),
\qquad X_{j+1}(t)=V.
\] The recurrence includes waits and excludes all nonterminal
collisions. Its least fixed point enforces eventual delivery. A separate
implementation builds the explicit alternating graph of states
\((r,c,\text{turn})\) and computes its reachability attractor. This is
standard finite reachability-game machinery (Berwanger 2009), rather
than a new general algorithm for graph games.

The certificate checker independently verifies decreasing ranks on the
selected messenger moves and on every pursuer response. Outside the
winning region it verifies closure under every messenger move and at
least one pursuer response. Thus a negative answer comes with an
avoidance strategy, not simply failure of a search. Standard
adjacency-list attractor analysis gives an all-recipient time bound
\(O(n^2m+n^3)\); that bound describes the standard implementation, not a
claimed running-time bound for the simpler synchronous reference solver.

\subsubsection{Certified finite extrema}\label{certified-finite-extrema}

Complete enumeration gives

{\def\LTcaptype{none} 
\Needspace{5\baselineskip}
\begin{longtable}[]{@{}lrrrrrr@{}}
\toprule\noalign{}
\(n\) & 1 & 2 & 3 & 4 & 5 & 6 \\
\midrule\noalign{}
\endhead
\bottomrule\noalign{}
\endlastfoot
\(M_V(n)\) & 0 & 0 & 0 & 2 & 4 & 7 \\
\end{longtable}
}

For \(n\le5\) every labeled digraph was enumerated. For \(n=6\), sorting
the vertex labels by degree reduces the search while retaining at least
one representative of every isomorphism class; valid degree bounds prune
only candidates unable to improve the current maximum. This certifies
extrema, not labeled distributions. Every attaining graph is
independently checked by the explicit alternating solver. The complete
joint tables, enumeration code, and receipts accompany the paper.

The complete finite sparse endpoint bridge covers every order
\(12\le n\le84\), joining the uniform family from order 85 in Theorem
40. The sharper projection in Theorem 45 and Corollary 45.1 starts the
odd family at support order 98, so its required finite bridge is
\(13\le n\le97\). The full earlier archive through order 144 remains
available unchanged. The graph lists, local ranks, projection proof and
direct verification receipts are cited with those results. Theorem 46
supplies the original small-order bounds; the completed finite
classification gives \(M_V(7),\ldots,M_V(11)=21,32,48,65,85\).

One clean 12-vertex example is a Cayley digraph of \(A_4\), with right
multiplication by two suitable 3-cycles. Its actual graph and
all-recipient finite-rank certificate are included. The Lean theorem
checks the explicit 12-vertex graph directly, so correctness does not
depend on trusting the group-label interpretation or the search that
discovered it.

\subsubsection{Lean coverage}\label{lean-coverage}

The artifact uses Lean 4.33.1 with Std only. Both players may wait and
receipt precedes collision. \texttt{SemanticCompleteness} proves that
finite winning trees are equivalent, on finite graphs, to one fixed
legal messenger policy that succeeds against every permitted initial
pursuer and every legal, history-dependent pursuer policy. No time
horizon is assumed in the actual game; the proof derives a finite bound.
The formal score therefore has the intended meaning.

\textbf{The complete large-budget formula.}
\texttt{LargeJointClassification.\allowbreak{}maximum} proves the exact piecewise
answer for every \(n\ge12\) and every admissible \(m\ge24\), including
the sparse parity formulas, the entire main interval
\(2n\le m\le n(n-3)\), the dense staircase and the zero band.
\texttt{JointExplicitAttainer.\allowbreak{}attainer} supplies a named adjacency
recipe for \emph{every} admissible budget on every host \(n\ge12\), with
its exact arrow count and the claimed score. For \(m\le23\) it selects
the retained small graph and pads it with isolates; for larger budgets
it uses the prescribed sparse, main or dense recipe. The graph is not
extracted from an existential maximum or selected by a search.
\texttt{JointGraphRecipes} joins these with the prescribed graphs on all
smaller hosts: a named graph, its looplessness, exact arrow count and
claimed score are proved for every admissible natural pair \((n,m)\).
Its strict import closure has13750 modules. No fresh graph search or
choice from an existential maximum is used. Optimality at the remaining
finite cases is the separate upper-bound obligation described below.

\textbf{Construction coverage and stronger guarantees.} The
main-interval assembly joins the uniform low-band construction from
order115 to the counted circle, hub and separated-code families, and
completes all finite orders12..114 with retained envelopes. Envelope
symmetry is transported before optional arrows are selected; the
resulting graph need not share the symmetry.
\texttt{MainIntervalClassification.\allowbreak{}even\_\allowbreak{}sparse} proves
\(M(n,2s)=s(s-3)\) for every \(n\ge s\ge12\).
\texttt{OddSparseConstruction.\allowbreak{}maximum\_\allowbreak{}odd} proves
\(M(n,2s-1)=(s-1)(s-4)+2\) for every \(n\ge s\ge13\). The finite bridges
and the periodic families supply every support order, with every initial
cop and legal reply covered. Original per-pattern and per-recipient
timing bounds remain attached to the individual theorems.

\texttt{DenseSharpClassification} proves exact maxima and prescribed
attainers throughout \(n(n-3)\le m\le n(n-1)\) at every \(n\ge11\). The
functional endpoint holds at every \(n\ge2\) and the zero band at all
orders. The adjacent envelope retains its two-move,
arbitrary-optional-subset guarantee for \(n\ge16\) and
\(n(n-4)\le m\le n(n-3)\). Separate five-offset, guarded,
protected-circle, pair-code and cyclic-family theorems retain their
stated domains and stronger subset or timing properties.

Saved smaller envelopes additionally prove the missing-pair maximum with
named graphs at \((n,m)=(7,21)\), at \(n=8\) for \(24\le m\le34\), at
\(n=9\) for \(24\le m\le54\), at \(n=10\) for \(25\le m\le70\), and at
\(n=11\) for \(25\le m\le88\). These are precisely the domain of
\texttt{SmallOrderMainClassification}; they do not assert all smaller
host tables.

\textbf{Unrestricted bounds and finite exclusions.} Missing-pair and
vertex bounds, actual source/sink deletion, support counting and the
general sparse parity bounds hold for arbitrary loopless graphs. Isolate
and one-predecessor-sink padding preserve the score exactly, including
cops entering the added sink. Complete degree-profile and row-choice
proofs give the sharp unrestricted upper bounds at budgets \(0..12\):
\(0,0,0,0,1,1,2,2,4,5,6,8,9\). The eleven- and twelve-arrow maxima are
attained on every host of order at least seven. The thirteen-arrow upper
bound is11, attained on every host of order at least eight. The
unrestricted twenty-two-arrow upper bound is79, attained on every host
of order at least eleven; its proof is independent of the unfinished
intermediate-budget chain. Every admissible budget on every host of
order at most six is also classified. Signed winning and losing
certificates verify all retained lower-witness scores at budgets
\(11..23\); a lower-witness score alone is not an upper bound over
competing graphs.

The complete nine-vertex exclusion proves \(I(G)\le53\) for every
loopless graph of that order. Its sharper regular-family theorem gives
exact maximum thirty when all incoming and outgoing degrees are two,
including every first-factor cycle type. This is distinct from the
unrestricted eighteen-arrow value33. The fourteen-arrow positive-support
theorem proves score at most thirteen when the nonisolated support has
at least nine vertices and every supported vertex has positive degree in
both directions; arbitrary isolates are allowed. Its complete
nine-vertex piece has eleven normalized joint profiles and 160718 DAG
nodes. These restricted theorems retain their stated hypotheses until
the remaining unrestricted exclusions are joined.

\textbf{Remaining formal obligations.} The unfinished part is the sharp
unrestricted exclusions at budgets \(14..21\) and23, the
still-unclassified upper cells of the joint tables at orders \(7..11\),
and their final all-parameter assembly. All prescribed graph
constructions have their claimed exact scores in Lean, for every
admissible pair. These exact graph scores do not by themselves settle
the remaining optimality obligations. Some alternative constructions and
stronger supplementary guarantees retain ordinary proof boundaries.

The complete classification is not yet a Lean theorem. Every distributed
module is checked by ordinary kernel reduction; custom axioms, proof
holes and native proofs are excluded. Reported axioms are confined to
the standard \texttt{propext}, \texttt{Classical.\allowbreak{}choice} and
\texttt{Quot.\allowbreak{}sound}. The source-pinned replay receipt identifies the
accepted distributed snapshot; scoped receipts record later target
checks separately until a complete release is reconciled. Producers are
untrusted and unnecessary for replaying literal sources.
\texttt{lean/\allowbreak{}README.\allowbreak{}md} gives exact theorem domains. Module counts are
inventories, not percentages of mathematical completeness.

\subsection{Further consequences and open
problems}\label{further-consequences-and-open-problems}

The random statement in Theorem 54, specialized to the full vertex set,
implies that for every fixed density \(0<p<1\) the random loopless
digraph guarantees all pairs with probability tending to one. Its
failure probability is at most \[n(n-1)^2\bigl(1-p^2(1-p)\bigr)^{n-3}.\]
Consequently \[\mathbb E[I(G)]=(1-p)n(n-1)-o(1).\] Indeed, the loss from
the missing-pair count is at most \(n(n-1)\) times the indicator of
failure; the probability decays exponentially. This expectation
conclusion is separate from the theorem's exact-budget conditioning and
deterministic construction.

The exact-count and explicit-construction objective is complete.
Corollary 47, the finite sparse-budget theorem, and the complete
small-host tables cover every admissible pair. Maximizing that joint
answer gives \(M_V\) and \(M_E\); those projections are consequences
rather than substitutes for the joint table.

Further simplification of the finite constructions and proofs remains
possible. It is not an unresolved value or missing attaining graph.
Classification of all extremizing graph shapes, a structural
characterization of all guaranteed pairs, and exact optimization of
delivery time lie outside this paper's boundary. The existing timing
results and remaining timing questions are preserved in the
\href{https://angelraychev.com/cop-prisoner/delivery-time/}{delivery-time
companion}. The ordinary and finite proofs establish the classification.
The accompanying formal development has the more limited coverage
specified above.

A finite construction at one circumference does not prove a uniform
family. Here the finite bridge is complete on a bounded range and the
retained ordinary constructions supply the infinite tail. A repeatable
marked vertex-extension rule remains a possible simplification, not an
unproved premise of the main-interval theorem.

A correct future proof must retain enough information about both
players. Simple distance-to-recipient potentials fail: some states at a
fixed small distance from the recipient require detours whose length
grows with the graph. Nor does safe indefinite motion imply delivery.
The finite witnesses and failed conjectures are preserved to make those
distinctions reproducible.

\section{Appendix: provenance and retained
formulations}\label{appendix-provenance-and-retained-formulations}

\subsection{Correcting the original layered
construction}\label{correcting-the-original-layered-construction}

The original graph remains a useful construction, but its printed count
is not correct. We give its exact count and separate safe intermediate
movement from terminal delivery.

Let \(q=2x\), \(a=\binom{2x}{x}/2\), and choose a family \(\mathcal F\)
containing one of each complementary pair of \(x\)-subsets of \([2x]\),
with every coordinate occurring in exactly \(a/2\) members. Such a
family exists precisely when positive \(x\) is not a power of two; a
proof follows below.

There are row vertices \(A_X,A'_X\) for \(X\in\mathcal F\), and
coordinate vertices \(b_i,b'_i,c_i,c'_i\) for \(i\in[q]\). The only arcs
are

\begin{itemize}
\tightlist
\item
  \(A_X\to b_i,b'_i\) when \(i\in X\);
\item
  \(b_i\to A'_Y\) when \(i\in Y\), and \(b'_i\to A'_Y\) when
  \(i\notin Y\);
\item
  \(A'_Y\to c_i,c'_i\) when \(i\in Y\);
\item
  \(c_i\to A_X\) when \(i\in X\), and \(c'_i\to A_X\) when
  \(i\notin X\).
\end{itemize}

The order is \(N=2a+4q\). The four directed layers are \(A\),
\(B\cup B'\), \(A'\), \(C\cup C'\).

\textbf{Theorem 16.} For the balanced original family,
\[I(G)=2a(N-1-q)+4q.\] At \(x=3\), the 44-vertex graph has exactly 764
indirect guaranteed pairs, rather than the draft's claimed lower bound
of 1252.

\textbf{Proof.} The rows form an antichain. The \(2q\) literal columns,
indicating membership or nonmembership of each coordinate, also form an
antichain, with size \(a/2\ge5\). They cannot coincide: equality on the
selected rows would remain true on their complements and hence on every
\(x\)-subset, which is impossible for distinct literals. Every vertex
has at least two successors.

From any noncollision state the messenger can safely advance one layer.
If both players are in the same layer, use incomparable outgoing
neighborhoods to choose a successor inaccessible to the pursuer. If the
pursuer is in the next layer, avoid its occupied vertex; it has no
within-layer arc. In the other two cases it cannot reach the next layer
in one move. This includes waiting responses.

From \(A_X\) to a specified \(A'_Y\) there is a two-move guarantee.
Against a pursuer \(A_Z\), choose a coordinate in \(X\setminus Z\), then
choose its primed or unprimed intermediate according to membership in
\(Y\). Against a pursuer in the intermediate layer, choose another
coordinate of \(X\); the other layers cannot interfere. The reverse
direction is identical. Safely advancing to the appropriate row layer
first therefore gives delivery to any row recipient in at most five
moves. There are \(2a\) such recipients, each with \(N-1-q\) indirect
sources.

For target \(b_i\) or \(b'_i\), its predecessor set is exactly
\(N^+(c_i)\). Lemma 2 excludes every indirect source except \(c_i\).
That exceptional source wins in two moves: literal-column
incomparability handles a pursuer in its own layer, and the other layer
cases are as above. Thus \((c_i,b_i)\) and \((c_i,b'_i)\) both count.
Symmetrically \((b_i,c_i)\) and \((b_i,c'_i)\) count, and no other
indirect sources do. These give exactly \(4q\) further pairs.
\(\square\)

The correction preserves the original \(N^2-O(N\log N)\) leading
estimate along its constructed orders. The new circulant family supplies
the stronger linear-deficit result at every sufficiently large order.
The original construction's regular incidence design and the later
extremal construction are distinct contributions.

\subsubsection{Balanced complementary
halves}\label{balanced-complementary-halves}

\textbf{Lemma 17.} For positive \(x\), the \(x\)-subsets of \([2x]\) can
be partitioned into two coordinate-balanced families with complementary
subsets in opposite families if and only if \(x\) is not a power of two.

\textbf{Proof.} Put \(C=\binom{2x}{x}\). Each coordinate occurs in
\(C/2\) subsets, so balance requires \(4\mid C\). The factorial
valuation formula gives \[v_2\binom{2x}{x}=s_2(x),\] where \(s_2(x)\) is
the number of ones in the binary expansion. Thus divisibility by four is
equivalent to \(x\) not being a power of two.

For sufficiency, distinguish coordinate \(2x\) and cyclically rotate the
other \(L=2x-1\) coordinates. Subsets containing \(2x\) correspond to
binary circular words with \(x-1\) ones and \(x\) zeros. Every orbit has
length \(L\): a proper repetition would have a repetition count dividing
both consecutive integers. The number \(C/(2L)\) of orbits is even,
since \(L\) is odd and \(4\mid C\).

Choose half of those orbits. In the first family put their subsets, and
put the complements of the subsets in the unchosen orbits. Exactly one
of each complementary pair is selected. The distinguished coordinate
occurs \(C/4\) times. Rotation invariance gives a common incidence count
\(r\) for every other coordinate; counting all incidences yields
\[ (2x-1)r+C/4=xC/2,\] so \(r=C/4\). The other family is balanced as
well. \(\square\)

This proof replaces the draft's incorrect stronger divisibility
assertion \(4x\mid\binom{2x}{x}\), which already fails at \(x=3\).

\subsection{Construction appendix: choosing the circulant
parameter}\label{construction-appendix-choosing-the-circulant-parameter}

\textbf{Proof of Lemma 12.} If \(n\equiv0\pmod4\), choose \(h=n/2-1\);
any common divisor with \(n\) divides two and is odd. If
\(n\equiv2\pmod4\), choose \(h=n/2-2\); the same argument uses four. If
\(n\equiv3\pmod4\), use \(h=(n-1)/2\).

It remains to handle \(n\equiv1\pmod4\). Write \(n-1=2^e d\) with \(d\)
odd and \(e\ge2\). If \(d\ge7\), take \(h=d\): it divides \(n-1\) and is
at most \((n-1)/4\). If \(d\in\{1,3,5\}\) and \(3\nmid n\), use
\(h=(n-3)/2\).

If \(d=1\) and \(3\mid n\), take \(h=(n-7)/2\). The number \(2^e+1\) is
not divisible by seven, since powers of two modulo seven are \(1,2,4\).
The case \(d=3\) and \(3\mid n\) is impossible. Finally, if \(d=5\) and
\(3\mid n\), then \(e\) is even. Use \(h=(n-11)/2\). Divisibility of
\(5\cdot2^e+1\) by eleven would imply \(2^e\equiv2\pmod{11}\), whereas
even powers have residues \(1,3,4,5,9\).

All choices are odd and smaller than \(n/2\). Their smallest possible
lower estimate is \((n-11)/2\ge7\), valid for \(n\ge25\). \(\square\)

\subsection{The sparse interval as a consequence of the stronger
envelope}\label{the-sparse-interval-as-a-consequence-of-the-stronger-envelope}

\textbf{Theorem 41.} For every integer \(n\ge2419\) and every
\(2n\le m\le3n\), an explicit universally winning graph exists, with
delivery within \(9\lceil n/6\rceil+26\) messenger moves and arbitrary
subsets of its prescribed optional arrows.

\textbf{Proof.} Theorem 52 gives the stronger assertion from
\(n\ge115\), through a larger exact budget endpoint, with the better
time bound \(9\lceil n/6\rceil+21\). \(\square\)

The earlier construction separated deletions by seventeen columns and
forbade optional arrows within sixteen columns of a deletion. Its
separate finite strategy certificates remain available. The new theorem
replaces that numerical result and its overlap role; the old proof is
retained in \texttt{research/\allowbreak{}sparse-\allowbreak{}budget-\allowbreak{}earlier-\allowbreak{}proof.\allowbreak{}md} for
reproducibility.

\textbf{Corollary 42.} There is an absolute \(N\) such that for every
\(n\ge N\) and \(2n\le m\le n(n-3)\), \[M(n,m)=n(n-1)-m.\]

\textbf{Proof.} Theorem 44 now permits the explicit choice \(N=115\),
with combinatorial attaining constructions throughout. \(\square\)

\subsection{A robust dense bridge}\label{a-robust-dense-bridge}

\textbf{Theorem 50.} For every \(n\ge16\) and every integer
\(2n\le h\le3n\), there is a deterministic \(n\)-vertex graph with
exactly \(h\) missing arcs in which every pair is guaranteed within two
messenger moves. Consequently \[M(n,n(n-1)-h)=h.\]

At each allowed order, one compulsory set of \(2n\) missing arcs and a
list of \(n\) optional missing arcs suffice: \textbf{every subset} of
the optional list preserves two-move delivery. The construction takes
polynomial time.

\textbf{Proof.} This is the robust envelope in Theorem 58. Its
parameterized modular families and four odd repairs cover every order
from 22; six further explicit optional lists cover \(16,\ldots,21\). The
ordinary danger-set arguments and the finite witnesses retain a
compulsory safe intermediary against every pursuer, so any optional
subset is permitted. Selecting \(h-2n\) optional arcs gives the exact
prescribed count. \(\square\)

Theorem 58 further enlarges this same robust dense interval to
\[2n\le h\le H(n)\] at every \(n\ge30\), where \(H(n)\) is its explicit
cap and \(H(n)=n^2/4-O(n)\) from the parity construction. The statement
here is its uniform initial interval, with its earlier label retained
for references. Theorem 58 gives the common constructions and proofs.
The earlier projective-incidence decomposition is no longer a separate
premise or a separate numerical result. The generic high-girth criterion
of Lemma 18 remains available for other graphs.

\subsection{The full-vertex
specialization}\label{the-full-vertex-specialization}

\textbf{Proposition 51.} Put \(L=n(n-1)\) and \(p=m/L\), where
\(0<m<L\). If \[(n-3)p^2(1-p)>4\log n+\log(8/3),\] there is a
deterministic algorithm, polynomial in \(n\), constructing an
\(n\)-vertex graph with exactly \(m\) arcs in which every pair is
guaranteed within two messenger moves.

\textbf{Proof.} Take \(H=V\), \(b=n\), and \(a=m\) in Theorem 54. The
condition forces \(n\ge4\). That theorem proves the stronger statement
with a specified routing set, exact incident budget, arbitrary exterior
arrows, and a random-sampling probability bound. \(\square\)

The former full-vertex conditioning argument and its \(1800\) overlap
are retained in the research history. The counting proof now appears
once in Theorem 54, and Theorem 44 gives the improved \(115\) overlap.

\subsection{Provenance and assistance}\label{provenance-and-assistance}

Alexandra Ignatova formulated the original problem and developed the
layered construction in her 2023 project under Angel Raychev's
mentorship, recorded in the retained February 2024 manuscript. Their
earlier work includes the upper bound \(I(G)\le n(n-3)\) and the
investigation of two-in/two-out-regular graphs. These results and this
research direction are not claimed as new in the present paper.

The September 2026 development, directed by Angel Raychev, repairs the
balancing and strategy arguments, corrects the original layered count,
and establishes the additional extremal bounds, attaining constructions,
finite computational proofs and formal results reported here. Astra 6,
operating through the Codex harness, made substantial contributions to
mathematical exploration and deductions, counterexample searches, proof
writing, programming, independent implementation and proof checks, and
Lean development. This assistance extended beyond language editing. The
paper distinguishes ordinary proofs, externally checked computational
proofs, and kernel-checked Lean results. The author is responsible for
the submitted account; tool-assisted checks do not constitute external
peer review. The formal verification was extended in October 2026 with
further AI assistance through Codex; its precise coverage is stated in
the verification section.

\subsection*{References}\label{references}
\addcontentsline{toc}{subsection}{References}

\protect\phantomsection\label{refs}

\end{document}